
\ifx\shlhetal\undefinedcontrolsequence\let\shlhetal\relax\fi

\input amstex
\expandafter\ifx\csname mathdefs.tex\endcsname\relax
  \expandafter\gdef\csname mathdefs.tex\endcsname{}
\else \message{Hey!  Apparently you were trying to
  \string\input{mathdefs.tex} twice.   This does not make sense.} 
\errmessage{Please edit your file (probably \jobname.tex) and remove
any duplicate ``\string\input'' lines}\endinput\fi




\catcode`\X=12\catcode`\@=11

\def\n@wcount{\alloc@0\count\countdef\insc@unt}
\def\n@wwrite{\alloc@7\write\chardef\sixt@@n}
\def\n@wread{\alloc@6\read\chardef\sixt@@n}
\def\r@s@t{\relax}\def\v@idline{\par}\def\@mputate#1/{#1}
\def\l@c@l#1X{\firstpart.#1}\def\gl@b@l#1X{#1}\def\t@d@l#1X{{}}

\def\crossrefs#1{\ifx\all#1\let\tr@ce=\all\else\def\tr@ce{#1,}\fi
   \n@wwrite\cit@tionsout\openout\cit@tionsout=\jobname.cit 
   \write\cit@tionsout{\tr@ce}\expandafter\setfl@gs\tr@ce,}
\def\setfl@gs#1,{\def\@{#1}\ifx\@\empty\let\next=\relax
   \else\let\next=\setfl@gs\expandafter\xdef
   \csname#1tr@cetrue\endcsname{}\fi\next}
\def\m@ketag#1#2{\expandafter\n@wcount\csname#2tagno\endcsname
     \csname#2tagno\endcsname=0\let\tail=\all\xdef\all{\tail#2,}
   \ifx#1\l@c@l\let\tail=\r@s@t\xdef\r@s@t{\csname#2tagno\endcsname=0\tail}\fi
   \expandafter\gdef\csname#2cite\endcsname##1{\expandafter
     \ifx\csname#2tag##1\endcsname\relax?\else\csname#2tag##1\endcsname\fi
     \expandafter\ifx\csname#2tr@cetrue\endcsname\relax\else
     \write\cit@tionsout{#2tag ##1 cited on page \folio.}\fi}
   \expandafter\gdef\csname#2page\endcsname##1{\expandafter
     \ifx\csname#2page##1\endcsname\relax?\else\csname#2page##1\endcsname\fi
     \expandafter\ifx\csname#2tr@cetrue\endcsname\relax\else
     \write\cit@tionsout{#2tag ##1 cited on page \folio.}\fi}
   \expandafter\gdef\csname#2tag\endcsname##1{\expandafter
      \ifx\csname#2check##1\endcsname\relax
      \expandafter\xdef\csname#2check##1\endcsname{}%
      \else\immediate\write16{Warning: #2tag ##1 used more than once.}\fi
      \multit@g{#1}{#2}##1/X%
      \write\t@gsout{#2tag ##1 assigned number \csname#2tag##1\endcsname\space
      on page \number\count0.}%
   \csname#2tag##1\endcsname}}

\def\multit@g#1#2#3/#4X{\def\t@mp{#4}\ifx\t@mp\empty%
      \global\advance\csname#2tagno\endcsname by 1 
      \expandafter\xdef\csname#2tag#3\endcsname
      {#1\number\csname#2tagno\endcsnameX}%
   \else\expandafter\ifx\csname#2last#3\endcsname\relax
      \expandafter\n@wcount\csname#2last#3\endcsname
      \global\advance\csname#2tagno\endcsname by 1 
      \expandafter\xdef\csname#2tag#3\endcsname
      {#1\number\csname#2tagno\endcsnameX}
      \write\t@gsout{#2tag #3 assigned number \csname#2tag#3\endcsname\space
      on page \number\count0.}\fi
   \global\advance\csname#2last#3\endcsname by 1
   \def\t@mp{\expandafter\xdef\csname#2tag#3/}%
   \expandafter\t@mp\@mputate#4\endcsname
   {\csname#2tag#3\endcsname\lastpart{\csname#2last#3\endcsname}}\fi}
\def\t@gs#1{\def\all{}\m@ketag#1e\m@ketag#1s\m@ketag\t@d@l p
\let\realscite\scite
\let\realstag\stag
   \m@ketag\gl@b@l r \n@wread\t@gsin
   \openin\t@gsin=\jobname.tgs \re@der \closein\t@gsin
   \n@wwrite\t@gsout\openout\t@gsout=\jobname.tgs }
\outer\def\localtags{\t@gs\l@c@l}
\outer\def\globaltags{\t@gs\gl@b@l}
\outer\def\newlocaltag#1{\m@ketag\l@c@l{#1}}
\outer\def\newglobaltag#1{\m@ketag\gl@b@l{#1}}

\newif\ifpr@ 
\def\m@kecs #1tag #2 assigned number #3 on page #4.%
   {\expandafter\gdef\csname#1tag#2\endcsname{#3}
   \expandafter\gdef\csname#1page#2\endcsname{#4}
   \ifpr@\expandafter\xdef\csname#1check#2\endcsname{}\fi}
\def\re@der{\ifeof\t@gsin\let\next=\relax\else
   \read\t@gsin to\t@gline\ifx\t@gline\v@idline\else
   \expandafter\m@kecs \t@gline\fi\let \next=\re@der\fi\next}
\def\pretags#1{\pr@true\pret@gs#1,,}
\def\pret@gs#1,{\def\@{#1}\ifx\@\empty\let\n@xtfile=\relax
   \else\let\n@xtfile=\pret@gs \openin\t@gsin=#1.tgs \message{#1} \re@der 
   \closein\t@gsin\fi \n@xtfile}

\newcount\sectno\sectno=0\newcount\subsectno\subsectno=0
\newif\ifultr@local \def\ultralocal{\ultr@localtrue}
\def\firstpart{\number\sectno}
\def\lastpart#1{\ifcase#1 \or a\or b\or c\or d\or e\or f\or g\or h\or 
   i\or k\or l\or m\or n\or o\or p\or q\or r\or s\or t\or u\or v\or w\or 
   x\or y\or z \fi}

\def\resetall{\global\advance\sectno by 1\subsectno=0
   \gdef\firstpart{\number\sectno}\r@s@t}
\def\resetsub{\global\advance\subsectno by 1
   \gdef\firstpart{\number\sectno.\number\subsectno}\r@s@t}
\def\newsection#1\par{\resetall\vskip0pt plus.3\vsize\penalty-250
   \vskip0pt plus-.3\vsize\bigskip\bigskip
   \message{#1}\leftline{\bf#1}\nobreak\bigskip}
\def\subsection#1\par{\ifultr@local\resetsub\fi
   \vskip0pt plus.2\vsize\penalty-250\vskip0pt plus-.2\vsize
   \bigskip\smallskip\message{#1}\leftline{\bf#1}\nobreak\medskip}


\newdimen\marginshift

\newdimen\margindelta
\newdimen\marginmax
\newdimen\marginmin

\def\margininit{       
\marginmax=3 true cm                  
				      
\margindelta=0.1 true cm              
\marginmin=0.1true cm                 
\marginshift=\marginmin
}    

\def\t@gsjj#1,{\def\@{#1}\ifx\@\empty\let\next=\relax\else\let\next=\t@gsjj
   \def\@@{p}\ifx\@\@@\else
   \expandafter\gdef\csname#1cite\endcsname##1{\citejj{##1}}
   \expandafter\gdef\csname#1page\endcsname##1{?}
   \expandafter\gdef\csname#1tag\endcsname##1{\tagjj{##1}}\fi\fi\next}
\newif\ifshowstuffinmargin
\showstuffinmarginfalse
\def\jjtags{\ifx\shlhetal\relax 
  \else
\ifx\shlhetal\undefinedcontrolseq
\else
\showstuffinmargintrue
\ifx\all\relax\else\expandafter\t@gsjj\all,\fi\fi \fi
}

\def\tagjj#1{\realstag{#1}\oldmginpar{\zeigen{#1}}}
\def\citejj#1{\rechnen{#1}\mginpar{\zeigen{#1}}}     

\def\rechnen#1{\expandafter\ifx\csname stag#1\endcsname\relax ??\else
                           \csname stag#1\endcsname\fi}

\newdimen\theight

\def\marginfont{\sevenrm}

\def\trymarginbox#1{\setbox0=\hbox{\marginfont\hskip\marginshift #1}%
		\global\marginshift\wd0 
		\global\advance\marginshift\margindelta}

\def \oldmginpar#1{%
\ifvmode\setbox0\hbox to \hsize{\hfill\rlap{\marginfont\quad#1}}%
\ht0 0cm
\dp0 0cm
\box0\vskip-\baselineskip
\else 
             \vadjust{\trymarginbox{#1}%
		\ifdim\marginshift>\marginmax \global\marginshift\marginmin
			\trymarginbox{#1}%
                \fi
             \theight=\ht0
             \advance\theight by \dp0    \advance\theight by \lineskip
             \kern -\theight \vbox to \theight{\rightline{\rlap{\box0}}%
\vss}}\fi}

\newdimen\upordown
\global\upordown=8pt
\font\tinyfont=cmtt8 
\def\mginpar#1{\smash{\hbox to 0cm{\kern-10pt\raise7pt\hbox{\tinyfont #1}\hss}}}
\def\mginpar#1{{\hbox to 0cm{\kern-10pt\raise\upordown\hbox{\tinyfont #1}\hss}}\global\upordown-\upordown}


\def\t@gsoff#1,{\def\@{#1}\ifx\@\empty\let\next=\relax\else\let\next=\t@gsoff
   \def\@@{p}\ifx\@\@@\else
   \expandafter\gdef\csname#1cite\endcsname##1{\zeigen{##1}}
   \expandafter\gdef\csname#1page\endcsname##1{?}
   \expandafter\gdef\csname#1tag\endcsname##1{\zeigen{##1}}\fi\fi\next}
\def\verbatimtags{\showstuffinmarginfalse
\ifx\all\relax\else\expandafter\t@gsoff\all,\fi}
\def\zeigen#1{\hbox{$\scriptstyle\langle$}#1\hbox{$\scriptstyle\rangle$}}


\def\margintag#1{\ifshowstuffinmargin\oldmginpar{\zeigen{#1}}\fi}

\def\(#1){\edef\dot@g{\ifmmode\ifinner(\hbox{\noexpand\etag{#1}})
   \else\noexpand\eqno(\hbox{\noexpand\etag{#1}})\fi
   \else(\noexpand\ecite{#1})\fi}\dot@g}

\newif\ifbr@ck
\def\eat#1{}
\def\[#1]{\br@cktrue[\br@cket#1'X]}
\def\br@cket#1'#2X{\def\temp{#2}\ifx\temp\empty\let\next\eat
   \else\let\next\br@cket\fi
   \ifbr@ck\br@ckfalse\br@ck@t#1,X\else\br@cktrue#1\fi\next#2X}
\def\br@ck@t#1,#2X{\def\temp{#2}\ifx\temp\empty\let\neext\eat
   \else\let\neext\br@ck@t\def\temp{,}\fi
   \def\teemp{#1}\ifx\teemp\empty\else\rcite{#1}\fi\temp\neext#2X}
\def\resetbr@cket{\gdef\[##1]{[\rtag{##1}]}}
\def\references{\resetbr@cket\newsection References\par}

\newtoks\symb@ls\newtoks\s@mb@ls\newtoks\p@gelist\n@wcount\ftn@mber
    \ftn@mber=1\newif\ifftn@mbers\ftn@mbersfalse\newif\ifbyp@ge\byp@gefalse
\def\defm@rk{\ifftn@mbers\n@mberm@rk\else\symb@lm@rk\fi}
\def\n@mberm@rk{\xdef\m@rk{{\the\ftn@mber}}%
    \global\advance\ftn@mber by 1 }
\def\rot@te#1{\let\temp=#1\global#1=\expandafter\r@t@te\the\temp,X}
\def\r@t@te#1,#2X{{#2#1}\xdef\m@rk{{#1}}}
\def\b@@st#1{{$^{#1}$}}\def\str@p#1{#1}
\def\symb@lm@rk{\ifbyp@ge\rot@te\p@gelist\ifnum\expandafter\str@p\m@rk=1 
    \s@mb@ls=\symb@ls\fi\write\f@nsout{\number\count0}\fi \rot@te\s@mb@ls}
\def\byp@ge{\byp@getrue\n@wwrite\f@nsin\openin\f@nsin=\jobname.fns 
    \n@wcount\currentp@ge\currentp@ge=0\p@gelist={0}
    \re@dfns\closein\f@nsin\rot@te\p@gelist
    \n@wread\f@nsout\openout\f@nsout=\jobname.fns }
\def\m@kelist#1X#2{{#1,#2}}
\def\re@dfns{\ifeof\f@nsin\let\next=\relax\else\read\f@nsin to \f@nline
    \ifx\f@nline\v@idline\else\let\t@mplist=\p@gelist
    \ifnum\currentp@ge=\f@nline
    \global\p@gelist=\expandafter\m@kelist\the\t@mplistX0
    \else\currentp@ge=\f@nline
    \global\p@gelist=\expandafter\m@kelist\the\t@mplistX1\fi\fi
    \let\next=\re@dfns\fi\next}
\def\symbols#1{\symb@ls={#1}\s@mb@ls=\symb@ls} 
\def\bigsymbol{\textstyle}
\symbols{\bigsymbol\ast,\dagger,\ddagger,\sharp,\flat,\natural,\star}
\def\ftnumbers{\ftn@mberstrue} \def\ftsymbols{\ftn@mbersfalse}
\def\paginal{\byp@ge} \def\resetftnumbers{\ftn@mber=1}
\def\ftnote#1{\defm@rk\expandafter\expandafter\expandafter\footnote
    \expandafter\b@@st\m@rk{#1}}

\long\def\jump#1\endjump{}
\def\ssum{\mathop{\lower .1em\hbox{$\textstyle\Sigma$}}\nolimits}

\def\qed{\nobreak\kern 1em \vrule height .5em width .5em depth 0em}
\def\newneq{\hbox{\rlap{\hbox to 1\wd9{\hss$=$\hss}}\raise .1em 
   \hbox to 1\wd9{\hss$\scriptscriptstyle/$\hss}}}
\def\subsetne{\setbox9 = \hbox{$\subset$}\mathrel{\hbox{\rlap
   {\lower .4em \newneq}\raise .13em \hbox{$\subset$}}}}
\def\supsetne{\setbox9 = \hbox{$\subset$}\mathrel{\hbox{\rlap
   {\lower .4em \newneq}\raise .13em \hbox{$\supset$}}}}

\def\vbar{\mathchoice{\vrule height6.3ptdepth-.5ptwidth.8pt\kern-.8pt}
   {\vrule height6.3ptdepth-.5ptwidth.8pt\kern-.8pt}
   {\vrule height4.1ptdepth-.35ptwidth.6pt\kern-.6pt}
   {\vrule height3.1ptdepth-.25ptwidth.5pt\kern-.5pt}}
\def\f@dge{\mathchoice{}{}{\mkern.5mu}{\mkern.8mu}}
\def\b@c#1#2{{\rm \mkern#2mu\vbar\mkern-#2mu#1}}
\def\b@b#1{{\rm I\mkern-3.5mu #1}}
\def\b@a#1#2{{\rm #1\mkern-#2mu\f@dge #1}}
\def\bb#1{{\count4=`#1 \advance\count4by-64 \ifcase\count4\or\b@a A{11.5}\or
   \b@b B\or\b@c C{5}\or\b@b D\or\b@b E\or\b@b F \or\b@c G{5}\or\b@b H\or
   \b@b I\or\b@c J{3}\or\b@b K\or\b@b L \or\b@b M\or\b@b N\or\b@c O{5} \or
   \b@b P\or\b@c Q{5}\or\b@b R\or\b@a S{8}\or\b@a T{10.5}\or\b@c U{5}\or
   \b@a V{12}\or\b@a W{16.5}\or\b@a X{11}\or\b@a Y{11.7}\or\b@a Z{7.5}\fi}}

\catcode`\X=11 \catcode`\@=12




\let\thischap\jobname

\def\partof#1{\csname returnthe#1part\endcsname}
\def\CHAPOF#1{\csname returnthe#1chap\endcsname}

\def\chapof#1{\CHAPOF{#1}}

\def\setchapter#1,#2,#3;{%
  \expandafter\def\csname returnthe#1part\endcsname{#2}%
  \expandafter\def\csname returnthe#1chap\endcsname{#3}%
}

\def\setprevious#1 #2 {%
  \expandafter\def\csname set#1page\endcsname{\input page-#2}
}


 \setchapter  E53,B,N;       \setprevious E53 null
 \setchapter  300z,B,A;       \setprevious 300z E53
 \setchapter  88r,B,I;       \setprevious 88r 300z
 \setchapter  600,B,II;       \setprevious  600 88r
 \setchapter  705,B,III;       \setprevious   705 600
 \setchapter  734,B,IV;        \setprevious   734 705
 \setchapter  300x,B,;      \setprevious   300x 734

 \setchapter 300a,A,V.A;      \setprevious 300a 88r
 \setchapter 300b,A,V.B;       \setprevious 300b 300a
 \setchapter 300c,A,V.C;       \setprevious 300c 300b
 \setchapter 300d,A,V.D;       \setprevious 300d 300c
 \setchapter 300e,A,V.E;       \setprevious 300e 300d
 \setchapter 300f,A,V.F;       \setprevious 300f 300e
 \setchapter 300g,A,V.G;       \setprevious 300g 300f

  \setchapter  E46,B,VI;      \setprevious    E46 734
  \setchapter  838,B,VII;      \setprevious   838 E46

\newwrite\pageout
\def\rememberpagenumber{\let\setpage\relax
\openout\pageout page-\jobname  \relax \write\pageout{\setpage\the\pageno.}}

\def\recallpagenumber{\csname set\jobname page\endcsname
\def\headmark##1{\rightheadtext{\chapof{\jobname}.##1}}\WRITETOC}
\def\setupchapter#1{\leftheadtext{\chapof{\jobname}. #1}}

\def\setpage#1.{\pageno=#1\relax\advance\pageno1\relax}

\def\cprefix#1{
\edef\theotherpart{\partof{#1}}\edef\theotherchap{\chapof{#1}}%
\ifx\theotherpart\thispart
   \ifx\theotherchap\thischap 
    \else 
     \theotherchap%
    \fi
   \else 
     \theotherchap\fi}

\def\sectioncite[#1]#2{%
     \cprefix{#2}#1}

\edef\thispart{\partof{\thischap}}
\edef\thischap{\chapof{\thischap}}

\def\lastpage of '#1' is #2.{\expandafter\def\csname lastpage#1\endcsname{#2}}

\def\yCITE[#1]#2{\cprefix{#2}.\scite{#2-#1}}

\newwrite\writetoc
\immediate\openout\writetoc \jobname.toc
\def\addcontents#1{\def\WRITETOC{\immediate\write\writetoc{\noexpand\tocentry{\chapof{\jobname}}{#1}{\number\pageno}}}}



\def\spuriousreset{}


\expandafter\ifx\csname citeadd.tex\endcsname\relax
\expandafter\gdef\csname citeadd.tex\endcsname{}
\else \message{Hey!  Apparently you were trying to
\string\input{citeadd.tex} twice.   This does not make sense.} 
\errmessage{Please edit your file (probably \jobname.tex) and remove
any duplicate ``\string\input'' lines}\endinput\fi

\sectno=-1   
\localtags
\jjtags
\NoBlackBoxes
\define\mr{\medskip\roster}
\define\sn{\smallskip\noindent}
\define\mn{\medskip\noindent}
\define\bn{\bigskip\noindent}
\define\ub{\underbar}
\define\wilog{\text{without loss of generality}}
\define\ermn{\endroster\medskip\noindent}

\define\dbca{\dsize\bigcap}
\define\rest{\restriction}
\define\dbcu{\dsize\bigcup}

\newbox\noforkbox \newdimen\forklinewidth
\forklinewidth=0.3pt   
\setbox0\hbox{$\textstyle\bigcup$}
\setbox1\hbox to \wd0{\hfil\vrule width \forklinewidth depth \dp0
                        height \ht0 \hfil}
\wd1=0 cm
\setbox\noforkbox\hbox{\box1\box0\relax}
\def\unionstick{\mathop{\copy\noforkbox}\limits}
\def\nonfork#1#2_#3{#1\unionstick_{\textstyle #3}#2}
\def\nonforkin#1#2_#3^#4{#1\unionstick_{\textstyle #3}^{\textstyle #4}#2}     
%
\setbox0\hbox{$\textstyle\bigcup$}
\setbox1\hbox to \wd0{\hfil{\sl /\/}\hfil}
\setbox2\hbox to \wd0{\hfil\vrule height \ht0 depth \dp0 width
                                \forklinewidth\hfil}
\wd1=0cm
\wd2=0cm
\newbox\doesforkbox
\setbox\doesforkbox\hbox{\box1\box0\relax}
\def\nunionstick{\mathop{\copy\doesforkbox}\limits}

\def\fork#1#2_#3{#1\nunionstick_{\textstyle #3}#2}
\def\forkin#1#2_#3^#4{#1\nunionstick_{\textstyle #3}^{\textstyle #4}#2}     

\define \nl{\newline}
\magnification=\magstep 1
\documentstyle{amsppt}

{    
\catcode`@11

\ifx\alicetwothousandloaded@\relax
  \endinput\else\global\let\alicetwothousandloaded@\relax\fi

\gdef\subjclass{\let\savedef@\subjclass
 \def\subjclass##1\endsubjclass{\let\subjclass\savedef@
   \toks@{\def\usualspace{{\rm\enspace}}\eightpoint}%
   \toks@@{##1\unskip.}%
   \edef\thesubjclass@{\the\toks@
     \frills@{{\noexpand\rm2000 {\noexpand\it Mathematics Subject
       Classification}.\noexpand\enspace}}%
     \the\toks@@}}%
  \nofrillscheck\subjclass}
} 


\expandafter\ifx\csname alice2jlem.tex\endcsname\relax
  \expandafter\xdef\csname alice2jlem.tex\endcsname{\the\catcode`@}
\else \message{Hey!  Apparently you were trying to
\string\input{alice2jlem.tex}  twice.   This does not make sense.}
\errmessage{Please edit your file (probably \jobname.tex) and remove
any duplicate ``\string\input'' lines}\endinput\fi

\expandafter\ifx\csname bib4plain.tex\endcsname\relax
  \expandafter\gdef\csname bib4plain.tex\endcsname{}
\else \message{Hey!  Apparently you were trying to \string\input
  bib4plain.tex twice.   This does not make sense.}
\errmessage{Please edit your file (probably \jobname.tex) and remove
any duplicate ``\string\input'' lines}\endinput\fi

\def\renewcommand{\newcommand}	       
\edef\cite{\the\catcode`@}%
\catcode`@ = 11
\let\@oldatcatcode = \cite
\chardef\@letter = 11
\chardef\@other = 12
%
%
%
%
\def\@innerdef#1#2{\edef#1{\expandafter\noexpand\csname #2\endcsname}}%
%
%
\@innerdef\@innernewcount{newcount}%
\@innerdef\@innernewdimen{newdimen}%
\@innerdef\@innernewif{newif}%
\@innerdef\@innernewwrite{newwrite}%
%
%
%
\def\@gobble#1{}%
%
%
%
\ifx\inputlineno\@undefined
   \let\@linenumber = \empty 
\else
   \def\@linenumber{\the\inputlineno:\space}%
\fi
%
%
%
\def\@futurenonspacelet#1{\def\cs{#1}%
   \afterassignment\@stepone\let\@nexttoken=
}%
\begingroup 
\def\\{\global\let\@stoken= }%
\\ 
\endgroup
\def\@stepone{\expandafter\futurelet\cs\@steptwo}%
\def\@steptwo{\expandafter\ifx\cs\@stoken\let\@@next=\@stepthree
   \else\let\@@next=\@nexttoken\fi \@@next}%
\def\@stepthree{\afterassignment\@stepone\let\@@next= }%
%
%
%
\def\@getoptionalarg#1{%
   \let\@optionaltemp = #1%
   \let\@optionalnext = \relax
   \@futurenonspacelet\@optionalnext\@bracketcheck
}%
%
%
\def\@bracketcheck{%
   \ifx [\@optionalnext
      \expandafter\@@getoptionalarg
   \else
      \let\@optionalarg = \empty
      \expandafter\@optionaltemp
   \fi
}%
\def\@@getoptionalarg[#1]{%
   \def\@optionalarg{#1}%
   \@optionaltemp
}%
%
%
%
\def\@nnil{\@nil}%
\def\@fornoop#1\@@#2#3{}%
\def\@for#1:=#2\do#3{%
   \edef\@fortmp{#2}%
   \ifx\@fortmp\empty \else
      \expandafter\@forloop#2,\@nil,\@nil\@@#1{#3}%
   \fi
}%
\def\@forloop#1,#2,#3\@@#4#5{\def#4{#1}\ifx #4\@nnil \else
       #5\def#4{#2}\ifx #4\@nnil \else#5\@iforloop #3\@@#4{#5}\fi\fi
}%
\def\@iforloop#1,#2\@@#3#4{\def#3{#1}\ifx #3\@nnil
       \let\@nextwhile=\@fornoop \else
      #4\relax\let\@nextwhile=\@iforloop\fi\@nextwhile#2\@@#3{#4}%
}%
%
%
%
\@innernewif\if@fileexists
\def\@testfileexistence{\@getoptionalarg\@finishtestfileexistence}%
\def\@finishtestfileexistence#1{%
   \begingroup
      \def\extension{#1}%
      \immediate\openin0 =
         \ifx\@optionalarg\empty\jobname\else\@optionalarg\fi
         \ifx\extension\empty \else .#1\fi
         \space
      \ifeof 0
         \global\@fileexistsfalse
      \else
         \global\@fileexiststrue
      \fi
      \immediate\closein0
   \endgroup
}%
%
%
%
%
\def\bibliographystyle#1{%
   \@readauxfile
   \@writeaux{\string\bibstyle{#1}}%
}%
\let\bibstyle = \@gobble
%
%
\let\bblfilebasename = \jobname
\def\bibliography#1{%
   \@readauxfile
   \@writeaux{\string\bibdata{#1}}%
   \@testfileexistence[\bblfilebasename]{bbl}%
   \if@fileexists
      \nobreak
      \@readbblfile
   \fi
}%
\let\bibdata = \@gobble
%
%
\def\nocite#1{%
   \@readauxfile
   \@writeaux{\string\citation{#1}}%
}%
\@innernewif\if@notfirstcitation
%
%
\def\cite{\@getoptionalarg\@cite}%
%
%
\def\@cite#1{%
   \let\@citenotetext = \@optionalarg
   \printcitestart
   \nocite{#1}%
   \@notfirstcitationfalse
   \@for \@citation :=#1\do
   {%
      \expandafter\@onecitation\@citation\@@
   }%
   \ifx\empty\@citenotetext\else
      \printcitenote{\@citenotetext}%
   \fi
   \printcitefinish
}%
\newif\ifweareinprivate
\weareinprivatetrue
\ifx\shlhetal\undefinedcontrolseq\weareinprivatefalse\fi
\ifx\shlhetal\relax\weareinprivatefalse\fi
\def\@onecitation#1\@@{%
   \if@notfirstcitation
      \printbetweencitations
   \fi
   \expandafter \ifx \csname\@citelabel{#1}\endcsname \relax
      \if@citewarning
         \message{\@linenumber Undefined citation `#1'.}%
      \fi
     \ifweareinprivate
      \expandafter\gdef\csname\@citelabel{#1}\endcsname{%
\strut 
\vadjust{\vskip-\dp\strutbox
\vbox to 0pt{\vss\parindent0cm \leftskip=\hsize 
\advance\leftskip3mm
\advance\hsize 4cm\strut\openup-4pt 
\rightskip 0cm plus 1cm minus 0.5cm ?  #1 ?\strut}}
         {\tt
            \escapechar = -1
            \nobreak\hskip0pt\pfeilsw
            \expandafter\string\csname#1\endcsname
             \pfeilso
            \nobreak\hskip0pt
         }%
      }%
     \else  
      \expandafter\gdef\csname\@citelabel{#1}\endcsname{%
            {\tt\expandafter\string\csname#1\endcsname}
      }%
     \fi  
   \fi
   \csname\@citelabel{#1}\endcsname
   \@notfirstcitationtrue
}%
%
%
\def\@citelabel#1{b@#1}%
%
%
\def\@citedef#1#2{\expandafter\gdef\csname\@citelabel{#1}\endcsname{#2}}%
%
%
%
\def\@readbblfile{%
   \ifx\@itemnum\@undefined
      \@innernewcount\@itemnum
   \fi
   \begingroup
      \def\begin##1##2{%
         \setbox0 = \hbox{\biblabelcontents{##2}}%
         \biblabelwidth = \wd0
      }%
      \def\end##1{}
      %
      %
      \@itemnum = 0
      \def\bibitem{\@getoptionalarg\@bibitem}%
      \def\@bibitem{%
         \ifx\@optionalarg\empty
            \expandafter\@numberedbibitem
         \else
            \expandafter\@alphabibitem
         \fi
      }%
      \def\@alphabibitem##1{%
         \expandafter \xdef\csname\@citelabel{##1}\endcsname {\@optionalarg}%
         \ifx\biblabelprecontents\@undefined
            \let\biblabelprecontents = \relax
         \fi
         \ifx\biblabelpostcontents\@undefined
            \let\biblabelpostcontents = \hss
         \fi
         \@finishbibitem{##1}%
      }%
      \def\@numberedbibitem##1{%
         \advance\@itemnum by 1
         \expandafter \xdef\csname\@citelabel{##1}\endcsname{\number\@itemnum}%
         \ifx\biblabelprecontents\@undefined
            \let\biblabelprecontents = \hss
         \fi
         \ifx\biblabelpostcontents\@undefined
            \let\biblabelpostcontents = \relax
         \fi
         \@finishbibitem{##1}%
      }%
      \def\@finishbibitem##1{%
         \biblabelprint{\csname\@citelabel{##1}\endcsname}%
         \@writeaux{\string\@citedef{##1}{\csname\@citelabel{##1}\endcsname}}%
         \ignorespaces
      }%
      %
      %
      \let\em = \bblem
      \let\newblock = \bblnewblock
      \let\sc = \bblsc
      \frenchspacing
      \clubpenalty = 4000 \widowpenalty = 4000
      \tolerance = 10000 \hfuzz = .5pt
      \everypar = {\hangindent = \biblabelwidth
                      \advance\hangindent by \biblabelextraspace}%
      \bblrm
      \parskip = 1.5ex plus .5ex minus .5ex
      \biblabelextraspace = .5em
      \bblhook
      \input \bblfilebasename.bbl
   \endgroup
}%
%
%
\@innernewdimen\biblabelwidth
\@innernewdimen\biblabelextraspace
%
%
%
\def\biblabelprint#1{%
   \noindent
   \hbox to \biblabelwidth{%
      \biblabelprecontents
      \biblabelcontents{#1}%
      \biblabelpostcontents
   }%
   \kern\biblabelextraspace
}%
%
%
%
\def\biblabelcontents#1{{\bblrm [#1]}}%
%
%
\def\bblrm{\rm}%
%
%
\def\bblem{\it}%
%
%
\def\bblsc{\ifx\@scfont\@undefined
              \font\@scfont = cmcsc10
           \fi
           \@scfont
}%
%
%
\def\bblnewblock{\hskip .11em plus .33em minus .07em }%
%
%
\let\bblhook = \empty
%
%
%
\def\printcitestart{[}
\def\printcitefinish{]}
\def\printbetweencitations{, }
\def\printcitenote#1{, #1}
%
%
%
\let\citation = \@gobble
%
%
%
\@innernewcount\@numparams
%
%
\def\newcommand#1{%
   \def\@commandname{#1}%
   \@getoptionalarg\@continuenewcommand
}%
%
%
\def\@continuenewcommand{%
   \@numparams = \ifx\@optionalarg\empty 0\else\@optionalarg \fi \relax
   \@newcommand
}%
%
%
\def\@newcommand#1{%
   \def\@startdef{\expandafter\edef\@commandname}%
   \ifnum\@numparams=0
      \let\@paramdef = \empty
   \else
      \ifnum\@numparams>9
         \errmessage{\the\@numparams\space is too many parameters}%
      \else
         \ifnum\@numparams<0
            \errmessage{\the\@numparams\space is too few parameters}%
         \else
            \edef\@paramdef{%
               \ifcase\@numparams
                  \empty  No arguments.
               \or ####1%
               \or ####1####2%
               \or ####1####2####3%
               \or ####1####2####3####4%
               \or ####1####2####3####4####5%
               \or ####1####2####3####4####5####6%
               \or ####1####2####3####4####5####6####7%
               \or ####1####2####3####4####5####6####7####8%
               \or ####1####2####3####4####5####6####7####8####9%
               \fi
            }%
         \fi
      \fi
   \fi
   \expandafter\@startdef\@paramdef{#1}%
}%
%
%
%
%
\def\@readauxfile{%
   \if@auxfiledone \else 
      \global\@auxfiledonetrue
      \@testfileexistence{aux}%
      \if@fileexists
         \begingroup
            \endlinechar = -1
            \catcode`@ = 11
            \input \jobname.aux
         \endgroup
      \else
         \message{\@undefinedmessage}%
         \global\@citewarningfalse
      \fi
      \immediate\openout\@auxfile = \jobname.aux
   \fi
}%
%
%
\newif\if@auxfiledone
\ifx\noauxfile\@undefined \else \@auxfiledonetrue\fi
%
%
%
%
\@innernewwrite\@auxfile
\def\@writeaux#1{\ifx\noauxfile\@undefined \write\@auxfile{#1}\fi}%
%
%
%
\ifx\@undefinedmessage\@undefined
   \def\@undefinedmessage{No .aux file; I won't give you warnings about
                          undefined citations.}%
\fi
%
%
\@innernewif\if@citewarning
\ifx\noauxfile\@undefined \@citewarningtrue\fi
%
%
%
\catcode`@ = \@oldatcatcode

\def\pfeilso{\leavevmode
            \vrule width 1pt height9pt depth 0pt\relax
           \vrule width 1pt height8.7pt depth 0pt\relax
           \vrule width 1pt height8.3pt depth 0pt\relax
           \vrule width 1pt height8.0pt depth 0pt\relax
           \vrule width 1pt height7.7pt depth 0pt\relax
            \vrule width 1pt height7.3pt depth 0pt\relax
            \vrule width 1pt height7.0pt depth 0pt\relax
            \vrule width 1pt height6.7pt depth 0pt\relax
            \vrule width 1pt height6.3pt depth 0pt\relax
            \vrule width 1pt height6.0pt depth 0pt\relax
            \vrule width 1pt height5.7pt depth 0pt\relax
            \vrule width 1pt height5.3pt depth 0pt\relax
            \vrule width 1pt height5.0pt depth 0pt\relax
            \vrule width 1pt height4.7pt depth 0pt\relax
            \vrule width 1pt height4.3pt depth 0pt\relax
            \vrule width 1pt height4.0pt depth 0pt\relax
            \vrule width 1pt height3.7pt depth 0pt\relax
            \vrule width 1pt height3.3pt depth 0pt\relax
            \vrule width 1pt height3.0pt depth 0pt\relax
            \vrule width 1pt height2.7pt depth 0pt\relax
            \vrule width 1pt height2.3pt depth 0pt\relax
            \vrule width 1pt height2.0pt depth 0pt\relax
            \vrule width 1pt height1.7pt depth 0pt\relax
            \vrule width 1pt height1.3pt depth 0pt\relax
            \vrule width 1pt height1.0pt depth 0pt\relax
            \vrule width 1pt height0.7pt depth 0pt\relax
            \vrule width 1pt height0.3pt depth 0pt\relax}

\def\pfeilsw{ \leavevmode 
            \vrule width 1pt height0.3pt depth 0pt\relax
            \vrule width 1pt height0.7pt depth 0pt\relax
            \vrule width 1pt height1.0pt depth 0pt\relax
            \vrule width 1pt height1.3pt depth 0pt\relax
            \vrule width 1pt height1.7pt depth 0pt\relax
            \vrule width 1pt height2.0pt depth 0pt\relax
            \vrule width 1pt height2.3pt depth 0pt\relax
            \vrule width 1pt height2.7pt depth 0pt\relax
            \vrule width 1pt height3.0pt depth 0pt\relax
            \vrule width 1pt height3.3pt depth 0pt\relax
            \vrule width 1pt height3.7pt depth 0pt\relax
            \vrule width 1pt height4.0pt depth 0pt\relax
            \vrule width 1pt height4.3pt depth 0pt\relax
            \vrule width 1pt height4.7pt depth 0pt\relax
            \vrule width 1pt height5.0pt depth 0pt\relax
            \vrule width 1pt height5.3pt depth 0pt\relax
            \vrule width 1pt height5.7pt depth 0pt\relax
            \vrule width 1pt height6.0pt depth 0pt\relax
            \vrule width 1pt height6.3pt depth 0pt\relax
            \vrule width 1pt height6.7pt depth 0pt\relax
            \vrule width 1pt height7.0pt depth 0pt\relax
            \vrule width 1pt height7.3pt depth 0pt\relax
            \vrule width 1pt height7.7pt depth 0pt\relax
            \vrule width 1pt height8.0pt depth 0pt\relax
            \vrule width 1pt height8.3pt depth 0pt\relax
            \vrule width 1pt height8.7pt depth 0pt\relax
            \vrule width 1pt height9pt depth 0pt\relax
      }


\def\widestnumber#1#2{}

\def\citewarning#1{\ifx\shlhetal\relax 
    \else
    \par{#1}\par
    \fi
}

\def\rm{\fam0 \tenrm}

\def\fakesubhead#1\endsubhead{\bigskip\noindent{\bf#1}\par}



%
%
%

%

\font\textrsfs=rsfs10
\font\scriptrsfs=rsfs7
\font\scriptscriptrsfs=rsfs5

\newfam\rsfsfam
\textfont\rsfsfam=\textrsfs
\scriptfont\rsfsfam=\scriptrsfs
\scriptscriptfont\rsfsfam=\scriptscriptrsfs

\edef\oldcatcodeofat{\the\catcode`\@}
\catcode`\@11

\def\Cal@@#1{\noaccents@ \fam \rsfsfam #1}

\catcode`\@\oldcatcodeofat


\expandafter\ifx \csname margininit\endcsname \relax\else\margininit\fi

\long\def\red#1\endred{}
\long\def\green#1\endgreen{}
\long\def\blue#1\endblue{}
\long\def\private#1\endprivate{}

\def\endred{ \unmatched endred! }
\def\endgreen{ \unmatched endgreen! }
\def\endblue{ \unmatched endblue! }
\def\endprivate{ \unmatched endprivate! }

\ifx\latexcolors\undefinedcs\def\latexcolors{}\fi

\def\emptycs{}
\def\evaluatelatexcolors{%
        \ifx\latexcolors\emptycs\else
        \expandafter\xxevaluate\latexcolors\xxfertig\evaluatelatexcolors\fi}
\def\xxevaluate#1,#2\xxfertig{\setupthiscolor{#1}%
        \def\latexcolors{#2}}


\font\smallfont=cmsl7
\def\rutgerscolor{\ifmmode\else\endgraf\fi\smallfont
\advance\leftskip0.5cm\relax}
\def\setupthiscolor#1{\edef\tmptmpcs{\noexpand\bgroup\noexpand\rutgerscolor
\noexpand\def\noexpand\currentcolor{#1}%
\noexpand}%
\expandafter\let\csname#1\endcsname\tmptmpcs
\def\tmptmpcs{\checkColorUnmatched{#1}\popthecolor}
\expandafter\let\csname end#1\endcsname\tmptmpcs}

\def\checkColorUnmatched#1{\def\expectcolor{#1}%
    \ifx\expectcolor\currentcolor   
    \else \edef\failhere{\noexpand\tryingToClose '\currentcolor' with end\expectcolor}\failhere\fi}

\def\currentcolor{???}

\def\popthecolor{\ifmmode\else\endgraf\fi\egroup}

\expandafter\def\csname#1\endcsname{}

\evaluatelatexcolors

 \let\outerhead\head
 \def\head{\innerhead}
 \let\innerhead\outerhead

 \let\outersubhead\subhead
 \def\subhead{\innersubhead}
 \let\innersubhead\outersubhead

 \let\outersubsubhead\subsubhead
 \def\subsubhead{\innersubsubhead}
 \let\innersubsubhead\outersubsubhead

 \let\outerproclaim\proclaim
 \def\proclaim{\innerproclaim}
 \let\innerproclaim\outerproclaim

 %
 %
 %
 %

\def\demo#1{\medskip\noindent{\it #1.\/}}
\def\enddemo{\smallskip}

\def\remark#1{\medskip\noindent{\it #1.\/}}
\def\endremark{\smallskip}

\pageheight{8.5truein}
\topmatter
\title{Strongly dependent theories} \endtitle
\author {Saharon Shelah \thanks {\null\newline 
This research was supported by the Israel Science Foundation. 
Publication 863. \null\newline
I would like to thank Alice Leonhardt for the beautiful typing. \null\newline
} \endthanks} \endauthor  

\affil{The Hebrew University of Jerusalem \\
Einstein Institute of Mathematics \\
Edmond J. Safra Campus, Givat Ram \\
Jerusalem 91904, Israel
 \medskip
 Department of Mathematics \\
 Hill Center-Busch Campus \\
  Rutgers, The State University of New Jersey \\
 110 Frelinghuysen Road \\
 Piscataway, NJ 08854-8019 USA} \endaffil

\abstract  We further investigate the class of models of 
a strongly dependent (first order complete)
theory $T$, continuing \cite{Sh:715}, \cite{Sh:783} and relatives.
Those are properties (= classes)
somewhat parallel to superstability among stable theory, though are
different from it even for stable theories.  We show equivalence of
some of their definitions, investigate relevant ranks and give some
examples, e.g. the first order theory of the $p$-adics is strongly dependent.
The most notable result is: if $|A| + |T| \le \mu,\bold I \subseteq
{\frak C}$ and $|\bold I| \ge \beth_{|T|^+}(\mu)$ then some $\bold J
\subseteq \bold I$ of cardinality 
$\mu^+$ is an indiscernible sequence over $A$.
\endabstract
\endtopmatter
\document

\newpage

\head {Annotated content} \endhead
 \spuriousreset
\bn
\S0 Introduction, pg.5-7
\bn
\S1 Strong dependent: basic variant, pg.8-28
\mr
\item "{${{}}$}"  [We define $\kappa_{\text{ict}}(T)$ and
strongly dependent ($=$ strongly$^1$ dependent 
$\equiv \kappa_{\text{ict}}(T)=\aleph_0$),
(\scite{dp1.2}), note preservation passing from $T$ to 
$T^{\text{eq}}$, preservation under interpretation 
(\scite{dp1.2.1}), equivalence of some versions of
``$\bar \varphi$ witness $\kappa < \kappa_{\text{ict}}(T)$"
(\scite{dp1.2.2}), and we deduce that without loss of generality 
$m=1$ in (\scite{dp1.2}).  An observation (\scite{dq.6})
 will help to prove the equivalence of some
variants.  To some extent, indiscernible sequences can replace an
element and this is noted in \scite{dp1.8}, \scite{dp1.9} dealing with
the variant $\kappa_{\text{icu}}(T)$.  We end with some examples, in particular
(as promised in \cite{Sh:783}) the first order theory of the $p$-adic
is strongly dependent and this holds for similar fields and for some
ordered abelian groups expanded by subgroups.  
Also there is a (natural) strongly stable not strongly$^2$ stable $T$.]
\endroster
\bn
\S2 Cutting indiscernible sequence and strongly$^\ell$, pg.29-42
\mr
\item "{${{}}$}"  [We give equivalent conditions to strongly dependent
by cutting indiscernibles (\scite{dp1.2.4}) and recall the parallel 
result for $T$
dependent.  Then we define $\kappa_{\text{ict},2}(T)$ (in
\scite{df2.3.5}) and show that it always almost is equal to
$\kappa_{\text{ict}}(T)$ in \scite{df2.4.8}.  The exceptional case is
``$T$ is strongly dependent but not strongly$^2$ dependent" for which we give
equivalent conditions (\scite{df2.3.5} and \scite{df2.2.9}.]
\endroster
\bn
\S3 Ranks, pg.43-53
\mr
\item "{${{}}$}"  [We define $M_0 \le_A M_1,M_0 \le_{A,p} M_2$ (in
\scite{dp1.1}) and observe some basic properties in \scite{dp1.1A}.
Then in \scite{dp1.3} for most $\ell=1,\dotsc,12$ we define
$<_\ell,<^\ell_{\text{at}},<^\ell_{\text{pr}},\le^\ell$, explicit
$\bar\Delta$-splitting and last but not least the ranks
dp-rk$_{\bar\Delta,\ell}({\frak x})$.  Easy properties are in
\scite{dp1.4}, the equivalence of ``rank is infinite", is $\ge
|T|^+,T$ is strongly dependent in \scite{dp1.4} and more basic
 properties in \scite{dp3.3}.  We then add more cases $(\ell > 12)$ to
 the main definition in order to deal with (verssion of) strongly 
dependency and  then have parallel claims.]
\endroster
\bn
\S4 Existence of indiscernibles, pg.54-57
\mr
\item "{${{}}$}"  [We prove that if $\mu \ge |A|+|T|$ and $a_\alpha
 \in {}^m{\frak C}$ for $\alpha < \beth_{\mu^+}$ \ub{then} for some $u
 \subseteq \beth_{\mu^+}$ of cardinality $\mu^+,\langle
 a_\alpha:\alpha \in u \rangle$ is indiscernible over $A$.]
\endroster
\bn
\S5 Concluding Remarks, pg.58-83
\mr
\item "{${{}}$}"  [We consider shortly several further relatives of strongly
 dependent.
\sn
\item "{${{}}$}"  $(A) \quad$ Ranks for dependent theories
\nl
We redefine explicitly $\bar\Delta$-splitting and
 dp-rk$_{\bar\Delta,\ell}$ for more cases, i.e. more $\ell$'s and for
the case of finite $\Delta_\ell$'s in a
 way fitting dependent $T$ (in \scite{dp1.5}), point out the
 basic equivalence (in \scite{dp1.5}), consider a variant
 (\scite{dp5.21}) and questions (\scite{dp1.5.A},\scite{dp2.53}).
\sn
\item "{${{}}$}"  $(B) \quad$ Minimal theories (or types)
\nl
We consider minimality, i.e., some candidates are 
parallel to $\aleph_0$-stable theories which are
 minimal.  It is hoped that some such definition will throw light on
the place of o-minimal theories.  We also consider
content minimality of types.
\sn
\item "{${{}}$}"  $(C) \quad$ Local ranks for super dependent and
indiscernibility
\nl
We deal with local ranks, giving a wide
 family parallel to superstable and then define some ranks parallel
 to those in \S3.
\sn
\item "{${{}}$}"  $(D) \quad$ Strongly$^2$ stable fields
\nl
We comment on strongly$^2$ dependent/stable fields.  In
particular for every infinite non-algebraically closed field
$K,\text{Th}({\frak K})$ is not strongly$^2$ stable.
\sn
\item "{${{}}$}"  $(E) \quad$ Strongly$^3$ dependent
\nl
We introduce strong$^{(3,*)}$ dependent/stable theories and remark on
them.  This is related to dimension
\sn
\item "{${{}}$}"  $(F) \quad$ Representability and strongly$_k$ dependent
\nl
We define and comment on representability and $\langle \bar b_t:t \in I\rangle$
being indiscernible for $I \in {\frak k}$. 
\sn
\item "{${{}}$}"  $(G) \quad$ strongly$_3$ stable and primely regular
types.  
\nl
We prove the density of primely regular types (for
strongly$_3$ stable $T$) and we comment how definable groups help.
\sn
\item "{${{}}$}"  $(H) \quad T$ is $n$-dependent
\nl
We consider strengthenings $n$-independent of ``$T$ is independent".]
\endroster
\newpage

\head {\S0 Introduction} \endhead  \resetall \sectno=0
 \spuriousreset
\bigskip

Our motivation is trying to solve the equations ``x/dependent =
superstable/stable" (e.g. among complete first order theories).  In
\cite[\S3]{Sh:783} mainly two approximate solutions are suggested:
strongly$^\ell$ dependent for $\ell=1,2$; here we try to investigate
them not relying on \cite[\S3]{Sh:783}.  
We define $\kappa_{\text{ict}}(T)$ generalizing $\kappa(T)$,
the definition has the form ``$\kappa < \kappa_{\text{ict}}(T)$ iff
there is a sequence $\langle \varphi_i(\bar x,\bar y_i):i <
\kappa\rangle$ of formulas such that ...".

Now $T$ is strongly dependent (= strongly$^1$ dependent) 
iff $\kappa_{\text{ict}}(T) = \aleph_0$;
prototypical examples are: the theory of dense linear orders, the
theory of real closed fields, the model completion of the theory of
trees (or trees with levels), and the theory of the $p$-adic fields
(and related fields with valuations).  (The last one is strongly$^1$ not
strongly$^2$ dependent, see \scite{dp0.16}.)

For $T$ superstable, if $\langle \bar a_t:t \in I\rangle$ is an
indiscernible set over $A$ and $C$ is finite then for some finite $I^*
\subseteq I,\langle \bar a_t:t \in I \backslash I^*\rangle$ is
indiscernible over $A \cup C$; moreover over $A \cup C \cup\{\bar
a_t:t \in I^*\}$.  In \S2 we investigate the parallel here, when $I$
is a linear order, complete for simplicity 
(as in \cite[\S3]{Sh:715} for dependent
theories).  But we get two versions: strongly$^\ell$ dependent $\ell=1,2$
according to whether we like to generalize the first version of the
statement above or the ``moreover".  

Next, in \S3, we define and investigate rank, not of types but of
related objects ${\frak x} = (p,M,A)$ where, e.g. $p \in \bold S^m(M \cup
A)$; but there are several variants.  For some of them we prove ``$T$
is strongly dependent iff the rank is always $< \infty$ iff the rank
is bounded by some $\gamma < |T|^+$".  We first deal with the ranks
related to ``strongly$^1$ dependent" and then for the ones related
to ``strongly$^2$ dependent".

Further evidence for those ranks being of interest is in \S4 where we
use them to get indiscernibles.  Recall that if $T$ is stable, $|A| \le
\lambda = \lambda^{|T|},a_\alpha \in {\frak C}$ for $\alpha < \mu :=
\lambda^+$ then for some stationary 
$S \subseteq \mu,\langle a_\alpha:\alpha \in
S\rangle$ is indiscernible over $A,|S|=\mu$, we can write this as
$\lambda \rightarrow (\lambda)^{< \omega}_{T,\mu}$.  We can get similar
theorems from set theoretic assumptions: e.g. $\mu$ a measurable
cardinal, very interesting and important but not for the present model
theoretic investigation.

We may wonder:  Can we classify first order theories by $\lambda
\rightarrow_T (\mu)_\kappa$, as was asked by Grossberg and the author
(see on this question \cite[2.9-2.20]{Sh:702}).  
A positive answer appears in \cite{Sh:197}, \ub{but} under a very strong
assumption on $T$:  not only $T$ is dependent but every subsets
$P_1,\dotsc,P_n$ of $|M|$ the theory Th$(M,P_1,\dotsc,P_n)$ is
dependent, i.e., being dependent is preserved by monadic expansions.

Here we prove that if $T$ is strongly stable and $\mu \ge |T|$ then
$\beth_{\mu^+} \rightarrow_T (\mu^+)^{< \omega}_{\mu^+}$.  We certainly
hope for a better result (using $\beth_n(|T|)$ for some fix $n$ or
even $(2^\mu)^+$ instead of
$\beth_{\mu^+}$) and weaker assumptions, say 
``$T$ is dependent" (or less) instead ``$T$ is strongly dependent".  
But still it seems worthwhile to prove  \scite{ind.1} particularly having
waited so long for something.

Let strongly$^\ell$ stable mean strongly$^\ell$ dependent and stable.
As it happens (for $T$), 
being superstable implies strong$^2$ stable implies
strong$^1$ stable but the inverses fail.  So strongly$^\ell$ dependent
does not really solve the equation we have started with.  However,
this is not necessarily bad, the notion 
``strongly$^\ell$ stable" seems interesting
in its own right; this applies to the further variants.

We give a ``simplest" example of a theory 
$T$ which is strongly$^1$ stable and not
strongly$^2$ stable in the end of \S1 as well as prove that the
(theories of the) $p$-adic field is strongly stable (for any prime
$p$) as well as similar enough fields.

In \S5 we comment on some further properties and ranks.  Such further
properties hopefully will be crucial in \cite{Sh:F705}, if it
materializes; it tries to deal mainly with
$K^{\text{or}}$-representable theories and contain other beginning
as well.  We comment on ranks parallel to those in \S3 suitable for
all dependent theories.

We further try to look at theories of fields.
Also we deal with the search for families of dependent theories $T$
which are unstable but ``minimal", much more well behaved.
For many years it seems quite bothering that we do not know how to
define o-minimality as naturally arising from a parallel to stability
theory rather than as an analog to minimal theories or generalizes examples
related to the theory of the field of the reals and its expansions.
Of course, the answer may be a somewhat larger class.  This motivates
Firstenberg-Shelah \cite{FiSh:E50} (on Th$(\Bbb R)$, specifically 
on ``perpendicularly is simple"), and some definitions 
in \S5.  Another approach to this question
is of Onshuus in his very illuminating works on th-forking
\cite{On0x1} and \cite{On0x2}.
\bn
A result from \cite[\S3,\S4]{Sh:783} used in \cite{FiSh:E50} says that
\proclaim{\stag{0.gr.1} Claim}    Assume $T$ is strongly$^2$
dependent
\mr
\item "{$(a)$}"  if $G$ is a definable group in ${\frak C}_T$ and $h$ is a
definable endomorphism of $G$ with finite kernels \ub{then} $h$ is almost
onto $G$, i.e., the index $(G:\text{\rm Rang}(h))$ is finite
\sn
\item "{$(b)$}"  it is not the case that: there are 
definable (with parameters) subset
$\varphi({\frak C},\bar a_1)$ of ${\frak C}$, an equivalence relation
$E_{\bar a_2} = E(x,y,\bar a_2)$ on $\varphi({\frak C},\bar a_1)$ with
infinitely many equivalence classes and
$\vartheta(x,y,z,\bar a_3)$ such that $E(c,c,\bar a_2) \Rightarrow
\vartheta(x,y,c,\bar a_3)$ is a one-to-one function from (a co-finite
subset of) $\varphi({\frak C},\bar a_1)$ into $c/E_{\bar a_2}$.
\endroster
\endproclaim
\bn
We continue investigating dependent theories in \cite{Sh:900},
\cite{Sh:877}, \cite{Sh:906}, more recently \cite{Sh:F931} and
Kaplan-Shelah and concerning definable 
groups in \cite{Sh:876}, \cite{Sh:886}.

We thank Moran Cohen, Itay Kaplan, Aviv Tatarski and a referee for
pointing out deficiencies.
\bigskip

\demo{\stag{0.0.7} Notation}  1) Let $\varphi^{\bold t}$ be $\varphi$
if $\bold t = 1$ or $\bold t =$ true \ub{and} 
$\neg \varphi$ if $\bold t=0$ or $\bold t=$
false.
\nl
2) $\bold S^\alpha(A,M)$ is the set of complete types over $A$ in $M$
(i.e. finitely satisfiable in $M$) in the free variables $\langle
x_i:i < \alpha \rangle$.
\enddemo
\newpage

\head {\S1 Strongly dependent: basic variant} \endhead  \resetall \sectno=1
 \spuriousreset
\bigskip

\demo{\stag{dp1.0} Convention}  1) $T$ is complete first order fixed.
\nl
2) ${\frak C} = {\frak C}_T$ a monster model for $T$.
\enddemo
\bn
Recall from \cite{Sh:783}:
\definition{\stag{dp1.2} Definition}  1) $\kappa_{\text{ict}}(T) =
\kappa_{\text{ict},1}(T)$ is
the minimal $\kappa$ such that for no $\bar \varphi = \langle
\varphi_i(\bar x,\bar y_i):i < \kappa \rangle$ is $\Gamma_\lambda =
\Gamma^{\bar \varphi}_\lambda$
consistent with $T$ for some ($\equiv$ every) $\lambda$ where $\ell
g(\bar x) = m,\ell g(\bar y^i_m) = \ell g(\bar y_i)$ and

$$
\Gamma_\lambda = \{\varphi_i(\bar x_\eta,\bar
y^i_\alpha)^{\text{if}(\eta(i)=\alpha)}:\eta \in{}^\kappa
\lambda,\alpha < \lambda \text{ and } i < \kappa\}.
$$
\mn
1A) We say that $\bar\varphi = \langle \varphi_i(\bar x,\bar y_i):i < \kappa
\rangle$ witness $\kappa < \kappa_{\text{ict}}(T)$ (with $m = \ell
g(\bar x)$) \ub{when} it is as in part (1).
\nl
2) $T$ is strongly dependent (or strongly$^1$ dependent) \ub{when}
$\kappa_{\text{ict}}(T) = \aleph_0$.
\enddefinition
\bn
Easy (or see \cite{Sh:783}):
\demo{\stag{dp.1.2N} Observation}  If $T$ is strongly dependent
\ub{then} $T$ is dependent.
\enddemo
\bigskip

\demo{\stag{dp1.2.1} Observation}  1) $\kappa_{\text{ict}}(T^{\text{eq}}) =
\kappa_{\text{ict}}(T)$.
\nl
2) If $T_\ell = \text{\rm Th}(M_\ell)$ for $\ell=1,2$ then
$\kappa_{\text{ict}}(T_1) \le \kappa_{\text{ict}}(T_2)$
when:
\mr
\item "{$(*)$}"  $M_1$ is (first order) interpretable in $M_2$.
\ermn
3) If $T' = \text{ Th}({\frak C},c)_{c \in A}$ then
$\kappa_{\text{ict}}(T') = \kappa_{\text{ict}}(T)$.
\nl
4) If $M$ is the disjoint sum of $M_1,M_2$ (or the product) and
Th$(M_1)$, Th$(M_2)$ are dependent then so is Th$(M)$; so $M_1,M_2,M$
has the same vocabulary.
\nl
5) In Definition \scite{dp1.2}, for some
   $\lambda,\Gamma^{\bar\varphi}_\lambda$ is consistent with $T$ iff
   for every $\lambda,\Gamma^{\bar\varphi}_\lambda$ is consistent with $T$.
\enddemo
\bigskip

\remark{Remark}  1) Concerning Part (4) for ``strongly dependent", see
Cohen-Shelah \cite{CoSh:E65}.
\endremark
\bigskip

\demo{Proof}  Easy.  \hfill$\square_{\scite{dp1.2.1}}$
\enddemo
\bigskip

\demo{\stag{dp1.2.2} Observation}   Let $\ell g(\bar x)=m;\bar
\varphi = \langle \varphi_i(\bar x,\bar y_i):i < \kappa \rangle$ and
let $\bar \varphi' = \langle \bar\varphi'_i(\bar x,\bar y'_i):i < \kappa
\rangle$ where $\varphi'_i(\bar x,\bar y'_i) = [\varphi_i(\bar x,\bar y^1_i) 
\wedge \neg \varphi_i(\bar x,\bar y^2_i)]$ and let $\bar\varphi'' = \langle
\varphi''_i(\bar x,\bar y''_i):i < \kappa \rangle$ where
$\bar\varphi''_i(\bar x,\bar y''_i) = [\bar\varphi_i(\bar x,\bar
y^1_i) \equiv \neg \varphi_i(\bar x,\bar y^2_i)]$.  \ub{Then}
$\circledast^1_{\bar\varphi} \Rightarrow \circledast^2_{\bar\varphi}
\Leftrightarrow \circledast^3_{\bar\varphi} 
\Leftrightarrow (\exists \eta \in {}^\kappa 2)
\circledast^2_{\bar\varphi^{[\eta]}} \Leftrightarrow (\exists \eta \in
{}^\kappa 2) \circledast^3_{\bar \varphi^{[\eta]}}$ and
$\circledast^\ell_{\bar\varphi} \Leftrightarrow 
\circledast^\ell_{{\bar\varphi}'} \Leftrightarrow
\circledast^\ell_{{\bar\varphi}''}$ for $\ell=2,3$ and
$\circledast^3_{\bar\varphi} \Leftrightarrow 
\circledast^1_{\bar\varphi'} \Leftrightarrow
\circledast^1_{\varphi''}$ where $\bar\varphi^{[\eta]} =
\langle \varphi_i(\bar x,\bar y_i)^{\eta(i)}:i < \kappa \rangle$ and
\mr
\item "{$\circledast^1_{\bar\varphi}$}"  $\bar \varphi$ witness $\kappa <
\kappa_{\text{ict}}(T)$
\sn
\item "{$\circledast^2_{\bar \varphi}$}"  we can find
$\langle \bar a^i_k:k < \omega,i < \kappa \rangle$
in ${\frak C}$ such that $\ell g(\bar a^i_k) = \ell g(\bar
y_i),\langle \bar a^i_k:k < \omega \rangle$ is indiscernible over
$\cup\{\bar a^j_k:j<\kappa,j \ne i,k <  \omega\}$ for each $i <
\kappa$ and $\{\varphi_i(\bar x,\bar a^i_0) \wedge \neg \varphi_i
(\bar x,\bar a^i_1):i < \kappa\}$ is consistent, i.e. finitely
satisfiable in ${\frak C}$
\sn
\item "{$\circledast^3_{\bar\varphi}$}"  like
$\circledast^2_{\bar\varphi}$ but in the end
\nl
$\{\varphi_i(\bar x,\bar a^i_0) \equiv \neg \varphi_i(\bar x,\bar
a^i_1):i < \kappa\}$ is consistent.
\endroster
\enddemo
\bigskip

\remark{\stag{dp1.2.7} Remark}  1) We could have added the 
indiscernibility condition to $\circledast^1_{\bar \varphi}$, i.e., to
\scite{dp1.2}(1) as this
variant is equivalent to $\circledast^1_{\bar \varphi}$.
\nl
2) Similarly we could have omitted the indiscernibility condition in
$\circledast^2_{\bar \varphi}$ but demand in the end: ``if $k_\ell <
\ell_i < \omega$ for $i < \kappa$ then $\{\varphi_i(\bar x,\bar
a^i_{k_i}) \wedge \neg \varphi_i(\bar x,a^i_{\ell_i}):i < \kappa\}$ is
consistent'' and get an equivalent condition.
\nl
3) Similarly we could have omitted the indiscernibility condition
in $\circledast^3_{\bar \varphi}$ but demand in the
end ``if $k_i < \ell_i <\omega$ for $i < \kappa$ then
$\{\varphi_i(\bar x,\bar a^i_{k_i}) \equiv \neg \varphi_i(\bar x,\bar
a^1_{\ell_i}):i < \kappa\}$ is consistent" and get an equivalent
condition.
\nl
4) We could add $\circledast^3_{\bar\varphi} \Leftrightarrow
\circledast^1_{\bar\varphi'}$. 
\nl
5) In $\circledast^2_{\bar\varphi},\circledast^3_{\bar\varphi}$ (and
the variants above) we can replace $\omega$ by any $\lambda$), see
\scite{dp1.2.3}).
\nl
6) What about $\circledast^2_{\bar\varphi} \Rightarrow
\circledast^1_{\bar\varphi}$?  We shall now describe a model whose theory is
a counterexample to this implication.  We define a model $M,\tau_M =
\{P,P_i,R_i:i < \kappa\},P$ a unary predicate, $P_i$ a unary
predicate, $R_i$ a binary predicate:
\mr
\item "{$(a)$}"  $|M|$ the universe of $M$ is 
$(\kappa \times \Bbb Q) \cup {}^\kappa \Bbb Q$
\sn
\item "{$(b)$}"  $P^M = {}^\kappa \Bbb Q$
\sn
\item "{$(c)$}"  $P^M_i = \{i\} \times \Bbb Q$
\sn
\item "{$(d)$}"  $R^M_i = \{(\eta,(i,q)):
\eta \in {}^\kappa \Bbb Q,q \in \Bbb Q \text{ and } 
\Bbb Q \models \eta(i) \ge q\}$
\sn
\item "{$(e)$}"  $\varphi_i(x,y) = P(x) \wedge P_i(y) \wedge 
R_i(x,y) \text{ for } i < \kappa$.
\ermn
Now
\mr
\item "{$(\alpha)$}"   Why (for Th$(\bar M)$) do we have
$\circledast^2_{\bar\varphi}$?
\ermn
For $i < \kappa,k < \omega$ let $a^i_k = (i,k) \in P^M_i$ recalling $\omega
\subseteq \Bbb Q$.
\nl
Easily $\langle a^i_k:k < \omega,i < \kappa\rangle$ are as required in
$\circledast^2_{\bar\varphi}$.
E.g. the unique $\eta \in {}^\kappa \Bbb Q$ realizing the type.  Also 
for each $i < \kappa$, the sequence $\langle a^i_k:k <
\omega\rangle$ is indiscernible over $\{a^j_m:j < \kappa,j \ne i$ and
$m < \omega\}$.
\nl
Why?  Because for every automorphism $\pi$ of the rational order
$(\Bbb Q,<)$, for the given $i < \kappa$ we can define a 
function $\hat \pi_i$ with domain $M$ by
\mr
\item "{$(*)_1$}"   for $j < \kappa$ and $q \in \Bbb Q$ we let
$\hat\pi_i((j,q))$ be $(j,q)$ if $j \ne i$ and $(j,\pi(q))$ if $j=i$
\sn
\item "{$(*)_2$}"   for $\eta,\nu \in {}^\chi \Bbb Q$ we have
$\hat\pi_i(\eta) = \nu$ iff $(\forall j < \kappa)(j
\ne i \Rightarrow \eta(j) = \nu(j))$ and $\nu(i) = \pi_i(\eta(i))$.
\ermn
So $\hat\pi_i$ is an automorphism of $M$ over $\dbcu_{j \ne i} P^M_j$
which includes the function $\{(a^i_q,a^i_{\pi(q)}):q \in \Bbb Q\}$
\mr
\item "{$(\beta)$}"   Why (for Th$(M))$, we do not have
$\circledast^1_{\bar\varphi}$?  
\nl
Because $M \models (\forall y_1,y_2)[P_i(y_1) \wedge P_i(y_2) 
\wedge y_1 \ne y_2 \rightarrow 
\dsize \bigvee^2_{\ell=1} (\forall x)(\varphi_i(x,y_\ell) \wedge P(x)
\rightarrow \varphi_i(x,y_{3-\ell}))]$.
\endroster
\endremark
\bigskip

\demo{Proof}  The following series of implications clearly suffices.
\mn
\ub{$\circledast^1_{\bar \varphi}$ implies
$\circledast^2_{\bar \varphi}$}

Why?  As $\circledast^1_{\bar \varphi}$, clearly for any $\lambda \ge
\aleph_0$ we can find $\bar a^i_\alpha \in {}^{\ell g(\bar y_i)}{\frak
C}$ for $i < \kappa,\alpha < \lambda$ and $\langle \bar c_\eta:\eta \in
{}^\omega \lambda \rangle,\bar c_\eta \in {}^{\ell g(\bar x)}{\frak
C}$ such that $\models \varphi_i[\bar c_\eta,\bar a^i_\alpha]$ iff
$\eta(i) = \alpha$.  By some applications of Ramsey theorem (or polarized
partition relations) \wilog \, $\langle \bar a^i_\alpha:\alpha <
\lambda \rangle$ is indiscernible over $\cup\{\bar a^j_\beta:j <
\kappa,j \ne i,\beta < \lambda\}$ for each $i < \omega$.  
Now those $\bar a^i_\alpha$'s
witness $\circledast^2_{\bar \varphi}$ as $\bar c_\eta$ witness the
consistency of the required type when $\eta \in {}^\kappa\{0\}$.
\bn
\ub{$\circledast^2_{\bar \varphi} \Rightarrow
\circledast^3_{\bar \varphi}$} (hence in particular
$\circledast^2_{{\bar \varphi}'} \Rightarrow
\circledast^3_{{\bar \varphi}'}$ and $\circledast^2_{{\bar \varphi}''} 
\Rightarrow \circledast^3_{{\bar \varphi}''}$).

Trivial; read the definitions.
\bn
\ub{$\circledast^3_{\bar \varphi} \Rightarrow
\circledast^2_{\bar \varphi}$} (hence in particular
$\circledast^3_{\bar \varphi'} \Rightarrow \circledast^3_{\bar
\varphi'}$ and $\circledast^2_{\bar \varphi''} \Rightarrow
\circledast^3_{\bar \varphi''}$).

By compactness, for the dense linear order $\Bbb R$ we can find
$\bar a^i_t$ for $i < \kappa,t \in \Bbb R$ such that for each $i <
\kappa$ the sequence
$\langle \bar a^i_t:t \in \Bbb R\rangle$ indiscernible over
$\cup\{\bar a^j_s:j \ne i,j < \kappa,s \in \Bbb R\}$ and for any $s_0
<_{\Bbb R} s_1$ the set $\{\varphi_i(\bar x,
\bar a^i_{s_0}) \equiv \neg \varphi_i(\bar x,\bar a^i_{s_1}):
i < \kappa\}$ is consistent, say realized by $\bar c = \bar c_{s_0,s_1}$.  Now
let $u = \{i < \kappa:{\frak C} \models \varphi_i[\bar c,\bar
a^i_{s_0}]\}$ and for $n < \omega$ define $\bar b^i_n$ as 
$\bar a^i_{s_0 + n(s_1-s_0)}$ if $i \in u$ and
as $\bar a^i_{s_1-n(s_1-s_0)}$ if $i \in \kappa \backslash u$.  Now
$\langle \bar b^i_n:n < \omega,i < \kappa\rangle$ exemplifies
$\circledast^2_{\bar \varphi}$.
\bn
\ub{$\circledast^2_{\bar \varphi}$ implies
$\circledast^1_{\bar \varphi'}$} (hence by the above
$\circledast^2_{\bar \varphi} \Rightarrow \circledast^2_{\bar
\varphi'}$ and $\circledast^3_{\bar \varphi} \Rightarrow
\circledast^3_{\bar \varphi'}$).

Let $\langle \bar a^i_\alpha:\alpha < \omega,i < \kappa \rangle$
witness $\circledast^2_{\bar \varphi}$ and $\bar c$ realizes
$\{\varphi_i(\bar x,a^i_0) \wedge \neg \varphi_i(\bar x,\bar a^i_1):i
< \kappa\}$.  Without loss of generality 
$a^i_t$ is well defined for every $t \in \Bbb Z$ not
just $t \in \omega$, (and $i < \kappa$), and $\langle a^i_t:t \in \Bbb
Z\rangle$ is an indiscernible sequence over $\{a^j_s:j \in \kappa
\backslash \{i\}$ and $s \in \Bbb Z\}$.  Also \wilog \, 
for each $i < \kappa,\langle
\bar a^i_\alpha:\alpha \in [2,\omega)\rangle$ as well as $\langle
a^i_{-1-n}:n \in \omega\rangle$ are indiscernible sequences 
over $\cup\{\bar a^j_t:j < \kappa,j \ne i$ and $t \in \Bbb Z\} 
\cup \{\bar c\}$.

For $t \in \Bbb Z,i < \kappa$ let $\bar b^i_t = \bar a^i_{2 t} \char 94
\bar a^i_{2 t+1}$, so ${\frak C} \models \varphi'_i[\bar c,\bar b^i_0]$ (as
this just means ${\frak C} \models \varphi_i(\bar c,\bar a^i_0) \wedge \neg
\varphi_i[\bar c,\bar a^i_1])$ and ${\frak C} \models \neg \varphi'_i[\bar
c,\bar b^i_s]$ when $s \in \Bbb Z \backslash \{0\}$ 
(as the sequences $\bar c \char 94 \bar a^i_{2s}$ 
and $\bar c \char 94 \bar a^i_{2s+1}$ realize the
same type).  So $\langle \bar b^i_\alpha:\alpha < \omega,i < \kappa
\rangle$ witness $\circledast^1_{\bar \varphi'}$.
\bn
\ub{$\circledast^3_{\bar \varphi'}$ implies
$\circledast^3_{{\bar\varphi}''}$}.

Read the definitions.
\bn
\ub{$\circledast^3_{{\bar\varphi}''}$ implies that for some $\eta \in
{}^\kappa 2$ we have $\circledast^1_{{\bar\varphi}^{[\eta]}}$}.

As in the proof of $\circledast^2_\varphi \Rightarrow
\circledast^1_{\bar \varphi'}$; but we elaborate:
let $\big< \langle \bar a^i_\alpha \char 94 \bar b^i_\alpha:\alpha <
\omega\rangle:i < \kappa\big>$ witness $\circledast^3_{\bar\varphi''}$
noting $\bar\varphi'' = \langle \varphi''_i(\bar x,\bar y^i_1,\bar
y^i_2):i < \kappa\rangle$ where 
$\ell g(\bar y^i_1) = \ell g(\bar y_i) = \ell g(\bar y^i_2)$.  Let $\bar
c$ realize $\{\varphi''_i(\bar x,\bar a^i_0,\bar b^i_0) \equiv \neg
\varphi''_i(\bar x,\bar a^i_1,\bar b^i_1):i < \kappa\}$.  Without loss
of generality for each $i < \kappa$ the sequence $\langle \bar
a^i_\alpha \char 94 \bar b^i_\alpha:2 \le \alpha < \omega\rangle$ is
indiscernible over $\cup\{\bar a^j_\alpha \char 94 \bar b^j_\alpha:j
\in \kappa \backslash \{i\}$ and $\alpha < \omega\} \cup \bar c$.

By this extra indiscernibility assumption for each $i < \kappa$ we can
find $\ell_0(i),\ell_1(i) \in \{0,1\}$ such that $n \ge 2 \Rightarrow
{\frak C} \models \varphi_i[\bar c,\bar a^i_n]^{\ell_0(i)} \wedge
\varphi_i[\bar c,\bar b^i_n]^{\ell_1(i)}$.  By the choice of $\bar c$ we have
${\frak C} \models \varphi''_i(\bar c,\bar a^i_0,\bar b^i_0)
\equiv \varphi''_i(\bar c,\bar a^i_1,\bar b^i_1)$, hence by the choice of
$\varphi''_i$, we \ub{cannot} have ${\frak C} \models \varphi_i[\bar
c,\bar a^i_0]^{\ell_0(i)} \wedge \varphi_i[\bar c,\bar a,\bar
b^i_0]^{\ell_1(i)} \wedge \varphi_i[\bar c,\bar a^i_1]^{\ell_0(i)}
\wedge \varphi_i[\bar c,\bar b^i_1]^{\ell_1(i)}$.

Hence there are $\ell_3(i),\ell_4(i) \in \{0,1\}$ such that
\mr
\item "{$\bullet$}"  $\ell_4(i) = 0 \Rightarrow {\frak C} \models
\varphi_i[\bar c,\bar a^i_{\ell_3(i)}]^{1-\ell_0(i)}$
\sn
\item "{$\bullet$}"  $\ell_4(i) = 1 \Rightarrow {\frak C} \models
\varphi_i[\bar c,\bar b^i_{\ell_3(i)}]^{1-\ell_1(i)}$.
\ermn
Lastly choose $\eta = \langle 1- \ell_{\ell_4(i)}(i):i <
\kappa\rangle$ and we choose $\langle \bar d^i_\alpha:\alpha <
\omega,i < \kappa\rangle$ as follows:
\mr
\item "{$\bullet$}"  if $\ell_4(i) = 0$ and $n=0$ then $\bar d^i_n =
\bar a^i_{\ell_3(i)}$
\sn
\item "{$\bullet$}"  if $\ell_4(i) = 0$ and $n>0$ then $\bar d^i_n =
\bar a^i_{1+n}$
\sn
\item "{$\bullet$}"  if $\ell_4(i) = 1$ and $n=0$ then $\bar d^i_n =
\bar b^i_{\ell_3(i)}$
\sn
\item "{$\bullet$}"  if $\ell_4(i) = 1$ and $n > 0$ then $\bar d^i_n =
\bar b^i_{1+n}$.
\ermn
Now check that $\langle \bar d^i_\alpha:\alpha < \omega$ and $i <
\kappa\rangle$ witness $\circledast^1_{\bar\varphi^{[\eta]}}$.
\bn
\ub{$\circledast^3_{\bar
\varphi^{[\eta]}},\circledast^3_{\bar\varphi}$ are equivalent where $\eta\in
{}^\kappa 2$}.

Why?  Because the formula $(\varphi_i(x,\bar a^i_0) \equiv \neg
\varphi_i(x,\bar a^i_1))$ is equivalent to $(\varphi_i(x,
a^i_0)^{\eta(i)} \equiv \neg \varphi_i(x,\bar a^i_1)^{\eta(i)}$.
\hfill$\square_{\scite{dp1.2.2}}$
\enddemo
\bigskip

\demo{\stag{dp1.2.3} Observation}  1) In Definition \scite{dp1.2} \wilog
\, $m (= \ell g(\bar x))$ is 1.
\nl
2) For any $\kappa$ we have: 
$\kappa < \kappa_{\text{ict}}(T)$ \ub{iff} for some
infinite linear order $I_i$ (for $i < \kappa$) and $\langle \bar
a^i_t:t \in I_i,i < \kappa \rangle$ such that $\langle \bar a^i_t:t
\in I_i\rangle$ is indiscernible over $\cup\{\bar a^j_s:s \in I_j$
and $j \ne i,j < \kappa\} \cup A$ and finite $C$, for $\kappa$
ordinals $i < \kappa$, the sequence $\langle \bar a^i_t:t \in
I_i\rangle$ is not indiscernible over $A \cup C$.
\nl
3) In \scite{dp1.2.2}, for any $\lambda (\ge \aleph_0)$ from the statement
   $\circledast^2_{\bar\varphi}$ we get an equivalent one if we
   replace $\omega$ by $\lambda$; similarly for $\circledast^3_{\bar\varphi}$. 
\enddemo
\bigskip

\demo{Proof}  1) For some $m$, 
there is $\bar \varphi = \langle \varphi_i(\bar x,
\bar y_i):i < \kappa \rangle,\ell g(\bar x)=m$ witnessing $\kappa <
\kappa_{\text{ict}}(T)$; \wilog \, $m$ is minimal.  Fixing $\bar\varphi$
 by \scite{dp1.2.2} we know that
$\circledast^2_{\bar \varphi}$ from that observation \scite{dp1.2.2} hold.
Let $\langle \bar a^i_\alpha:i <
\kappa,\alpha < \lambda \rangle$ exemplify $\circledast^2_{\bar
\varphi}$ with $\lambda$ instead $\omega$ 
and let $\bar c = \langle c_i:i < m\rangle$ realize
$\{\varphi_i(\bar x,\bar a^i_0) \wedge \neg \varphi_i(\bar x,\bar
a^i_1):i < \kappa\}$.
\bn
\ub{Case 1}:  For some $u \subseteq \kappa,|u| < \kappa$ for every
$i \in \kappa \backslash u$ the sequence 
$\langle \bar a^i_\alpha:\alpha < \lambda
\rangle$ is an indiscernible sequence over $\cup\{\bar a^j_\beta:j \in
\kappa \backslash u \backslash \{i\}\} \cup \{c_{m-1}\}$.

In this case for $i \in \kappa \backslash u$
let $\psi_i(\bar x',\bar y'_i) := \varphi_i(\bar x
\restriction (m-1),\langle x_{m-1}\rangle \char 94 \bar y_i)$ and
$\bar \psi = \langle \psi_i(\bar x',\bar y'_i):i \in \kappa
\backslash u \rangle$ and $\bar b^i_\alpha = \langle c_{m-1} \rangle
\char 94 \bar a^i_\alpha$ for $\alpha < \lambda,i \in \kappa
\backslash u$ and $\bar\varphi = \langle \psi_i(\bar x',\bar y'_i):i
\in \kappa \backslash u\rangle$.  
Now $\langle \bar b^i_\alpha:\alpha < \lambda,i \in
\kappa \backslash u \rangle$ witness that (abusing our notation)
$\circledast^2_{\bar\psi}$ holds (the 
consistency exemplified by $\bar c \restriction (m-1))$, hence (in the
notation of \scite{dp1.2.2}) $\circledast^1_{\bar\psi^{[\eta]}}$ holds
for some $\eta \in {}^{\kappa \backslash u}2$ 
contradiction to the minimality of $m$.
\bn
\ub{Case 2}:  Not Case 1.

We choose $v_\zeta$ by induction on $\zeta < \kappa$ such that
\mr
\item "{$\bigotimes_\zeta$}"  $(a) \quad v_\zeta \subseteq \kappa
\backslash \cup\{v_\varepsilon:\varepsilon < \zeta\}$
\sn
\item "{${{}}$}"  $(b) \quad v_\zeta$ is finite
\sn
\item "{${{}}$}"  $(c) \quad$ for some $i \in v_\zeta,\langle
\bar a^i_\alpha:\alpha < \lambda \rangle$ is not indiscernible over
\nl

\hskip25pt $\cup\{\bar a^j_\beta:j \in v_\zeta \backslash \{i\},\beta <
\lambda\} \cup \{c_{m-1}\}$
\sn
\item "{${{}}$}"  $(d) \quad$ under $(a)+(b)+(c),|v_\zeta|$ is
minimal.
\ermn
In the induction step, the set $u_\zeta =
\cup\{v_\varepsilon:\varepsilon < \zeta\}$ cannot exemplify case 1, so
for some ordinal $i(\zeta) \in \kappa \backslash u_\zeta$ the sequence
$\langle \bar a^{i(\zeta)}_\alpha:\alpha < \lambda\rangle$ is 
not indiscernible over $\cup\{\bar a^j_\beta:j \in \kappa
\backslash u_\zeta \backslash \{i(\zeta)\}$ and $\beta < \lambda\}
\cup \{c_{m-1}\}$, so by the finite character of indiscernibility, there is a
finite $v \subseteq \kappa \backslash u_\zeta
\backslash \{i(\zeta)\}$ such that $\langle \bar
a^{i(\zeta)}_\alpha:\alpha < \lambda \rangle$ is not indiscernible
over $\cup\{\bar a^j_\beta:j \in v,\beta < \lambda\} \cup \{c_{m-1}\}$.
So $v' = \{i(\zeta)\} \cup v$ satisfies $(a)+(b)+(c)$ hence some
finite $v_\zeta \subseteq \kappa \backslash u_\zeta$ satisfies clauses
$(a),(b),(c)$ and $(d)$. 

Having carried the induction let $i_*(\zeta) \in v_\zeta$ 
exemplify clause (c).  We can find a sequence
$\bar d_\zeta$ from $\cup\{\bar a^j_\beta:j \in v_\zeta \backslash
\{i_*(\zeta)\}$ and $\beta < \lambda\}$ such that $\langle \bar
a^{i_*(\zeta)}_\alpha:\alpha < \lambda \rangle$ is not indiscernible
over $\langle c_{m-1}\rangle \char 94 \bar d_\zeta$.

Also we can find $n(\zeta) < \omega$ and ordinals
$\beta_{\zeta,\ell,0} < \beta_{\zeta,\ell,1} < \ldots <
\beta_{\zeta,\ell,n(\zeta)-1} < \lambda$ for $\ell=1,2$ such that the sequences
$\bar d \char 94 \bar a^{i_*(\zeta)}_{\beta_\zeta,1,0} 
\char 94 \ldots \char 94 \bar a^{i_*(\zeta)}_{\beta_\zeta,1,n(\zeta)-1}$ 
and $\bar d \char 94 \bar a^{i_*(\zeta)}_{\beta_\zeta,2,0} 
\char 94 \ldots \char 94 \bar
a^{i_*(\zeta)}_{\beta_\zeta,2,n(\zeta)-1}$ realize different
types over $c_{m-1}$.  

Now we consider $\bar a^{i_*(\zeta)}_\beta \char 94 \ldots 
\char 94 \bar a^{i_*(\zeta)}_{\beta + n(\zeta)-1}$ where $\beta := \text{
max}\{\beta_{\zeta,1,n(\zeta)-1} + 1,\beta_{\zeta,2,n(\zeta)-1} +1\}$,
so renaming without loss of generality 
$\beta_{\zeta,1,n(\zeta)-1} < \beta_{\zeta,2,0}$.  Omitting some
$a^{i_*(\zeta)}_\beta$'s \wilog \, $\beta_{\beta_\zeta,1,m} =
m,\beta_{\zeta,2,m} = n(\zeta) +m$ for $m < n(\zeta)$.  Now we define
$\bar b^\zeta_\beta := \bar d_\zeta \char 94 \bar
a^{i_*(\zeta)}_{n(\zeta)\beta} \char 94 \ldots \char 94 \bar
a^{i_*(\zeta)}_{n(\zeta)\beta + n(\zeta)-1}$ for $\beta <
\lambda,\zeta < \kappa$.

By the indiscernibility of $\langle \bar a^{i_\zeta(*)}_\gamma:\gamma
< \lambda\rangle$ over $\bar d_\zeta \cup \bigcup\{\bar a^j_\beta:j
\in \kappa \backslash v_\zeta,\beta < \lambda\} \subseteq
\cup\{a^j_\beta:j \in \kappa \backslash \{i_\zeta(*)\},\beta <
\lambda\}$ we can deduce that 
$\langle\bar b^\zeta_\beta:\beta < \lambda\rangle$
is an indiscernible sequence over $\cup\{\bar b^\varepsilon_\beta:
\varepsilon \in \kappa \backslash\{\zeta\}$ and
$\beta < \lambda\}$.  But by an earlier sentence $\bar b^\zeta_0,\bar
b^\zeta_1$ realizes different types over $c_{m-1}$ so we can choose
$\varphi'_\zeta(x,\bar y_\zeta)$ such that ${\frak C} \models
\varphi'_\zeta(c_{m-1},\bar b^\zeta_0) \wedge \neg \varphi'_i
(c_{m-1},\bar b^i_1)$ for $i < \kappa$.

So $\langle \bar b^\zeta_\alpha:\alpha < \omega,\zeta < \kappa\rangle$ and
$\bar \varphi' = \langle \varphi'_\zeta(x,\bar y_\zeta):
\zeta < \kappa\rangle$ satisfy the demands
on $\langle \bar a^i_k:k < \omega,i < \kappa\rangle,\langle
\varphi_i(x,\bar y_i):i < \kappa\rangle$ in $\circledast^2_{\bar\varphi}$
\ub{for $m=1$} (by \scite{dp1.2.2}'s notation), 
so by \scite{dp1.2.2} also $\circledast^1_{\bar\varphi^{[\eta]}}$
 holds for some $\eta \in {}^\kappa 2$ so we are done.
\nl
2) Implicit in the proof of part (1) (and see case 1 in the proof of
\scite{dp1.2.4}).
\nl
3) Trivial.
\nl
${{}}$  \hfill$\square_{\scite{dp1.2.3}}$
\enddemo
\bn
A relative of $\kappa_{\text{ict}}(T)$ is
\definition{\stag{dp1.8} Definition}  1) $\kappa_{\text{icu}}(T) =
\kappa_{\text{icu},1}(T)$ is
the minimal $\kappa$ such that for no $m < \omega$ and $\bar\varphi =
\langle \varphi_i(\bar x_i,\bar y_i):i < \kappa \rangle$ with $\ell
g(\bar x^i) = m \times n_i$ can we find $\bar a^i_\alpha \in 
{}^{\ell g(\bar y_i)}{\frak C}$ for $\alpha < \lambda,i < \kappa$ and $\bar
c_{\eta,n} \in {}^m {\frak C}$ for $\eta \in {}^\kappa \lambda$ such that:
\mr
\item "{$(a)$}"  $\langle \bar c_{\eta,n}:n < \omega \rangle$ is an
indiscernible sequence over $\cup\{\bar a^i_\alpha:\alpha < \lambda,i
< \kappa\}$
\sn
\item "{$(b)$}"  for each $\eta \in {}^\kappa \lambda$ and $i <
\kappa$ we have ${\frak C} \models \varphi_i(\bar c_{\eta,0} \char 94
\ldots \char 94 \bar c_{\eta,n_i-1},\bar
a^i_\alpha)^{\text{if}(\alpha=\eta(i))}$.
\ermn
2) If $\bar \varphi$ is as in (1) then we say that it witnesses
$\kappa < \kappa_{\text{icu}}(T)$.
\nl
3) $T$ is strongly$^{1,*}$ dependent if $\kappa_{\text{icu}}(T) = \aleph_0$.   
\enddefinition
\bigskip

\proclaim{\stag{dp1.9} Claim}  1) $\kappa_{\text{icu}}(T) \ge
\kappa_{\text{ict}}(T)$. 
\nl
2) If {\rm cf}$(\kappa) > \aleph_0$ \ub{then} $\kappa_{\text{icu}}(T) >
\kappa \Leftrightarrow \kappa_{\text{ict}}(T) > \kappa$.
\nl
3) The parallels of \scite{dp1.2.1}, \scite{dp1.2.2},
\scite{dp1.2.3}(2) hold\footnote{and of course more than
\scite{dp1.2.3}(2), using an indiscernible sequence of $m_*$-tuples,
for any $m_* < \omega$.}.
\endproclaim
\bigskip

\demo{Proof}  1) Trivial.
\nl
2) As in the proof of \scite{dp1.2.3}.
\nl
3) Similar.  \hfill$\square_{\scite{dp1.9}}$
\enddemo
\bn
\centerline{$* \qquad * \qquad *$}
\bn
To translate statement on several indiscernible sequences to one (e.g.
in \scite{dp1.2.4}) notes:
\demo{\stag{dq.6} Observation}  Assume that for each $\alpha <
\kappa,I_\alpha$ is an infinite linear order,
the sequence $\langle \bar a_t:t \in I_\alpha
\rangle$ is indiscernible over $A \cup \cup\{\bar a_t:t \in I_\beta$
and $\beta \in \kappa \backslash \{\alpha\}\}$ (and for notational
simplicity $\langle I_\alpha:\alpha < \kappa \rangle$ are pairwise
disjoint) and let $I = \Sigma\{I_\alpha:\alpha < \kappa\},t \in
I_\alpha \Rightarrow \ell g(\bar a_t) = \zeta(\alpha)$ and lastly
for $\alpha \le \kappa$ we let
$\xi(\alpha) = \Sigma\{\zeta(\beta):\beta < \alpha\}$.

\ub{Then} there is $\langle \bar b_t:t \in I \rangle$ such that
\mr
\item "{$(a)$}"   $\ell g(\bar b_t) = \xi(\kappa)$
\sn
\item "{$(b)$}"   $\langle \bar b_t:t \in I \rangle$ is an
indiscernible sequence over $A$
\sn
\item "{$(c)$}"   $t \in I_\alpha \Rightarrow \bar a_t = \bar b_t
\restriction [\xi_\alpha,\xi_\alpha + \zeta_\alpha)$
\sn
\item "{$(d)$}"   if $C \subseteq {\frak C}$ and ${\Cal P}$ is a set
of cuts of $I$ such that [$J$ is a convex subset of $I$ not divided by
any member of ${\Cal P} \Rightarrow \langle \bar b_t:t \in J\rangle$
is indiscernible over $A \cup C$] \ub{then} we can find $\langle {\Cal
P}_\alpha:\alpha < \kappa\rangle,{\Cal P}_\alpha$ is a set of cuts of
$I_\alpha$ such that 
$\Sigma\{|{\Cal P}_\alpha|:\alpha < \kappa\} =
|{\Cal P}|$ and if $\alpha < \kappa,J$ is a convex subset of
$I_\alpha$ not divided by any member of ${\Cal P}_\alpha$ then
$\langle \bar a_t:t \in J\rangle$ is indiscernible over $A \cup C$
\sn
\item "{$(e)$}"  if $C \subseteq {\frak C}$ and ${\Cal P}$ is a set of
cuts of $I$ such that [$J$ is a convex subset of $I$ not divided by any
member of ${\Cal P} \Rightarrow \langle \bar b_t:t \in J\rangle$ is
indiscernible over $A \cup C \cup\{b_s:s \in I \backslash J\}$]
\ub{then} we can find $\langle {\Cal P}_\alpha:\alpha < \kappa
\rangle,{\Cal P}_\alpha$ is a set of cuts of $I_\alpha$ such that
$\Sigma\{|{\Cal P}_\alpha|:\alpha < \kappa\} = |{\Cal P}|$ and if $\alpha <
\kappa,J$ is a convex subset of $I_\alpha$ not divided by any member
of ${\Cal P}_\alpha$ then $\langle \bar a_t:t \in J\rangle$ is
indiscernible over $A \cup C \cup\{\bar a_t:t \in I \backslash J\}$
\sn
\item "{$(f)$}"  moreover in clause (d),(e) we can choose ${\Cal
P}_\alpha$ as the set of non-trivial cuts of $I_\alpha$
induced by ${\Cal P}$, i.e.$\{(J' \cap I_\alpha,
J'' \cap I_\alpha):(J',J'') \in {\Cal P}\} \backslash
\{(I_\alpha,\emptyset),(\emptyset,I_\alpha)\}$.
\endroster
\enddemo
\bigskip

\demo{Proof}  Straightforward.  E.g.,

Without loss of generality $\langle
I_\alpha:\alpha < \kappa \rangle$ are pairwise disjoint and let $I =
\Sigma\{I_\alpha:\alpha < \kappa\}$.  We can find $\bar b^\alpha_t \in
{}^{\zeta(\alpha)} {\frak C}$ for $t \in I,\alpha < \kappa$ such
that: if $n < \omega,\alpha_0 < \ldots < \alpha_{n-1} <
\kappa,t^\ell_0 <_I \ldots <_I t^\ell_{k_\ell-1}$ and $s^\ell_0
<_{I_{\alpha_\ell}} \ldots <_{I_{a_\ell}} s^\ell_{k_\ell-1}$
for $\ell < n$
then the sequence $(\bar b^{\alpha_0}_{t^0_0} \char 94 \ldots
\char 94 \bar b^{\alpha_0}_{t^0_{k_0-1}}) \char 94 \ldots \char 94
(\bar b^{\alpha_{n-1}}_{t^{n-1}_0} \char 94 \ldots \char 94
\bar b^{\alpha_{n-1}}_{t^{n-1}_{k_{n-1}-1}})$ realizes the same type
as the sequence $(\bar a^\alpha_{s^0_0} \char 94 \ldots \char 94 \bar
a^{\alpha_0}_{s^0_{k_n-1}}) \char 94 \ldots \char 94 (\bar
a^{\alpha_{n-1}}_{s^{n-1}_0} \char 94 \ldots \char 94 \bar
a^{\alpha_{n-1}}_{s^{n-1}_{k_{n-1}-1}})$; this is 
possible by compactness.  Using an
automorphism of ${\frak C}$ \wilog \, $t \in I_\alpha \Rightarrow \bar
b^\alpha_t = \bar a^\alpha_t$.  Now for $t \in I$ let $\bar a^*_t$ be $(\bar
a^0_t \char 94 \bar a^1_t \char 94 \ldots \char 94 \bar a^1_\alpha
\ldots)_{\alpha < \kappa}$.
\bn
Clauses (a)+(b)+(c) trivially hold and clauses (d),(e),(f) follows.
  \hfill$\square_{\scite{dq.6}}$ 
\enddemo
\bn
\centerline{$* \qquad * \qquad *$}
\bn
In the following we consider ``natural'' examples which are strongly
dependent, see more in \scite{dq1.8}.
\proclaim{\stag{dw1.10} Claim}  1) Assume $T$ is a complete first order
theory of an ordered abelian group expanded by some individual
constants and some unary predicates $P_i(i < i(*))$ which are
subgroups and $T$ has elimination of quantifiers.
\nl
$T$ is strongly dependent \ub{iff} we cannot find $i_n < i(*)$ and $\iota_n
\in \Bbb Z \backslash \{0\}$ for $n < \omega$ such that:
\mr
\item "{$(*)$}"  we can find $b_{n,\ell} \in {\frak C}$ for $n,\ell <
\omega$ such that
{\roster
\itemitem{ $(a)$ }  $\ell_1 < \ell_2 \Rightarrow \iota_n(b_{n,\ell_2} -
b_{n,\ell_1}) \notin P^{\frak C}_{i_n}$
\sn
\itemitem{ $(b)$ }  for every $\eta \in {}^\omega \omega$ there is
$c_\eta$ such that $c_\eta - b_{n,\eta(n)} \in P^{\frak C}_{i_n}$
for $n < \omega$.
\endroster}
\ermn
2) Let $M$ be $(\Bbb Z,+,-,0,1,<,P_n)$ where $P_n = \{\text{na}:a \in
\Bbb Z\}$ so we know that $T = \text{\rm Th}(M)$ has elimination of
quantifiers.  \ub{Then} $T$ is strongly dependent hence 
{\rm Th}$(\Bbb Z,+,-,0,<)$ is strongly dependent.
\endproclaim
\bigskip

\remark{\stag{dw1.11} Remark}   1) This generalizes the parallel theorem 
for stable abelian groups. \nl
2) Note, if $G$ is the ordered abelian group with sets of elements
$\Bbb Z[x]$, addition of $\Bbb Z[x]$ and $p(x) > 0$ iff the leading
coefficient is $>0$, in $\Bbb Z,P_n$ as above (so definable), then
Th$(G)$ is not strongly dependent using $P_n$ for $n$ prime.
\nl
2) On elimination of quantifiers for ordered abelian groups, see
   Gurevich \cite{Gu77}.
\endremark
\bigskip

\demo{Proof}  1)  The main point is the if direction.
We use the criterion from \scite{dp1.2.4}(2),(4) below.  
So let $\langle \bar a_t:t \in I \rangle$ be an
infinite indiscernible sequence and $c \in {\frak C}$ (with $\bar a_t$
not necessarily finite).  Without loss
of generality ${\frak C} \models ``c>0"$ and 
$\bar a_t = \langle a_{t,\alpha}:\alpha <
\alpha^*\rangle$ list the members of $M_t$, a model and even a
$|T|^+$-saturated model, (see \scite{dp1.2.4}(4)) and let 
$p_t = \text{ tp}(c,M_t)$.

Note that
\mr
\item "{$(*)_1$}"  if $a_{s,i} = a_{t,j}$ and $s \ne t$ then $\langle
a_{r,i}:r \in I\rangle$ is constant.
\ermn
Obviously \wilog \, $c \notin \cup\{M_t:t \in I\}$
but ${\frak C}$ is torsion free (as an abelian group because it is
ordered) hence
\mr
\item "{$(*)_2$}"  $\iota \in \Bbb Z
\backslash \{0\} \Rightarrow \iota c \notin \cup\{M_t:t \in I\}$ 
\sn
\item "{$(*)_3$}"  for $t \in I,a \in M_t$ and $\iota \in 
\Bbb Z \backslash \{0\}$
 let $\eta^\iota_a \in {}^{i(*)+1}2$
be such that $[\eta^\iota_a(i(*))=1 \Leftrightarrow \iota c > a]$ and 
for $i<i(*),[\eta^\iota_a(i) = 1 \Leftrightarrow \iota c - a 
\in P^{\frak C}_i]$
\sn
\item "{$(*)_4$}"   for $t \in I$ and $a \in M_t$ 
let $p_a := \dbcu_{\iota \in \Bbb Z \backslash \{0\}}
(p^\iota_a \cup q^\iota_a$) where
\footnote{recall that $\varphi^1 = \varphi,\varphi^0 = \neg \varphi$}
$p^\iota_a(x) := \{\iota x \ne a,(\iota x>a)^{\eta^\iota_a(i(*))}\}$ and
$q^\iota_a(x) := \{P_i(\iota x-a)^{\eta^\iota_a(i)}:i < i(*)\}$.
\ermn
Now 
\mr
\item "{$\boxdot_0$}"  for $\iota \in \Bbb Z \backslash \{0\}$ and
$\alpha < \alpha^*$ let $I^\iota_\alpha = \{t \in I:a_{t,\alpha} <
\iota c\}$
\sn
\item "{$\boxdot_1$}"  $\langle u_{-1},u_0,u_1\rangle$ is a partition of
$\alpha^*$ where
{\roster
\itemitem{ $(a)$ }   $u_{-1} = \{\alpha < \alpha^*$: for every $s <_I t$ we
have ${\frak C} \models a_{t,\alpha} < a_{s,\alpha}\}$
\sn
\itemitem{ $(b)$ }  $u_0 = \{\alpha < \alpha^*$: for every 
$s <_I t$ we have ${\frak C} \models a_{s,\alpha} =
a_{t,\alpha}\}$
\sn
\itemitem{ $(c)$ }   $u_1 = \{\alpha < \alpha^*$: for every $s <_I t$
we have ${\frak C} \models a_{s,\alpha} < a_{t,\alpha}\}$
\endroster}
\item "{$\boxdot_2$}"  if $\iota \in \Bbb Z \backslash \{0\}$ 
\ub{then} 
{\roster
\itemitem{ $(a)$ }  $I^\iota_\alpha$ is an initial segment of $I$ when
$\alpha \in u_1$
\sn
\itemitem{ $(b)$ }  $I^\iota_\alpha$ is an end segment of $I$ when
$\alpha \in u_{-1}$
\sn
\itemitem{ $(c)$ }  $I^\iota_\alpha \in \{\emptyset,I\}$ when $\alpha
\in u_0$
\sn
\itemitem{ $(d)$ }  $\{I^\iota_\alpha:\alpha \in u_1\}
\backslash \{\emptyset,I\}$ has at most 2 members.
\endroster}
\ermn
[Why?  Recall $<^{\frak C}$ is a linear order.  So for each 
$\iota \in \Bbb Z \backslash \{0\},\alpha \in u_1$ by the definition
of $u_1$ the set 
$I^\iota_\alpha := \{t \in I:a_{t,\alpha} < \iota c\}$ 
is an initial segment of $I$, also $t \in I \backslash I^\iota_\alpha
\Rightarrow \iota c <^{\frak C} a_{t,\alpha}$ as $c \notin \cup\{M_s:s
\in I\}$ by $(*)_2$.  

Now suppose $\alpha,\beta \in u_1$ and $|I^\iota_\beta
\backslash I^\iota_\alpha| > 1$ and $I^\iota_\alpha,I^\iota_\beta
\notin \{\emptyset,I\}$ then choose $t_1 <_I t_2$ from $I^\iota_\beta
\backslash I^\iota_\alpha$ and $t_0 \in I^\iota_\alpha,t_3 \in I \backslash
I^\iota_\beta$.  As $I^\iota_\alpha,I^\iota_\beta$ are initial
segments and $t_0 <_I t_1 <_I t_2 <_I t_3$, necessarily
${\frak C} \models ``a_{t_0,\alpha} < \iota c < a_{t_1,\alpha}
\wedge a_{t_2,\beta} < \iota c < a_{t_3,\beta}"$.  If
$a_{t_1,\alpha} \le^{\frak C} a_{t_2,\beta}$ we can deduce a
contradiction $({\frak C} \models ``\iota c < a_{t_1,\alpha} \le a_{t_2,\beta}
< \iota c"$). Otherwise by the indiscernibility of the sequence $\langle
(a_{t,\alpha},a_{t,\beta}):t \in I\rangle$ we get ${\frak C} \models
a_{t_3,\beta} < a_{t_0,\alpha}$ and a similar contradiction.  So 
$|I^\iota_\beta \backslash I^\iota_\alpha| \le 1$.  

So $I^\iota_\alpha,I^\iota_\beta \notin \{\emptyset,I\} \Rightarrow
|I^\iota_\beta \backslash I^\iota_\alpha| \le 1$ and by 
symmetry $|I^\iota_\alpha \backslash I^\iota_\beta| \le 1$. 
So $|\{I^\iota_\alpha:\alpha \in u_1\} \backslash \{\emptyset,I\}| 
\le 2$, i.e. clause (d) of $\boxdot_2$ holds; the other clauses should
be clear.]  

Now clearly
\mr
\item "{$\boxdot_3$}"  if $\alpha,\beta < \alpha(*),\iota \in \Bbb Z
\backslash \{0\}$ and $a_{t,\alpha} = -a_{t,\beta}$ 
(for some equivalently for every $t \in I$) \ub{then}:
{\roster
\itemitem{ $(a)$ }  $(\alpha \in u_1) \equiv (\beta \in u_{-1})$
\sn
\itemitem{ $(b)$ }  $((\iota c) < a_{t,\alpha}) \equiv (a_{t,\beta} <
((-\iota)c))$ recalling $\iota c,(-\iota) c \notin \dbcu_{t \in I}
M_t$
\sn
\itemitem{ $(c)$ }  $I^\iota_\alpha = I \backslash I^\iota_\beta$.
\endroster}
\ermn
Also
\mr
\item "{$\boxdot_4$}"   if $\iota_1,\iota_2$ are from
$\{1,2,\ldots\}$ and $\iota_1 a_{t,\alpha} = \iota_2 a_{t,\beta}$ \ub{then}
{\roster
\itemitem{ $(a)$ }  $[\alpha \in u_{-1} \equiv \beta \in
u_{-1}],[\alpha \in u_0 \equiv \beta \in u_0]$ and $[\alpha \in u_1
\equiv \beta \in u_1]$
\sn
\itemitem{ $(b)$ }  $(t \in I^{\iota_2}_\alpha) \Leftrightarrow 
(t \in I^{\iota_1}_\beta)$ hence $I^{\iota_2}_\alpha = I^{\iota_1}_\beta$.
\endroster}
\ermn
[Why?  Clause (a) is obvious.  For clause (b) note that $t \in
I^{\iota_2}_\alpha \Leftrightarrow a_{t,\alpha} < \iota_2 c
\Leftrightarrow \iota_1 a_{t,\alpha} < \iota_1(\iota_2 c)
\Leftrightarrow \iota_2 a_{t,\beta} < \iota_2(\iota_1 c)
\Leftrightarrow a_{t,\beta} < \iota_1 c \Leftrightarrow t \in
I^{\iota_1}_\beta$.] 

By symmetry, i.e. by $\boxdot_3$ clearly
\mr
\item "{$\boxdot_5$}"  the statement (c),(d) in $\boxdot_2$
holds for $\alpha \in u_{-1}$.
\ermn
Obviously
\mr
\item "{$\boxdot_6$}"  if $\alpha \in u_0$ then $I^\iota_\alpha \in
\{\emptyset,I\}$.
\ermn
Together
\mr
\item "{$\boxdot_7$}"  $\{I^\iota_\alpha:\alpha < \alpha^*$ 
and $\iota \in \Bbb Z \backslash \{0\}\} \backslash \{\emptyset,I\}$ 
hence has $\le 4$ members.
\ermn
Hence
\mr
\item "{$\circledast_0$}"   there are initial segments $J_\ell$ of $I$
for $\ell < \ell(*) \le 4$ such that: if $s,t$ belongs to $I$ and
$\ell < \ell(*) \Rightarrow [s \in J_\ell \equiv t \in J_\ell]$ 
\ub{then} 
$\eta^\iota_{a_{t,\alpha}}(i(*)) = \eta^\iota_{a_{s,\alpha}}(i(*))$.
\ermn
[Why?  By the above and the 
definition of $\eta^\iota_{a_{t,\alpha}}(i(*))$ we are done.]
\mr
\item "{$\circledast_1$}"  for each $t \in I$ we have 
$\cup\{p_a(x):a \in M_t\} \vdash p_t(x)$.
\nl
[Why?  Use the elimination of quantifiers and the closure properties
of $M_t$.  That is, every formula in $p_t(x)$ is equivalent to a
Boolean combination of quantifier free formulas.  So it suffices to
deal with the cases $\varphi(x,\bar a) \in p_t(x)$ which is atomic or
negation of atomic and $x$ appear.  
As for $b_1,b_2 \in {\frak C}$ exactly one of the possibilities
$b_1 < b_2,b_1 = b_2,b_2 < b_1$ holds and by symmetry, it suffices to deal with
$\sigma_1(x,\bar a) > \sigma_2(x,\bar a),\sigma_1(x,\bar a) =
\sigma_2(x,\bar a),P_i(\sigma(x,\bar a)),\neg P_i(\sigma(x,\bar a))$
where $\sigma(x,\bar y),\sigma_1(x,\bar y),\sigma_2(x,\bar y)$ are 
terms in $\Bbb L(\tau_T)$.
As we can substract, it suffices to deal with $\sigma(x,\bar a)
> 0,\sigma(x,\bar a) = 0,P_i(\sigma(x,\bar a)),\neg P_i(\sigma(x,\bar
a))$.  By linear algebra as $M_t$ is closed under the operations,
without loss of generality  
$\sigma(x,\bar a) = \iota x - a_{t,\alpha}$ for some 
$\iota \in \Bbb Z$ and $\alpha < \alpha^*$, and \wilog \, 
$\iota \ne 0$.  The case $\varphi(x) = (\iota x - a_{t,\alpha} = 0) 
\in p(x)$ implies $c \in M_t$ (as $M$ is torsion free) 
which we assume does
not hold.  In the case $\varphi(x,\bar a) = (\iota x - a_{t,\alpha} > 0)$ use
$p^\iota_{a_{t,\alpha}}(x)$, in the case $\varphi(x,\bar a) = P_i(\iota x -
a_{t,\alpha})$ or $\varphi(x,\bar a) = \neg P_i(\iota x - a_{t,\alpha})$ 
use $q^\iota_{a_{t,\alpha}}(x)$ for
$\eta^\iota_{a_{t,\alpha}}(i)$.]
\sn
\item "{$\circledast_2$}"  if $\iota \in \Bbb Z \backslash \{0\},
n < \omega$ and $a_0,\dotsc,a_{n-1} \in
M_t$ then for some $a \in M_t$ we have $\ell < n \wedge i < i(*) \wedge
\eta^\iota_{a_\ell}(i) =1 \Rightarrow \eta^\iota_a(i)=1$ 
\nl
[Why?  Let $a' \in M_t$ realize $p_t \restriction
\{a_0,\dotsc,a_{n-1}\}$, exist as $M_t$ was chosen as 
$|T|^+$-saturated; less is necessary.
Now $\iota c - a_\ell \in P^{\frak C}_i \Rightarrow \iota a' - a_\ell 
\in P^{\frak C}_i \Rightarrow (\iota c - \iota a') = 
((\iota c - a_\ell) - (\iota a' -
a_\ell)) \in P^{\frak C}_i$ and let $a := \iota a'$.]
\sn
\item "{$\circledast_3$}"  assume $\iota \in \Bbb Z \backslash \{0\},
i < i(*),\alpha < \alpha^*,s_1 <_I s_2$ 
and $t \in I \backslash \{s_1,s_2\}$ then:
{\roster
\itemitem{ $(a)$ }  if $\eta^\iota_{a_{s_1,\alpha}}(i) = 1$ and
$\eta^\iota_{a_{s_2,\alpha}}(i)=0$ then $\eta^\iota_{a_{t,\alpha}}(i)=0$
\sn
\itemitem{ $(b)$ }  if $\eta^\iota_{a_{s_1,\alpha}}(i)=0$ and
$\eta^\iota_{a_{s_2,\alpha}}(i)=1$ then $\eta^\iota_{a_{t,\alpha}}(i)=0$.
\endroster}
[Why?  As we can invert the order of $I$ it is enough to prove clause
(a).  By the choice of $a \mapsto \eta^\iota_a$ we have $\iota c-a_{s_1,\alpha}
\in P^{\frak C}_i,\iota c-a_{s_2,\alpha} \notin P^{\frak C}_i$
hence $a_{s_1,\alpha} - a_{s_2,\alpha} \notin P^{\frak C}_i$ hence also
$a_{s_2,\alpha} - a_{s_1,\alpha} \notin P^{\frak C}_i$.
\nl
By the indiscernibility we have $a_{t,\alpha} -
a_{s_1,\alpha} \notin P^{\frak C}_i$ and as $\iota c-a_{s_1,\alpha} \in
P^{\frak C}_i$ we can deduce $\iota c-a_{t,\alpha} \notin P^{\frak C}_i$
hence $\eta^\iota_{a_{t,\alpha}}(i)=0$.  So we are done.]
\sn
\item "{$\circledast_4$}"  for each $\iota \in \Bbb Z \backslash \{0\},
i<i(*)$ and $\alpha < \alpha^*$ the set $I^\iota_{i,\alpha} := 
\{t:\eta^\iota_{a_{t,\alpha}}(i)=1\}$ is $\emptyset,I$ or a
singleton.
\nl
[Why?  By $\circledast_3$.]
\sn
\item "{$\circledast_5$}"  if $I_* = \cup\{I^\iota_{i,\alpha}:\iota \in \Bbb Z
\backslash \{0\},i < i(*),\alpha < \alpha^*$ and $I^\iota_{i,\alpha}$
is a singleton$\}$ is infinite then (possibly inverting $I$)
we can find $t_n \in I$
and $\beta_n < \alpha^*,\iota_n \in \Bbb Z \backslash \{0\}$ and 
$i_n < i(*)$ for $n < \omega$ such that
{\roster
\itemitem{ $(a)$ }   $t \in I$ then $[\iota_n c-a_{t,\beta_n} \in 
P^{\frak C}_{i_n}] \Leftrightarrow t = t_n$ for every $n < \omega$
\sn
\itemitem{ $(b)$ }  $\langle a_{t,\beta_n} - a_{s,\beta_n}:s \ne t \in
I \rangle$ are pairwise not equal mod $P^{\frak C}_{i_n}$
\sn
\itemitem{ $(c)$ }  $t_n < t_{n+1}$ for $n < \omega$.
\endroster}
\ermn
[Why?  Should be clear.]
\mr
\item "{$\circledast_6$}"  if $I_* = \cup\{I^\iota_{i,\alpha}:\iota \in \Bbb Z
\backslash \{0\},\alpha < \alpha^*,i < i(*)$ and $I^\iota_{i,\alpha}$ is
a singleton$\}$ is finite and $J_\ell (\ell < \ell(*) \le 6)$ are as in
$\circledast_0$, \ub{then} 
tp$(\bar a_s,\{c\}) = \text{\rm tp}(\bar a_t,\{c\})$ 
whenever $(s,t \in I \backslash I_*)
\wedge \dsize \bigwedge_{\ell < \ell(*)}(s \in J_\ell \equiv t \in J_\ell)$
recalling $\bar a_t$ list the elements of $M_t$.
\ermn
[Why?  By $\circledast_4$ and $\circledast_1$ (and $\circledast_0$)
recalling the choice of $p_a$ in $(*)_4$.]

Assume $c,\langle \bar a_t:t \in I\rangle$ exemplify $T$ is not
strongly dependent then $I_*$ cannot be finite (by $\circledast_6$)
hence $I_*$ is infinite so by $\circledast_5$ 
we can find $\langle(t_n,\beta_n,\iota_n,i_n):n < \omega\rangle$ as
there.

That is, for $n < \omega,\ell < \omega$ let $b_{n,\ell} :=
a_{t_\ell,\beta_n}$.  So
\mr
\item "{$\circledast_7$}"  $\iota_n c - b_{n,\ell} \in P^{\frak
C}_{i_n}$ \ub{iff} $\iota_n c - a_{t_\ell,\beta_n} \in P^{\frak
C}_{i_n}$ iff $t_\ell = t_n$ iff $\ell=n$
\sn
\item "{$\circledast_8$}"  if $\ell_1 < \ell_2$ then $b_{n,\ell_2} -
b_{n,\ell_2} \notin P^{\frak C}_{i_n}$.
\ermn
[Why?  By clause (b) of $\circledast_5$.]

Now
\mr
\item "{$\circledast_9$}"  if $\eta \in {}^\omega \omega$ is
increasing then there $c_\eta \in {\frak C}$ such that
\nl
$n < \omega \Rightarrow \iota_n c_\eta - b_{n,\eta(n)} \in P^{\frak
C}_{i_n}$.
\ermn
[Why?  As $\langle \bar a_t:t \in I\rangle$ is an indiscernible
sequence, there is an automorphism $f = f_\eta$ of ${\frak C}$ which
maps $\bar a_{t_n}$ to $\bar a_{t_{\eta(n)}}$ for $t \in I$ so
$f_\eta(b_{\eta,n}) = b_{n,\eta(n)}$.  Hence $c_\eta = f_\eta(c)$
satisfies $n < \omega \Rightarrow \iota_n f(c) - b_{\eta,\eta(n)}  \in
P^{\frak C}_{i_n}$.]

Now $\langle b_{n,\ell}:n,\ell < \omega\rangle$ almost satisfies
$(*)$ of \scite{dw1.10}.  Clause (a) holds by $\circledast_8$ and
clause (b) holds for all increasing $\eta \in {}^\omega \omega$.  By
compactness we can find $\langle \bar b'_{n,\ell}:n,\ell <
\omega\rangle$ satisfying (a) + (b) of $(*)$ of \scite{dw1.10}.
\nl
[Why?  Let $\Gamma = \{P_{i_n}(\iota_n x_\eta - y_{n,\eta(n)}):\eta \in
{}^\omega \omega,n < \omega\} \cup \{\neg P_{i_n}(\iota_n x_{n,\ell_1}
- \iota_n x_{n,\ell_2}):n < \omega,\ell_1 < \ell_2 < \omega\}$.  If
$\Gamma$ is satisfied in ${\frak C}$ we are done, otherwise there is a
finite inconsistent $\Gamma' \subseteq \Gamma$, let $n_*$ be such that:
if $y_{n,\ell}$ appear in $\Gamma'$ then $n,\ell < n_*$.  But the
assignment $y_{n,\ell} \mapsto b_{n n_* + \ell}$ for $n < n_*,\ell <
n_*$ exemplified that $\Gamma'$ is realized, so we have proved half of
the claim.  The other direction should be clear, too.]
\nl
2) The first assertion (on $T$) holds by part (1); the second holds as
the set of terms $\{0,1,2,\dotsc,n-1\}$ is provably a set of
representatives for $\Bbb Z/P_n$ which is finite.
\hfill$\square_{\scite{dw1.10}}$ 
\enddemo
\bn
\margintag{dw1.12}\ub{\stag{dw1.12} Example}:   Th$(M)$ is not strongly stable when $M$ satisfies:
\mr
\item "{$(a)$}"  has universe ${}^\omega \Bbb Q$
\sn
\item "{$(b)$}"  is an abelian group as a power of $(\Bbb Q,+)$,
\sn
\item "{$(c)$}"  $P^M_n = \{f \in M:f(n)=0\}$, a subgroup.
\endroster
\bn
We now consider the $p$-adic fields and more generally valued fields.
\definition{\stag{dt.7} Definition}  1) We define a valued field $M$ as
one in the Denef-Pas language, i.e., a model $M$ such that:
\mr
\item "{$(a)$}"  the elements of $M$ are of three sorts:
{\roster
\itemitem{ $(\alpha)$ }  the field $P^M_0$ which (as usual) we call
$K^M$, so $K = K^M$ is the field of $M$ and has universe $P^M_0$ so we
have appropriate individual constants (for $0,1$), and the field operations
(including the inverse which is partial)
\sn
\itemitem{ $(\beta)$ }  the residue field $P^M_1$ which (as usual) is
called $k^M$, so $k=k^M$ is a field with universe $P^M_1$ so with the
appropriate $0,1$ and field operations
\sn
\itemitem{ $(\gamma)$ }  the valuation ordered abelian group $P^M_2$
which (as usual) we call $\Gamma^M$, so $\Gamma = \Gamma^M$ is 
an ordered abelian group with universe $P^M_2$ so with $0$, addition,
subtraction and the order
\endroster}
\item "{$(b)$}"  the functions (and individual constants)
of $K^M,k^M,\Gamma^M$ and the order of $\Gamma^N$ (actually mentioned
in clause (a))
\sn
\item "{$(c)$}"  val$^M:K^M \rightarrow \Gamma^M$, the valuation
\sn
\item "{$(d)$}"  ac$^M:K^M \rightarrow k^M$, the function giving the
``leading coefficient" (when as in natural cases the members of $K$
are power series)
\sn
\item "{$(e)$}"  of course, the sentences saying that the following hold:
{\roster
\itemitem{ $(\alpha)$ }  $\Gamma^M$ is an ordered abelian group
\sn
\itemitem{ $(\beta)$ }  $k$ is a field
\sn
\itemitem{ $(\gamma)$ }  $K$ is a field
\sn
\itemitem{ $(\delta)$ }  val,ac satisfies the natural demands.
\endroster}
\ermn
1A) Above we replace ``language" by $\omega$-language \ub{when}: in clause
$(b)$, i.e. $(a)(\gamma)$, $\Gamma^M$ has $1_\Gamma$ 
(the minimal positive elements) and we replace (d) by
\mr
\item "{$(d)^-_\omega$}"  ac$^M_n:K^M \rightarrow k^M$ satisfies: 
$\dsize \bigwedge_{\ell < n} \text{ ac}^M_\ell(x) = \text{ ac}^M_k(y)
\Rightarrow \text{ val}^M(x-y) > \text{ val}^m(x) + n$.
\ermn
2) We say that such $M$ (or Th$(M)$) has elimination of the field
quantifier \ub{when}: every first order formula (in the language of
Th$(M)$) is equivalent to a Boolean combination of atomic formulas,
formulas about $k^M$ (i.e., all variable, free and bounded vary on
$P^M_1$) and formulas about $\Gamma^M$; note this definition requires
clause (d) in part (1).
\enddefinition
\bn
It is well known that (on \scite{dt.8},\scite{dt.9} see, e.g. \cite{Pa90}, 
\cite{CLR06}).
\proclaim{\stag{dt.8} Claim}  1) Assume $\Gamma$ is a divisible ordered
abelian group and $k$ is a perfect field of characteristic zero.  
Let $K$ be the fixed power series for $(\Gamma,k)$,
i.e. $\{f:f \in {}^\Gamma k$ and {\rm supp}$(f)$ is well ordered$\}$ where
{\rm supp}$(f) = \{s \in \Gamma:f(s) \ne 0_k\}$.  \ub{Then} the 
model defined by $(K,\Gamma,k)$ has elimination of the field quantifiers.
\nl
2) For $p$ prime, we can consider the
$p$-adic field as a valued field in the Denef-Pas $\omega$-language 
and its first order theory has elmination of the field quantifiers 
(this version of the $p$-adics and the original one are (first-order)
bi-interpretable; note that the field $k$ here is finite
and formulas speaking on $\Gamma$ which is the
ordered abelian group $\Bbb Z$ are well understood). 
\endproclaim
\bn
We will actually be interested only in valuation fields $M$ with
elimination of the field quantifiers.  It is well known that
\proclaim{\stag{dt.9} Claim}  Assume ${\frak C} = {\frak C}_T$ is a
(monster,i.e. quite saturated) valued field in the Denef-Pas
language (or in the $\omega$-language) with elimination of
the field quantifiers.  If $M \prec {\frak C}$ \ub{then}
\mr
\item "{$(a)$}"  it satisfies the cellular decomposition of Denef
which implies
\footnote{note: $p \in \bold S^1(A,M),A \subseteq M$ is a little more
complicated}:
\nl
if $p \in \bold S^1(M)$ and $P_0(x) \in p$ \ub{then} $p$ is equivalent to
\nl
$p^{[*]} := \cup\{p^{[*]}_c:c \in P^M_0\}$ where $p^{[*]}_c =
p^{[*,1]}_c \cup p^{[*,2]}_c$ and
\nl
$p^{[*,1]}_c = \{\varphi(\text{\rm val}(x-c),\bar d) \in p:
\varphi(x,\bar y)$ is a formula speaking on
$\Gamma^M$ only so $\bar d \subseteq \Gamma^M,c \in P^M_0\}$ and 
\nl
$p^{[*,2]}_c = \{\varphi(\text{\rm ac}(x-c),\bar d) \in p:
\varphi$ speaks on $k^M$ only$\}$ but for the $\omega$-language we
should allow 
$\varphi(\text{\rm ac}_0(x-c),\dotsc,\text{\rm ac}_n(x-c),\bar d)$ for
some $n < \omega$
\sn
\item "{$(b)$}"  if $p \in \bold S^1(M),P_0(x) \in p$ and 
$c_1,c_2 \in P^M_0$ and
{\rm val}$^M(x-c_1) <^{\Gamma^M} \text{\rm val}^M(x-c_2)$ belongs to $p(x)$
\ub{then} $p^{[*]}_{c_2}(x) \vdash p^{[*]}_{c_1}(x)$ and even $\{\text{\rm
val}(x-c_1) < \text{\rm val}(x-c_2)\} \vdash p^{[*]}_{c_1}(x)$
\sn
\item "{$(c)$}"  for $\bar c \in {}^{\omega >}(k^M)$, the type
\text{\rm tp}$(\bar c,\emptyset,k^M)$ determines 
{\rm tp}$(\bar c,\emptyset,M)$ and similarly for $\Gamma^M$.
\endroster
\endproclaim
\bigskip

\proclaim{\stag{dp0.16} Claim}  1) The first order theory $T$ of the
$p$-adic field is strongly dependent.
\nl
2) For any theory $T$ of a valued field $\Bbb F$ which has elimination
of the field quantifier we have:
\nl
$T$ is strongly dependent iff the theory of the valued ordered group 
and the theory of the residue fields of $\Bbb F$ are strongly
dependent.
\nl
3) Like (2) when we use the $\omega$-language and we assume $k^M$ is finite.
\endproclaim
\bigskip

\remark{\stag{dp0.16.19} Remark}  1) In \scite{dp0.16} we really get that $T$ is
strongly dependent over the residue field + the valuation ordered 
abelian group.
\nl
2)  We had asked in a preliminary version of [Sh:783,\S3]:
show that the theory of the $p$-adic field is strongly dependent.
Udi Hrushovski has noted that the criterion (St)$_2$
presented there (and repeated in \scite{0.gr.1} here from
\cite[3.10=ss.6]{Sh:783}) apply so $T$ is not 
strongly$^2$ dependent.   Namely take the
following equivalence relation $E$ on $\Bbb Z_p$:val$(x-y) \ge
\text{\rm val}(c)$, where $c$ is some fixed element with infinite valuation.
Given $x$, the map $y \mapsto (x + cy)$ is a bijection between $\Bbb
Z_p$ and the class $x/E$.
\nl
3) By \cite[\S3]{Sh:783}, the theory of real closed fields,
i.e. Th$(\Bbb R)$ is strongly dependent.
Onshuus shows that also the theory of the field of the reals is not
   strongly$^2$ dependent (e.g. though Claim \cite[3.10=ss.6]{Sh:783}
does not apply but its proof works using pairwise not too near $\bar b$'s,
in general just an uncountable set of $\bar b$'s). 
\nl
4) See more in \S5.
\endremark
\bn
Of course,
\demo{\stag{dp0.16.23} Observation}  1) For 
a field $K$, Th$(K)$ being strongly dependent is
preserved by finite extensions in the field theoretic sense by
\scite{dp1.2.1}(2).
\nl
2) In \scite{dp0.16}, if 
we use the $\omega$-language and $k^N$ is infinite, the theory
is not strongly dependent.
\enddemo
\bigskip

\demo{Proof}  1) Recall that by \scite{dw1.10}(2), the theory
of the valued group (which is an ordered abelian group) 
is strongly dependent, and this trivially holds for
the residue field being finite.  So by \scite{dt.8}(2) we can apply part
(3).
\nl
2) We consider the models of $T$ as having three sorts:
$P^M_0$ the field, $P^M_1$ the ordered abelian group (like value of
valuations) and $P^M_2$ the residue field.
  
Let
\mr
\item "{$\boxdot_1$}"  $(a) \quad I$ be an infinite linear order, \wilog \,
complete and 
\nl

\hskip25pt dense (and with no extremal members),
\sn
\item "{${{}}$}"  $(b) \quad \langle \bar a_t:t \in I 
\rangle$ an indiscernible sequence, $\bar a_t \in 
{}^\alpha{\frak C}$  and let $c \in {\frak C}$
\nl

\hskip25pt  (a singleton!)
\ermn
and we shall prove
\mr
\item "{$\boxdot_2$}"  for some finite $J \subseteq I$ we have: if $s,t \in I
\backslash J$ and $(\forall s \in J)(r <_I s \equiv r <_I t)$ then
$\bar a_s,\bar a_t$ realizes the same type over $\{c\}$.
\ermn
This suffices by \scite{dp1.2.4} and as there by \scite{dp1.2.4}(4) \wilog 
\mr
\item "{$\boxdot_3$}"  $\bar a_t =
\langle a_{t,i}:i < \alpha \rangle$ list the elements of an elementary
submodel $M_t$ of ${\frak C} = {\frak C}_T$ (we may assume $M_t$ is
$\aleph_1$-saturated; alternatively we could have assumed that it is quite
complete).
\ermn
Easily it follows that it suffices to prove (by the
L.S.T. argument but not used)
\mr
\item "{$\boxdot'_2$}"  for every countable $u \subseteq \alpha$ there
is a finite $J \subseteq I$ which is O.K. for $\langle \bar a_t
\restriction u:t \in I\rangle$.
\ermn
Let $\bold f_{t,s}$ be the mapping $a_{s,i} \mapsto a_{t,i}$ for $i <
\alpha$; clearly it is an isomorphism from $M_s$ onto $M_s$.

Now
\mr
\item "{$\boxdot_4$}"  $p_t = \text{\rm tp}(c,M_t)$
so $(p_t)^{[*]}_a$ for $a \in M_t$ is well defined in \scite{dt.9}(a). 
\ermn
The case $P_2(x) \in \dbca_t p_t$ is easy and the case $P_1(x) \in \dbca_t
p_t$ is easy, too, by an assumption (and clause (c) of \scite{dt.9}), so 
we can assume $P_0(x) \in \dbca_t p_t(x)$.

Let ${\Cal U} = \{i < \alpha:a_{s,i} \in P^{\frak C}_0$ for every
($\equiv$ some) $s \in I\}$.

Now for every $i \in {\Cal U}$
\mr
\item "{$(*)^1_i$}"  the function $(s,t) \mapsto
\text{\rm val}^{\frak C}(a_{t,i} - a_{s,i})$ for $s <_I t$ 
satisfies one of the following
\sn
\ub{Case $(a)^1_i$}:  it is constant
\sn
\ub{Case $(b)^1_i$}: it depends just on $s$ and is a strictly monotonic
(increasing, by $<_\Gamma$) function of $s$
\sn 
\ub{Case $(c)^1_i$}:  it depends just on $t$ and is a strictly monotonic
(decreasing, by $<_\Gamma$) function of $t$.
\ermn
[Why?  This follows by inspection.]
\nl
For $\ell = -1,0,1$ let ${\Cal U}_\ell := \{i \in {\Cal U}$: if
$\ell = 0,1,-1$ then case $(a)^1_i,(b)^1_i,(c)^1_i$ 
respectively of $(*)^1_i$
holds$\}$ so $\langle {\Cal U}_{-1},{\Cal U}_0,{\Cal U}_1\rangle$ is a
partition of ${\Cal U}$.

For $i,j \in {\Cal U}_1$ we shall prove more than 
\mr
\item "{$(*)^2_{i,j}$}"   we have $i,j \in {\Cal U}_1$ and
the function $(s,t) \mapsto
\text{\rm val}^{\frak C}(a_{t,j} - a_{s,i})$ for $s <_I t$
satisfies one of the following:
\sn
\ub{Case $(a)^2_{i,j}$}:  val$^{\frak C}(a_{t,j} - a_{s,i})$ is constant
\sn
\ub{Case $(b)^2_{i,j}$}: val$^{\frak C}(a_{t,j} - a_{s,i})$ depends only on $s$
and is a monotonic (increasing) function of $s$ and is equal to
val$^{\frak C}(a_{s_1,i} - a_{s,i})$ when $s <_I s_1$
\sn 
\ub{Case $(c)^2_{i,j}$}:  val$^{\frak C}(a_{t,j} - a_{s,i})$ 
depends only on $t$
and is a monotonic (increasing) function of $t$ and is equal to
val$^{\frak C}(a_{t,j} - a_{t_1,j})$ when $t <_I t_1$.
\ermn
[Why $(*)^2_{i,j}$ holds?  In this case we give full checking.

First, assume: for some (equivalently every) $t \in I$ the sequence
$\langle \text{val}^{\frak C}(a_{t,j} - a_{s,i}):s$ satisfies $s <_I
t\rangle$ is $<_\Gamma$-decreasing with $s$ recalling that we have
assumed $I$ is a linear order with neither first nor last element.  
Choose $s_1 <_I s_2 <_I t$ so by the present assumption 
we have val$^{\frak C}(a_{t,j} - a_{s_2,i}) 
<_\Gamma \text{ val}^{\frak C}(a_{t,j} -
a_{s_1,i})$ hence val$^{\frak C}((a_{t,j} - a_{s_2,i}) - (a_{t,j} -
a_{s_1,i})) = \text{ val}^{\frak C}(a_{t,j} - a_{s_2,i})$ which means
val$^{\frak C}(a_{t,j} - a_{s_2,i}) = \text{ val}^{\frak C}
(-(a_{s_2,i} - a_{s_1,i})) = \text{ val}^{\frak C}(a_{s_2,i} -
a_{s_1,i})$.  So in the right side $t$ does not appear, in the left
side $s_1$ does not appear, hence by the equality the left side,
$\text{val}^{\frak C}(a_{t,j} - a_{s_2,i})$, does not 
depend on $t$ and the right side, val$^{\frak C}(a_{s_2,i} - a_{s_1,i})$
 does not depend on $s_2$ but as $i \in{\Cal U}_1$ it does not depend
on $s_1$.  Together by the indiscernibility for $s <_I t$ we have 
val$^{\frak C}(a_{t,i} - a_{s,i})$ is
constant, i.e. case $(a)^2_{i,j}$ holds.  So we can from now on assume:
for each $t \in I$ the sequence 
$\langle\text{val}^{\frak C}(a_{t,j}-a_{s,i}):s$ satisfies
$s <_I t\rangle$ is constant or for each $t \in I$ it is
$<_\Gamma$-increasing with $s$.

Second, assume: for some (equivalently every) $s \in I$ the sequence
$\langle\text{val}^{\frak C}(a_{t,j}-a_{s,i}):t$ satisfies $s <_I
t\rangle$ is $<_\Gamma$-decreasing with $t$.  As in ``first" we can show
that case $(a)^2_{i,j}$ holds.  So from now on we can assume that for
every $s \in I$ the sequence 
$\langle\text{val}^{\frak C}(a_{t,j}-a_{s,i}):t$ satisfies $s <_I
t\rangle$ is constant or for every $s \in I$ the sequence is
$<_\Gamma$-increasing with $s$.

Third, assume: for some (equivalently every) $t \in I$ the sequence
$\langle\text{val}^{\frak C}(a_{t,j}-a_{s,i}):s$ satisfies $s <_I
t\rangle$ is constant.  This implies that $s <_I t \Rightarrow
\text{val}^{\frak C}(a_{t,j}-a_{s,i}) = e_t$ for some $\bar e = \langle
e_t:t \in I\rangle$.  If for some (equivalently every) $s \in I$ the sequence
$\langle\text{val}^{\frak C}(a_{t,j}-a_{s,i}):t$ satisfies $s <_I
t\rangle$ is constant then clearly case $(a)^2_{i,j}$ holds so we can
assume this fails so by the end of ``second" this sequence is
$<_\Gamma$-increasing hence $\langle e_t:t \in I\rangle$ is
$<_\Gamma$-increasing.  So most of the requirements in case $(c)^2_{i,j}$
holds; still we have to show that $t <_I t_1 \Rightarrow
\text{val}(a_{t,j}-a_{t_1,j}) = e_t$.

Let $s <_I t <_I t_1$, we know that $e_{t} <_\Gamma e_{t_1}$, which means
that $\text{val}^{\frak C}(a_{t,j}-a_{s,i}) <_\Gamma
\text{val}^{\frak C}(a_{t_1,j}-a_{s,i})$.  This implies that
$\text{val}^{\frak C}((a_{t,j}-a_{s,i}) - (a_{t_1,j}-a_{s,i})) = 
\text{ val}^{\frak C}(a_{t,j}-a_{s,i})$ which means that
$\text{val}^{\frak C}(a_{t,j}-a_{t_1,j}) = \text{val}^{\frak C}(a_{t,j}
- a_{s,i}) = e_t$ as required; so case $(c)^2_{i,j}$ and we are done
(if ``Third..." holds).

Fourth, assume that for some (equivalently every) $s \in I$ the sequence
$\langle\text{val}^{\frak C}(a_{t,j}-a_{s,i}):t$ satisfies $s <_I
t\rangle$ is constant, then we proceed as in ``third" getting
case $(b)^2_{i,j}$ instead of case $(c)^2_{i,j}$.

So assume that none of the above occurs, hence for every (equivalently
some) $t \in I$ the sequence $\langle\text{val}^{\frak
C}(a_{t,j}-a_{s,i}):s$ satisfies $s <_I t\rangle$ is
$<_\Gamma$-increasing (with $s$, 
by ``first"...and ``third" above) and for every
(equivalently some) $s \in I$ the sequence 
$\langle\text{val}^{\frak C}(a_{t,j}-a_{s,i}):t$ satisfies $s <_I
t\rangle$ is $<_\Gamma$-increasing (with $t$, 
by ``second" and ``fourth" above).

Hence we have $s <_I t_1 <_I t_2 \Rightarrow 
\text{val}^{\frak C}(a_{t_1,j}-a_{s,i}) <_\Gamma \text{ val}^{\frak
C}(a_{t_2,j}-a_{s,i}) \Rightarrow \text{ val}^{\frak C}(a_{t_1,j}
-a_{s,i}) = \text{val}^{\frak C}((a_{t_2,j}-a_{s,i}) -
(a_{t_1,j} - a_{s,i})) = \text{ val}(a_{t_2,j} - a_{t_1,j})$ hence
$\text{val}^{\frak C}(a_{t_1,j}-a_{s,i})$ does not depend on $s$ as
$s$ does not appear on the left side, but, see above, it 
is $<_\Gamma$-increasing
with $s$, contradiction.  So we have finished proving $(*)^2_{i,j}$.]
\mr
\item "{$(*)^3_i$}"   for each $i \in {\Cal U}_1$, for 
some $t^*_i \in \{-\infty\} \cup I \cup\{+ \infty\}$ we have:
{\roster
\itemitem{ $(a)^3_i$ }  $\text{\rm val}^{\frak C}(c-a_{s,i}) =
\text{\rm val}^{\frak C}(a_{t,i} - a_{s,i})$ when $s <_I t$ and 
$s \in I_{<t^*_i}$ 
\sn
\itemitem{ $(b)^3_i$ }  $\langle\text{\rm val}^{\frak C}(c-a_{s,i}):s
\in I_{>t^*_i}\rangle$ is constant and if $r \in I_{> t^*_i}$
and $s <_I t$ are from $I_{> t^*_i}$ then 
$\text{val}^{\frak C}(c-a_{r,i}) <_\Gamma 
\text{val}^{\frak C}(a_{t,i}-a_{s,i})$

\sn
\itemitem{ $(c)^3_i$ }  $\text{\rm ac}^{\frak C}(c-a_{s,i}) = 
\text{\rm ac}^{\frak C}(a_{t,i} - a_{s,i})$ when $s <_I t$ and 
$s \in I_{<t^*_i}$
\sn
\itemitem{ $(d)^3_i$ }  $\langle\text{\rm ac}^{\frak C}(c-a_{s,i}):s \in 
I_{>t^*_i}\rangle$ is constant.
\endroster}
\ermn
[Why?  Recall the definition of ${\Cal U}_1$ which appeared just after
$(*)^1_i$ recalling that we are assuming $I$ is a complete linear
order, see $\boxdot_1(a)$.]
\mr
\item "{$(*)_4$}"  the set $J_1 = \{t^*_i:i \in {\Cal U}_1\}$
has at most one member in $I$.
\ermn
[Why?  Otherwise we can find $i,j$ from ${\Cal U}_1$ 
such that $t^*_i \ne t^*_j$ are from $I$.  
Now apply $(*)^2_{i,j} + (*)^3_i + (*)^3_j$.]

So \wilog \, 
\mr
\item "{$(*)_5$}"  $J_1$ is empty.
\ermn
[Why?  If not let $J_0 = \{t_*\}$ and we can  get it is enough to
prove the claim for $I_{<t_*}$ and for $I_{>t_*}$.]

Now  
\mr
\item "{$\boxplus_1$}"  if $i \in {\Cal U}_1$ and
$t^*_i = \infty$ \ub{then}
{\roster
\itemitem{ $(a)$ }  for every $s_0 <_I s_1 <_I s_2 <_I s_3$ we have
\sn
\itemitem{ $(b)$ }  $\{\text{\rm val}^{\frak C}(x - a_{s_3,i}) > 
\text{\rm val}^{\frak C}(a_{s_2,i} a- a_{s_1,i})\} \vdash
p^{[*]}_{a_{s_0,i}}$ and
\sn
\itemitem{ $(c)$ }   $c$ satisfies the formula in the left side;
on $p^{[*]}_{a_{s_0,j}}$, see $\boxdot_3$.
\endroster}
\ermn
[Why?  By clause (b) of \scite{dt.9} and $(*)^3_i + (*)^3_j$ and reflect.]

Hence
\mr
\item "{$\boxplus_2$}"   if ${\Cal W}_1 = \{i \in {\Cal U}_1:t^*_i =
\infty\}$ \ub{then} $\boxtimes_{{\Cal W}_1}$ 
\nl
where for ${\Cal W} \subseteq {\Cal U}$ we let
{\roster
\itemitem{ $\boxtimes_{\Cal W}$ }  if $s <_I t$ then
$\boxtimes^{s,t}_{\Cal W}$ where for ${\Cal U}' \subseteq {\Cal U}$:
\sn
\itemitem{ $\boxtimes^{s,t}_{{\Cal U}'}$ }   ${\Cal U}'
\subseteq \alpha,s,t \in I$ and 
$\bold f_{t,s}$ maps $\cup\{p^{[*]}_{a_{s,i}}:i \in {\Cal U}'\}$ onto
$\cup\{p^{[*]}_{a_{t,i}}:i \in {\Cal U}'\}$.
\endroster}
\ermn
[Why?  Should be clear as $J_1 = \emptyset$ and the indiscernibility
of $\langle \bar a_t:t \in I\rangle$ and 
$\boxplus_1$.]
\mr
\item "{$\boxplus_3$}"   assume that: for every 
$i \in {\Cal U}_1$ satisfying $t^*_i = - \infty$, there is 
$j \in {\Cal U}_1$ such that $t^*_j = -\infty$ and 
$s,t \in I \Rightarrow \text{\rm val}^{\frak C}(c - a_{t,j}) > \text{\rm
val}^{\frak C}(c-a_{s,i})$.  \ub{Then}:
{\roster
\itemitem{ $\odot_3$ }  if $s_0 <_I s_1 <_I s_2$ then 
$\{\text{\rm val}^{\frak C}(x-a_{s_2,j}) >
\text{\rm val}^{\frak C}\{(c-a_{s_1,i})\} \vdash
p^{[*]}_{a_{s_0,i}}$ and the formula on the left is satisfied
by $c$.
\endroster}
\ermn
[Why?  Should be clear.]

Hence
\mr
\item "{$\boxplus_4$}"   if the assumption of $\boxplus_3$ holds then
$\boxtimes_{{\Cal W}_2}$ holds for ${\Cal W}_2 = \{i \in {\Cal U}_1:t^*_i =
- \infty\}$.
\ermn
[Why?  As in $\boxplus_2$.]

Consider the assumption
\mr
\item "{$\boxplus_5$}"  the hypothesis of $\boxplus_3$ fails and
let $j(*) \in {\Cal U}_1$ exemplify this (so in particular $t^*_{j(*)} =
- \infty$).  Let ${\Cal W}_3 = \{i \in {\Cal U}_1:t^*_i = - \infty$ and
val$^{\frak C}(c-a_{s,j(*)}) > \text{\rm val}^{\frak C}(c-a_{t,i})$
for any $s,t \in I\}$ and ${\Cal W}_4 = \{i \in {\Cal U}_1:t^*_i = - \infty$
and $i \notin {\Cal W}_3\}$ so $j(*) \in {\Cal W}_4$
\sn
\item "{$\boxplus_6$}"  if $\boxplus_5$ then
$\boxtimes_{{\Cal W}_3}$.
\ermn
[Why?  Similarly to the proof of $\boxplus_2$.]
\mr
\item "{$\boxplus_7$}"  if $\boxplus_5$ then
{\roster
\itemitem{ $(a)$ }  $\langle \text{\rm val}^{\frak C}(c-a_{s,j}):s \in
I$ and $j \in {\Cal W}_4\rangle$ is constant
\sn
\itemitem{ $(b)$ }   $\text{\rm val}^{\frak C}(c-a_{r,j(*)}) 
<_\Gamma \text{ \rm val}^{\frak C}(a_{t,i}-a_{s,i})$ hence 
$(p_s)^{[*]}_{a_{s,j(*)}} \vdash
(p_s)^{[*]}_{a_{s,i}}$ when $i \in {\Cal W}_4$ and $s <_I t \wedge r \in I$
\sn
\itemitem{ $(c)$ }  for some finite $J_1 \subseteq I$ we have: if $s,t
\in J \backslash J_1$ and $(\forall r \in J_1)(s <_I s \equiv r <_I
t))$ then tp$(\text{\rm val}^{\frak C}(c-a_{s,j(*)}),M_s) = 
\bold f_{s,t}(\text{\rm tp}(\text{\rm val}^{\frak C}(c-a_{t,j(*)}),M_t))$
\sn
\itemitem{ $(d)$ }  for some finite $J_2 \subseteq I$ we have: if $s,t
\in I \backslash J_2$ and $(\forall r \in J_r)(r <_I s \equiv r <_I t)$ 
then tp$(\text{\rm ac}^{\frak C}(c-a_{s,j(*)}),M_s) = \bold
f_{s,t}(\text{\rm tp}(\text{\rm ac}^{\frak C}(c-a_{t,j(*)}),M_t)$
\sn
\itemitem{ $(e)$ }  for some finite $J_3 \subseteq I$ we
have: if $s,t \in I \backslash J_3$ and $(\forall r \in J)(r <_I s
\equiv r <_I t$, then $\boxtimes^{s,t}_{{\Cal W}_4}$
\endroster}
[Why?  Let $i \in {\Cal W}_4$; so $i \in {\Cal W}_2$, hence $i \in
{\Cal U}_1$, which means that case $(b)^1_i$ of
$(*)^1_i$ holds, so for each $t \in I$ the sequence
$\langle\text{val}^{\frak C}(a_{t,i} - a_{s,i}):s$ satisfies $s <_I
t\rangle$ is $<_\Gamma$-increasing.  Also as $i \in{\Cal W}_2$ clearly
$t^*_i= - \infty$ hence by $(*)^3_i (b)^3_i$ we have
$\langle\text{val}^{\frak C}(c-a_{s,i}):s \in I\rangle$ is constant,
call it $e_i$.  All this apply to $j(*)$, too.  Now as $i \in {\Cal
W}_4$ we know that for some $s_1,t_1 \in I$ we have $\text{val}^{\frak
C}(c-a_{s_1,j(*)}) \le_\Gamma \text{ val}^{\frak C}(c-a_{t_1,i})$,
i.e. $e_{j(*)} \le_\Gamma e_i$.  By the choice of $j(*)$ for every $j
\in {\Cal U}_1$ such that $t^*_j = - \infty$, i.e. for every $j \in
{\Cal W}_2$ for some (equivalently every) $s,t \in I$ we have
val$^{\frak C}(c-a_{s,j}) \le \text{ val}^{\frak C}(c-a_{t,j(*)})$.
In particular this holds for $j=i$, hence for some
$s_2,t_2 \in I$ we have $\text{val}^{\frak C}(c-a_{s_2,i}) \le
\text{ val}^{\frak C}(c-a_{t_2,j(*)})$, i.e $e_i \le_\Gamma e_{j(*)}$
so together with the previous sentence, $e_i = e_{j(*)}$, so clause
(a) of $\boxplus_7$ holds.  Also, the first phrase in clause (b) is 
easy (using $(*)^3_i(b)^3_i$ second phrase); the second phrase of (b)
follows because $e_i = e_{j(*)}$.
For clause (c) note that Th$(\Gamma^M)$ is strongly stable, for 
clause (d) note that Th$(k^M)$ is strongly dependent.  
\nl
Lastly, for clause (e) combine the earlier clauses.]
\sn
\item "{$\boxplus_8$}"  for some finite $J \subseteq I$, if $s,t \in I
\backslash J$ and $(\forall r \in J)(r <_I s \equiv r <_I t)$
\ub{then} $\boxtimes^{s,t}_{{\Cal U}_1}$
\nl
[Why?  If the hypothesis of $\boxplus_3$ holds let $J = \emptyset$ and
if it fails (so $\boxplus_5,\boxplus_6,\boxplus_7$ apply), let 
$J$ be as in $\boxplus_7(d),(e)$, so it partitions $I$ to
finitely many intervals.  It is enough to prove
$\boxtimes^{s,t}_{\Cal W}$ for several 
${\Cal W} \subseteq {\Cal U}_1$ which covers ${\Cal U}_1$.  Now by $\boxplus_2$
this holds for ${\Cal W}_1 = \{i \in {\Cal U}_1:t^*_i = \infty\}$.  If the
assumption of $\boxplus_3$ holds we get the same for ${\Cal W}_2$ by
$\boxplus_4$ and if it fails we get it for ${\Cal W}_3$ by
$\boxplus_6$ and for ${\Cal W}_4$ by $\boxplus_7(e)$ and the choice of
$J$.  Using ${\Cal U}_1 = {\Cal W}_1 \cup {\Cal W}_2,{\Cal W}_2 =
{\Cal W}_3 \cup {\Cal W}_4$ we are done.]
\ermn
As we can replace $I$ by its inverse
\mr
\item "{$\boxplus_9$}"  for some finite $J \subseteq I$ if $s,t \in I
\backslash J$ and $(\forall r)(r <_I s \equiv r <_I t)$ then
$\boxtimes^{s,t}_{{\Cal U}_{-1}}$.
\ermn
So we are left with ${\Cal U}_0$.
For $i \in {\Cal U}_0$ let $e_{0,i} = \text{ val}(a_{t,i} - a_{s,i})$
for $s <_I t$, well defined by the definition of ${\Cal U}_0$.  Let
${\Cal W}_5 := \{i \in {\Cal U}_0$: for every (equivalently some) $s
\ne t\in I$, val$^{\frak C}(c-a_{s,i}) < \text{ val}(a_{t,i} - a_{s,i})\}$ and
let ${\Cal W}_6 := {\Cal U}_0 \backslash {\Cal W}_5$.

Obviously
\mr
\item "{$\boxplus_{10}$}"  we have $\boxtimes_{{\Cal W}_5}$.
\ermn
Easily
\mr
\item "{$\boxplus_{11}$}"   if $i,j \in {\Cal W}_6$ then case
$(a)^2_{i,j}$ of $(*)^2_{i,j}$ holds.
\ermn
[Why?  By $(*)^2_{i,j}$ and as $i,j \in {\Cal W}_6 \Rightarrow
(*)^1_i(a) + (*)^1_i(b)$.]
\mr
\item "{$\boxplus_{12}$}"   if $i,j \in {\Cal W}_6$ and $s \ne t \in
I$ then val$^{\frak C}(a_{t,j} - a_{s,i}) = e_{0,i}$.
\ermn
[Why?  As ${\Cal W}_6 = {\Cal U}_0 \backslash {\Cal W}_5$.]

Hence
\mr
\item "{$\boxplus_{13}$}"  $\langle e_{0,i}:i \in {\Cal W}_6\rangle$
is constant, call the constant value 
$e_*$ so $s \ne t \in I \wedge i,j \in {\Cal W}_6
\Rightarrow \text{ val}^{\frak C}(a_{t,j} - a_{s,i}) = e_*$.
\ermn
Easily
\mr
\item "{$\boxplus_{14}$}"  for every $i \in {\Cal W}_6$ the set
$I_{i,c} := \{s \in I:\text{val}^{\frak C}(c-a_{s,i}) > e_*\}$ has at most one
member
\sn
\item "{$\boxplus_{15}$}"   let ${\Cal W}_7 := \{i \in {\Cal
W}_6:I_{i,c} \ne \emptyset\}$ and let $\{t^{**}_i\} = I_{i,c}$ for $i
\in {\Cal W}_7$
\sn 
\item "{$\boxplus_{16}$}"  if $i,j \in {\Cal W}_7$ then $t^{**}_i =
t^{**}_j$.
\ermn
[Why?  Otherwise let $t \in I$ be such that $t^{**}_i < t \wedge t^{**}_j < t$,
now val$^{\frak C}(c-a_{t^{**}_i,j}) > \text{ val}^{\frak C}(a_{t,i} 
- a_{t^{**}_i,i}) = e_*$ and val$^{\frak C}(c - a_{t^{**}_j,j}) 
> \text{ val}^{\frak C}(a_{t,j} -
a_{t^{**}_j,j}) = e_*$ hence $e_* < \text{ val}^{\frak C}
((c - a_{t^{**}_i,i}) -
(c - a_{t^{**}_j,j})) = \text{ val}^{\frak C}(a_{t^{**}_j,j} - a_{t^{**}_i,i})$
but the last one is $e_*$ by $\boxplus_{12}$, contradiction.]
\mr
\item "{$\boxplus_{17}$}"  \wilog \, ${\Cal W}_7 = \emptyset$.
\ermn
[Why?  E.g. as otherwise we can prove separately for $I_{<t^{**}_i}$ and
for $I_{>t^{**}_i}$ for any $i \in {\Cal W}_7$.]
\mr 
\item "{$\boxplus_{18}$}"  if $i,j \in {\Cal W}_6$ and $s \ne t \in I$
\ub{then} ac$^{\frak C}(c - a_{t,j}) - \text{ ac}^{\frak C}
(c - a_{s,i}) = \text{ ac}^{\frak C}(a_{s,i} - a_{t,j})$.
\ermn
[Why?  As val$^{\frak C}(c - a_{t,j}),\text{val}^{\frak C}(c - a_{s,i})$ and 
$\text{val}^{\frak C}(c_{s,i} - (c_{t,j})$ are all equal to $e_*$.]

The rest should be clear.
\nl
3) For the $\omega$-language: the proof is similar.
\hfill$\square_{\scite{dp0.16}}$
\enddemo
\newpage

\head {\S2 Cutting indiscernible sequence and strongly$^+$ dependent} \endhead  \resetall \sectno=2
 \spuriousreset
\bigskip

\demo{\stag{dp1.2.4} Observation}  1) The following conditions on $T$
are equivalent, for $\alpha \ge \omega$
\mr
\item "{$(a)$}"  $T$ is strongly dependent, i.e., $\aleph_0 =
\kappa_{\text{ict}}(T)$ 
\sn
\item "{$(b)_\alpha$}"  if $I$ is an infinite linear order, $\bar a_t
\in {}^\alpha{\frak C}$ for $t \in I,\bold I = 
\langle \bar a_t:t \in I\rangle$ is
an indiscernible sequence and $C \subseteq {\frak C}$ is finite, 
\ub{then} there is a convex equivalence relation $E$ on $I$ with
finitely many equivalence classes such
that $s E t \Rightarrow \text{ tp}(\bar a_s,C) = 
\text{ tp}(\bar a_t,C)$
\sn
\item "{$(c)_\alpha$}"  if $\bold I = \langle \bar a_t:t \in I\rangle$ is as
above and $C \subseteq {\frak C}$ is finite \ub{then} there is a convex
equivalence relation $E$ on $I$ with finitely many equivalence classes
such that: if $s \in I$ then $\langle
\bar a_t:t \in (s/E) \rangle$ is an indiscernible sequence over $C$.
\ermn
2) We can add to the list in (1)
\mr
\item "{$(b)'_\alpha$}"   like $(b)_\alpha$ but $C$ a singleton
\sn
\item "{$(c)'_\alpha$}"  like $(c)_\alpha$ but the set $C$ is a
singleton.
\ermn
3) We can in part (1),(2) clauses
$(c)_\alpha,(b)_\alpha,(b)'_\alpha,(c)'_\alpha$ restrict ourselves to
well order $I$.
\nl
4) In parts (1),(2),(3), given $\kappa = \kappa^{< \theta},\theta > 
|T|$, in clauses $(b)_\kappa,(c)_\kappa$ and their
parallels we can add that ``$\bar a_\alpha$ is the universe of a
$\theta$-saturated model"; moreover we allow $\bold I$ to be:
\mr
\widestnumber\item{$(iii)$}
\item "{$(i)$}"  $\bold I = \langle \bar a_u:u \in
[I]^{< \aleph_0} \rangle$ is indiscernible over $A$ (see Definition
\scite{dw5.1}(2))
\sn
\item "{$(ii)$}"   $\bar a_{\{t\}} = \bar a_t$,
\sn
\item "{$(iii)$}"   each $\bar a_t$ is the universe of a $\theta$-saturated
model
\sn
\item "{$(iv)$}"   for some infinite linear orders $I_{-1},I_1$ and
some $\bold I' = \langle \bar a'_u:u \in [I_{-1} + I
+ I_1]^{< \aleph_0} \rangle$ indiscernible over $A = \text{\rm
Rang}(\bar a_\emptyset)$ we have:
{\roster
\itemitem{ $(\alpha)$ }  $u \in [I]^{< \aleph_0} \Rightarrow \bar a'_u
= \bar a_u$
\sn
\itemitem{ $(\beta)$ }   for every $B \subseteq A$ of
cardinality $< \theta$, every subtype of the type of $\langle \bar
a_u:u \in [I_{-1} + I_1]^{< \aleph_0}\rangle$ over $\langle \bar a_u:u \in
[I]^{< \aleph_0}\rangle$ of cardinality $< \theta$ is realized in $A$
(we can use only $A$ and $\langle \bar a_t:t \in I\rangle$, of course).
\endroster}
\endroster
\enddemo
\bigskip

\remark{Remark}  1) Note that \scite{df2.4.8} below says more for the
cases $\kappa_{\text{ict}}(T) > \aleph_0$ so no point to deal with it
here.
\nl
2) We can in \scite{dp1.2.4} add in
$(b)_\alpha,(c)_\alpha,(b_\alpha)',(c_\alpha)'$ ``over a fixed $A$"
by \scite{dp1.2.1}(3).
\nl
3) By \scite{dq.6} we can translate this to the case of a family of
indiscernible sequences.
\endremark
\bigskip

\demo{Proof}  1) Let $\kappa = \omega$ (to serve in the proof of a subsequence
observation).
\bn
\ub{$\neg(a) \Rightarrow \neg(b)_\alpha$}

Let $\lambda > \aleph_0$, as in the proof of \scite{dp1.2.2}, because
we are assuming $\neg(a)$, there are 
$\bar\varphi = \langle \varphi_i(\bar x,\bar y_i):i < \kappa\rangle$ and
$\langle \bar a^i_\alpha:i < \omega,\alpha <
\lambda \rangle$ witnessing $\circledast^2_{\bar\varphi}$ from there.

For $\alpha < \lambda$ let 
$\bar a^*_\alpha \in {\frak C}$ be the concatenation of $\langle \bar
a^i_\alpha:i < \kappa \rangle$, possibly with repetitions so it has length
$\kappa$.

Let $\eta = \langle \omega n:n < \omega \rangle$ and $\bar b^*$
realizes $\{\varphi_n(x,\bar a^n_{\omega n}) \wedge \neg
\varphi_n(x,\bar a^n_{\omega n+1}):n < \omega\}$.

So for each $n$, tp$(\bar a^n_{\omega n},\bar b^*) \ne \text{
tp}(a^n_{\omega n+1},\bar b^*)$ hence tp$(\bar a^*_{\omega n},\bar
b^*) \ne \text{ tp}(\bar a^*_{\omega n+1},\bar b^*)$.  So any convex
equivalence relation on $\lambda$ as required (i.e. such that
$\alpha E \beta \Rightarrow \text{ tp}(\bar a^*_\alpha,\bar b^*) =
\text{ tp}(\bar a^*_\beta,\bar b^*))$ satisfies $n < \omega
\Rightarrow \neg(\omega n)E(\omega n+1)$; it certainly shows $\neg(b)_\alpha$.
\bn
\ub{$\neg(b)_\alpha \Rightarrow \neg(c)_\alpha$}

Trivial.
\bn
\ub{$\neg(c)_\alpha \Rightarrow \neg(a)$}

Let $\langle \bar a_t:t \in I\rangle$ and $C$ exemplify
$\neg(c)_\alpha$, and assume toward contradiction that $(a)$ holds.
Without loss of generality $I$ is a dense linear
order (so with neither first nor last element) and is complete and let
$\bar c$ list $C$.
\nl
So
\mr
\item "{$(*)$}"  for no convex equivalence relation $E$ on $I$ with
finitely many equivalence classes do we have $s \in I \Rightarrow \langle \bar
a_t:t \in (s/E) \rangle$ is an indiscernible sequence over $C$.
\ermn
We now choose $(E_n,I_n,\Delta_n,J_n)$ by induction on $n$ such that
\mr
\item "{$\circledast$}"  $(a) \quad E_n$ is a convex equivalence
 relation on $I$ such that each equivalence 
\nl

\hskip25pt class is dense (so with no extreme member!) or is a singleton
\sn
\item "{${{}}$}"  $(b) \quad \Delta_n$ is a finite set of formulas
(each of the form $\varphi(\bar x_0,\dotsc,\bar x_{m-1},\bar y)$,
\nl

\hskip25pt $\ell g(\bar x_\ell) = 
\alpha$, for some $m,\ell g(\bar y) = \ell g(\bar c)$)
\sn
\item "{${{}}$}"  $(c) \quad I_0 = I,E_0$ is the equality, $\Delta_0 =
\emptyset$ 
\sn
\item "{${{}}$}"  $(d) \quad I_{n+1}$ is one of the equivalence
classes of $E_n$ and is infinite
\sn
\item "{${{}}$}"  $(e) \quad \Delta_{n+1}$ is a finite set of formulas
such that $\langle \bar a_t:t \in I_{n+1} \rangle$ is not \nl

\hskip25pt $\Delta_{n+1}$-indiscernible over $C$
\sn
\item "{${{}}$}"  $(f) \quad E_{n+1} \restriction I_{n+1}$ is a convex
equivalence relation with finitely many 
\nl

\hskip25pt classes,  each dense (no
extreme member) or singleton, if $J$ 
\nl

\hskip25pt is an infinite equivalence class of $E_{n+1} \restriction I_{n+1}$ 
\nl

\hskip25pt then $\langle \bar a_t:t \in J\rangle$ is
$\Delta_{n+1}$-indiscernible over $C$ and 
\nl

\hskip25pt $|I_{n+1}/E_{n+1}|$ is minimal under those conditions
\sn
\item "{${{}}$}"  $(g) \quad E_{n+1} \restriction (I \backslash
I_{n+1}) = E_n \restriction (I \backslash I_{n+1})$, so $E_{n+1}$
refines $E_n$
\sn
\item "{${{}}$}"  $(h) \quad$ we choose $(\Delta_{n+1},E_{n+1})$ such that
if possible $I_{n+1}/E_{n+1}$ has
\nl

\hskip25pt $\ge 4$ members. 
\ermn
There is no problem to carry the induction as $T$ is dependent (see
\scite{df1.2.5}(1) below which says more or see
\cite[3.4+Def 3.3]{Sh:715}).
\nl
For $n>0,E_n \restriction I_n$ is an equivalence relation on
$I_n$ with finitely many equivalence classes, each convex; so as $I$
is a complete linear order clearly
\mr
\item "{$(*)_1$}"  for each $n>0$ there are $t^n_1 <_I \ldots <
t^n_{k(n)-1}$ from $I_n$ such that $s_1 \in I_n \wedge s_2 \in I_n
\Rightarrow [s_1 E_n s_2 \equiv (\forall k)(s_1 < t^n_k \equiv s_2 <
t^n_k \wedge s_1 > t^n_k \equiv s_2 > t^n_k)]$.
\ermn
As $n > 0 \Rightarrow E_n \ne E_{n-1}$ clearly
\mr
\item "{$(*)_2$}"  $k(n) \ge 2$ and $|I_n/E_n| = 2k(n) - 1$
\sn
\item "{$(*)_3$}"  $\{I_{n,\ell}:\ell < k(n)\} \cup \{\{t^n_\ell\}:0
< \ell < k(n)\}$ are the equivalence classes of $E_n \restriction
I_n$; 
\nl
where
\sn
\item "{$(*)_4$}"  for non-zero $n < \omega,\ell < k(i)$ we define
$I_{n,\ell}$:
\nl

if $0 < \ell < k(n)-1$ then $I_{n,\ell} =
(t^n_\ell,t^n_{\ell+1})_{I_n}$
\nl

if $0 = \ell$ then $I_{n,\ell} = (-\infty,t^n_\ell)_{I_n}$
\nl

if $\ell = k(n)-1$ then $I_{n,\ell} = (t^n_\ell,\infty)_{I_n}$.
\ermn
As (see end of clause (f))) we cannot omit any $t^n_\ell(\ell < k(n))$
and transitivity of equality of types clearly
\mr
\item "{$(*)_5$}"   for each $\ell < k(n)-1$ for some $m$ and $\varphi =
\varphi(x_0,\dotsc,\bar x_{m-1},\bar y) \in \Delta_n$ there are $s_0
<_I \ldots <_I s_{m-1}$ from $I_{n,\ell}$ and $s'_0 <_I \ldots <_I
s'_{m-1}$ from $I_{n,\ell} \cup \{t^n_{\ell+1}\} \cup I_{n,\ell+1}$ such
that ${\frak C} \models \varphi[\bar a_{s_0},\dotsc,\bar c] \equiv
\neg \varphi[a_{s'_0},\dotsc,\bar c]$.
\ermn
Hence easily
\mr
\item "{$(*)_6$}"  $J \in \{I_{n,\ell}:\ell < k(n)\}$ iff $J$ is a
maximal open interval of $I_n$ such that $\langle \bar a_t:t \in
J\rangle$ is $\Delta_n$-indiscernible over $C$.
\ermn
By clause (h) and $(*)_6$
\mr
\item "{$(*)_7$}"  if $k(n) < 4$ and $\ell < k(n)$ then $\langle a_t:t \in
I_{n,\ell}\rangle$ is an indiscernible sequence over $C$
\nl
hence
\item "{$(*)_8$}"  if $k(n) < 4$ then for no $m>n$ do we have $I_m
\subseteq I_n$.
\ermn
Note that 
\mr
\item "{$(*)_9$}"  $m < n \Rightarrow I_n \subset I_m \vee I_n \cap
I_m =0$.
\endroster
\bn
\ub{Case 1}:  There is an infinite $u \subseteq \omega$ such that $\langle
I_n:n < \omega \rangle$ are pairwise disjoint.

For each $n$ we can find $\bar c_n \in {}^{\omega >}C$ and $k_n <
\omega$ (no connection to $k(n)$ from above!) and
$\varphi(\bar x_0,\dotsc,\bar x_{k_n-1},\bar y) \in \Delta_n$ such that
$\langle \bar a_t:t \in I_n \rangle$ is not $\varphi_n(\bar
x_0,\dotsc,\bar x_{k_n-1},\bar c)$-indiscernible (so $\ell g(\bar
x_\ell) = \alpha$).  So we can find 
$t^\ell_0 < \ldots < t^\ell_{k_n-1}$ in $I_n$
for $\ell =1,2$ such that $\models \varphi_n[\bar a_{t^\ell_0},
\dotsc,\bar a_{t^\ell_{k_n-1}},\bar c_n]^{\text{if}(\ell=2)}$.  By
minor changes in $\Delta_n,\varphi_n$, \wilog \, $\bar c_n$ is without
repetitions hence without loss of generality 
$n < \omega \Rightarrow \bar c_n = \bar c_*$.

Without loss of generality $\Delta_n$ is closed under negation and
without loss of generality 
$t^1_{k_n-1} <_I t^2_0$.  We can choose $t^m_k \in I_n (m < \omega,m
\notin \{1,2\},k < k_n)$ such that for every $m < \omega,k < k_n$ we have
$t^m_k <_I t^m_{k+1},t^m_{k_n-1} <_I
t^{m+1}_0$; let $\bar a^*_{n,m} = \bar a_{t^m_0} \char 94 \ldots \char
94 \bar a_{t^m_{k_n-1}}$ and let $\bar x = \langle x_i:i < \ell g(\bar
c_*)\rangle$.  So for every $\eta \in {}^\omega \omega$ the type
$\{\neg \varphi_n(\bar a^*_{n,\eta(n)},\bar x) \wedge \varphi_n(\bar
a^*_{n,\eta(n)+1},\bar x):n < \omega\}$ is consistent.  This is enough
for showing $\kappa_{\text{ict}}(T) > \aleph_0$.
\bn
\ub{Case 2}:  There is an infinite $u \subseteq \omega$ such that
$\langle I_n:n \in u \rangle$ is decreasing.

For each $n \in u,E_n \restriction I_n$ has an infinite equivalence
class $J_n$ (so $J_n \subseteq I_n$) such that $n < m \wedge \{n,m\}
\subseteq u \Rightarrow I_m \subseteq J_n$.  By $(*)_8$ clearly for
each $n \in u,k(n) \ge 4$ hence we can find $\ell(n) <k(n)$ such that $I'_n
= (I_{n,\ell(n)} \cup \{t^n_{\ell,n}\} \cup I_{n,\ell(n)+1})$ is
disjoint to $J_m$.  Now $\langle I'_n:n \in u\rangle$
are pairwise disjoint and we continue as in Case 1.  

By Ramsey theorem at least one of the two cases occurs so we are done.
\nl
2) By induction on $|C|$.  
\nl
3),4)  Easy by now.  \hfill$\square_{\scite{dp1.2.4}}$
\enddemo
\bn
Recall
\demo{\stag{df1.2.5} Observation}  1) Assume that $T$ is 
dependent, $\langle \bar
a_t:t \in I \rangle$ is an indiscernible sequence, $\Delta$ a finite set
of formulas, $C \subseteq {\frak C}$ finite.  \ub{Then} for some
convex equivalence relation $E$ on $I$ with finitely many equivalence
classes, each equivalence class in an infinite open convex set or is 
a singleton such that
for every $s \in I,\langle \bar a_t:t \in s/E\rangle$ is an
$\Delta$-indiscernible sequence over $\cup\{\bar a_t:t \in I\backslash (s/E)\}
\cup C$.
\nl
2) If $I$ is dense and complete 
there is the least fine such $E$.  In fact for $J$ an open convex
subset of $I$ we have: $J$ is an $E$-equivalence class \ub{iff} $J$
is a maximal open convex subset of $I$ such that $\langle \bar a_t:t \in
J\rangle$ is $\Delta$-indiscernible over $C \cup \cup\{\bar a_t:t \in
I \backslash J\}$. 
\nl
3) Assume if $I$ is dense (no extreme elements) and complete then
there are $t_1 <_I \ldots < t_{k-1}$ such that stipulating $t_0 = -
\infty,t_k = \infty,I_\ell = (t_\ell,t_{\ell+1})_I$, we have
\mr
\item "{$(a)$}"  $\langle \bar a_t:t \in I_\ell\rangle$ is
indiscernible over $C$
\sn
\item "{$(b)$}"  if $\ell \in \{1,\dotsc,k-1\})$ and $t^-_\ell <_I
t_\ell <_I t^+_\ell$, then $\langle a_t:t \in
(t^-_\ell,t^+_\ell)_I\rangle$ is not $\Delta$-indiscernible over $C$.
\endroster
\enddemo
\bigskip

\demo{Proof}  1) See clause (b) of \cite[Claim 3.2]{Sh:715}.
\nl
2),3)   Done inside
the proof of \scite{dp1.2.4} and see proof of \scite{df2.2.9}.
\hfill$\square_{\scite{df1.2.5}}$
\enddemo
\bigskip

\definition{\stag{df2.3.5} Definition}  1) We say that $\bar\varphi =
\langle \varphi_i(\bar x,\bar y_i):i < \kappa \rangle$ witness $\kappa
< \kappa_{\text{ict},2}(T)$ \ub{when} there are a sequence $\langle \bar
a_{i,\alpha}:\alpha < \lambda,i < \kappa \rangle$ and $\langle \bar
b_i:i < \kappa \rangle$ such that
\mr
\item "{$(a)$}"  $\langle \bar a_{i,\alpha}:\alpha < \lambda \rangle$
is an indiscernible sequence over $\cup\{\bar a_{j,\beta}:j \in \kappa
\backslash \{i\}$ and $\beta < \lambda\}$ for each $i < \kappa$
\sn
\item "{$(b)$}"  $\bar b_i \subseteq \cup\{\bar
a_{j,\alpha}:j<i,\alpha < \lambda\}$
\sn
\item "{$(c)$}"  $p = \{\varphi_i(\bar x,\bar a_{i,0} \char 94 
\bar b_i),\neg \varphi_i(\bar x,\bar a_{i,1} \char 94 \bar b_i):i < \kappa\}$
is consistent (= finitely satisfiable in ${\frak C}$).
\ermn
2) $\kappa_{\text{ict},2}(T)$ is
the first $\kappa$ such that there is no witness for $\kappa <
\kappa_{\text{ict},2}(T)$.  
\nl
3) $T$ is strongly$^2$ dependent (or strongly$^+$ dependent) if
$\kappa_{\text{ict},2}(T)=\aleph_0$. 
\nl
4) $T$ is strongly$^2$ stable if it is strongly$^2$ dependent and
stable.
\enddefinition
\bigskip

\demo{\stag{df2.3.27} Observation}  If $M$ is a valued field in the
sense of Definition and $|\Gamma^M| > 1$ \ub{then} $T := \text{
Th}(M)$ is not strongly$^2$ dependent.
\enddemo
\bigskip

\demo{Proof}  Let $a \in \Gamma^M$ be positive, $\varphi_0(x,a) :=
(\text{val}(x) \ge a),E(x,y,a) := (\text{ val}(x,y) \ge 2a)$ and
$F(x,y) = x^2+y$ (squaring in $K^M$).  Now for $b \in
\varphi_0(M,\bar a)$, the funciton $F(-,b)$ is $(\le 2)$-to-1 function from
$\varphi_0(M,a)$ to $b/E$.  So we can apply \cite[\S4]{Sh:783}.

Alternatively let $a_n \in \Gamma^M,a_n <_{\Gamma^M} a_{n+1}$ for $n <
\omega$ be such that there are $b_{n,\alpha} \in K^M$ for $\alpha <
\omega$ such that $\alpha < \beta < \omega \Rightarrow \text{
val}^M(b_{n,\alpha}-b_{n,\beta}) > a_n \le a_n$ and val$(b_n,\alpha)$.
Without loss of generality for eac $n < \omega$ the sequence
$\langle b_{n,\alpha}:\alpha < \omega\rangle$ is indiscernible over
$\{b_{n_1,\alpha_1}:n_1 \in \omega \backslash \{n\},\alpha <
\omega\} \cup \{a_{n_1}:n_1 < \omega\}$.  Now for $\eta \in {}^\omega
\omega$ clearly $p_\eta = \{\text{val}(x-\Sigma\{a_{m,\eta(m)}:m <
n\}) > a_n:n < \omega\}$, it is consistent, and we have an example.
\hfill$\square_{\scite{df2.3.27}}$ 
\enddemo
\bn
Note that the definition of strongly$^2$ dependent here (in
\scite{df2.3.5}) is equivalent to the one in \cite[3.7]{Sh:783}(1) by
$(a) \Leftrightarrow (e)$ of Claim \scite{dq.7} below.
\nl
The following example shows that there is a difference even among the
stable $T$.
\nl
\margintag{dq1.8}\ub{\stag{dq1.8} Example}:   There is a strongly$^1$ stable not
strongly$^2$ stable $T$ (see Definition \scite{df2.3.5}).
\bigskip

\demo{Proof}   Fix $\lambda$ large enough.
Let $\Bbb F$ be a field, let $V$ be a vector
space over $\Bbb F$ of infinite dimension, let $\langle V_n:n < \omega
\rangle$ be a decreasing sequence of subspaces of $V$ with $V_n/V_{n+1}$
having infinite dimension $\lambda$ and $V_0 = V$ and $V_\omega =
\cap\{V_n:n < \omega\}$ have dimension $\lambda$.  Let $\langle
x^n_\alpha + V_{n+1}:\alpha < \lambda\rangle$ be a basis of
$V_n/V_{n+1}$ and let $\langle x^{\omega,i}_\alpha:i \in \Bbb Z$ and
$\alpha < \lambda\rangle$ be a basis of $V_\omega$.  
Let $M = M_\lambda$ be the following model:
\mr
\item "{$(a)$}"  universe: $V$
\sn
\item "{$(b)$}"  individual constants: $0^V$
\sn
\item "{$(c)$}"  the vector space operations: 
$x+y,x-y$ and cx for $c \in \Bbb F$
\sn
\item "{$(d)$}"  functions: $F^M_1$, a linear unary 
function: $F^M_1(x^n_\alpha) =
x^{n+1}_\alpha,F^M_1(x^{\omega,i}_\alpha) = x^{\omega,i+1}_\alpha$
\sn
\item "{$(e)$}"  $F^M_2$, a linear unary function:
\nl
$F^M_2(x^0_\alpha) = x^0_\alpha,F^M_2(x^{n+1}_\alpha) =
x^n_\alpha$ and $F^M_2(x^{\omega,i}_\alpha) = x^{\omega,i-1}_\alpha$
\sn
\item "{$(f)$}"  predicates: $P^M_n = V_n$ so $P_n$ unary
\ermn
Now
\mr
\item "{$(*)_0$}"  for any models $M_1,M_2$ of Th$(M_\lambda)$ with uncountable
$\cap\{P^{M_\ell}_n:n < \omega\}$ for $\ell=1,2$, the set ${\Cal F}$
exemplify $M_1,M_2$ are $\Bbb L_{\infty,\aleph_0}$-equivalent where:

${\Cal F}$ is the family of partial isomorphisms $f$ from $M_1$ into
$M_2$ such that 

for some $n,\langle N_i:\kappa > 0,i = \omega\rangle$
we have:
{\roster
\itemitem{ $(a)$ }  Dom$(f) = \bigoplus_{i<n} N_i \oplus N_\omega$
\sn
\itemitem{ $(b)$ }  $N_i \subseteq P^{M_1}_i$ is a subspace when $i < n
\vee i = \omega$
\sn
\itemitem{ $(c)$ }  $N_i$ is of finite dimension
\sn
\itemitem{ $(d)$ }  $N_i \cap P^{M_1}_{i+1} = \{0\}$ for $i < n$
\sn
\itemitem{ $(e)$ }  similar conditions on 
$N'_i = f(N_i)$ for $i < n \vee i = \omega$
\endroster}
\item "{$(*)_1$}"  $T = \text{ Th}(M_\lambda)$ has elimination of quantifiers
\nl
[Why?  Easy.]
\ermn
Hence
\mr
\item "{$(*)_2$}"  $T$ does not depend on $\lambda$
\sn
\item "{$(*)_3$}"   $T$ is stable.
\ermn
[Why?  As if $N_1$ is $\aleph_1$-saturated, $N_1 \prec N_2$ then
$\{\text{tp}(a,N_1,N_2):a \in {\frak C}\}$ has cardinality $\le
\|N_1\|^{\aleph_0}$ by $(*)_2$.]

Now
\mr
\item "{$(*)_4$}"  $T$ is not strongly$^2$ dependent.
\nl
[Why?  By \scite{0.gr.1}. Alternatively, define 
a term $\sigma_n(y)$ by induction on 
$n:\sigma_0(y) = y,\sigma_{n+1}(y) = F_1(\sigma_n(y))$, and 
for $\eta \in {}^\omega \lambda$ increasing let   

$$
\align
p_\eta(y) = \{&P_1(y - \sigma_0(x^0_{\eta(0))}),
P_2(y - \sigma_0(x^0_{\eta(0)}) - \sigma_1(x^1_{\eta(1)})),\ldots, \\
  &P_n(y -\Sigma\{\sigma_\ell(x^\ell_{\eta(\ell)}):\ell < n\}),\ldots\}.
\endalign
$$
\mn
Clearly each $p_\eta$ is finitely satisfiable in $M_\lambda$.  Easily
this proves that $T$ is not strongly$^2$ stable $I$.]
\ermn
So it remains to prove
\mr
\item "{$(*)_5$}"  $T$ is strongly stable.
\ermn
Why this holds?  We work in ${\frak C} = {\frak C}_T$.
Let $\lambda \ge (2^\kappa)^+$ be large enough 
and $\kappa = \kappa^{\aleph_0}$.  We
shall prove $\kappa_{\text{ict}}(T) = \aleph_0$ by the variant of
$(b)'_\omega$ from \scite{dp1.2.4}(3), this suffices.
Let $\langle \bar
a_\alpha:\alpha < \lambda \rangle$ be an indiscernible sequence over a set
$A$ such that $\ell g(\bar a_\alpha) \le \kappa$.  By
\scite{dq.6} \wilog \, each $\bar a_\alpha$ enumerate the set of elements of an
elementary submodel $N_\alpha$ of ${\frak C}$ which include $A$ and is
$\aleph_1$-saturated.
\nl
Without loss of generality ($I \cap \Bbb Z = \emptyset$ and):
\mr
\item "{$\boxdot_1$}"  for some $\bar a'_n (n \in \Bbb Z),A \supseteq
c \ell(A' \cup \cup\{\bar a'_i:i \in \Bbb Z\})$ and 
$\langle \bar a'_n:n < 0\rangle \char 94 \langle \bar a_\alpha:
\alpha < \lambda\rangle \char 94 \langle \bar a'_n:n \ge 0 \rangle$ is 
an indiscernible sequence over $A'$ and $\langle \bar a_\alpha:\alpha < \lambda
\rangle \char 94 \langle A \rangle$ is linearly 
independent over $A',A$ is the universe of $N,N$ is
$\aleph_1$-saturated and $N \cap N_\alpha$ is $\aleph_1$-saturated
(and does not depend on $\alpha$)
\ermn
Hence by $(*)_0$
\mr
\item "{$\boxdot_2$}"  $(a) \quad \alpha \ne \beta 
\wedge a_{\alpha,i} = a_{\beta,j}
\Rightarrow a_{\alpha,i} = a_{\beta,i} \in A$
\sn
\item "{${{}}$}"  $(b) \quad$ if $u \subseteq \lambda$ then
cl$(\cup\{\bar a_\alpha:\alpha \in u\} \cup A\})$ is $\prec {\frak C}$
\sn
\item "{${{}}$}"  $(c) \quad$ if $u \subseteq \lambda$ is finite we
get an $\aleph_1$-saturated model (not really used).
\ermn
(We can use the stronger \scite{dp1.2.4}(4)).
Easily
\mr
\item "{$\boxdot_3$}"  if $a \in N_\alpha,b \in c
\ell(\cup\{N_\beta:\beta < \alpha\} \cup A)$ \ub{then}:
{\roster
\itemitem{ $(a)$ }  $a=b \Rightarrow a \in A$
\sn
\itemitem{ $(b)$ }  $a-b \in P^{\frak C}_n \Rightarrow (\exists c \in A)(a-c
\in P^{\frak C}_n \wedge b-c \in P^{\frak C}_n)$.
\endroster}
\ermn
[Why?  Let $b = \sigma^{\frak C}(\bar a_{\beta_0},\dotsc,\bar
a_{\beta_{m-1}},\bar a),\bar a \in {}^{\omega >}A,\sigma$ a term,
$\beta_0 < \beta_1 < \ldots < \beta_{m-1},\bar c \in {}^{\omega
>} A$, then for every $k < \omega$ large enough $b' := \sigma^{\frak
C}(a'_k,\bar a'_{k+1},\ldots,\bar a_{k+m-1},\bar a)$ belongs to $A$
and in Case (a): $a = b \Rightarrow a = b'$ and in case (b): $a-b \in
P^{\frak C}_n \Rightarrow a-b' \in P^{\frak C}_n$.]
\mr
\item "{$\boxdot_4$}"  if $a_\ell \in c \ell(\cup\{N_\alpha:\alpha \in 
u_\ell\} \cup A)$ and $u_\ell \subseteq \lambda$ for $\ell=1,2$ \ub{then}:
{\roster
\itemitem{ $(a)$ }  if $a_1 = a_2$ then for some $b \in c \ell
(\cup\{N_\alpha:\alpha \in u_1 \cap u_2\} \cup A)$ we have $a_1-b =
a_2-b \in A$
\sn
\itemitem{ $(b)$ }  if $a_1-a_2 \in P^{\frak C}_n$ then for some $b
\in c \ell(\{N_\alpha:\alpha \in u_1 \cap u_2\} \cup A)$ and $c \in A$
we have $a_2-b-c \in P^{\frak C}_n$ and $a_2-b-c \in P^{\frak C}_n$.
\endroster}
\ermn
[Why?  Similarly to $\boxdot_2$.]

Now let $c \in {\frak C}$, the proof splits to cases.
\bn
\ub{Case 1}:  $c \in c \ell(\cup\{\bar a_\beta:\beta < \lambda\} \cup A)$.

So for some finite $u \subseteq \lambda,c \in c\ell(\cup\{\bar
a_\beta:\beta \in u\})$, easily $\langle \bar a_\beta:\beta \in
\lambda \backslash u\rangle$ is an indiscernible set over $A \cup \{c\}$, and
we are done.
\bn
\ub{Case 2}:  For some finite $u \subseteq \lambda$, for every $n$ for
some $c_n \in c \ell(\cup\{\bar a_\beta:\beta \in u\} \cup A)$ we have $c-c_n
\in P^M_n$ (but not case 1).

Clearly $u$ is as required.  (In fact, easily 
$c \ell(\{\bar a_\beta:\beta \in u\} \cup A)$ is $\aleph_1$-saturated
(as $u$ is finite, by $\boxdot_2(c)$) 
hence there is $c^* \in c \ell(\cup\{a_\beta:\beta
\in u\} \cup A)$ such that $n < \omega \Rightarrow c^* - c_n \in P^M_n$.
\bn
\ub{Case 3}:  Neither case 1 nor case 2 (less is needed).

Let $n(1) < \omega$ be maximal such that for some $c_{n(1)} \in A$ we have
$c - c_{n(1)} \in P^M_{n(1)}$ (for $n=0$ every $c' \in A$ is O.K.; by
not Case 2 such $n(1)$ exists).
\bn
\ub{Subcase 3A}:  There is $n(2) \in (n(1),\omega)$ and $c_{n(2)} \in
c  \ell(\{\bar a_\beta:\beta < \lambda\} \cup A)$ such that $c - c_{n(2)} \in
P^M_{n(2)}$.  

Let $u$ be a finite subset of $\lambda$ such that
$c_{n(2)} \in c \ell(\{\bar a_\beta:\beta \in u\} \cup A)$, now $u$ is as
required (by $\boxdot_3 + \boxdot_4$ above).
\bn
\ub{Subcase 3B}:  Not subcase 3A.

Choose $u=\emptyset$ works because neither case 1,
nor case 2 hold with $u = \emptyset$.  \hfill$\square_{\scite{dq1.8}}$
\enddemo
\bigskip

\remark{\stag{dq1n.8.7} Remark}  We can prove a claim parallel to
\scite{dw1.10}, i.e. replacing strong dependent by strongly$^2$ dependent. 
\endremark
\bigskip

\proclaim{\stag{df2.3.51} Claim}  1) $\kappa_{\text{ict},2}(T^{\text{eq}}) =
\kappa_{\text{ict},2}(T)$.
\nl
2) If $T_\ell = \text{\rm Th}(M_\ell)$ for $\ell=1,2$ then
$\kappa_{\text{ict},2}(T_1) \ge \kappa_{\text{ict},2}(T_2)$
when:
\mr
\item "{$(*)$}"  $M_1$ is (first order) interpretable in $M_2$.
\ermn
3) If $T' = \text{\rm Th}({\frak C},c)_{c \in A}$ then
$\kappa_{\text{ict},2}(T') = \kappa_{\text{ict},2}(T)$.
\nl
4) If $M$ is the disjoint sum of $M_1,M_2$ (or the product) and
{\rm Th}$(M_1)$, {\rm Th}$(M_2)$ are strongly$^2$ dependent then 
so is {\rm Th}$(M)$.
\endproclaim
\bigskip

\demo{Proof}  Similar to \scite{dw1.10}. \hfill$\square_{\scite{df2.3.51}}$
\enddemo
\bn
Now $\kappa_{\text{ict}}(T)$ is very close to being equal to
$\kappa_{\text{ict},2}(T)$. 
\proclaim{\stag{df2.4.8} Claim}  1) If $\kappa = \kappa_{\text{ict},2}(T)
\ne \kappa_{\text{ict}}(T)$ \ub{then}:
\mr
\item "{$(a)$}"  $\kappa_{\text{ict},2}(T) = \aleph_1 \wedge
\kappa_{\text{ict}}(T) = \aleph_0$
\sn
\item "{$(b)$}"  there is an indiscernible sequence $\langle \bar
a_t:t \in I \rangle$ with $\bar a_t \in {}^\omega{\frak C}$ and $c \in
{\frak C},I$ is dense complete for clarity such that 
{\roster
\itemitem{ $(*)$ }  for no finite $u \subseteq I$ do we have: if $J$
is a convex subset of $I$ disjoint to $u$ then 
$\langle \bar a_t:t \in J \rangle$ is 
indiscernible over $\cup\{\bar a_t:t \in I \backslash J\} \cup \{c\}$.
\endroster}
\ermn
2) If $T$ is strongly$^+$ dependent then $T$ is strongly dependent.
\nl
3) In the definition of $\kappa_{\text{ict},2}(T)$, \wilog \,\, $m=1$.
\endproclaim
\bigskip

\demo{Proof}  1) We use Observation \scite{dp1.2.2}.
Obviously $\kappa_{\text{ict}}(T) \le
\kappa_{\text{ict},2}(T)$, the rest is proved together with 
\scite{df2.2.9} below.
\nl
2) Easy.  
\nl
3) Similar to the proof of \scite{dp1.2.3} or better to use
\scite{df2.2.9}(1),(2).   \hfill$\square_{\scite{df2.4.8}}$
\enddemo
\bigskip

\proclaim{\stag{dq.7} Claim}   The following conditions on $T$ are equivalent:
\mr
\item "{$(a)$}"  $\kappa_{\text{ict},2}(T) > \aleph_0$
\sn
\item "{$(b)$}"  we can find $A$ and an indiscernible sequence
$\langle \bar a_t:t \in I\rangle$ over $A$ satisfying 
$\bar a_t \in {}^\omega{\frak C}$ and $t_n \in I$ increasing with $n$ 
and $\bar c \in {}^{\omega >}{\frak C}$ such that for every $n$
\nl

$t_n <_I t \Rightarrow \text{\rm tp}(\bar a_{t_n},A \cup \bar c \cup
\{\bar a_{t_m}:m<n\}) \ne \text{\rm tp}(\bar a_t,A \cup \bar c \cup
\{\bar a_{t_m}:m<n\})$
\sn
\item "{$(c)$}"  similarly to (b) but 
$t_n <_I t \Rightarrow \text{\rm tp}(\bar a_{t_m},A \cup \bar c \cup
\{\bar a_s:s <_I t_n\}) \ne \text{\rm tp}(a_t,A \cup \bar c \cup
\{\bar a_s:s <_I t_n\})$
\sn
\item "{$(d)$}"  we can find $A$ and a sequence $\langle \bar a^n_t:t
\in I_n \rangle,I_n$ an infinite order such that $\langle \bar a^n_t:t
\in I_n\rangle$ is indiscernible over $A \cup \cup\{\bar a^m_t:m \ne
n,m < \omega,t \in I_n\}$ and for some 
$\bar c \in {}^{\omega >}{\frak C}$ for each $n,\langle 
\bar a^n_t:t \in I_n\rangle$ is not indiscernible over $A \cup \bar c \cup
\cup\{\bar a^m_t:t \in I_m,m <n\}$
\sn
\item "{$(e)$}"  we can find a sequence $\langle \varphi_n(x,\bar
y_n,\dotsc,\bar y_0):n < \omega \rangle$ and $\langle \bar
a^n_\alpha:\alpha < \lambda,n < \omega\rangle$ such that: for every
$\eta \in {}^\omega \lambda$ the set $p_\eta = \{\varphi_n(\bar x,\bar
a^n_\alpha,\bar a^{n-1}_{\eta(n-1)},\dotsc,\bar
a^9_{\eta(0)})^{\text{if}(\alpha=\eta(n))}:n < \omega,\alpha <
\lambda\}$ is consistent.
\endroster
\endproclaim
\bigskip

\demo{Proof}  Should be clear from the proof of \scite{dp1.2.4} (more
\scite{df2.3.5}).  \hfill$\square_{\scite{dq.7}}$ 
\enddemo
\bigskip

\demo{\stag{df2.2.9} Observation}  1) For any $\kappa$ and $\zeta \ge
\kappa$ we have $(d) \Leftrightarrow (c)_\zeta \Rightarrow (b)_\zeta
\Leftrightarrow (a)$; if in addition we assume 
$\neg(\aleph_0 = \kappa_{\text{ict}}(T)
< \kappa = \aleph_1 = \kappa_{\text{ict},2}(T))$  \ub{then} we have 
also $(c)_\zeta \Leftrightarrow (b)_\zeta$ so all the
following conditions on $T$ are equivalent;
\mr
\item "{$(a)$}"  $\kappa \ge \kappa_{\text{ict}}(T)$
\sn
\item "{$(b)_\zeta$}"  if $\langle \bar a_t:t \in I\rangle$ is an
indiscernible sequence, $I$ a linear order, 
$\bar a_t \in {}^\zeta{\frak C}$ and $C \subseteq {\frak C}$
is finite \ub{then} for some set ${\Cal P}$ of $< \kappa$ initial
segments of $I$ we have:
{\roster
\itemitem{ $(*)$ }  if $s,t \in I$ and $(\forall J \in {\Cal P})(s \in
J \equiv t \in J)$ then $\bar a_s,\bar a_t$ realizes the same type
over ${\frak C}$ (if $I$ is complete this means: for some $J \subseteq
I$ of cardinality $< \kappa$, if $s,t \in I$ realizes the same
quantifier free type over $J$ in $I$ then $\bar a_s,\bar a_t$
realizes the same type over $C$)
\endroster}  
\item "{$(c)_\zeta$}"  like (b) but strengthening the conclusion to: if
$n < \omega,s_0 <_I \ldots <_I s_{n-1},t_0 <_I \ldots <_I t_n$ and
$(\forall \ell < n)(\forall k <n)(\forall J \in {\Cal P})[s_\ell \in J
= t_k \in J]$ then $\bar a_{s_0} \char 94 \ldots \char 94 \bar
a_{t_{n-1}}$ and $\bar a_{t_0} \char 94 \ldots \char 94 \bar a_{t_{n-1}}$
realize the same type over $C$
\sn
\item "{$(d)$}"   $\kappa \ge \kappa_{\text{ict},2}(T)$.
\ermn
2) We can in clause $(b)_\zeta,(c)_\zeta$ add $|C|=1$ and/or demand
$I$ is well ordered (for the last use \scite{dq.6}).
\enddemo
\bigskip

\demo{Proof}  We shall prove various implications which together
obviously suffice (for \scite{df2.2.9} and \scite{df2.4.8}(1) and
\scite{df2.4.8}(3)).   
\sn
\ub{$\neg(a) \Rightarrow \neg(b)_\zeta$}

Let $\lambda \ge \kappa$.  As in the proof of \scite{dp1.2.2} there
are $\bar\varphi = \langle \varphi_i(\bar x,\bar y_i):i < \kappa
\rangle,m = \ell g(\bar x)$ and $\langle \bar a^i_\alpha:i <
\kappa,\alpha < \lambda \rangle$ exemplifying
$\circledast^2_{\bar\varphi}$ from \scite{dp1.2.2}, so necessarily 
$\bar a^\ell_\alpha$ is non-empty.  Without loss of generality $\ell
g(\bar a^0_\alpha) \le \omega$ and \wilog \, $\zeta \ge \omega^2$.
Let $\bar a^*_\alpha \in {}^\zeta{\frak C}$ be $\bar a^0_\alpha
\char 94 \bar a^1_\alpha \char 94 \ldots \char 94 \bar a'_\alpha$
were $\bar a'_\alpha$ has length $\zeta - \dsize \Sigma_{\ell <
\kappa} \ell g(\bar a^i_\alpha)$ and is constantly the first member of
$\bar a^0_\alpha$.  Let $\bar c$ realize $p = \{\varphi_i(\bar x,\bar
a_{2i}) \wedge \neg \varphi_i(\bar x,\bar a_{2i+1}):i < \kappa\}$.

Easily $\bar c$ (or pedantically Rang$(\bar c))$ and $\langle \bar
a^*_\alpha:\alpha < \lambda \rangle$ exemplifies $\neg(b)_\zeta$.
\bn
\ub{$(a) \Rightarrow (b)_\zeta$}.

If $\kappa = \aleph_0$, this holds by \scite{dp1.2.4}(1); in general,
this holds by the proof of \scite{dp1.2.4}(1) and this is why there we
use $\kappa$.
\bn
\ub{$\neg(b)_\zeta \Rightarrow \neg(c)_\zeta$}

Obvious.
\bn
\ub{$\neg(a) \Rightarrow \neg(d)$}

The witness for $\neg(a)$ is a witness for $\neg(d)$.
\bn
\ub{$\neg(d) \Rightarrow \neg(c)_\zeta$}

Let $\langle \varphi_i(\bar x,\bar y_i):i < \kappa \rangle$ witness
$\neg(d)$, i.e., witness $\kappa < \kappa_{\text{ict},2}(T)$, so there
are $\langle \bar a_{i,\alpha}:\alpha < \lambda,i < \kappa \rangle$
and $\langle \bar b_i:i < \kappa \rangle$ satisfying clauses
(a),(b),(c) of Definition \scite{df2.3.5}.  By Observation \scite{dq.6}
we can find an indiscernible sequence $\langle \bar a^*_\alpha:\alpha
< \lambda \times \kappa\rangle,\ell g(\bar a^*_\alpha) = \zeta_\kappa$
where $\zeta_j := \Sigma\{\ell g(\bar y_i):i<j\}$ such that $i < \kappa
\wedge \alpha < \lambda \Rightarrow \bar a^*_i \restriction
[\zeta_i,\zeta_{i+1}) = \bar a^i_\alpha$.  Now $\langle \bar
a^*_\alpha:\alpha < \lambda \times \kappa \rangle,\bar c$ witness
$\neg(c)_{\zeta_\kappa}$, because if ${\Cal P}$ is as required in
$(c)_{\zeta_\kappa}$ then easily $(\forall i < \kappa)(\exists J
\in {\Cal P})(J \cap [\lambda i,\lambda i + \lambda] \notin
\{\emptyset,[\lambda i,\lambda i + \lambda\}$, hence $|{\Cal P}| \ge
\kappa$.  Now clearly $\zeta_\kappa \le \zeta$ hence repeating the first
element $(\zeta - \kappa)$ times we get $\langle \bar b^i_\alpha:\alpha
< \lambda \kappa \rangle$ which together with $\bar c$ exemplify
$\neg(c)_\zeta$.
\medskip

It is enough to prove
\mr
\item "{$(*)$}"  assume $\neg(c)_\zeta$ then 
{\roster
\itemitem{ $(i)$ }  $\neg(d)$
\sn
\itemitem{ $(ii)$ }   $\neg(a)$ except possibly when (a) + (b) of
\scite{df2.4.8}(1) holds, in particular $\aleph_0 =
\kappa_{\text{ict}}(T) < \kappa = \aleph_1 = \kappa_{\text{ict},2}(T)$.
\endroster}
\ermn
Toward this we can assume that
\mr
\item "{$\boxtimes$}"  $T$ is dependent and $C,\langle \bar a_t:t \in
I\rangle$ form a witness to $\neg(c)_\zeta$.
\ermn
Let $\bar c$ list $C$ without repetitions and \wilog \,
$I$ is a dense complete linear order (so with no extreme elements).
Let $\ell g(\bar x_\ell) = \zeta$ for $\ell < \omega$ be pairwise
disjoint with no repetitions, of course, $\ell g(\bar y) =
\ell g(\bar c) < \omega$ (pairwise disjoint) and let $\bar\varphi = 
\langle \varphi_i = 
\varphi_i(\bar x_0,\dotsc,\bar x_{n(i)-1},\bar y):i < |T|\rangle$
list all such formulas in $\Bbb L(\tau_T)$.  For each $i < |T|$ by
\scite{df1.2.5}(1),(2) there are $m(i) < \omega$ and $t_{i,1} <_I \ldots <_I
t_{i,m(i)-1}$ as there and $m(i)$ is minimal, so 
stipulating $t_{i,0} = - \infty,t_{i,m(i)} = \infty$ we have:
\mr
\item "{$(*)_1$}"  if $s'_0 <_I \ldots <_I s'_{m(i)-1}$ and 
$s''_0 <_I \ldots <_I s''_{m(i)-1}$ and $s'_\ell,s''_\ell$ realize 
the same quantifier free
type over $\{t_{i,1},\ldots,t_{i,m(i)-1}\}$ in the linear order $I$
for each $\ell < m(i)$ \ub{then}
${\frak C} \models ``\varphi_i[\bar a_{s'_0},\dotsc,
\bar a_{s'_{m(i)-1}},\bar c] \equiv
\varphi_i[\bar a_{s''_0},\dotsc,\bar a_{s''_{m(i)-1}} \bar c]"$.
\ermn
For each $i < |T|$, for each $\ell \in \{1,\dotsc,m(i)\}$ we can
choose $\omega_{\ell,i}$ and find $w_{i,\ell}$ such that
\mr
\item "{$(*)_2$}"  $(a) \quad w_{i,\ell} \subseteq I \backslash
\{t_{i,\ell}\}$
\sn
\item "{${{}}$}"  $(b) \quad w_{i,\ell}$ is finite
\sn
\item "{${{}}$}"  $(c) \quad$ if $s_1 < t_{i,\ell(i)} < s_2$ then 
$\langle \bar a_t:t \in (s_1,s_2)_I\rangle$ is not
$\{\varphi_i\}$-indiscernible
\nl

\hskip25pt  over $C \cup \{\bar a_t:t \in w_i\}$.
\ermn
If the set $\{t_{i,k}:i < |T|,k=1,\dotsc,m(i)-1\}$ 
has cardinality $< \kappa$ we are done, so assume that
\mr
\item "{$(*)_3$}"  $\{t_{i,\ell}:i < |T|$ and $\ell \in [1,m(i)]\}$
has cardinality $\ge \kappa$.
\ermn
\ub{Case 1}:  $\kappa > \aleph_0$ (so we have to prove $\neg(a)$).

By Hajnal free subset theorem and by $(*)_3$ there is 
$u_0 \subseteq |T|$ of order type $\kappa$ such
that $i \in u_0 \Rightarrow \{t_{i,\ell}:\ell=1,\dotsc,m(i)-1\}
\nsubseteq \{t_{j,\ell}:j \in u_0 \cap i$ and $\ell=1,\dotsc,m(j)-1\}$.

There are $u \subseteq u_\ell$ of
cardinality $\kappa$ and a sequence 
$\langle \ell(i):i \in u \rangle,0 < \ell(i)
<m(i)$ such that $\langle t_{i,\ell(i)}:i \in u \rangle$ is with no
repetitions and disjoint to $\{t_{i,\ell}:i \in u$ and $\ell \ne
\ell(i)\} \cup \bigcup\{w_{i,\ell(i)}:i \in u\}$.  We shall now prove
$\kappa < \kappa_{\text{ict}}(T)$, this gives $\neg(a),\neg(d)$ so it
suffices.   

Clearly by \scite{dp1.2.2} it suffices to show
($\lambda$ any cardinality $\ge \aleph_0$)
\mr
\item "{$\boxdot_u$}"  there are $\bar a^i_\alpha \in {}^\zeta{\frak
C}$ for $i \in u,\alpha < \lambda$ and set $A$ such that
{\roster
\itemitem{ $(a)$ }  $\langle \bar a^i_\alpha:\alpha < \lambda \rangle$
is an indiscernible sequence over $\cup\{\bar a^j_\beta:j \in u,j \ne
i,\alpha < \lambda\} \cup A$
\sn
\itemitem{ $(b)$ }  $\langle \bar a^i_\alpha:\alpha < \lambda\rangle$
is not $\{\varphi_i\}$-indiscernible over $A \cup \bar c$.
\endroster}
\ermn
By compactness it suffices to prove $\boxdot_v$ for any finite $v
\subseteq u$ and $\lambda = \aleph_0$; also we can replace $\lambda$ by
any infinite linear order.  

We can find $\langle(s_{1,i},s_{2,i}):i \in v\rangle$ such that
\mr
\item "{$(*)_4$}"  $s_{1,i} <_I t_{i,\ell(i)} <_I s_{2,i}$ (for $i \in
v$)
\sn
\item "{$(*)_5$}"  $(s_{1,i},s_{2,i})_I$ is disjoint to
$\cup\{(s_{1,j},s_{2,j}):j \in v \backslash \{i\}\} \cup
\bigcup\{w_{j,\ell(j)} \in v\}$.
\ermn
So $\left< \langle a^j_t:t \in (s_{1,j},s_{2,j})_I\rangle:j \in v
\right>$ and choosing 
$A = \cup\{\bar a_t:t \in w_{i,\ell(i)},i \in v\}$ are as required above.
So we are done.    
\enddemo
\bn
\ub{Case 2}:  $\kappa = \aleph_0$ so we have to prove $\neg(d)$ and
clause (ii) of $(*)$ and (for proving part (2) of the present \scite{df2.2.9})
that \wilog \, $|C| = 1$.

We can find $A$ and $u$
\mr
\item "{$\boxdot^1$}"  $(a) \quad A \subseteq C$
\sn
\item "{${{}}$}"  $(b) \quad u \subseteq I$ is finite
\sn
\item "{${{}}$}"  $(c) \quad$ if $n < \omega$ and $t^\ell_0 <_I \ldots
<_I t^\ell_{n-1}$ for $\ell=1,2$ and 
\nl

\hskip25pt $(\forall k<n)(\forall s \in
u)(t^1_k=s \equiv t^2_k = s \wedge t^1_k <_I s \equiv t^2_k <_I s)$ then 
\nl

\hskip25pt $\bar a_{t^1_0} \char 94 \ldots \char 94 \bar a_{t^1_{n-1}},\bar
a_{t^2_0} \char 94 \ldots \char 94 \bar a_{t^2_{n-1}}$ realize the
same type over $A$
\sn
\item "{${{}}$}"  $(d) \quad$ if $A',u'$ satisfies (a)+(b)+(c) then
$|A'| \le |A|$.
\ermn
This is possible because $C$ is finite and the empty set satisfies
 clauses (a),(b),(c) for $A$ by our present assumption
$A \ne C$, so let $c \in C \backslash A$.  Now we try to
choose $(i_k,\ell_k,w_k)$ by induction on $k < \omega$
\mr
\item "{$\circledast$}"  $(a) \quad i_k < \kappa$
\sn
\item "{${{}}$}"  $(b) \quad 1 \le \ell_k \le m(i)-1$
\sn
\item "{${{}}$}"  $(c) \quad t_{i_k,\ell_k} \in I \backslash \omega$
\sn
\item "{${{}}$}"  $(d) \quad w_k \supseteq u \cup w_0 \cup \ldots
\cup w_{k-1} \cup \{t_{i_0,k_0},\dotsc,t_{i_{k-1},\ell_{k-1}}\}$
\sn
\item "{${{}}$}"  $(e) \quad w_k \subseteq I \backslash
\{t_{i_k,\ell_k}\}$ is finite
\sn
\item "{${{}}$}"  $(f) \quad$ if $s' <_I t_{i_k,\ell_k} <_I s''$ then
$\langle \bar a_t:t \in (s',s'')_I\rangle$ is not indiscernible over
\nl

\hskip25pt $\{\bar a_s:s \in w_k\} \cup \bar c$.
\ermn
If we are stuck in $k$ then $w_{k-1} \in [I]^{< \aleph_0}$ when
$k>0,u$ when $k=0$ show that $\langle \bar a_t:t \in I\rangle,A \cup
\{c\}$ contradict the choice of $A$ and so 
witness $\neg(c)_\zeta$.  If we succeed then we prove as in Case 1
that $\kappa_{\text{ict}}(\text{Th}({\frak C},a)_{a \in A}) >
\aleph_0$ so by \scite{dp1.2.1} we get $\kappa_{\text{ict}}(T) >
\aleph_0$.
\nl
${{}}$    \hfill$\square_{\scite{df2.2.9}}$
\bigskip

\demo{\stag{df2.2.35} Conclusion}  $T$ is strongly$^2$ dependent by
Definition \scite{df2.3.5} iff $T$ is strongly$^2$ dependent by
\cite[\S3,3.7]{Sh:783} which means
we say $T$ is strongly$^2$ (or strongly$^+$) dependent when:
if $\langle \bar{\bold a}_t:t \in I\rangle$ is an indiscernible
sequence over $A$, $t \in I \Rightarrow \ell g(\bar{\bold a}_t) = \alpha$
and $\bar b \in {}^{\omega >}({\frak C})$ \ub{then} we
can divide $I$ to finitely many convex sets $\langle I_\ell:\ell <
k\rangle$ such that for each $\ell$ the sequence $\langle \bar{\bold
a}_t:t \in I_\ell \rangle$ is an indiscernible sequence over 
$\{\bar a_s:s \in I \backslash I_\ell\} \cup A \cup \bar b$.
\enddemo
\bn
\centerline{$* \qquad * \qquad *$}
\bn
\ub{Discussion}:  Now we define ``$T$ is strongly$^{2,*}$ dependent",
parallely to \scite{dp1.8}, \scite{dp1.9} from the end of \S1.
\bigskip

\definition{\stag{dq2.8.4} Definition}  1) $\kappa_{\text{icu},2}(T)$ is
the minimal $\kappa$ such that for no $m < \omega$ and $\bar\varphi =
\langle \varphi_i(\bar x_i,\bar y_i):i < \kappa \rangle$ with $\ell
g(\bar x^i) = m \times n_i$ can we find $\bar a^i_\alpha \in {}^{\ell
g(\bar y_i)}{\frak C}$ for $\alpha < \lambda,i < \kappa$ and $\bar
c_{\eta,n} \in {}^m {\frak C}$ for $\eta \in {}^\kappa \lambda$ such that:
\mr
\item "{$(a)$}"  $\langle \bar c_{\eta,n}:n < \omega \rangle$ is an
indiscernible sequence over $\cup\{\bar a^i_\alpha:\alpha < \lambda,i
< \kappa\}$
\sn
\item "{$(b)$}"  for each $\eta \in {}^\kappa \lambda$ and $i <
\kappa$ we have ${\frak C} \models \varphi_i(\bar c_{\eta,0} \char 94
\ldots \bar c_{\eta,n_i-1},\bar a^i_\alpha)^{\text{if}(\alpha=\eta(i))}$.
\ermn
2) If $\bar \varphi$ is as in (1) then we say that it witnesses
$\kappa < \kappa_{\text{icu},2}(T)$.
\nl
3) $T$ is strongly$^{1,*}$ dependent if $\kappa_{\text{jcu}}(T) =
 \aleph_0$.   
\enddefinition
\bigskip

\proclaim{\stag{dq.2.9.4} Claim}  1) $\kappa_{\text{icu},2}(T) \le
\kappa_{\text{ict},2}(T)$. 
\nl
2) If {\rm cf}$(\kappa) > \aleph_0$ then $\kappa_{\text{icu},2}(T) >
\kappa \Leftrightarrow \kappa_{\text{ict},2}(T) > \kappa$.
\nl
3) The parallel of \scite{dp1.2.1}, \scite{dp1.2.2}, \scite{dp1.2.3}(2) holds.
\endproclaim
\newpage

\head {\S3 Ranks \\
\S (3A) Rank for strongly dependent $T$} \endhead  \resetall 
 \spuriousreset
\bn
\margintag{df1.0}\ub{\stag{df1.0} Explanation/Thesis}:  
\mr
\item "{$(a)$}"  For stable theories we normally consider not just a
model $M$ (and say a type in it, but all its elementary extensions; we
analyze them together
\sn
\item "{$(b)$}"   for dependent theories we should be more liberal,
allowing to replace $M$ by $N^{[a]}$ when $M \prec N \prec N_1,\bar a
\in {}^{\ell g(\bar a)}(N_1)$, ($N^{[\bar a]}$ is the sum of $N$ by
restrictions of relation in $N_1$ definable with parameters from $\bar a$
\sn
\item "{$(c)$}"  this motivates some of the ranks below).
\ermn
Such ranks relate to strongly$^1$ dependent, they have relatives for
strongly$^2$ dependent.

Note that we can represent the ${\frak x} \in K'_{\ell,m}$ (and ranks)
close to \cite[\S1]{Sh:783} particularly $\ell = 9$.
\bigskip

\definition{\stag{dp1.1} Definition}  1) Let $M_0 \le_A M_1$ for
$M_0,M_1 \prec {\frak C}$ and $A \subseteq {\frak C}$ means that:
\mr
\item "{$(a)$}"  $M_0 \subseteq M_1$ (equivalently $M_0 \prec M_1$)
\sn
\item "{$(b)$}"  for every $\bar b \in M_1$, the type tp$(\bar b,M_0
\cup A)$ is f.g. (= finitely satisfiable) in $M_0$.
\ermn
2) Let $M_0 \le_{A,p} M_1$ for $M_0,M_1 \prec {\frak C},A \subseteq
{\frak C}$ and $p \in \bold S^{< \omega}(M_1 \cup A)$ or just $p$ is
a $(< \omega)$-type over $M_1 \cup A$ means that
\mr
\item "{$(a)$}"  $M_0 \subseteq M_1$
\sn
\item "{$(b)$}"  if  $\bar b \in M_1,\bar c \in M_0,\bar a_1 \in A,\bar a_2 \in
A,{\frak C} \models \varphi_1[\bar b,\bar a_1,\bar c]$ and 
$\varphi_2(\bar x,\bar b,\bar a_2,\bar c) \in p$ 
or is just a (finite) conjunction of members of $p$ (e.g. empty)
\ub{then} for some $\bar b' \in M_0$ we have ${\frak C} \models
\varphi_1[\bar b'_1,\bar a_1,\bar c]$ and 
$\varphi_2(\bar x,\bar b',\bar a_2,\bar c) \in p$ or just is a finite
conjunction of members of $p$.
\endroster
\enddefinition
\bigskip

\demo{\stag{dp1.1A} Observation}  1) $M_0 \le_{A,p} M_1$ 
implies $M_0 \le_A M_1$. 
\nl
2) If $p = \text{ tp}(\bar b,M_1 \cup A) \in \bold S^m(M_1 \cup A)$
\ub{then} $M_0 \le_{A,p} M_1$ iff $M_1 \le_{A \cup \bar b} M_2$.
\nl
3) If $M_0 \le_A M_1 \le_A M_2$ then $M_0 \le_A M_2$.
\nl
4) If $M_0 \le_{A,p \restriction (M_1 \cup A)} M_1 \le_{A,p} M_2$ then $M_0
\le_{A,p} M_2$. \nl
5) If the sequences $\langle M_{1,\alpha}:\alpha \le \delta
\rangle,\langle A_\alpha:\alpha \le \delta\rangle$ are increasing
continuous, $\delta$ a limit ordinal and $M_0 \le_{A_\alpha}
M_{1,\alpha}$ for $\alpha < \delta$ then $M_0 \le_{A_\delta}
M_{1,\delta}$.  Similarly using $<_{A_\alpha,p_\alpha}$.
\nl
6) If $M_1 \subseteq M_2$ and $p$ is an $m$-type over $M_1 \cup A$
then $M_1 \le_A M_2 \Leftrightarrow M_1 \le_{A,p} M_2$.
\enddemo
\bigskip

\demo{Proof}  Easy.
\enddemo
\bn
\margintag{dp2.17}\ub{\stag{dp2.17} Discussion}:  1) Note that the 
ranks defined below are related to \cite[\S1]{Sh:783}.  An alternative
presentation (for $\ell \in \{3,6,9,12\}$) is that we define $M_A$ as
$(M,a)_{a \in A}$ and $T_A = \text{ Th}({\frak C},a)_{a \in A}$ and we
consider $p \in \bold S(M_A)$ and in the definition of ranks to extend
$A$ and $p$ we use appropriate $q \in \bold S(N_B),M_A \prec N_A,A
\subseteq B$.  Originally we prsent here many variants, but now we
present only two $(\ell = 8,9)$, retaining the others in \S (5A).   
\nl
2) We may change the definition, each time retaining from $p$ 
only one formula with little change in the claims.
\nl
3) We can define ${\frak x} \in K_{\ell,m}$ such that it has also
$N^{\frak x}$ where $M^{\frak x} \subseteq N^{\frak x}
(\prec {\frak C}_T)$ and:
\mr
\item "{$(A)$}"   change the definition of 
${\frak x} \le^\ell_{\text{at}} {\frak y}$ to:
{\roster
\itemitem{ $(a)$ }  $N^{\frak y} \subseteq N^{\frak x}$
\sn
\itemitem{ $(b)$ }  $A^{\frak x} \subseteq A^{\frak y} \subseteq
A^{\frak x} \cup N^{\frak x}$
\sn
\itemitem{ $(c)$ }   $M^{\frak x} \subseteq M^{\frak y} \subseteq
N^{\frak x}$
\sn
\itemitem{ $(d)$ }  $p^{\frak y} \subseteq p^{\frak x}$
\endroster}
\sn
\item "{$(B)$}"  change ``${\frak y}$ explicitly $\bar\Delta$-split
$\ell$-strongly over ${\frak x}$" according to and replacing in
$(e),(e)' \, p^{{\frak x}'}$ by $p^{\frak x}$
\sn
\item "{$(C)$}"  dp-rk$^m_{\bar\Delta,\ell}$ is changed accordingly.
\ermn
So now dp-rk$^m_{\bar\Delta}$ may be any ordinal so
\scite{dp1.4} may fail, but the result in \S3 becomes stronger
covering also some model of non-strongly dependent.
\bigskip

\definition{\stag{dp1.3} Definition}  1) For $\ell=8,9$ let 

$$
\align
K_{m,\ell} = \bigl\{{\frak x}:&{\frak x} = (p,M,A),M \text{ a model } \prec 
{\frak C}_T,A \subseteq {\frak C}_T, \\
  &p \in \bold S^m(M \cup A) \text{ and if } \ell=9 \text{ then} \\
  &p \text{ is finitely satisfiable in } M\bigr\}.
\endalign
$$
\mn
If $m=1$ we may omit it.

For ${\frak x} \in K_{m,\ell}$ let 
${\frak x} = (p^{\frak x},M^{\frak x},A^{\frak x}) = (p[{\frak
x}],M[{\frak x}],A[{\frak x}])$ and $m = m({\frak x})$ recalling 
$p^{\frak x}$ is an $m$-type.
\nl
2) For ${\frak x} \in K_{m,\ell}$ let $N_{\frak x}$ be $M^{\frak x}$
expanded by $R_{\varphi(\bar x,\bar y,\bar a)} = \{\bar b \in {}^{\ell g(\bar
y)}M:\varphi(\bar x,\bar b,\bar a) \in p\}$ for $\varphi(\bar x,\bar y,\bar z)
\in \Bbb L(\tau_T),\bar a \in {}^{\ell g(\bar z)}A$ and 
$R_{\varphi(\bar y,\bar a)} = \{\bar b \in {}^{\ell g(\bar y)}M:{\frak
C} \models \varphi[\bar b,\bar a]\}$ for $\varphi(\bar y,\bar z) \in
\Bbb L(\tau_T),\bar a \in {}^{\ell g(\bar y)}{\frak C}$; let $\tau_{\frak x} =
\tau_{N_{\frak x}}$.  
\nl
2A) In part (1) and (2): if 
we omit $p$ we mean $p = \text{ tp}(<>,M \cup
A)$ so we can write $N_{A}$, a $\tau_A$-model so in this case 
$p = \{\varphi(\bar b,\bar a):
\bar b \in M,\bar a \in M$ and ${\frak C} \models \varphi[\bar b,\bar a]\}$.
\nl
3) For ${\frak x},{\frak y} \in K_{m,\ell}$ let
\mr
\item "{$(\alpha)$}"  ${\frak x} \le^\ell_{\text{pr}} {\frak y}$ means
that ${\frak x},{\frak y} \in K_{m,\ell}$ and
{\roster
\itemitem{ $(a)$ }  $A^{\frak x} = A^{\frak y}$
\sn 
\itemitem{ $(b)$ }  $M^{\frak x} \le_{A[{\frak x}]} M^{\frak y}$
\sn
\itemitem{ $(c)$ }  $p^{\frak x} \subseteq p^{\frak y}$
\sn
\itemitem{ $(d)$ }  $M^{\frak x} \le_{A[{\frak x}],
p[{\frak y}]} M^{\frak y}$ 
\endroster}
\sn
\item "{$(\beta)$}"  ${\frak x} \le^\ell {\frak y}$ means that 
for some $n$ and 
$\langle {\frak x}_k:k \le n \rangle,{\frak x}_k \le^\ell_{\text{at}}
{\frak x}_{k+1}$ for $k < n$ and $({\frak x},{\frak y}) = ({\frak x}_0,{\frak
x}_n)$ 
\nl
where
\sn
\item "{$(\gamma)$}"  ${\frak x} \le^\ell_{\text{at}} {\frak y}$ iff
$({\frak x},{\frak y} \in K_{m,\ell}$ and) for some ${\frak x}' \in
K_{m,\ell}$ we have
{\roster
\itemitem{ $(a)$ }  ${\frak x} \le^\ell_{\text{pr}} {\frak x}'$
\sn
\itemitem{ $(b)$ }   $A^{\frak x} \subseteq A^{\frak y} \subseteq
A^{\frak x} \cup M^{{\frak x}'}$
\sn
\itemitem{ $(c)$ }  $M^{\frak y} \subseteq M^{{\frak x}'}$
\sn
\itemitem{ $(d)$ }  $p^{\frak y} \supseteq p^{{\frak x}'}
\restriction (A^{\frak x} \cup M^{\frak y})$ so $\ell \in \{1,4\}
\Rightarrow p^{\frak y} = p^{{\frak x}'} \restriction M^{\frak y}$ and
$p^{\frak y} = p^{{\frak x}'}
\restriction (M^{\frak y} \cup A^{\frak y})$.
\endroster}
\ermn
4) For ${\frak x},{\frak y} \in K_{m,\ell}$ we say that ${\frak y}$ 
explicitly $\bar\Delta$-splits $\ell$-strongly over ${\frak x}$
\ub{when}: $\bar\Delta = (\Delta_1,\Delta_2),\Delta_1,\Delta_2 \subseteq
\Bbb L(\tau_T)$ and for some ${\frak x}'$ and 
$\varphi(\bar x,\bar y) \in \Delta_2$ we have clauses 
(a),(b),(c),(d) of part (3)$(\gamma)$ and
\mr
\item "{$(e)$}"   there are $\bar{\bold a}$ such that
{\roster
\itemitem{ $(\alpha)$ }  $\bar{\bold a} = \langle \bar a_i:i < \omega
+ 1 \rangle$ is $\Delta_1$-indiscernible over $A^{\frak x} \cup M^{\frak y}$
\sn
\itemitem{ $(\beta)$ }  $A^{\frak y} \backslash A^{\frak x} = \cup
\{\bar a_i:i < \omega\}$; yes $\omega$ not $\omega + 1$! 
(note that $``A^{\frak y} \backslash A^{\frak x} ="$ and 
not ``$A^{\frak y} \backslash A^{\frak x} \supseteq"$ as
we use it in $(e)(\gamma)$ in the proof of \scite{dp1.4})
\sn
\itemitem{ $(\gamma)$ }  $\bar a_i \in M^{{\frak x}'}$ for $i < \omega
+1$ and $\bar b \in {}^{\omega >}(A^{\frak x})$
\sn
\itemitem{ $(\delta)$ }  $\varphi(\bar x,\bar a_k \char 94 \bar b) 
\wedge \neg \varphi(\bar x,\bar a_\omega \char 94 \bar b)$ belongs
to $p^{{\frak x}'}$ for $k < \omega$.
\endroster}
\ermn
5) We define dp-rk$^m_{\bar\Delta,\ell}:K_{m,\ell} \rightarrow \text{ Ord }
 \cup\{\infty\}$ by 
\mr
\item "{$(a)$}"  $\text{dp-rk}^m_{\bar\Delta,\ell}({\frak x}) 
\ge 0 \text{ always}$
\sn
\item "{$(b)$}"  $\text{dp-rk}^m_{\bar\Delta,\ell}({\frak x}) \ge 
\alpha +1 \text{ \ub{iff} there is } {\frak y} \in K_{m,\ell} 
\text{ which explicitly } \bar\Delta$-splits $\ell$-strongly over ${\frak x} 
\text{ and dp-rk}_{\bar\Delta,\ell}({\frak y}) \ge \alpha$
\sn
\item "{$(c)$}"   $\text{dp-rk}^m_{\bar\Delta,\ell}({\frak x}) \ge \delta 
\text{ iff dp-rk}^m_{\bar\Delta,\ell}({\frak x})
\ge \alpha \text{ for every } \alpha < \delta$ when $\delta$ is
a limit ordinal.
\ermn
Clearly well defined.  We may omit $m$ from dp-rk as ${\frak x}$
determines it.
\nl
6) Let dp-rk$^m_{\bar\Delta,\ell}(T) = 
\cup\{\text{dp-rk}_{\bar\Delta,\ell}({\frak x}):{\frak x} \in
K_{m,\ell}\}$; if $m=1$ we may omit it.
\nl
7) If $\Delta_1 = \Delta_2 = \Delta$ we may write $\Delta$ instead of
$(\Delta_1,\Delta_2)$. 
If $\Delta = \Bbb L(\tau_T)$ then we may omit it.
\enddefinition
\bigskip

\remark{Remark}  There are obvious monotonicity and inequalities.
\endremark
\bigskip

\demo{\stag{dp1.3A} Observation}  1) $\le^\ell_{\text{pr}}$ is a
partial order on $K_\ell$.
\nl
2) $K_{m,9} \subseteq K_{m,8}$.
\nl
3) ${\frak x} \le^8_{\text{pr}} {\frak y} \Leftrightarrow 
{\frak x} \le^9_{\text{pr}} {\frak y}$.
\nl
4) ${\frak x} \le^8_{\text{at}} {\frak y} \Leftrightarrow 
{\frak x} \le^9_{\text{at}} {\frak y}$.
\nl
5) ${\frak y}$ explicitly $\bar\Delta$-splits $8$-strongly over
${\frak x}$ iff ${\frak y}$ explicitly $\bar\Delta$-splits
$9$-strongly over ${\frak x}$.
\nl
6) If ${\frak x} \in K_{m,9}$ then 
dp-rk$^m_{\bar\Delta,9}({\frak x})
\le \text{\rm dp-rk}^m_{\bar\Delta,8}({\frak x})$.
\nl
7) If $\bar a \in {}^m {\frak C}$ and ${\frak y} = (\text{\rm tp}(\bar
a,M \cup A),M,A)$ and ${\frak x} = (\text{\rm tp}(\bar a,M \cup
A),M,A)$ \ub{then} ${\frak x} \in K_{m,8}$.
\nl
8) In part (7) if tp$(\bar a,M \cup A)$ is finitely satisfiable in
$M$ then also ${\frak y} \in K_{m,9}$.
\nl
9) If ${\frak x} \in K_{m,\ell}$ and $\kappa > \aleph_0$ then there is
${\frak y} \in K_{m,\ell}$ such that ${\frak x} \le^\ell_{\text{pr}}
{\frak y}$ and $M^{\frak y}$ is $\kappa$-saturated, moreover $M^{\frak
y}_{A[{\frak y}],p[{\frak y}]}$ is $\kappa$-saturated (hence in
Definition \scite{dp1.1}(4) \wilog \, $M^{{\frak x}'}$ is 
$(|M^{\frak x} \cup A^{\frak x}|^+)$-saturated).
\enddemo
\bigskip

\demo{Proof}  Easy.
\enddemo
\bigskip

\proclaim{\stag{dp1.4} Claim}  1) For each $\ell=8,9$ we have
{\rm dp-rk}$_\ell(T) = \infty$ iff {\rm dp-rk}$_\ell(T) 
\ge |T|^+$ iff $\kappa_{\text{ict}}(T) > \aleph_0$.
\nl
2) For each $m \in [1,\omega)$, similarly using 
{\rm dp-rk}$^m_\ell(T)$, hence the properties do not depend on such $m$.
\endproclaim
\bigskip

\remark{\stag{dp1.4.3} Remark}  In the implications in the proof we
allow more cases of $\ell$.
\endremark
\bigskip

\demo{Proof}  Part (2) has the same proof as part (1) when we recall
\scite{dp1.2.3}.
\mn
\ub{$\kappa_{\text{ict}}(T) > \aleph_0$ implies
dp-rk$_\ell(T)=\infty$}

By the assumption there is a sequence 
$\bar \varphi = \langle \varphi_n(x,\bar y_n):n < \omega
\rangle$ exemplifying $\aleph_0 < \kappa_{\text{ict}}(T)$.  
Let $\lambda > \aleph_0$ and $I$ be $\lambda \times \Bbb Z$ ordered
lexicographically and let $I_\alpha =\{\alpha\} \times \Bbb Z$ and
$I_{\ge \alpha} = [\alpha,\lambda) \times \Bbb Z$.  
As in \scite{dp1.2.2} by 
Ramsey theorem and compactness we can
find $\langle \bar a^n_t:t \in I \rangle$ (in ${\frak C}_T$) such that
\mr
\item "{$\circledast$}"  $(a) \quad \ell g(\bar a^n_t) = \ell
g(\bar y_n)$
\sn
\item "{${{}}$}"  $(b) \quad \langle \bar a^n_t:t \in I
\rangle$ is an indiscernible sequence over $\cup\{\bar a^m_t:m <
\omega,m \ne n$ \nl

\hskip25pt and $t \in I\}$
\sn
\item "{${{}}$}"  $(c) \quad$ for every $\eta \in {}^\omega I,
p_\eta = \{\varphi_n(x,\bar a^n_t)^{\text{if}(\eta(n)=t)}:
n < \omega,t \in I\}$ is \nl

\hskip25pt consistent (i.e., finitely satisfiable in ${\frak C}$).
\ermn
Choose a complete $T_1 \supseteq T$ with Skolem functions
and $M^* \models T_1$ expanding ${\frak C}$ be such that 
in it $\langle \bar a^n_\alpha:t \in I,n <
\omega \rangle$ satisfies $\circledast$ also in $M^*$; exists by
Ramsey theorem.  Let $M^*_n$ be the
Skolem hull in M$^*$ of $\cup\{\bar a^m_t:m<n,t \in I_1\} \cup
\{\bar a^m_t:m \in [n,\omega)$ and $t \in I\}$ and let 
$M_n = M^*_n \restriction \tau(T)$.
So we have $M_n \prec {\frak C}$ which includes 
$\{\bar a^m_t:t \in I,m \in [n,\omega)\}$ such that $M_{n+1}
\prec M_n$ and $\langle \bar a^n_t:t \in I_{\ge 2}\rangle$ is
an indiscernible sequence over $M_{n+1} \cup \{\bar a^m_t:m <n,t
\in I\}$ hence $\langle a^n_t:t \in I_2\rangle$ is an indiscernible
sequence over $M_{n+1} \cup A_m$; the indiscernibility holds even in
$M^*$ where $A_n = \{\bar a^m_t:m<n$ and $t \in I_1\}$.
We delay the case $\ell=9$.  Let 
$\eta \in {}^\omega I$ be chosen as: $\langle (2,i):i < \omega \rangle$.
Let $p \in \bold S(M_0)$ be such that it includes $p_\eta$.  

Lastly, let ${\frak x}_n  = {\frak x}'_n = (p_n,M_n,A_n)$ 
where $p_n = p \restriction (A_n \cup M_n)$.  By
\scite{dp1.3A}(7) clearly ${\frak x}_n \in K_\ell$.

It is enough to show that dp-rk$_\ell({\frak x}_n) < \infty
\Rightarrow \text{ dp-rk}_\ell({\frak x}_n) > 
\text{ dp-rk}_\ell({\frak x}_{n+1})$ as by the ordinals being well ordered this
implies that dp-rk$_\ell({\frak x}_n)=\infty$ for every $n$.  By
Definition \scite{dp1.3}(5) clause (b) it is enough to show 
(fixing $n < \omega$) that ${\frak x}_{n+1}$ explicitly split $\ell$-strongly 
over ${\frak x}_n$ using $\langle \bar a^n_{(1,i)}:i < \omega \rangle
\char 94 \langle \bar a^n_{(2,n)}\rangle$.
To show this, see Definition \scite{dp1.3}(4) we use ${\frak
x}'_n := {\frak x}_n$, 
clearly ${\frak x}_n \le^\ell_{\text{pr}} {\frak x}'_n$ as
${\frak x}_n = {\frak x}'_n \in K_\ell$ so clause (a) of Definition
\scite{dp1.3}(3)$(\gamma)$ holds.  Also $A^{{\frak x}_n} \subseteq
A^{{\frak x}_{n+1}} \subseteq A^{{\frak x}_n} \cup M^{{\frak x}'_n}$
as $A^{{\frak x}_{n+1}} = A^{{\frak x}_n} \cup \{\bar
a^n_t:t \in I_1\}$ and $\cup\{\bar a^n_t:t \in I_1\} 
\subseteq M^{{\frak x}_n}$ so clause
(b) of Definition \scite{dp1.3}(3)$(\gamma)$ holds.  Also $M^{{\frak
x}_{n+1}} \subseteq M^{{\frak x}'_n}$ and $p^{{\frak x}_{n+1}}
\supseteq p^{{\frak x}'_n} \restriction (A^{{\frak x}_n} \cup
M^{{\frak x}_{n+1}})$ holds trivially so also clause (c),(d) of
Definition \scite{dp1.3}(3)$(\gamma)$ holds.

Lastly, $\varphi_n(x,\bar a^n_{(1,i)})$ 
for $i < \omega,\neg \varphi_n(x,\bar a_{(2,n)})$ belongs 
to $p_\eta$ hence to $p^{{\frak x}_{n+1}}$
hence by renaming also clause (e) from Definition \scite{dp1.3}(4) holds.
So we are done. 

We are left with the case $\ell=9$.  For the proof above to work we
need just that $p(\in \bold S(M_0))$ satisfies $n < \omega \Rightarrow
p \restriction (M_n\cup A_n)$ is finitely satisfiable in $M_n$.
Toward this \wilog \, for each $n$ there is a function symbol $F_n \in
\tau(M^*)$ such that: if $\eta \in {}^n I$ then $c_\eta :=
F^{M^*}_n(\bar a^0_{\eta(0)},\dotsc,\bar a^{n-1}_{\eta(n-1)})$ realizes
$\{\varphi_m(x,\bar a^m_t)^{\text{if}(t=\eta(m))}:m<n$ and
$\alpha < \lambda\}$, so $F_n$ has arity $\Sigma\{\ell g(y_m):m<n\}$.

Let $D$ be a uniform ultrafilter on $\omega$ and let $c_\omega \in
{\frak C}$ realize $p^* = \{\psi(x,\bar b):\bar b \subseteq
M_0,\psi(x,\bar y) \in \Bbb L(\tau_{M^*})$ and $\{n:{\frak C}
\models \psi(c_{\eta \restriction n},\bar b)\} \in D$, so clearly $p =
\text{ tp}(c_\omega,M_0,{\frak C}) \in \bold S(M_0)$ extends
$\{\varphi_n(x,\bar a^m_t)^{\text{if}(t=\eta(n))}:n< \omega$
and $t \in I\}$.  So we have just to check that $p_n = p
\restriction (A_n \cup M_n)$ is finitely satisfiable in $M_n$, so let
$\vartheta(\bar x,\bar b) \in p_n$, so we can find $k(*) < \omega
(\subseteq \Bbb Z)$ such that $\bar b$ is included in the 
Skolem hull $M^*_{n,k(*)}$ of
$\cup\{\bar a^m_{(1,a)}:m < n$ and $a \in \Bbb Z \wedge a < k(*)\} 
\cup\{\bar a^m_t:m \in [n,\omega),t \in I\}$ inside $M^*$.

Let $\nu \in {}^\omega \lambda$ be defined by

$$
\nu(m) = \eta(m) \text{ for } m \in [n,\omega)
$$

$$
\nu(m) = (1,k(*) -n + m) \text{ for } m < n.
$$
\mn
By the indiscernibility:
\mr
\item "{$(*)_1$}"  for every $n,{\frak C} \models \psi(c_{\eta
\restriction n},\bar b) \equiv \psi(c_{\nu \restriction n},\bar b)$
\ermn
and by the choice of $p$
\mr
\item "{$(*)_2$}"   $\{n:{\frak C} \models 
\psi(c_{\eta \restriction n},\bar b)\}$ is infinite but clearly
\sn
\item "{$(*)_3$}"   $c_{\eta \restriction m} \in M_n$ for $m <
\omega$.
\ermn
Together we are done.
\mn
\ub{dp-rk$_\ell(T)=\infty$ implies dp-rk$_\ell(T) \ge |T|^+$}:

Trivial. 
\mn
\ub{dp-rk$_\ell(T) \ge |T|^+ \Rightarrow \kappa_{\text{ict}}(T) > \aleph_0$}:

We choose by induction on $n$ sequences $\bar \varphi^n$ and $\langle
{\frak x}^n_\alpha:\alpha < |T|^+ \rangle,\langle \bar{\bold
a}^n_\alpha:\alpha < |T|^+ \rangle$ such that:
\mr
\item "{$\circledast_n$}"  $(a) \quad \bar\varphi^n = \langle
\varphi_m(x,\bar y_m):m < n \rangle$; that is $\bar\varphi^n = \langle
\varphi^n_m(x,\bar y^n_m):m <n\rangle$ and 
\nl

\hskip20pt $\varphi^n_m(x,\bar y^n_m)
= \varphi^{n+1}_m(x,\bar y^{m+1}_m)$ for $m<n$ 
so we call it $\varphi_m(x,\bar y_m)$
\sn
\item "{${{}}$}"  $(b) \quad {\frak x}^n_\alpha \in K_\ell$ and
dp-rk$_\ell({\frak x}^n_\alpha) \ge \alpha$
\sn
\item "{${{}}$}"  $(c) \quad \bar{\bold a}^n_\alpha = \langle 
\bar a^{n,m}_{\alpha,k}:k < \omega,m<n \rangle$ where the sequence
$\bar a^{n,m}_{\alpha,k}$ is from $A^{{\frak x}^n_\alpha}$.
\sn
\item "{${{}}$}"  $(d) \quad$ for each $\alpha < |T|^+$ and $m<n$ the
sequence $\langle \bar a^{n,m}_{\alpha,k}:k < 
\omega \rangle$ is \nl

\hskip20pt indiscernible over 
$\cup\{\bar a^{n,i}_{\alpha,k}:i<n,i \ne m,k<\omega\} 
\cup M^{{\frak x}^n_\alpha} \cup A^n_\alpha$
\sn
\item "{${{}}$}"  $(e) \quad$ we have
$\bar b^{n,m}_\alpha \subseteq A^{{\frak x}^n_\alpha} = 
\cup\{\bar a^{n,i}_{\alpha,k}:i < m,k < \omega\} \cup A^n_\alpha$ 
\nl

\hskip35pt for $m<n$ such that:
\nl

\hskip35pt if $\eta \in {}^n \omega$ and 
$m<n \Rightarrow \bar b^{n,m}_\alpha \subseteq 
\cup\{\bar a^{n,i}_{\alpha,k}:i<m,k < \eta(i)\} \cup A^n_\alpha$
\nl

\hskip35pt \ub{then} 
$(p^{{\frak x}^n_\alpha} \restriction M^{{\frak x}^n_\alpha}) 
\cup\{\varphi_m(\bar a^{n,m}_{\alpha,\eta(m)},\bar
b^{n,m}_\alpha) \wedge \neg \varphi_m(\bar x,\bar
a^{n,m}_{\alpha,\eta(m)+1},\bar b^{n,m}_\alpha):$
\nl

\hskip35pt $m<n\}$ is finitely satisfiable in ${\frak C}$.
\ermn
For $n=0$ this is trivial by the assumption rk-dp$_\ell(T) \ge |T|^+$
see Definition \scite{dp1.3}(6) (and \scite{dp1.3}(7)).

For $n+1$, for every $\alpha < |T|^+$, (as rk-dp$_\ell
({\frak x}^n_{\alpha +1}) > \alpha$ by Definition
\scite{dp1.3}(5)) we can find ${\frak z}^n_\alpha,{\frak y}^n_\alpha,
\varphi^n_\alpha(x,\bar y^n_\alpha),\langle
\bar a^{n,*}_{\alpha,k}:k < \omega \rangle$ such that Definition
\scite{dp1.3}(4) is satisfied with $({\frak x}^n_{\alpha +1},
{\frak z}^n_\alpha,{\frak y}^n_\alpha,\varphi^n_\alpha(x,\bar y_\alpha),\langle
\bar a^{n,*}_{\alpha,k}:k < \omega \rangle)$ here standing for
$({\frak x},{\frak x}',{\frak y},\varphi(x,\bar y),\langle \bar a_k:k
< \omega \rangle)$ there such that rk-dp$_\ell({\frak y}^n_\alpha) \ge
\alpha$ and we also have $\bar a^{n,*}_{\alpha,\omega},
\bar b^{n,*}_\alpha$ here standing for $\bar a_\omega,\bar b$ there.  
So for some formula $\varphi_n(x,\bar y_n)$ the set $S_n
= \{\alpha < |T|^+:\varphi^n_\alpha(x,\bar y^n_\alpha) =
\varphi_n(x,\bar  y_n)\}$ is unbounded in $|T|^+$, 
so $\bar \varphi^{n+1}$ is well defined so clause $(a)$ of 
$\circledast_{n+1}$ holds.

For $\alpha < |T|^+$ let $\beta_n(\alpha) = \text{ Min}(S_n \backslash
\alpha)$ and let ${\frak x}^{n+1}_\alpha = {\frak y}^n_{\beta(\alpha)}$ so
clause (b) of $\circledast_{n+1}$ holds.  Let 
$\langle \bar a^{n+1,m}_{\alpha,k}:k < \omega \rangle$ be $\langle
\bar a^{n,m}_{\beta(\alpha)+1,k}:k < \omega \rangle$ if $m<n$ and $\langle
\bar a^{n,*}_{\beta(\alpha),k}:k < \omega\rangle$ if $m=n$ and let
$A^{n+1}_\alpha = A^n_{\beta(\alpha)}$, so clauses (c)
+ (d) from $\circledast_{n+1}$ holds.  Also 
we let $\bar b^{n+1,m}_\alpha$ is $\bar b^{n,m}_{\beta(\alpha)}$ if
$m<n$ and is $\bar b^{n,*}_{\beta(\alpha)}$ if $m=n$.
Next we check clause (e) of $\circledast_{n+1}$.

Let $\eta \in {}^{n+1}\omega$ be as required in sub-clause $(\gamma)$
of clause (e) of $\circledast_{n+1}$.  By the induction hypothesis $(p^{{\frak
x}^n_{\alpha +1}} \restriction M^{{\frak x}^n_{\alpha +1}}) \cup
\{\varphi(x,\bar a^{n,m}_{\alpha,\eta(m)}),\bar b^{n,m}_\alpha) \wedge
\neg \varphi(x,\bar a^{n,m}_{\alpha,\eta(m)+1},\bar
b^{n,m}_\alpha):m<n\}$ is finitely satisfiable in ${\frak C}$.  

By clause (d) of \scite{dp1.3}(3)$(\alpha)$ it follows that
$(p^{{\frak z}^n_\alpha} \restriction M^{{\frak z}^n_\alpha}) \cup
\{\varphi(x,\bar a^{n,m}_{\alpha +1,\eta(m)}),b^{n,m}_\alpha) \wedge
\neg \varphi(x,\bar a^{n,m}_{\alpha +1,\eta(m)+1}):m<\omega\}$ is
finitely satisfiable in ${\frak C}$ (i.e. we use $M^{{\frak
x}^n_{\alpha +1}} \le_{A[{\frak z}^n_\alpha],p[{\frak z}^n_\alpha] \restriction
M[{\frak z}^n_\alpha]} M^{{\frak z}^n_\alpha}$ which suffice; we use
freely the indiscernibility).

Hence, by
monotonicity for each $k<\omega$ and using other names, the set
$(p^{{\frak z}^n_\alpha} \restriction (M^{{\frak y}^n_\alpha}
\cup\{\bar a^{n+1,m}_{\alpha,k}:k \le \eta(n)$ or $k=\omega\} \cup
A^n_{\alpha +1}) \cup
\{\varphi(\bar x,\bar a^{n+1,m}_{\alpha,\eta(m)},\bar b^{m,n}_\alpha)
\wedge \neg \varphi(x,\bar a^{n+1,m}_{\alpha,\eta(m)+1};\bar
b^{n,m}_\alpha):m<n\}$ is finitely satisfiable in ${\frak C}$.

Similarly
$(p^{{\frak z}^n_\alpha} \restriction (M^{{\frak y}^n_\alpha})
\cup\{\varphi(x,\bar a^{n+1,m}_{\alpha,\eta(n)},\bar b^{n+1,n}_\alpha)
\wedge \neg \varphi(x,\bar a^{n+1,n}_{\alpha,\omega})\} \cup
\{\varphi(x,\bar a^{n+1,n}_{\alpha,\eta(m)},\bar b^{n+1,m}_\alpha)
\wedge \neg \varphi(x,\bar a^{n+1,m}_{\alpha,\eta(m+1)},\bar
b^{n+1,m}_\alpha):m<n\}$ is finitely satisfiable in ${\frak C}$.

But $\bar a^{n+1,m}_{\alpha,\omega},\bar a^{n+1,n}_{\alpha,\eta(n)+1}$
realizes the same type over a set including all the 
relevant elements so we can above
replace the first $(\bar a^{n+1,n}_{\alpha,\omega})$ by the second
$(\bar a^{n+1,n}_{\alpha,\eta(m)+1})$ so we are done proving
clause $(e)$ of $\circledast_{n+1}$.
\medskip

Having carried the induction 
it suffices to show that $\bar \varphi = \langle
\varphi_n(x,\bar y_n):n < \omega \rangle$ exemplifies that
$\kappa_{\text{ict}}(T) > \aleph_0$; for this 
it suffices to prove the assertion
$\circledast^2_{\bar\varphi}$ from \scite{dp1.2.2}(1).  
By compactness it suffices for each $n$ to 
find $\langle \bar a^{n,m}_k:k < \omega
\rangle$ for $m<n$ in ${\frak C}$ such that $\ell g(\bar a^{n,m}_k) =
\ell g(\bar y_n),\langle \bar a^{n,m}_k:k < \omega\rangle$ is
indiscernible over $\cup\{\bar a^{n,i}_k:k < \omega,i<n,i \ne m\}$ for
each $m<n$ and ${\frak C} \models (\exists x)[\dsize \bigwedge_{m<n}
(\varphi(x,\bar a^{n,m}_0) \wedge \neg \varphi(x,\bar a^{n,m}_1)]$.

We choose $\bar a^{n,m}_k = \bar
a^{n,m}_{\alpha,k(*)+k} \char 94 \bar b^{n,m}_\alpha$ where $k(*)$
is large enough such that $\cup\{\bar b^{n,m}_\alpha:m<n\} \subseteq
\cup \{\bar a^{n,m}_{\alpha,k}:m<n$ and $k<k(*)\}$ and let $\alpha=0$;
clearly we are done.  \hfill$\square_{\scite{dp1.4}}$
\enddemo
\bigskip

\demo{\stag{dp3.3} Observation}  1) If ${\frak x} \in K_\ell$ and $|T|
+ |A^{\frak x}| \le 
\mu < \|M^{\frak x}\|$ \ub{then} for some $M_0 \prec M^{\frak x}$
we have $\|M_0\| = \mu$ and for every ${\frak y} \le^\ell_{\text{pr}}
{\frak x}$ satisfying $M_0 \subseteq M^{\frak y}$  we have $\text{dp-rk}_\ell
({\frak y}) = \text{ dp-rk}_\ell({\frak x})$. 
\nl
1A) If dp-rk$_\ell({\frak x}) < \infty$ then it is $< |T|^+$.
Similarly dp-rk$_\ell(T)$, (with $(2^{|T|})^+$ this is easier).
\nl
1B) If dp-rk$_{\bar\Delta,\ell}({\frak x}) < \infty$ then it is $<
|\Delta_1 \cup \Delta_2|^+ + \aleph_0$.
\nl
2) If ${\frak x} \le^\ell_{\text{pr}} {\frak y}$ then 
dp-rk$_\ell({\frak x}) \ge \text{ dp-rk}_\ell({\frak y})$.
\nl
3) If ${\frak x} \le^\ell_{\text{pr}} {\frak y}$ and 
${\frak z}$ explicitly splits
$\ell$-strongly over ${\frak y}$ \ub{then} ${\frak z}$ explicitly
splits $\ell$-strongly over ${\frak x}$.
\nl
4) The previous parts hold for $m>1$, too.
\enddemo
\bigskip

\demo{Proof}  1) We do not need a really close look at the rank for
this.  First, fix an ordinal $\zeta$.

We can choose a vocabulary $\tau_{\zeta,\alpha,m}$ of cardinality $|A|+
|\zeta| + |T|$ such that:
\mr
\item "{$\circledast_1$}"   for any set $A$ fixing a sequence
$\bar{\bold a} = \langle a_\beta:\beta < \alpha \rangle$ listing the
elements of $A,M \prec {\frak C}$ and $p \in \bold S^m(M 
\cup\{a_\beta:\beta < \alpha\}),M_{A,p}$ or more exactly
$M_{\bar{\bold a},p}$ is a $\tau_{\zeta,\alpha,m}$-model;
\ermn
we let
\mr
\item "{$\circledast_2$}"  $(a) \quad$ ds$(\zeta) = \{\eta:\eta$ a
decreasing sequence of ordinals $< \zeta\}$
\sn
\item "{${{}}$}"  $(b) \quad \Gamma_\zeta = \{u:u$ is a subset of
ds$(\zeta)$ closed under initial segments$\}$ and 
\nl

\hskip25pt $\Gamma_\infty = \cup\{\Gamma_\zeta:\zeta$ an ordinal$\}$ 
\sn
\item "{${{}}$}"  $(c) \quad$ for $u \in \Gamma_\zeta$ let $\Xi^m_u =
\{\bar\varphi:\bar \varphi$ has the form $\langle \varphi_n(\bar
x,\bar y_\eta):\eta \in u \rangle$ where 
\nl

\hskip25pt $\bar x = \langle x_\ell:\ell
< m \rangle,\varphi_\eta(\bar x,\bar y_n) \in \Bbb L(\tau_T)\}$ and
\sn
\item "{$\circledast_3$}"  there are functions $\Phi_{\alpha,m}$ for
$m < \omega,\alpha$ an ordinal, satisfying
{\roster
\itemitem{ $(a)$ }  if $u \in \Gamma_\infty,\alpha \in \text{
Ord}$ and $\bar \varphi \in \Xi^m_u$, \ub{then} 
$\Phi_{\alpha,m}(u)$ is a set of first order sentences
\sn
\itemitem{ $(b)$ }  $\Phi_{\alpha,m}(u)$ is a set of first order
sentences
\sn
\itemitem{ $(c)$ }  if ${\frak x} \in K_{m,\ell}$ and $\bar{\bold a} =
\langle a_\beta:\beta < \alpha\rangle$ list $A^{\frak x}$ then
dp-rk$_\ell({\frak x}) \ge \zeta$ iff Th$(M_{\bar{\bold a},p[{\frak
x}]}) \cup \Phi_{\alpha,m}(\bar \varphi)$ is consistent for some
$\bar\varphi \in \Xi^m_{\text{ds}(\zeta)}$
\sn
\itemitem{ $(d)$ }  if $\bar\varphi,\bar\psi$ are isomorphic (see
below) then $\Phi_{\alpha,m}(\bar\varphi)$ is consistent iff
$\Phi_{\alpha,m}(\bar\psi)$ is;
\nl
where
\endroster}
\item "{$\circledast_4$}"  $\bar\varphi = \langle
\varphi_n(\bar x,\bar y_\eta):\eta \in u \rangle,\bar\psi = \langle
\psi_\eta(\bar x,\bar z_\eta):\eta \in v\rangle$ are isomorphic
when there is a one to one mapping function $h$ from $u$ onto $v$
preserving lengths, being initial segments and its negation such that
$\varphi_\eta(\bar x,\bar y_\eta) =\psi_{h(\eta)}(\bar x,\bar
z_{h(\eta)})$ for $\eta \in u$.
\ermn
[Why $\circledast_3$? just reflects on the definition.]

Now if $\zeta = \text{\rm dp-rk}_\ell({\frak x})$ is $\le \mu$
(e.g. $\zeta < |T|^+$) part (1) should be clear.  In the renaming case
if $\mu \ge |T|^+$ by (1A) we are done and otherwise use the implicit
characterization of ``$\infty = \text{\rm dp-rk}_\ell({\frak x})$".
\nl
1A) Now the proof is similar to the third part of the proof of
\scite{dp1.4}(1) and is standard.  We choose
by induction on $n$ a formula $\varphi_n(\bar x,\bar y_n) < |T|^+$ for
some decreasing sequence
$\eta^*_{m,\alpha}$ of ordinals $>\alpha$ of length $n$, we have
\mr
\item "{$\bigodot$}"  $\Phi_{n,\alpha}(\bar\varphi^n)$ is consistent
with Th$(M^{{\frak x}^n_\alpha}_{\bar{\bold a}[{\frak
x}^n_\alpha],p[{\frak x}^n_\alpha]})$ where
Dom$(\bar\varphi^{n,\alpha}) = \{\eta^*_{n,\alpha} \restriction
\ell:\ell \le n\}$ and $\varphi^{n,\alpha}_{\eta_{n,\alpha} \restriction
\ell}(\bar x,\bar y^{n,\alpha}_{\eta_{n,\alpha} \restriction \ell}) =
\varphi_\ell(\bar x,\bar y_\ell)$ for $\ell < n$.
\ermn
The induction should be clear and clearly is enough.
\nl
1B) Similarly.
\nl
2) We prove by induction on the ordinal $\zeta$ that
dp-rk$_\ell({\frak y}) \ge \zeta \Rightarrow \text{ dp-rk}_\ell({\frak
x}) \ge \zeta$.  For $\zeta=0$ this is trivial and for $\zeta$ a limit
ordinal this is obvious.  For $\zeta$ successor order let $\zeta =
\xi+1$ so there is ${\frak z} \in K_\ell$ which explicitly splits
$\ell$-strongly over ${\frak y}$ by part (3) and the definition of
dp-rk$_\ell$ we are done.
\nl
3) Easy as $\le^{\text{pr}}_\ell$ is transitive.
\nl
4) Similarly.   \hfill$\square_{\scite{dp3.3}}$
\enddemo
\bn
\centerline {$* \qquad * \qquad *$}
\bn
\ub{\S (3B) Ranks for strongly$^+$ dependent $T$}:

We now deal with a relative of Definition \scite{dp1.3} relevant for
``strongly$^+$ dependent".
\definition{\stag{dr.6} Definition}  1) For $\ell \in
\{14,15\}$ we define $K_{m,\ell} = K_{m,\ell-6}$ (and if
$m=1$ we may omit it and $\le^\ell_{\text{pr}} =
\le^{\ell-6}_{\text{pr}},\le^\ell_{\text{at}} =
\le^{\ell-6}_{\text{at}},\le^\ell = \le^{\ell -6}$.
\nl
2) For ${\frak x},{\frak y} \in K_{m,\ell}$ we say that ${\frak y}$ explicitly
$\bar\Delta$-split $\ell$-strongly over ${\frak x}$ when: $\bar\Delta
= (\Delta_1,\Delta_2),\Delta_1,\Delta_2 \subseteq \Bbb L(\tau_T)$ and
for some ${\frak x}'$ and $\varphi(\bar x,\bar y) \in \Delta_2$ with
$\ell g(\bar x) = m$ we
have clauses (a),(b),(c),(d) of clause $(\gamma)$ of Definition
\scite{dp1.3}(3) and
\mr
\item "{$(e)''$}"  there are $\bar b,\bar{\bold a}$ such that
{\roster
\itemitem{ $(\alpha)$ }  $\bar{\bold a} = \langle \bar a_i:i < \omega
\rangle$ is $\Delta_1$-indiscernible over $A^{\frak x} \cup M^{\frak
y}$
\sn
\itemitem{ $(\beta)$ }   $A^{\frak y} \supseteq A^{\frak x} \cup\{\bar
a_i:i < \omega\}$
\sn
\itemitem{ $(\gamma)$ }  $\bar b \subseteq A^{\frak x}$ and $\bar a_i \in
M^{\frak x}$ for $i < \omega$
\sn
\itemitem{ $(\delta)$ }  $\varphi(\bar x,\bar a_0 \char 94 \bar b)
\wedge \neg \varphi(\bar x,\bar a_1 \char 94 \bar b) \in p^{{\frak x}'}$.
\endroster}
\ermn
3) dp-rk$^m_\ell(T) = \cup\{\text{\rm dp-rk}_\ell({\frak x}) +1:
{\frak x} \in K_\ell\}$.
\nl
4) If $\Delta_1 = \Delta = \Delta_2$ we may write $\Delta$ instead of
$\bar\Delta$ and if $\Delta = \Bbb L(\tau_T)$ we may omit
$\Delta$.  Lastly,  if $m=1$ we may omit it.
\enddefinition
\bn
Similarly to \scite{dp1.3A}.
\demo{\stag{dr.7} Observation}  1) If ${\frak x},{\frak y} \in K_{15}$ then
``${\frak y}$ explicitly $\bar\Delta$-split $15$-strongly 
over ${\frak x}$" iff
``${\frak y}$ explicitly $\bar\Delta$-split $14$-strongly over ${\frak x}$".
\nl
2) If ${\frak x} \in K_{m,15}$ then
dp-rk$^m_{\bar\Delta,15}({\frak x}) \le 
\text{\rm dp-rk}^m_{\bar\Delta,14}({\frak x})$.
\enddemo
\bigskip

\demo{Proof}  Easy by the definition.
\enddemo
\bigskip

\proclaim{\stag{dr.8} Claim}  1) For $\ell=14$ we have {\rm dp-rk}$_\ell(T) =
\infty$ iff {\rm dp-rk}$_\ell(T) \ge |T|^+$ iff $\kappa_{\text{ict},2}(T) >
\aleph_0$.
\nl
2) For each $m \in [1,\omega)$ similarly using {\rm dp-rk}$^m_\ell(T)$.
\nl
3) The parallel of \scite{dp3.3} holds (for $\ell = 14,15$).
\endproclaim
\bigskip

\demo{Proof}  1) \ub{$\kappa_{\text{ict},2}(T) > \aleph_0$ implies
dp-rk$_\ell(T)=\infty$}.

As in the proof of \scite{dp1.4}.
\mn
\ub{dp-rk$_\ell(T) = \infty \Rightarrow \text{\rm dp-rk}_\ell(T) \ge
|T|^+$} for any $\ell$.

Trivial.
\mn
\ub{dp-rk$_\ell(T) \ge |T|^+ \Rightarrow \kappa_{\text{ict},2}(T)$}.

We repeat the proof of the parallel statement in \scite{dp1.4}, and we
choose $\bar b$ but not $\bar a^{n+1,n}_{\alpha,\omega}$. 
\nl
2) By part (1) and \scite{df2.4.8}(3).
\nl
3) Similar proof.  \hfill$\square_{\scite{dr.8}}$
\enddemo
\newpage

\head {\S4 Existence of indiscernibles} \endhead  \resetall \sectno=4
 \spuriousreset
\bn
Now we arrive to our main result.
\proclaim{\stag{ind.1} Theorem}  1) Assume
\mr
\item "{$(a)$}"  $\ell \in \{8,9\}$
\sn
\item "{$(b)$}"  $\infty > \zeta(*) = { \text{\rm dp-rk\/}}^m_{\ell}(T)$ 
so $\zeta(*)<|T|^+$
\sn
\item "{$(c)$}"  $\lambda_* = \beth_{2\times(\zeta(*)+1)}(\mu)$
\sn
\item "{$(d)$}"   $\bar a_\alpha \in {\frak C}_T$ for $\alpha <
\lambda^+_*,\ell g(\bar a_\alpha) = m$
\sn
\item "{$(e)$}"   $A \subseteq {\frak C}_T,|A| + |T| \le \mu$.
\ermn
\ub{Then} for some $u \in [\lambda^+_*]^{\mu^+}$, the sequence
$\langle \bar a_\alpha:\alpha \in u \rangle$ is an indiscernible
sequence over $A$.
\nl
2) If $T$ is strongly dependent, then for some $\zeta(*) < |T|^+$ part
(1) holds, i.e., if clauses (c),(d),(e) from there holds then the
conclusion there holds.
\endproclaim
\bigskip

\remark{\stag{ind.3} Remark}  0) This works for $\ell = 14,15,17,18$,
too, see \S (5A).
\nl
1) A theorem in this direction is natural as small
dp-rk points to definability and if the relevent types increases with
the index and are definable say over the first model then it follows
that the sequence is indiscernible.
\nl
2) The $\beth_{2 \times(\zeta +1)}(\mu)$ is more than needed, we can
use $\lambda^+_{\zeta(*)}$ where we define $\lambda_\zeta = \mu +
\Sigma\{(2^{\lambda_\xi})^+:\xi < \zeta\}$ by induction on $\zeta$.
\nl
3) We may like to have a one-model version of this theorem.  This will
be dealt with elsewhere.
\endremark
\bigskip

\demo{Proof}  Clearly ${\frak x} \in K_{m,\ell} 
\Rightarrow p^{\frak x} \in \bold S^m (A^{\frak x} \cup 
M^{\frak x})$ and we shall use clause $(e)$ of Definition \scite{dp1.3}(4).

By \scite{dp1.3A}(6), it 
is enough to prove this for $\ell=9$, but the
 proof is somewhat simpler for $\ell=8$, so we carry the proof for
 $\ell=8$ but say what more is needed for $\ell=9$.
We prove by induction on the ordinal $\zeta$ that
(note that the $M_\alpha$'s are increasing but not necessarily the
$p_\alpha$'s; this is not an essential point as by decreasing
somewhat the cardinals we can regain it):
\mr
\item "{$(*)_\zeta$}"   if the sequence ${\bold I}
= \langle \bar a_\alpha:\alpha < \lambda^+ \rangle$ satisfies
$\boxtimes_\zeta$ below then for some $u \in [\lambda^+]^{\mu^+}$ the
sequence $\langle \bar a_\alpha:\alpha \in u \rangle$ is an
indiscernible sequence over $A$ where (below, the 2 is an overkill, in
particular for successor of successor, but for
limit $\zeta$ we ``catch our tail"):
\sn
\item "{$\boxtimes_\zeta$}"  there are $\lambda,B,\bar M,\bar p$ such that
{\roster
\itemitem{ $(a)$ }   $\lambda = \lambda^{\beth_{2(\xi +1)}(\mu)}$ for
every $\xi < \zeta$
\sn
\itemitem{ $(b)$ }  $\bar M = \langle M_\alpha:\alpha < \lambda^+ \rangle$
and $M_\alpha \prec {\frak C}_T$ is increasing continuous (with $\alpha$)
\sn
\itemitem{ $(c)$ }   $M_\alpha$ has cardinality $\le \lambda$
\sn
\itemitem{ $(d)$ }  $\bar a_\alpha \in {}^m(M_{\alpha +1})$ for
$\alpha < \lambda^+$ 
\sn
\itemitem{ $(e)$ }  $p_\alpha = \text{ tp}(\bar a_\alpha,M_\alpha \cup A
\cup B)$
\sn
\itemitem{ $(f)$ }  $B \subseteq {\frak C},|B| \le \aleph_0$
\sn
\itemitem{ $(g)$ }   ${\frak x}_\alpha = (p_\alpha,M_\alpha,A \cup B)$
belongs to $K_{m,\ell}$ and satisfies dp-rk$^m_\ell({\frak x}_\alpha) 
< \zeta$.
\endroster}
\ermn
Why is this enough?  We apply $(*)$ for the case $\zeta
= \zeta(*)$ so $\lambda = \lambda_*$ and we choose $M_\alpha \prec
{\frak C}$ of cardinality $\lambda$ by induction on $\alpha <
\lambda^+$ such that $M_\alpha$ is increasing continuous, $\{\bar
a_\beta:\beta < \alpha\} \subseteq M_\alpha$.

If $\ell=8$ fine; if $\ell=9$ it seemed that we have a problem with clause
(g).  That is in checking ${\frak x}_\alpha \in K_{n,\ell}$ we have to
show that ``$p_\alpha$ is finitely satisfiable in $M_\alpha$".  But
this is not a serious one: in this case note that for some club $E$ 
of $\lambda^+$, for every $\alpha \in E$ the type 
we have tp$(a_\alpha,M_\alpha \cup A \cup B)$
is finitely satisfiable in $M_\alpha$.  So letting $M'_\alpha =
M_{\alpha'},a'_\alpha = \bar a_{\alpha'}$ when $\alpha < \lambda^+,\alpha' \in
E$ and otp$(C \cap \alpha') = \alpha$ and similarly $p'_\alpha = \text{
tp}(\bar a_{\alpha'},M_\alpha,{\frak C})$ we can use $\langle
(a'_\alpha,M'_\alpha,p'_\alpha):\alpha < \lambda^+\rangle$ so 
we are done.

So let us carry the induction; arriving to $\zeta$ we let 
$\theta_\ell = \beth_{2 \times \zeta + \ell}(\mu)$, for $\ell <
3$; note that $\theta^{\theta_\ell}_{\ell +1} = \theta_\ell$ and
$\lambda^{\theta_2} = \lambda$.  Let $\chi$ be large enough and
let ${\frak B} \prec ({\Cal H}
(\chi),\in,<^*_\chi)$ be of cardinality $\lambda$ such that ${\frak C},\bar M,
\bar p,\bar{\bold a},B,A$ belongs to ${\frak B}$ and $\lambda +1
\subseteq {\frak B}$ and $Y \subseteq {\frak B} \wedge |Y| \le
\theta_2 \wedge \lambda^{|Y|} = X  \Rightarrow Y \in {\frak B}$.  Let
$\delta(*) = {\frak B} \cap \lambda^+$ so \wilog \, cf$(\delta(*))$
satisfies $\lambda^{\text{cf}(\delta(*))} > \lambda$.
Let $\zeta^* = \text{ dp-rk}(p_{\delta(*)},M_{\delta(*)},A \cup B)$ and
$\theta = \theta_1$, hence $\lambda = \lambda^{\theta^+}$. 
We try by induction on $\varepsilon \le \theta^+ + \theta^+$
to choose $(N_{\alpha_\varepsilon},\alpha_\varepsilon)$ such that
\mr
\item "{$\circledast_\varepsilon$}"  $(a) \quad \alpha_\varepsilon <
\delta(*)$ is increasing with $\varepsilon$
\sn
\item "{${{}}$}"  $(b) \quad N_\varepsilon <_{A \cup
B,p_{\alpha(*)}} M_{\delta(*)}$ is increasing continuous with $\varepsilon$
\sn
\item "{${{}}$}"  $(c) \quad N_\varepsilon$ has cardinality $\theta$
\sn
\item "{${{}}$}"  $(d) \quad \xi < \varepsilon \Rightarrow
a_{\alpha_\xi} \in N_{\alpha_\varepsilon}$
\sn
\item "{${{}}$}"  $(e) \quad \bar a_{\alpha_\varepsilon}$ realizes
$p_{\delta(*)} \restriction (N_{\alpha_\varepsilon} \cup A \cup B)$ 
\sn
\item "{${{}}$}"  $(f) \quad$ if $p_{\alpha(*)}$ splits over
$N_\varepsilon \cup A \cup B$ then $p_{\delta(*)} \restriction 
(N_{\alpha_{\varepsilon +1}} \cup A \cup B)$ splits over \nl

\hskip35pt $N_\varepsilon \cup A \cup B$
\sn
\item "{${{}}$}"  $(g) \quad (p_{\alpha_\varepsilon} \restriction
(N_{\alpha_\varepsilon} \cup A \cup B),N_{\alpha_\varepsilon},A \cup
B) <_{\text{pr}} (p_{\delta(*)},M_{\delta(*)},A \cup B)$ and they
\nl

\hskip25pt  (have to) have the same dp-rk
\sn 
\item "{${{}}$}"  $(h) \quad N_\varepsilon \subseteq 
M_{\alpha_\varepsilon}$ (but not used).
\ermn
Clearly we can carry the definition.  Now the proof splits to two cases.
\bn
\ub{Case 1}:  For $\xi = \theta^+,p_{\alpha(*)}$ does not
split over $N_{\alpha_\xi} \cup A \cup B$.  

By clause (e) of $\circledast_\varepsilon$ 
clearly $\varepsilon \in [\xi,\xi + \theta^+) \Rightarrow \text{
tp}(\bar a_{\alpha_\varepsilon},N_\varepsilon \cup A \cup B)$
does not split over $N_{\alpha_\xi} \cup A \cup B$ and
increases with $\varepsilon$.  As $\langle N_{\xi + \varepsilon}:
\varepsilon < \theta \rangle$ is increasing and $\bar
a_{\alpha_\varepsilon} \in N_{\varepsilon +1}$ it follows that
tp$(\bar a_{\alpha_\varepsilon},N_{\theta^+} \cup \{\bar
a_\beta:\beta \in [\theta^+,\varepsilon)\} \cup A \cup B\}$ does not
split over $N_{\theta^+_1} \cup A \cup B$.  Hence by
\cite[I,\S2]{Sh:c} that the sequence $\langle \bar a_{\alpha_j}:j \in
[\xi,\xi +\theta^+)\rangle$ is an indiscernible sequence 
over $N_{\alpha_\xi} \cup A \cup B$ so as $M^+ \le \theta^+$ we are done.
\mn
\ub{Case 2}:  For $\xi = \theta^+,p_{\delta(*)}$ splits over
$N_{\alpha_\xi} \cup A \cup B$.

So we can find $\varphi(x,\bar y) \in \Bbb L(\tau_T)$ and
$\bar b,\bar c \in {}^{\ell g(\bar y)}(M_{\delta(*)} \cup A \cup B)$
realizing the same type over $N_{\alpha_\xi} \cup A \cup B$ and
$\varphi(\bar x,\bar b),\neg \varphi(\bar x,\bar c) \in
p_{\delta(*)}$.  So \wilog \, $\bar b = \bar b' \char 94 \bar d,
\bar c = \bar c' \char 94 \bar d$ where 
$\bar d \in {}^{\omega >}(A \cup B)$
and $\bar b',\bar c' \in {}^{m(*)}(M_{\delta(*)})$ for some
$m(*)$.  As $N_{\alpha_\xi} <_{A \cup B} M_{\delta(*)}$ 
(see clause $(b)$ of $\circledast_\xi$)
clearly there is $D$, an ultrafilter on ${}^{m(*)}(N_\xi)$ such that 
Av$(N_\xi \cup A \cup B,D) = \text{ tp}(\bar b',N_\xi \cup A 
\cup B) = \text{ tp}(\bar c',N_\xi \cup A \cup B)$.

Without loss of generality $\{\bar b'' \in
{}^{m(*)}(N_{\alpha_\xi}):\neg \varphi(\bar x,\bar b'',\bar d) \in
p_{\delta(*)}\}$ belongs to $D$, as otherwise we can replace $\varphi,\bar
b',\bar c'$ by $\neg \varphi,\bar c',\bar b'$.

Let $M_* = (M_{\delta(*)})_{A \cup B \cup \{\bar a_{\delta(*)}\}}$ and let
$M^+ \prec {\frak C}$ be such that $M_{\delta(*)} \subseteq M^+$ and
moreover $(M_*)_{A \cup B \cup \{\bar a_{\delta(*)}\}} 
\prec M^+_{A \cup B \cup \{\bar a_{\delta(*)}\}}$ and the latter is
$\lambda^+$-saturated.  Clearly
letting $p^+_\delta = (\text{tp}(\bar a_{\delta(*)},M^+ \cup A \cup
B)$ and ${\frak x}^+_{\delta(*)} = (p^+_{\delta(*)},M^+_{\delta(*)},A \cup B)$
we have ${\frak x}_{\delta(*)} \le_{\text{pr}} {\frak x}^+_{\delta(*)}$.
Note that $\varepsilon < \xi \Rightarrow (p_{\alpha_\varepsilon}
\restriction (N_{\alpha_\varepsilon} \cup A \cup
B),N_{\alpha_\varepsilon},A \cup B) \le_{\text{pr}} {\frak x}_{\delta(*)}$.

We can find $\langle \bar b_\alpha:\alpha < \omega +\omega \rangle$
such that $\bar b_\alpha \in {}^{m(*)}(M^+)$ realizes
Av$(N_{\alpha_\xi} \cup A \cup B \cup \{\bar b_\beta:\beta <
\alpha\},D)$ and \wilog \, $\bar b_\omega = \bar b'$.

We would like to apply the induction hypothesis to $\zeta' = \text{
dp-rk}({\frak x}_{\delta(*)})$, so let
\mr
\item "{$\boxdot$}"  $(a) \quad \lambda' = \theta$
\sn
\item "{${{}}$}"  $(b) \quad a'_\varepsilon = a_{\alpha_\varepsilon}$
for $\varepsilon < \theta^+$
\sn
\item "{${{}}$}"  $(c) \quad M'_\varepsilon = N_\varepsilon$
\sn
\item "{${{}}$}"  $(d) \quad p'_\varepsilon = \text{ tp}(\bar
a_{\alpha_\varepsilon},N_\varepsilon)$
\sn
\item "{${{}}$}"  $(e) \quad B' = B \cup \{\bar b_\alpha:\alpha <
\omega + \omega\}$
\sn
\item "{${{}}$}"  $(f) \quad A'=A$.
\ermn
We can apply the induction hypothesis to $\zeta'$, i.e., use
$(*)_{\zeta'}$ for some $u' \subseteq \theta^+$ of cardinality $\mu^+$
the sequence $\langle a'_\varepsilon:\varepsilon \in u'\rangle$ is
indisernible over $A$, hence the set $u :=
\{\alpha_\varepsilon:\varepsilon \in u'\}$ has cardinality $\mu^+$ and
the sequence $\langle a_\alpha:\alpha \in u\rangle$ is indiscernible
over $A$ so we are done.

But we have to check that the demands from $\boxtimes_{\zeta'}$ holds (for
$\theta^+$) $\bar M' = \langle M'_\varepsilon:\varepsilon <
\theta^+\rangle,\bar p' = \langle p'_\varepsilon:\varepsilon <
\theta^+\rangle$.
\bn
\ub{Clause (a)}:  As $\theta = \beth_{2 \times \zeta^*+1}(\mu)$
clearly for every $\xi < \zeta^*$ we have $\theta = \theta^{\beth_{2
\times (\xi+1)}}$ hence $\theta = \theta^{\beth_{2 \times (\xi +1)}}$.
\bn
\ub{Clause (b)}:  By $\circledast_\varepsilon(b),\bar M$ is increasing
continuous.
\bn 
\ub{Clause (c)}:  By $\circledast_\varepsilon(c)$.
\bn
\ub{Clause (d)}:  By $\circledast_\varepsilon(d)$.
\bn
\ub{Clause (e)}:  By the choice of $p'_\varepsilon$.
\bn
\ub{Clause (f)}:  By the choice of $B'$.
\bn
\ub{Clause (g)}:  Clearly ${\frak x}'_\varepsilon \in K_{m,\ell}$, but
why do we have dp-rk$({\frak x}'_\varepsilon) < \zeta^*$?  This is
equivalent to dp-rk$({\frak x}'_\varepsilon) < \text{ dp-rk}({\frak
x}_{\delta(*)})$.

Recall ${\frak x}_{\delta(*)} \le_{\text{pr}} {\frak x}^+_{\delta(*)}$
and ${\frak x}'_\varepsilon$ explicitly split $\ell$-strongly over
${\frak x}_{\delta(*)}$, hence by the definition of dp-rk we get
dp-rk$({\frak x}'_\varepsilon) < \text{ dp-rk}({\frak
x}_{\delta(*)})$.

What about the finitely satisfiable of $p'$ when
$\ell=9$? for some club $E$ of
$\theta^+,\varepsilon \in E \Rightarrow 
\text{ tp}(\bar a_{\alpha_\varepsilon},
N_{\alpha_\varepsilon} \cup A \cup B')$ is
finitely satisfiability in $N_{\alpha_\varepsilon}$.
\nl
2) By \scite{dp1.4}, dp-rk$^m_\ell(T) < |T|^+$ for $\ell =8$, so we can
apply part (1).  \hfill$\square_{\scite{ind.1}}$
\enddemo 
\newpage

\head {\S5 Concluding Remarks} \endhead  \resetall \sectno=5
 \spuriousreset
\bn
We comment on some things here which we intend to continue elsewhere
so the various parts ((A),(B),...) are not so connected.
\sn
(A) \quad \ub{Ranks for dependent theories}:
\bn
We note some generalizations of \S3, so Definition \scite{dp1.3} is replaced
by
\definition{\stag{e3.4} Definition}  1) For $\ell=1,2,3,4,5,
6,8,9,11,12$ (but not 7,10), let 

$$
\align
K_{m,\ell} = \bigl\{{\frak x}:&{\frak x} = (p,M,A),M \text{ a model } \prec 
{\frak C}_T,A \subseteq {\frak C}_T, \\
  &\text{if } \ell \in \{1,4\} \text{ then } p \in 
\bold S^m(M),\text{ if } \ell \notin \{1,4\} \text{ then} \\
  &p \in \bold S^m(M \cup A) \text{ and if } \ell=3,6,9,12 \text{ then} \\
  &p \text{ is finitely satisfiable in } M\bigr\}.
\endalign
$$
\mn
If $m=1$ we may omit it.

For ${\frak x} \in K_{m,\ell}$ let 
${\frak x} = (p^{\frak x},M^{\frak x},A^{\frak x}) = (p[{\frak
x}],M[{\frak x}],A[{\frak x}])$ and $m = m({\frak x})$ recalling 
$p^{\frak x}$ is an $m$-type.
\nl
2) For ${\frak x} \in K_{m,\ell}$ let $N_{\frak x}$ be $M$ expanded by
$R_{\varphi(\bar x,\bar y,\bar a)} = \{\bar b \in {}^{\ell g(\bar
y)}M:\varphi(\bar x,\bar b,\bar a) \in p\}$ for $\varphi(\bar x,\bar y,\bar z)
\in \Bbb L(\tau_T),\bar a \in {}^{\ell g(\bar z)}A$ and $\ell = 1,4
\Rightarrow \bar a = <>$ and
$R_{\varphi(\bar y,\bar a)} = \{\bar b \in {}^{\ell g(\bar y)}M:{\frak
C} \models \varphi[\bar b,\bar a]\}$ for $\varphi(\bar y,\bar z) \in
\Bbb L(\tau_T),\bar a \in {}^{\ell g(\bar y)}{\frak C}$; let $\tau_{\frak x} =
\tau_{N_{\frak x}}$.
\nl
2A)  If we omit $p$ we mean $p = \text{ tp}(<>,M \cup
A)$ so we can write $N_{A}$, a $\tau_A$-model so in this case 
$p = \{\varphi(\bar b,\bar a):
\bar b \in M,\bar a \in M$ and ${\frak C} \models \varphi[\bar b,\bar a]\}$.
\nl
3) For ${\frak x},{\frak y} \in K_{m,\ell}$ let
\mr
\item "{$(\alpha)$}"  ${\frak x} \le^\ell_{\text{pr}} {\frak y}$ means
that ${\frak x},{\frak y} \in K_{m,\ell}$ and
{\roster
\itemitem{ $(a)$ }  $A^{\frak x} = A^{\frak y}$
\sn 
\itemitem{ $(b)$ }  $M^{\frak x} \le_{A[{\frak x}]} M^{\frak y}$
\sn
\itemitem{ $(c)$ }  $p^{\frak x} \subseteq p^{\frak y}$
\sn
\itemitem{ $(d)$ }  if $\ell = 1,2,3,8,9$ then $M^{\frak x} 
\le_{A[{\frak x}],p[{\frak y}]} M^{\frak y}$ (for $\ell=1$ this
follows from clause (b))
\endroster}
\sn
\item "{$(\beta)$}"  ${\frak x} \le^\ell {\frak y}$ means that 
for some $n$ and 
$\langle {\frak x}_k:k \le n \rangle,{\frak x}_k \le^\ell_{\text{at}}
{\frak x}_{k+1}$ for $k < n$ and $({\frak x},{\frak y}) = ({\frak x}_0,{\frak
x}_n)$ 
\nl
where
\sn
\item "{$(\gamma)$}"  ${\frak x} \le^\ell_{\text{at}} {\frak y}$ iff
$({\frak x},{\frak y} \in K_{m,\ell}$ and) for some ${\frak x}' \in
K_{m,\ell}$ we have
{\roster
\itemitem{ $(a)$ }  ${\frak x} \le^\ell_{\text{pr}} {\frak x}'$
\sn
\itemitem{ $(b)$ }   $A^{\frak x} \subseteq A^{\frak y} \subseteq
A^{\frak x} \cup M^{{\frak x}'}$
\sn
\itemitem{ $(c)$ }  $M^{\frak y} \subseteq M^{{\frak x}'}$
\sn
\itemitem{ $(d)$ }  $p^{\frak y} \supseteq p^{{\frak x}'}
\restriction (A^{\frak x} \cup M^{\frak y})$ so $\ell \in \{1,4\}
\Rightarrow p^{\frak y} = p^{{\frak x}'} \restriction M^{\frak y}$ and
$\ell \notin \{1,4\} \Rightarrow p^{\frak y} = p^{{\frak x}'}
\restriction (M^{\frak y} \cup A^{\frak y})$.
\endroster}
\ermn
4) For ${\frak x},{\frak y} \in K_{m,\ell}$ we say that ${\frak y}$ 
explicitly $\bar\Delta$-splits $\ell$-strongly over ${\frak x}$
\ub{when}: $\bar\Delta = (\Delta_1,\Delta_2),\Delta_1,\Delta_2 \subseteq
\Bbb L(\tau_T)$ and for some ${\frak x}'$ and 
$\varphi(\bar x,\bar y) \in \Delta_2$ we have clauses 
(a),(b),(c),(d) of part (3)$(\gamma)$ and
\mr
\item "{$(e)$}"    when $\ell \in \{1,2,3,4,5,6\}$, in 
$A^{\frak y}$ there is a $\Delta_1$-indiscernible sequence 
$\langle \bar a_k:k < \omega\rangle$ over $A^{\frak x} \cup M^{\frak
y}$ such that $\bar a_k \in {}^{\omega >}(M^{{\frak x}'})$ and 
$\varphi(\bar x,\bar a_0),\neg \varphi(\bar x,\bar a_1) \in p^{{\frak
x}'}$ and $\bar a_k \subseteq A^{\frak y}$ for $k < \omega$
\sn
\item "{$(e)'$}"  when $\ell = 8,9,11,12$ there are $\bar b,
\bar{\bold a}$ such that
{\roster
\itemitem{ $(\alpha)$ }  $\bar{\bold a} = \langle \bar a_i:i < \omega
+ 1 \rangle$ is $\Delta_1$-indiscernible over $A^{\frak x} \cup M^{\frak y}$
\sn
\itemitem{ $(\beta)$ }  $A^{\frak y} \backslash A^{\frak x} = 
\{\bar a_i:i < \omega\}$; yes $\omega$ not $\omega + 1$! 
(note that $``A^{\frak x} ="$ and not ``$A^{\frak y} \backslash
A^{\frak x} \supseteq"$ as
we use it in $(e)(\gamma)$ in the proof of \scite{dp1.4})
\sn
\itemitem{ $(\gamma)$ }  $\bar b \subseteq A^{\frak x}$ and $\bar a_i 
\in M^{{\frak x}'}$ for $i < \omega +1$
\sn
\itemitem{ $(\delta)$ }  $\varphi(\bar x,\bar a_k \char 94 \bar b) 
\wedge \neg \varphi(\bar x,\bar a_\omega \char 94 \bar b))$ belongs
\footnote{this explains why $\ell = 7,10$ are missing}
to $p^{{\frak x}'}$ for $k < \omega$.
\endroster}
\ermn
5) We define dp-rk$^m_{\bar\Delta,\ell}:K_{m,\ell} \rightarrow \text{ Ord }
 \cup\{\infty\}$ by 
\mr
\item "{$(a)$}"  $\text{dp-rk}^m_{\bar\Delta,\ell}({\frak x}) 
\ge 0 \text{ always}$
\sn
\item "{$(b)$}"  $\text{dp-rk}^m_{\bar\Delta,\ell}({\frak x}) \ge 
\alpha +1 \text{ \ub{iff} there is } {\frak y} \in K_{m,\ell} 
\text{ which explicitly } \bar\Delta$-splits $\ell$-strongly over ${\frak x} 
\text{ and dp-rk}_{\bar\Delta,\ell}({\frak y}) \ge \alpha$
\sn
\item "{$(c)$}"   $\text{dp-rk}^m_{\bar\Delta,\ell}({\frak x}) \ge \delta 
\text{ iff dp-rk}^m_{\bar\Delta,\ell}({\frak x})
\ge \alpha \text{ for every } \alpha < \delta$ when $\delta$ is
a limit ordinal.
\ermn
Clearly well defined.  We may omit $m$ from dp-rk as ${\frak x}$
determines it.
\nl
6) Let dp-rk$^m_{\bar\Delta,\ell}(T) = 
\cup\{\text{dp-rk}_{\bar\Delta,\ell}({\frak x}):{\frak x} \in
K_{m,\ell}\}$; if $m=1$ we may omit it.
\nl
7) If $\Delta_1 = \Delta_2 = \Delta$ we may write $\Delta$ instead of
$(\Delta_1,\Delta_2)$. 
If $\Delta = \Bbb L(\tau_T)$ then we may omit it.
\nl
8) For ${\frak x} \in K_{m,\ell}$ let ${\frak x}^{[*]} = (p^{\frak x}
\restriction M^{\frak x},M^{\frak x},A^{\frak x})$.
\enddefinition
\bn
So Observation \scite{dp1.3A} is replaced by
\demo{\stag{e3.5} Observation}  1) $\le^\ell_{\text{pr}}$ is a
partial order on $K_\ell$.
\nl
2) $K_{m,\ell(1)} \subseteq K_{m,\ell(2)}$ when $\ell(1),\ell(2) \in
\{1,2,3,4,5,6,8,9,11,12\}$ and $\ell(1) \in \{1,4\} \Leftrightarrow
\ell(2) \in \{1,4\}$ and $\ell(2) \in \{3,6,9,12\} \Rightarrow \ell(1)
\in \{3,6,9,12\}$.
\nl
2A) $K_{m,\ell(1)} \subseteq \{{\frak x}^{[*]}:{\frak x} \in
K_{m,\ell(2)}\}$ when $\ell(1) \in \{1,4\},\ell(2) \in
\{1,\dotsc,6,8,9,11,12\}$.
\nl
2B) In (2A) equality holds if $x(\ell(1),\ell(2)) \in
\{(1,2),(1,3),(4,5),(4,6)\}$. 
\nl
3) ${\frak x} \le^{\ell(1)}_{\text{pr}} {\frak y} \Rightarrow 
{\frak x} \le^{\ell(2)}_{\text{pr}} {\frak y}$ when $(\ell(1),\ell(2))$ is
as in (2) and $\ell(2) \in \{2,3,8,9\} \Rightarrow \ell(1) \in
\{2,3,8,9\}$.
\nl
3B) ${\frak x} \le^{\ell(1)}_{\text{pr}} {\frak y} \Rightarrow {\frak
x}^{[*]} \le^{\ell(1)}_{\text{pr}} {\frak y}^{[*]}$ when the pair
$(\ell(1),\ell(2))$ is as in (2B).
\nl
4) ${\frak x} \le^{\ell(1)}_{\text{at}} 
{\frak y} \Rightarrow {\frak x} \le^{\ell(2)}_{\text{at}} {\frak y}$
when $(\ell(1),\ell(2))$ are as in part (3) (hence (2)).
\nl
4B) ${\frak x} \le^{\ell(1)}_{\text{at}} {\frak y} \Rightarrow 
{\frak x}^{[*]} \le^{\ell(2)}_{\text{at}} {\frak y}$
if $(\ell(1),\ell(2))$ are as in part (2A).
\nl
5) ${\frak y}$ explicitly $\bar\Delta$-splits $\ell(1)$-strongly over
${\frak x}$ implies ${\frak y}$ explicitly $\bar\Delta$-splits
$\ell(2)$-strongly over ${\frak x}$ when the pair
$(\ell(1),\ell(2))$ is as in parts (2),(3) and $\ell(1) \in \{1,2,3,4,5,6\}
\Leftrightarrow \ell(2) \in \{1,2,3,4,5,6\}$.
\nl
6) Assume $(\ell(1),\ell(2))$ is as in parts (2),(3),(5).  If ${\frak x}
\in K_{m,\ell(1)}$ then dp-rk$^m_{\bar\Delta,\ell(1)}({\frak x})
\le \text{\rm dp-rk}^m_{\bar\Delta,\ell(2)}({\frak x})$; i.e.,
\nl
$\{\ell(1),\ell(2)) \in \{(3,2),(2,5),(3,5),(6,5),(3,6)\} \cup
\{(9,8),(8,11),(9,11),(12,11),(9,12)\}$. 
\nl
7) Assume $\bar a \in {}^m {\frak C}$ and ${\frak y} = (\text{\rm tp}(\bar
a,M \cup A),M,A)$ and ${\frak x} = (\text{\rm tp}(\bar a,M \cup A),M,A)$.  Then
\mr
\item "{$(a)$}"  ${\frak x}^{[*]} = {\frak y}^{[*]}$
\sn
\item "{$(b)$}"  ${\frak x} \in K_{m,1} \cap K_{m,4}$
\sn
\item "{$(c)$}"  ${\frak y} \in K_{m,2} \cap K_{m,5}
\cap K_{m,8} \cap K_{m,11}$
\sn
\item "{$(d)$}"  if tp$(\bar a,M \cup A)$ is finitely satisfiable in
$M$ then also ${\frak y} \in K_{m,3} \cap K_{m,6} \cap K_{m,9} \cap
K_{m,12}$.
\ermn
8) If ${\frak x} \in K_{m,\ell(2)}$ then
dp-rk$_{\ell^m(2)}({\frak x}^{[*]}) \le \text{\rm dp-rk}_{\ell^m(2)}
({\frak x})$ when the pair $(\ell(1),\ell(2))$ is as in part (2A).
\nl
9) If ${\frak x} \in K_{m,\ell}$ and $\kappa > \aleph_0$ then there is
${\frak y} \in K_{m,\ell}$ such that ${\frak x} \le^\ell_{\text{pr}}
{\frak y}$ and $M^{\frak y}$ is $\kappa$-saturated, moreover $M^{\frak
y}_{A[{\frak y}],p[{\frak y}]}$ is $\kappa$-saturated (hence in
Definition \scite{dp1.1}(4) \wilog \, $M^{{\frak x}'}$ is $(|M^{\frak
x} \cup A^{\frak x}|^+)$-saturated).
\enddemo
\bigskip

\proclaim{\stag{e3.6} Claim}  In \scite{dp1.4} we can allow
$\ell=1,2,5$ (in addition to $\ell=8,9$).
\endproclaim
\bigskip

\demo{Proof}  Similar but:
\mn
\ub{$\kappa_{\text{ict}}(T) > \aleph_0$} implies dp-rk$_\ell(T) = D$
when $\ell \in \{1,2,4,5,8,9,11,12\}$:
\mr
\item "{$(A)$}"  Let $A_n = \cup\{\bar a^m_t:
m <n,t \in I_2\}$ if $\ell < 7$ and if $\ell > 7,
A_n = \{\bar a^m_t:m<n$ and $t \in I_1\}$.
\sn
\item "{$(B)$}"  ``${\frak x}_{n+1}$ explicitly split $\ell$-strongly
over ${\frak x}_n$" using $\langle \bar a^n_{(2,n+i)}:i <
\omega\rangle$ if $\ell < 7$ and $\langle a^n_{(1,i)}:i <
\omega\rangle \char 94 \langle \bar a^n_{2,n}\rangle$ if $\ell > 7$.
\sn
\item "{$(C)$}"  Similarly in ``Lastly...":  
Lastly, if $\ell < 7,\varphi_n(x,\bar a^n_{(1,n)}),\neg \varphi_n(x,\bar
a^n_{(1,n+1)})$ belongs to $p^{{\frak x}'_n}$ and even $p^{{\frak
x}_{n+1}}$ and if $\ell > 7,\varphi_n(x,\bar a^n_{(1,n)})$ 
for $n < \omega$,
$\neg \varphi_n(x,\bar a_{(2,n)})$ belongs to $p_\eta$ hence to 
$p^{{\frak x}_{n+1}}$
hence by renaming also clause (e) or (e)$^-$ from
Definition \scite{dp1.3}(4) holds.
So we are done. 
\endroster
\bn
\ub{$\text{dp-rk}_\ell(T) \ge |T|^+ \Rightarrow \kappa_{\text{ict}}(T)
> \aleph_0$} when $\ell=1,2,3,5,6,8,9$
\mr
\item "{$(D)$}"  In $\circledast_n(e)$ we use
\sn
\item "{$(E)$}"  $(\alpha) \quad$ if $\ell \in \{2,3,5,6\}$
and $m<n,k < \omega$ then
$\varphi_m(x,\bar a^{n,m}_{\alpha,k}) \in p^{{\frak x}^n_\alpha}$
\nl

\hskip35pt $\Leftrightarrow k =0$ hence
$\neg \varphi_m(x,\bar a^{n,m}_{\alpha,k})
\in p^{{\frak x}^n_\alpha}$
\nl

\hskip35pt $\Leftrightarrow k \ne 0$ for $k < 2$
\sn
\item "{${{}}$}"  $\qquad (\beta) \quad$ if $\ell=1$ then $p^{{\frak
x}^n_\alpha} \cup \{\varphi_m(x,\bar a^{n,m}_{\alpha,k})
^{\text{if}(k=0)}:m<n,k < 2\}$ is 
\nl

\hskip35pt consistent
\sn
\item "{${{}}$}"  $\qquad (\gamma) \quad$ if $\ell=8,9$ we also have
$\bar b^{n,m}_\alpha \subseteq A^{{\frak x}^n_\alpha} = 
\cup\{\bar a^{n,i}_{\alpha,k}:i < m,k < \omega\} \cup A^n_\alpha$ 
\nl

\hskip35pt for $m<n$ such that: if $\eta \in {}^n \omega$ and 
$m<n \Rightarrow \bar b^{n,m}_\alpha \subseteq$
\nl

\hskip35pt  $\cup\{\bar a^{n,i}_{\alpha,k}:i<m,k < \eta(i)\} \cup
A^n_\alpha$ then
\nl

\hskip35pt  $(p^{{\frak x}^n_\alpha} \restriction M^{{\frak x}^n_\alpha}) 
\cup\{\varphi_m(\bar a^{n,m}_{\alpha,\eta(m)},\bar
b^{n,m}_\alpha) \wedge \neg \varphi_m(\bar x,\bar
a^{n,m}_{\alpha,\eta(m)+1},\bar b^{n,m}_\alpha):$
\nl

\hskip35pt $m<n\}$ is finitely satisfiable in ${\frak C}$.
\sn
\item "{$(F)$}"  In checking clause (e) of $\circledast_{n+1}$
\sn
\ub{Case $\ell=1$}:  We know that $p^{{\frak x}^n_{\alpha +1}} \cup
\{\varphi_m(x,\bar a^{n,m}_{\alpha,k})^{\text{if}(k=0)}:m<n$ and $k<2\}$ is
consistent.  As ${\frak x}^n_{\alpha +1} \le^\ell_{\text{pr}} {\frak
z}^n_\alpha$ by clause $(\alpha)(d)$ of Definition \scite{dp1.3}(3) we
know that $q^{n+1}_\alpha := p^{{\frak z}^n_\alpha} \cup \{\varphi_m(x,\bar
a^{n,m}_{\alpha +1,k})^{\text{if}(k=0)}:m<n$ and $k<2\}$ is
consistent.
But $\varphi_n(x,\bar a^{n+1,m}_{\alpha,k}) = \varphi_n(x,\bar
a^{n,m}_{\alpha +1,k}) \in q^{n+1}_\alpha$ for $k<2,m<n$ and
$\varphi_n(x,\bar a^{n+1},m_{\alpha,k})^{\text{if}(k=0)} =
\varphi_n(x,\bar a^{n,*}_{\alpha,k})^{\text{if}(k=0)} \in
q^{n+1}_\alpha$ and $p^{{\frak x}^{n+1}_\alpha} \subseteq p^{{\frak
z}^n_\alpha} \subseteq q^{n+1}_\alpha$ hence 
$p^{{\frak x}^{n+1}_\alpha} \cup\{\varphi(x,\bar
a^{n,m}_{\alpha,k})^{\text{if}(k=0)}:m\le n$ and $k<2\}$ being a
subset of $q^{n+1}_\alpha$ is 
consistent, as required (this argument does not work for $\ell=4$).
\sn
\ub{Case 2}:  $\ell \in \{2,3,5,6\}$.
\nl
Straight.
\sn
\ub{Case 3}:  $\ell \in \{8,9\}$.
\endroster
\enddemo
\bn
As before
\demo{\stag{e3.8} Observation}  Like \scite{dp3.3} for
$\ell=1,2,3,4,5,6,8,9,11,12$. 
\enddemo
\bigskip

\definition{\stag{e3.9} Definition}  In Definition \scite{dr.6} we
allow $\ell = 17,18$.
\enddefinition
\bigskip

\demo{\stag{e3.10} Observation}  1) If ``${\frak y}$ explicitly
$\bar\Delta$-split $\ell(1)$-strongly over ${\frak x}$" then
``${\frak y}$ explicitly $\bar\Delta$-split $\ell(2)$-strongly over
${\frak x}$" when $(\ell(1),\ell(2)) \in
\{(15,14),(14,17),(18,17),(15,18)\} \cup \{(\ell,\ell +12):\ell = 2,3,5,6\}$.
\nl
2) If ${\frak x} \in K_{m,\ell(1)}$ then
dp-rk$^m_{\bar\Delta,\ell(1)}({\frak x}) \le 
\text{\rm dp-rk}^m_{\bar\Delta,\ell(2)}({\frak x})$ when
$(\ell(1),\ell(2))$ is as above.
\enddemo
\bigskip

\demo{Proof}  Easy by the definition.
\enddemo
\bigskip

\proclaim{\stag{e3.11} Claim}  1) In \scite{dr.8}(3) we allow $\ell =
17,18$.
\nl
2) $``\text{\rm dp-rk}_\ell(T) \ge |T|^+ \Rightarrow
   \kappa_{\text{ict}}(T) \ge \aleph_1$" we allow $\ell =14,15,17,18$.
\endproclaim
\bigskip

\proclaim{\stag{d4.1} Theorem}  In \scite{ind.1} we can allow
\mr
\item "{$(a)$}"  $\ell \in \{8,9,11,12\}$ and even $\ell \in
   \{14,15,17,18\}$.
\endroster
\endproclaim
\bigskip

\demo{Proof}  Similar to \scite{ind.1}.  \hfill$\square_{\scite{d4.1}}$
\enddemo
\bn
We can try to use ranks as in \S3 for $T$ which are just dependent.
In this case it is natural to revise the definition of the rank to
make it more ``finitary", say in Definition \scite{dp1.3}(4), clause
(e),(e)$'$ replace $\langle \bar a_k:k < \omega \rangle$ by a finite
long enough sequence.
\nl
Meanwhile just note that
\proclaim{\stag{dp1.5} Claim}  Let $\ell=1,2,3,5,6$ [and even $\ell =
14,15,17,18$].   For any finite
$\Delta \subseteq \Bbb L(\tau_T)$ we have: for every finite
$\Delta_1$, {\rm rk}$_{\Delta_1,\Delta,\ell}(T) = \infty$ \ub{iff} for
every finite $\Delta_1$, {\rm rk}$_{\Delta_1,\Delta,\ell}(T) \ge
\omega$ \ub{iff} some $\varphi(x,\bar y) \in \Delta$
has the independence property. 
\endproclaim
\bigskip

\demo{Proof}  Similar proof to \scite{dp1.4}, \scite{e3.6}.

Let $\langle \bar a_\alpha:\alpha < \omega \rangle
\subseteq M$ be indiscernible.

Let $\varphi(\bar x,\bar a_0),\neg \varphi(\bar x,\bar a_1) \in p$
exemplify ``$p$ splits strongly over $A_\varepsilon =
\cup\{M_{\alpha_\varepsilon}:\zeta < \varepsilon\} \cup A \cup B$ so
tp$(\bar a_0,A_\varepsilon) = \text{ tp}(\bar a_1,A_\varepsilon)$.
Let $A^+ = A \cup \bar a_0 \cup a_1$ and we find $u \subseteq
\{\alpha_\varepsilon:\varepsilon < \theta^+_1\}$ as required
\mr
\item "{$(*)$}"  there is $N^+ \prec M,\|N^*\| \le \theta$ such that
$N^* \prec N \prec M \Rightarrow \text{ dp-rk}(A,p \restriction (N^*
\cup A),N^*) = \text{ dp-rk}(A,p,M)$.  \hfill$\square_{\scite{dp1.5}}$
\endroster
\enddemo
\bn
\margintag{dp1.5.A}\ub{\stag{dp1.5.A} Question}:  1) Can such local ranks help us prove
some weak versions of ``every $p \in \bold S_\varphi(M)$ is
definable"?  (Of course, the first problem is to define such ``weak
definability"; see \cite[\S1]{Sh:783}).
\nl
2) Does this help for indiscernible sequences?
\bigskip

\definition{\stag{dp5.21} Definition}   We define
$K^x_{m,\ell}$ and 
$\text{dx-rk}^m_{\bar\Delta,\ell}$ for $x = \{p,c,q\}$ as follows:
\mr
\item "{$(A)$}"   \ub{for $x = p$}: as in Definition
\scite{dp1.3}(4),(5), \scite{e3.4}(4),(5)
\sn 
\item "{$(B)$}"  \ub{for $x=c$}:  as in Definition
\scite{dp1.3}(4),(5), \scite{e3.4}(4),(5) but
we demand that in clause (e),(e)$'$ of part (4) that $\{\varphi(\bar
x,\bar b_n):n < \omega\}$ is contradictory
\sn
\item "{$(C)$}"  \ub{for $x = q$}:  as in Definition
\scite{dp1.3}(4),(5), \scite{e3.4}(4),(5) but
clauses (e),(e)$'$ of part (4) we have $\bar a_\alpha$ from $A^{\frak
y}$ for $\alpha < \omega + \omega$ such that
\nl
$\{\varphi(x,a_\alpha)^{\text{if}(\alpha < \omega)}:\alpha < \omega + 
\omega\} \subseteq p^{{\frak x}'}$ and
in $(e')$ we have $\bar a_n$ from $A^{\frak y}$ and $\bold a_{\omega
+n}$ from $M^{{\frak x}'}$.  In details:
\sn
\item "{$(e)$}"    when $\ell \in \{1,2,3,4,5,6\}$, in 
$A^{\frak y}$ there is a $\Delta_1$-indiscernible sequence 
$\langle \bar a_k:k < \omega\rangle$ over $A^{\frak x} \cup M^{\frak
y}$ such that $\bar a_k \in {}^{\omega >}(M^{{\frak x}'})$ for $\alpha
< \omega$ and $\varphi(\bar x,\bar a_k),
\neg \varphi(\bar x,\bar a_{\omega +k}) \in 
p^{{\frak x}'}$ and $\bar a_k,\bar a_{\omega +k} 
\subseteq A^{\frak y}$ for $k < \omega$
\sn
\item "{$(e)'$}"  when $\ell = 8,9,11,12$ there are $\bar b,
\bar{\bold a}$ such that
{\roster
\itemitem{ $(\alpha)$ }  $\bar{\bold a} = \langle \bar a_i:i < \omega 
+ \omega \rangle$ is $\Delta_1$-indiscernible over $A^{\frak x} 
\cup M^{\frak y}$
\sn
\itemitem{ $(\beta)$ }  $A^{\frak y} \supseteq A^{\frak x} \cup \{\bar
a_i:i < \omega + \omega\}$; 
\sn
\itemitem{ $(\gamma)$ }  $\bar b \subseteq A^{\frak x}$ and $\bar a_i 
\in M^{{\frak x}'}$ for $i < \omega + \omega$
\sn
\itemitem{ $(\delta)$ }  $\varphi(\bar x,\bar a_k \char 94 \bar b) 
\wedge \neg \varphi(\bar x,\bar a_\omega \char 94 \bar b)$ belongs
\footnote{this explains why $\ell = 7,10$ are missing}
to $p^{{\frak x}'}$ for $k < \omega$.]]
\endroster}
\endroster
\enddefinition
\bn
\margintag{dp2.53}\ub{\stag{dp2.53} Question}:  Does Definition \scite{dp5.21} help
concerning question \scite{dp1.5.A}?
\bn
\margintag{dp5.22}\ub{\stag{dp5.22} Discussion}:  We can immitate \S3 with dc-rk or dq-rk
instead of dp-rk and use appropriate relatives of
$\kappa_{\text{ict}}(T)$.  But compare with \S4.
\bn
\centerline {$* \qquad * \qquad *$}
\bn
(B) \quad \ub{Minimal theories (or types)}:

It is natural to look for the parallel of minimal theories (see end of
the introduction).
\nl
A subsequent work of E. Firstenberg and the author \cite{FiSh:E50},
using \cite{Sh:757}, (see better \cite{Sh:E63})
considered a generalization of ``uni-dimensional stable $T$".
The generalization says (see \scite{dp0.13}(1))
\definition{\stag{dp3.2A} Definition}  1) $T$ is 
uni-dp-dimensional when: ($T$ is a dependent theory and)
every infinite non-trivial indiscernible sequences
of singletons $\bold I,\bold J$ have finite distance, see below.
\nl
2) (From \cite{Sh:93}) for indiscernibles sequences $\bold I,\bold J$
over $A$ we say that they are immediate $A$-neighbours 
\ub{if} $\bold I + \bold J$ is an indiscernible sequence or 
$\bold J + \bold I$ is an indiscernible sequence.  They have
distance $\le n$ if there are $\bold I_0,\dotsc,\bold I_n$ such that
$\bold I = \bold I_0,\bold J = \bold I_n$ and $\bold I_\ell,\bold
I_{\ell +1}$ are immediate $A$-neighbors (so indiscernible over $A$)
for $\ell < n$.  They are neighbors \footnote{we may prefer the local
version: for every finite $\Delta \subseteq \Bbb L(\tau_T)$ and finite
$A' \subseteq A$ (or $A' = A$) there are $\bold I',\bold J'$ realizing
the $\Delta$-type over $A'$ of $\bold I,\bold J$ respectively such
that $\bold I',\bold J'$ are (infinite) indiscernible sequences over
$A'$ (or $A$) and has distance over $A'$.} 
if they have distance $\le n$ for some $n$.
\nl
3) If $\bold I$ is an infinite indiscernible sequence over $A$ then
$\bold C_A(\bold I) = \cup\{\bold I':\bold I',\bold I$ have distance
$< \omega\}$.
\enddefinition
\bn
\margintag{dp3.2E}\ub{\stag{dp3.2E} Problem}:  1) Does uni-dp-dimensional theories have a
dimension theory? 
\nl
2) Can we characterize them?
\nl
3) If $p \in  \bold S^m(A)$, is there an indiscernible sequence $\bold
I \subseteq p({\frak C})$ based on $A$?, i.e. such that $\{F(\bold
C_A(\bold I)):F$ an automorphism of ${\frak C}$ over $A\}$ has
cardinality $< {\frak C}$ (equivalently $\le 2^{|T|+|A|}$) as is the case
for simple theories.
\bn
We can try another generalization.
\definition{\stag{dp0.14} Definition}  $T$ is dp$^\ell$-minimal when
dp-rk$^\ell({\frak x}) \le 1$ for every ${\frak x} \in K_\ell$,
i.e. $K_{m,\ell}$ for $m=1$.
\enddefinition
\bigskip

\demo{\stag{dp0.9} Hypothesis}  (till \scite{dp0.13.3})  Let $\ell$ be
as in Definition \scite{dp1.3}, \scite{e3.4}.
\enddemo
\bigskip

\remark{\stag{dp0.14.6} Remark}  For this property, $T$ and $T^{\text{eq}}$ may
differ.  Probably if we add only finitely many sorts, the ``finite
rank, i.e., dp-rk$^\ell({\frak x}) < n^* < \omega$ for every 
${\frak x} \in K_\ell$" is preserved.
\endremark
\bigskip

\demo{\stag{dp0.15} Observation}  $T$ is dp$^\ell$-minimal \ub{when}: for
every infinite indiscernible sequence $\langle \bar a_t:t \in
I\rangle,I$ complete, $\bar a_t \in {}^\alpha{\frak C}$ and element $c
\in {\frak C}$ there is $\{t\} \subseteq I$ as in \scite{dp1.2.4} (i.e., a
singleton or the empty set if you like) when $\ell \le 12$, and as in
\scite{dq.7} when $\ell \in \{14,\ldots\}$.
\enddemo
\bigskip

\demo{Proof}  Should be clear.  \hfill$\square_{\scite{dp0.15}}$
\enddemo
\bigskip

\proclaim{\stag{dp0.11} Claim}  1) For $\ell = 1,2$ we say
$T$ is {\rm dp}$^\ell$-minimal \ub{when}: 
there are no $\langle \bar a^i_n:n < \omega\rangle$ and 
$\varphi_i(x,\bar y_i)$ such that
\mr
\item "{$(a)$}"  for $i=1,2$, $\langle \bar a^i_n:n < \omega
\rangle$ is an indiscernible sequence over 
$\cup\{\bar a^{3 - i}_n:n < \omega\}$
\sn
\item "{$(b)$}"  for some $b \in {\frak C}$ we have
\nl
$\models \varphi_1(b,\bar a^1_0) \wedge \neg \varphi_2(b,\bar a^1_1)
\wedge \varphi_2(b,\bar a^2_0) \wedge \neg \varphi_2(b,\bar a^2_1)$.
\ermn
2) Similarly for {\rm rk-dp}$^\ell({\frak x}) \le n (< \omega)$.
\endproclaim
\bigskip

\demo{Proof}  Straight.
\enddemo
\bn
\margintag{dp0.12}\ub{\stag{dp0.12} Problem}: 1) Are dp$^\ell$-minimal theories $T$
similar to o-minimal theories?
\nl
2)  Characterize the dp$^\ell$-minimal theories of fields.
\nl
3) What are the implications between ``dp$^\ell$-minimal" for the
various $\ell$.
\nl
4) Above also for uni-dp-dimensionality.
\bigskip

\proclaim{\stag{dp0.13} Claim}  1) For $\ell=1,2$ the theory $T$,
{\rm Th}$(\Bbb R)$, the theory of real closed field is
uni-dp$^\ell$-dimensional; similarly for any o-minimal theory.
\nl
2) {\rm Th}$(\Bbb R)$ is {\rm dp}$^\ell$-minimal for $\ell=1,2$, similarly
for any o-minimal theory.
\nl
3) For prime $p$, the first order theory of the $p$-adic field is 
{\rm dp}$^1$-minimal. 
\endproclaim
\bigskip

\demo{Proof}  1) As in \cite{FiSh:E50}.
\nl
2) Repeat the proof in \cite[3.3]{Sh:783}(6).
\nl
3) By the proof of \scite{dp0.16}.  \hfill$\square_{\scite{dp0.13}}$
\enddemo
\bigskip

\remark{\stag{dp0.13.3} Remark}  If $T$ is a theory of valued fields
with elimination of field quantifier, see Definition
\scite{dt.7}(1),(2), and $k^{{\frak C}_T}$ is infinite this
fails.  But, if $\Gamma^{{\frak C}_T},k^{{\frak C}_T}$ are
dp$^1$-minimal \ub{then} the dp-rk for $T$ are $\le 2$.
\endremark
\bn
Another direction is:
\definition{\stag{dp5.16} Definition}  1) We say that a type $p(\bar
x)$ is content minimal \ub{when}:
\mr
\item "{$(a)$}"  $p(\bar x)$ is not algebraic
\sn
\item "{$(b)$}"  if $q(\bar x)$ extends $p(\bar x)$ and is not
algebraic then $\Phi_{q(\bar x)} = \Phi_{p(\bar x)}$, see below.
\ermn
2) $\Phi_{p(\bar x)} = \{\varphi(\bar x_0,\dotsc,\bar x_{n-1}):
\cup\{p(\bar x_\ell):\ell < n\} \cup \{\varphi(\bar x_1,\dotsc,\bar x_n)\}$ is
consistent, (see \cite{Sh:93}).
\enddefinition
\bn
\margintag{dp5.16.3}\ub{\stag{dp5.16.3} Question}:  Can we define reasonable dimension for
such types, at least for $T$ dependent or even strongly dependent?
\bn
\centerline{$* \qquad * \qquad *$}
\bn
(C) \quad \ub{Local ranks for super dependent and indiscernibles}:

Note that the original motivation of introducing ``strongly dependent"
in \cite{Sh:783} was to solve the equation: X/dependent =
superstable/stable.  However (the various variants) of strongly
dependent, when restricted to the family of stable theories, gives
classes which seem to me interesting but are not the class of
superstable $T$.  So the original question remains open.  Now 
returning to the search for ``super-dependent" we may consider another
generalization of superstable.
\definition{\stag{dp5.16.8} Definition}  1) We define
lc-rk$^m(p,\lambda) = \text{\rm lc-rk}^m_0(p,\lambda)$ for 
types $p$ which belongs to $\bold S^m_\Delta(A)$ for
some $A(\subseteq {\frak C})$ and finite $\Delta(\subseteq \Bbb L(\tau_T))$.

It is an ordinal or infinity and
\mr
\item "{$(a)$}"  lc-rk$^m(p,\lambda) \ge 0$ always
\sn
\item "{$(b)$}"  lc-rk$^m(p,\lambda) \ge \alpha = \beta +1$ \ub{iff} every
$\mu < \lambda$ there are finite $\Delta_1 \supseteq \Delta$ and
pairwise distinct $q_i \in \bold S^m_{\Delta_1}(A)$ extending $p$ such that
$i < 1 + \mu \Rightarrow \text{\rm lc-rk}^m(q_i,\lambda) \ge \beta$
\sn
\item "{$(c)$}"  lc-rk$^m(p,\lambda) \ge \delta,\delta$ a limit ordinal
\ub{iff} lc-rk$^m(p) \ge \alpha$ for every $\alpha < \delta$.
\ermn
2) For $p \in \bold S^m(A)$ let\footnote{Easily, if $\Delta_1
   \subseteq \Delta_2 \subseteq \Bbb L(\tau_T)$ are finite and $p_2
   \in \bold S^m_{\Delta_2}(A)$ and $p_1 = p_2 \rest \Delta_1$
  \ub{then} k-rk$^m(p_1) \ge \text{ lc-rk}^m(p_2)$.  So
   lc-rk$^m(p,\lambda)$ is well defined.}
 lc-rk$^m(p,\lambda)$ be min$\{\text{lc-rk}^m
(p,\lambda) \restriction \Delta:\Delta \subseteq \Bbb L(\tau_T)$
finite$\}$.
\nl
3) Let lc-rk$^m(T,\lambda) = \cup\{\text{lc-rk}^m(p,\lambda)+1:p \in
\bold S^m(A),A \subset {\frak C}\}$.
\nl
4) If we omit $\lambda$ we mean $\lambda = |T|^{++}$.
\enddefinition
\bn
\margintag{dw5.10}\ub{\stag{dw5.10} Discussion}:  There are other variants and they are
naturally connected to the existence of indiscernibles (for subsets of
${}^m{\frak C}$, concerning subsets of ${}^{|T|}{\frak C}$), probably
representability is also relevant (\cite{Sh:F705}).
\bigskip

\proclaim{\stag{dp5.16.9} Claim}  1) The following conditions on $T$
are equivalent (for all $\lambda > |T|^+$):
\mr
\item "{$(a)_\lambda$}"  for every $A$ and $p \in \bold S^m_\Delta(A)$ we have
{\rm lc-rk}$^m(p,\lambda) < \infty$
\sn
\item "{$(b)_\lambda$}"  for some $\alpha^* <|T|^+$ for every $A$ and $p \in
\bold S^m_\Delta(A)$ we have $\text{\rm lc-rk}^m(p,\lambda) < \alpha^*$
\sn
\item "{$(c)_\lambda$}"  there is no increasing chain $\langle
\Delta_n:n < \omega \rangle$ of finite subsets of $\Bbb L(\tau_T)$ and
$A$ and $\langle p_\eta:\eta \in {}^{\omega >}\lambda\rangle$ such
that $p_\eta \in \bold S^m_{\Delta_{\ell g(\eta)}}(A)$ and $\nu
\triangleleft \eta \Rightarrow p_\nu \subseteq p_\eta$ and if $\eta_1
\ne \eta_2$ are from ${}^n\lambda$ then $p_{\eta_1} \ne p_{\eta_2}$
\sn
\item "{$(c)_{\aleph_0}$}"  like $(c)_\lambda$ with $\langle
p_\eta:\eta \in {}^{\omega >}\omega \rangle$.
\ermn
2) Similarly restricting ourselves to $A = |M|$.
\endproclaim
\bigskip

\demo{Proof}  Easy.  \hfill$\square_{\scite{dp5.16.9}}$
\enddemo
\bn
Closely related is
\definition{\stag{dp5.16.10} Definition}   1) We define
$\text{\rm lc}_1-\text{rk}^m(p,\lambda)$ for types $p \in \bold S^m(A)$ for $A
\subseteq {\frak C}$ as an ordinal or infinitely by:
\mr
\item "{$(a)$}"  $\text{\rm lc}_1-\text{\rm rk}^m(p,\lambda)\ge 0$ always
\sn
\item "{$(b)$}"  $\text{\rm lc}_1-\text{\rm rk}^m(p,\lambda) \ge \alpha = \beta
+1$ \ub{iff} for every $\mu < \lambda$ and finite $\Delta \subseteq
\Bbb L(\tau_T)$ we can find pairwise distinct $q_i \in \bold S^m(A)$
for $i < 1 + \mu$ such that $p \restriction \Delta \subseteq q_i$ and
lc$_1-\text{\rm rk}^m(q_i,\lambda) \ge \beta$
\sn
\item "{$(c)$}"  $\text{\rm lc}_1-\text{\rm rk}^m(p,\lambda) \ge \delta$ for
$\delta$ a limit ordinal iff $\text{\rm lc}_1-\text{\rm rk}^m(p) 
\ge \alpha$ for every $\alpha < \delta$.
\ermn
2) If $\lambda = \beth_2(|T|)^{++}$ we may omit it.
\enddefinition
\bigskip

\proclaim{\stag{dp5.16.11} Claim}  1) The following conditions on $T$
are equivalent when $\lambda > \beth_2(|T|)^{++}$
\mr
\item "{$(a)_\lambda$}"   for every $A$ and $p \in \bold S^m(A)$ we have
$\text{\rm lc}_1-\text{\rm rk}^m(p,\lambda) < \infty$
\sn
\item "{$(b)_\lambda$}"  for some $\alpha^* < (2^{|T|})^+$ for every $A$ and
$p \in \bold S^m(A)$ we have lc$_1-\text{\rm rk}^m(p,\lambda) <
\alpha^*$
\sn
\item "{$(c)$}"  for no $A$ do we have a non-empty set ${\bold P}
\subseteq \bold S^m(A)$ such that for every $p \in {\bold P}$ and
finite $\Delta \subseteq \Bbb L(\tau_T)$ for some finite $\Delta_1$
the set $\{q \restriction \Delta_1:q \in {\bold P}$ and $q
\restriction \Delta = p \restriction \Delta\}$ has cardinality $>
\beth_2(|T|)^{++}$
\sn
\item "{$(d)_\lambda$}"  letting $\Xi = \cup\{\Xi_n:n < \omega\},\Xi_n
= \{\bar\Lambda:\bar\Lambda$ is a sequence of length $n$ of finite
sets of formulas $\varphi(\bar x,\bar y),\ell g(\bar x) = m\}$ there
is $\langle\Delta_{\bar\Lambda}:\bar\Lambda \in
\Xi\rangle$ where $\Delta_{\bar\Lambda}$ is a finite set of formulas
such that: for every $\lambda$ we can find $A$ and $\langle
p_{\bar\Lambda,\eta}:\bar\Lambda \in \Xi$ and $\eta \in {}^{\ell
g(\bar\Lambda)}\lambda\rangle$ such that:
{\roster
\itemitem{ $(\alpha)$ }  $p_{\bar\Lambda,\bar\eta} \in \bold S^m(A)$
\sn
\itemitem{ $(\beta)$ }   if $\bar\Lambda \in \Xi_n,\eta \in {}^n
\lambda$ and $\bar\Lambda' = \bar\Lambda \char 94 \langle
\Lambda_n\rangle \in \Xi_{n+1}$, then $p_{\bar\Lambda',\eta \char 94
<\alpha>} \restriction \bar\Lambda_n = p_{\bar\Lambda,\eta}
\restriction \Lambda_n$ for $\alpha < \lambda$ and $\langle
p_{\bar\Lambda',\eta \char 94 <\alpha>} \restriction
\Delta_{\bar\Lambda'}:\alpha < \lambda\rangle$ are pairwise distinct
\endroster}
\item "{$(e)_\lambda$}"  for some
$\langle\Delta_{\bar\Lambda}:\bar\Lambda \in \Xi\rangle$ as above the set
 $T \cup \Gamma_\lambda$ is consistent where $\Gamma$ is non-empty and:
{\roster
\itemitem{ $(\alpha)$ }  if $\bar\Lambda = \Xi_{n+1},\eta \in
{}^{n+1}\lambda$ and $\varphi(\bar x,\bar y) \in \Lambda_n$ then
$(\forall \bar y)[\dsize \bigwedge_{\ell < \ell g(\bar y)} P(y_\ell)
\rightarrow (\varphi(\bar x_{\bar\Lambda,\eta},\bar y) \equiv
\varphi(\bar x_{\bar\Lambda \restriction n,\eta \restriction n},\bar
y))]$
\sn
\itemitem{ $(\beta)$ }   if $\bar\Lambda \in \Xi_{n+1} \eta \in {}^n
\lambda$ and $\alpha < \beta < \lambda$, then $\dsize
\bigvee_{\varphi(x,\bar y)} \in \Delta_{\bar\Lambda}(\exists \bar
y)(\dsize \bigwedge_{\ell < \ell g(\bar y)} P(y_\ell) \wedge 
(\varphi(x_{\bar\Lambda,\eta \char 94<\alpha>}:\bar y) \equiv \neg
\varphi(\bar x_{\bar\Lambda,\eta \char 94 <\beta>},\bar y))$.
\endroster}
\ermn
2) Similarly restricting ourselves to the cases $A = |M|$, i.e. $A$ is
   the universe of some $M \prec {\frak C}$.
\endproclaim
\bigskip

\demo{Proof}  Similar.  \hfill$\square_{\scite{dp5.16.11}}$
\enddemo
\bigskip

\definition{\stag{dp5.16.12} Definition}  1) We define 
$\text{\rm lc}_2-\text{\rm rk}^m(p,\lambda)$, lc$_3-\text{\rm rk}^m(p,\lambda)$
like lc$_0-\text{\rm rk}^m(p,\lambda),\text{\rm lc}_1-\text{\rm rk}^m
(p,\lambda)$ respectively replacing ``$\Delta \subseteq \Bbb L(\tau_T)$ is
finite" by ``$\Delta \subseteq \Bbb L(\tau_T)$ and arity$(\Delta) <
\omega"$ where.
\nl
2) arity$(\varphi) =$ the number of free variables of $\varphi$,
arity$(\Delta) = \sup\{\text{arity}(\varphi):\varphi \in \Delta\}$ (if
we use the objects $\varphi(\bar x)$ we may use 
arity$(\varphi(\bar x)) = \ell g(\bar x))$.
\enddefinition
\bigskip

\proclaim{\stag{dp5.16.13} Claim}  The parallel of \scite{dp5.16.11}
for Definition \scite{dp5.16.12}.
\endproclaim
\bigskip

\remark{Remark}  Particularly the rank $\text{\rm lc}_3-\text{\rm rk}^m$ seems 
related to the existence of indiscernibility, i.e.
\endremark
\bn
\margintag{dp5.16.23}\ub{\stag{dp5.16.23} Conjecture}:  1) Assume, lc$_\ell$-rk$^m(T) <
\infty$ for some $\ell \le 3$.
We can prove (in ZFC!) that for
every cardinal $\mu$ for some $\lambda$ we have $\lambda \rightarrow
(\mu)_T$.
\nl
2) Moreover $\lambda$ is not too large, say is $\langle \beth_{\omega
+1}(\mu + |T|)$ (or just $< \beth_{(2^\mu)^+}$).
\bn
\centerline{$* \qquad * \qquad *$}
\bn
\head {(D) \quad Strongly$^2$ stable fields} \endhead
 \spuriousreset
\bigskip

A reasonable aim is to generalize the characterization of the
 superstable complete theories of fields.  Macintyre \cite{Ma71}
 proved that every infinite field whose first order theory is
 $\aleph_0$-stable, is algebraically closed.  Cherlin \cite{Ch78} proves
 that every infinite division ring whose first order theory in
 superstable is commutative, i.e. is a field so algebraically closed.  
Cherlin-Shelah \cite{ChSh:115}
prove ``any superstable theory Th$(K),K$ an infinite field is the theory of
algebraically closed fields" (and is true even for
division rings).  More generally we would like to replace stable by dependent
 and/or superstable by strongly dependent or at least strongly$^2$
 stable (or other variant).
\nl 
Of course, for strongly dependent we 
should allow at least the following cases (in addition to the
 algebraically closed fields): the first
order theory of the real field (not problematic as is the only one
with finite non-trivial Galois groups), the 
$p$-adic field for any prime $p$ and the first order theories covered by
 \scite{dp0.16}(2), i.e. Th$(K^{\Bbb F})$ for such $\Bbb F$.
\bn
So
\demo{\stag{fd.1} Conjecture}
\mr
\item "{$(a)$}"  if $K$ is an infinite field and $T = \text{ Th}(K)$ is
strongly$^2$ dependent (i.e., $\kappa_{\text{ict},2}(T) = \aleph_0$)
\ub{then} $K$ is an algebraically closed field (not strongly!!)
\sn
\item "{$(b)$}"  similarly for division rings
\sn
\item "{$(c)$}"  if $K$ is an infinite field and $T = \text{ Th}(K)$ is
strongly$^1$ dependent \ub{then} $K$ is finite or algebraically closed
or real closed or elementary equivalent to $K^{\Bbb F}$ for some $\Bbb
 F$ as in \scite{dp0.16}(2) (like the $p$-adics) or a finite algebraic
 extension of such a field
\sn
\item "{$(d)$}"    similarly to (c) for division rings.
\endroster
\enddemo
\bn
Of course it is even better to answer \scite{fd.2A}(1):
\nl
\margintag{fd.2A}\ub{\stag{fd.2A} Question}:  1) Characterize the fields with dependent
first order theory.
\nl
2) At least ``strongly dependent" (or another variant see (E),(F)
below).
\nl
3) Suppose $M$ is an ordered field and $T
= \text{ Th}(M)$ is dependent (or strongly dependent).  Can we
characterize?
\bigskip

\remark{Remark}  But we do not know this even for stability.
\nl
So adopting strongly dependent as our context we look what we can do.
\endremark
\bigskip

\proclaim{\stag{fd.2} Claim}    For a dependent $T$ and group $G$
interpreted in the monster model ${\frak C}$ of $T$; for every
$\varphi(x,\bar y) \in \Bbb L(\tau_T)$ there is $n_\varphi < \omega$
such that if $\alpha$ is finite $\langle \bar a_i:i < \alpha\rangle$
is such that $G \cap \varphi({\frak C},\bar a_i)$ is a subgroup of $G$
\ub{then} their intersection is the intersection of
some $\le n_\varphi$ of them.
\endproclaim
\bigskip

\remark{Remark}   If $T$ is stable this holds also for infinite
$\alpha$ by the Baldwin-Saxl \cite{BaSx76} theorem.
\endremark
\bigskip

\demo{Proof}  See \cite{Sh:F917}.
\enddemo
\bigskip

\proclaim{\stag{fd.3} Claim}  If the complete theory $T$ is
strongly$^2$ dependent \ub{then} ``finite kernel implies almost surjectivity"
which means that if in ${\frak C},G$ is a definable group, $\pi$ a
definable homomorphism from $G$ into $G$ with finite kernel then
$(G:\text{\rm Rang}(\pi))$ is finite.
\endproclaim
\bigskip

\demo{Proof}  By a general result from \cite[3.8=tex.ss.4.5]{Sh:783}
quoted here as \scite{0.gr.1}.  \hfill$\square_{\scite{fd.3}}$
\enddemo
\bigskip

\proclaim{\stag{fd.4} Claim}  Being strongly$^\ell$ dependent is
preserved under interpretation.
\endproclaim
\bigskip

\demo{Proof}  By \scite{dp1.2.1}, \scite{df2.3.51}.
\hfill$\square_{\scite{fd.4}}$ 
\enddemo
\bn
Hence the proof in \cite{ChSh:115} works ``except" the part on
``translating the connectivity", which rely on ranks not available here.  

However, if $T$ is stable this is fine hence we deduce
that we have
\demo{\stag{fd.5} Conclusion}  If $K$ is an infinite field and Th$(K)$ is
strongly$^2$ stable \ub{then} $T$ is algebraically closed.
\enddemo
\bigskip

\proclaim{\stag{fd.6.17} Claim}  Let $p$ be a prime.  $T$ is not
strongly dependent if $T$ is the theory of differentially closed
fields of characteristic $p$ \ub{or} $T$ is the theory of some
separably closed fields of characteristic $p$ which is not
algebraically closed.
\endproclaim
\bigskip

\demo{Proof}  The second case implies the first because if $\tau_1
\subseteq \tau_1,T_2$ a complete $\Bbb L(\tau_2)$-theory which is
strongly dependent then so is $T_1 = T_2 \cap \Bbb L(\tau_1)$. So let
$M$ be a $\aleph_1$-saturated separably closed field of characteristic
$p$ which is not algebraically closed.  Let $\varphi_n(x) = (\exists
y)(y^{p^n} = x)$ and $p_*(x) = \{\varphi_n(x):n < \omega\}$
and let $x E_n y$ mean $\varphi_n(x-y)$, so $E^M_n$ is an equivalent relation.

Let $\langle a_\alpha:\alpha < \omega_1\rangle$ be an indiscernible set
such that $\alpha < \beta < \omega_1 \Rightarrow a_\beta - a_\alpha
\notin \varphi_1(M)$.

Let $\psi_n(x,y_0,y_1,\dotsc,y_{n-1},z) = (\exists z)[\varphi_n(z) \wedge x
= y_0 + y^p_1 + \ldots + y^{p^{n-1}}_{n-1} + z]$.

Now by our understanding of Th$(M)$
\mr
\item "{$\circledast$}"  $(a) \quad$ if $b_\ell \in M$ 
for $\ell < n$ then $M \models (\exists x)\psi_n(x,b_0,\dotsc,b_{n-1})$
\sn
\item "{${{}}$}"  $(b) \quad$ in $M$ we have $\psi_{n+1}(x,y_0,\dotsc,y_n)
\vdash \psi_n(x,y_0,\dotsc,y_{n-1})$
\sn
\item "{${{}}$}"  $(c) \quad$ in $M$ we have, if 
$\psi_n(b,a_{\alpha_0},\dotsc,a_{\alpha_{n-1}})
\wedge \psi_n(b,a_{\beta_0},\dotsc,a_{\beta_{n-1}})$
\nl

\hskip25pt  then $\dsize \bigwedge_{\ell < n} \alpha_\ell = \beta_\ell$.
\ermn
[Why?  Clause (a) holds because if $b_\ell \in M$ for $\ell < n$
then $a = b_0 + b^p_1 + \ldots + b^{p^{n-1}}_{n-1}$ exemplifies
``$\exists x"$.  Clause (b) holds as if $M \models
\psi_{n+1}[a,b_0,\dotsc,b_{n-1},b_n]$ as witnessed by $z \mapsto d$,
then $M \models \psi_n[a,b_0,\dotsc,b_{n-1}]$ as witnessed by $z
\mapsto d + b^{p^n}_n$ which $\in \varphi_n(M)$ as $\varphi_n(M)$ is
closed under addition and $d \in \varphi_n(M)$ by $d \in
\varphi_{n+1}(M) \subseteq \varphi_n(M)$ and $b^{p^n}_n \in \varphi_n(M)$ as
$b_n$ witnesses it.  Lastly, to prove clause (c) assume that for $\ell=1,2$
we have $d^\ell = d^{p^n}_\ell \in \varphi_n(M),
b = a_{\alpha_0} + a^p_{\alpha_1} + a^{p^2}_{\alpha_2} + \ldots
+ a^{p^{n-1}}_{\alpha_{n-1}} + d^{p^n}_2$ and 
$b = a_{\beta_0} + a^p_{\beta_1} + a^{p^2}_{\beta_2} + \ldots + 
a^{p^{n-1}}_{\beta_{n-1}} + d^{p^n}_2$.  We prove this by induction on
$n$.  For $n=0$ this is trivial, $n=m+1$ substituting, etc., we get
$a_{\alpha_0} - a_{\beta_0} = (a^p_{\beta_1} - a^p_{\alpha_1}) + \ldots
+ (a^{p^{n-1}}_{\beta_{n-1}} - a^{p^{n-1}}_{\alpha_{n-1}}) + (d^{p^n}_2
- d^{p^n}_1) \in \varphi_1(M)$, so by an assumption on $\langle
a_\gamma:\gamma < \omega_1\rangle$ it follows that $\alpha_0 =
\beta_0$.  As there are unique $p$-th roots 
the original equation implies $a_{\alpha_1} +
a^p_{\alpha_2} + \ldots + a^{p^{n-2}}_{\alpha_{n-2}} + d^{p^n}_1 =
a_{\beta_1} + a^p_{\beta_2} + \ldots + a^{p^{n-2}}_{\beta_{n-2}} +
d^{p^n}_2$, and we use the induction hypothesis.] 

So together:
\mr
\item "{$\odot$}"  for every $\eta \in {}^\omega(\omega_1)$, there is
$b_\eta \in M$ such that
{\roster
\itemitem{ $(\alpha)$ }  $M \vdash
\psi_n(b_\eta,a_{\eta(0)},\dotsc,a_{\eta(n-1)})$ hence
\sn
\itemitem{ $(\beta)$ }  if $n < \omega,\nu \in {}^n(\omega_1),\nu
\ne \eta \restriction n$ then $M \models \neg
\psi_n(b_\eta,a_{\nu(0)},\dotsc,a_{\nu(m-1)})$.
\endroster}
\ermn
This suffices.  \hfill$\square_{\scite{fd.6.17}}$
\enddemo
\bn
\centerline{$* \qquad * \qquad *$}
\bn
(E) \quad \ub{On strongly$^3$ dependent}:

It is still not clear which versions of strong dependent (or
stable) will be most interesting.  Another reasonable version is
strongly$^3$ dependent and see more below.  
It has parallel properties and is natural.  Hopefully at least some
of those versions allows us to generalize weight (see \cite[V,\S3]{Sh:c}); 
we intend to return to it elsewhere.  
\nl
Meanwhile note:
\definition{\stag{ds5.1} Definition}  1) $T$ is strongly$^3$ dependent
if $\kappa_{\text{ict},3}(T) = \aleph_0$ (see below).
\nl
2) $\kappa_{\text{ict},3}(T)$ is the first $\kappa$ such that the
following \footnote{we may consider replacing $\delta$ by a linear
order and ask for $< \kappa$ cuts} holds:

if $\gamma$ is an ordinal,
$\bar a_\alpha \in {}^\gamma(M_{\alpha+1})$ for $\alpha < \delta,\langle
\bar a_\alpha:\alpha \in [\beta,\delta)\rangle$ is an indiscernible
sequence over $M_\beta$ for $\beta < \delta$ and $\beta_1 < \beta_2
\Rightarrow M_{\beta_1} \prec M_{\beta_2} \prec {\frak C}$ and $\bar c \in
{}^{\omega >}{\frak C}$ and cf$(\delta) \ge \kappa$ \ub{then} for some
$\beta < \kappa,\langle \bar a_\alpha:\alpha \in [\beta,\delta)\rangle$ is an
indiscernible sequence over $M_\beta \cup \bar c$.
\nl
3) We say $T$ is strongly$^\ell$ stable if $T$ is strongly$^\ell$ dependent
and is stable.
\nl
4) We define $\kappa_{\text{ict},3,*}(T)$ and strongly$^{3,*}$
dependent and strongly$^{3,*}$ stable as in the parallel 
cases (see Definition \scite{dp1.8},
\scite{dq2.8.4}), i.e., above we replace $\bar c$ by $\langle \bar
c_n:n < \omega \rangle$ indiscernible over $\cup\{M_\beta:\beta < \delta\}$. 
\enddefinition
\bigskip

\proclaim{\stag{dp5.4} Claim}    1) If $T$ is strongly$^{\ell+1}$
dependent \ub{then} $T$ is strongly$^\ell$ dependent for $\ell=1,2$.
\nl
2) $T$ is strongly$^\ell$ dependent iff $T^{\text{eq}}$ is; moreover
$\kappa_{\text{ict},\ell}(T) = \kappa_{\text{ict},\ell}(T^{\text{eq}})$.
\nl
3) If $T_1$ is interpretable in $T_2$ \ub{then}
 $\kappa_{\text{ict},\ell}(T_1) \le \kappa_{\text{ict},\ell}(T_2)$.
\nl
4) If $T_2 = \text{\rm Th}({\frak B}_{M,MA})$, see \cite[\S1]{Sh:783} 
and $T_1 =  \text{\rm Th}(M)$ \ub{then}
$\kappa_{\text{ict},\ell}(T_2) = \kappa_{\text{ict},\ell}(T_1)$.
\nl
5) $T$ is not strongly$^3$ dependent \ub{iff} we can find $\bar \varphi
= \langle \varphi_n(\bar x_0,\bar x_1,\bar y_n):
n < \omega \rangle,m = \ell g(\bar x_0))$ 
and for any infinite linear order $I$ we can find an
indiscernible sequence $\langle \bar a_t,\bar b_\eta:t \in I,\eta \in
{}^{\omega >} I$ increasing$\rangle$, see Definition \scite{dw5.1}
below such
that for any increasing sequence $\eta \in {}^\omega I$, the set
$\{\varphi_n(\bar x_0,\bar a_s,\bar b_{\eta \restriction n})^{\text{if}
(s=\eta(n))}:n < \omega$ and $\eta(n-1) <_I s \in I$ if $n>0\}$ of formulas is 
consistent (or use just $s =
\eta(n),\eta(n)+1$ or $\eta(n) \le_I s$, does not matter).
\nl
6) The parallel of parts (1)-(5) hold with strongly$^{3,*}$ instead of
strongly$^3$.  In particular, (parallel to part (5)), we have
$T$ is not strongly$^{3,*}$ dependent \ub{iff} we can find $\bar \varphi
= \langle \varphi_n(\bar x_0,\dotsc,\bar x_{k(n)},\bar y_n):
n < \omega \rangle,m = 
\ell g(\bar x))$ and for any infinite linear order $I$ we can find an
indiscernible sequence $\langle \bar a_t,\bar b_{\eta,t}:t \in I,\eta \in
{}^{\omega >} I$ increasing$\rangle$, see \scite{dw5.1} such
that for any increasing $\eta \in {}^\omega I,\{\varphi(\bar x_0,\bar a_s,
\bar b_{\eta \restriction n})^{\text{if}(s=\eta(n))}:n <
\omega$ and $\eta(n-1) <_I s$ if $n>0\} \cup 
\{\psi(\bar x_{i_0},\dotsc,\bar x_{i_{m-1}},\bar c) = 
\psi(\bar x_{j_0},\dotsc,\bar x_{j_{m-1}},\bar c):m < \omega,i_0 < \ldots <
 i_{m-1} < \omega,j_0 < \ldots < j_{m-1} < \omega$ and $\bar c
 \subseteq \cup\{\bar a_s,b_\rho:s \in I,\rho \in 
{}^{\omega >}I\text{ increasing}\}\}$ is consistent.
\endproclaim
\bigskip

\demo{Proof}  1)-4).  Easy.
\nl
5),6) As in \cite{Sh:897}.  \hfill$\square_{\scite{dp5.4}}$
\enddemo
\bn
Recall this 
definition applies to stable $T$ (i.e. Definition \scite{ds5.1}(3)).
\demo{\stag{dp5.11} Observation}  The theory $T$ is strongly$^3$
stable \ub{iff}: $T$ is stable and we cannot find 
$\langle M_n:n < \omega \rangle,\bar c \in 
{}^{\omega >}{\frak C}$ and $\bar{\bold a}_n \in
{}^\omega(M_{n+1})$ such that:
\mr
\item "{$(a)$}"  $M_n$ is $\bold F^a_\kappa$-saturated
\sn
\item "{$(b)$}"  $M_{n+1}$ is $\bold F^a_\kappa$-prime over $M_n \cup
\bar{\bold a}_n$
\sn
\item "{$(c)$}"  tp$(\bar{\bold a}_n,M_n)$ does not fork over $M_0$
\sn
\item "{$(d)$}"  tp$(\bar c,M_n \cup \bar{\bold a}_n)$ forks over
$M_n$.
\endroster
\enddemo
\bigskip

\demo{Proof}  Easy.  \hfill$\square_{\scite{dp5.11}}$
\enddemo
\bn
\margintag{dp5.11.4}\ub{\stag{dp5.11.4} Conjecture}  For strongly$^3$ stable $T$ we 
have dimension theory (including weight)
close to the one for superstable theories (as in \cite[V]{Sh:c}), we may
try to deal with it in \cite{Sh:839}; related to \S5(G) below.
\bn
(F) \quad \ub{Representability and strongly$_4$ dependent}:

In \cite{Sh:897} we deal with $T$ being fat or lean.  We say a class
$K$ of models is fat when for every ordinal $\alpha$ there are a
regular cardinal $\lambda$ and non-isomorphic models $M,N \in
K_\lambda$ which are EF$^+_{\alpha,\lambda}$-equivalent where
EF$^+_{\alpha,\lambda}$ is a strong version of ``the isomorphism
player has a winning strategy in a strong version of the
Ehrenfuecht-Fr\"asse game of length $\lambda$".  We prove there, that
consistently if $T$ is not strongly stable and $T_1
\supseteq T$, then PC$(T_1,T)$ is fat
(in a work in preparation \cite{Sh:F918} we show 
that it suffices to assume ``$T$ is not strongly$_4$-stable"; see below).

In \cite{Sh:F705}, a work under preparation, we 
shall deal with representability.  The weakest form (for
${\frak k}$ a class of index models, e.g. linears order) is:  an
e.g. first order $T$ is weakly ${\frak k}$-represented when for every
model $M$ of $T$ and say finite set $\Delta \subseteq \Bbb L(\tau_T)$ we
can find an index model $I \in {\frak k}$ and sequence
$\langle \bar a_t:t \in I\rangle$ of finite 
sequences from $M^{\frak C}$ (or just
singletons) which is $\Delta$-indiscernible, i.e., see below, such that
$|M| \subseteq \{a_t:t \in I\}$.

This is a parallel to stable and superstable when we play with
essentially the arity of the functions of ${\frak k}$ and the size of
$\Delta$'s considered.  
The thesis is that $T$ is stable \ub{iff} it, essentially can
be represented for essentially ${\frak k}$ the class of sets and parallel
representability for ${\frak k}$ derived for order characterize
versions of the class of dependent theories. We also 
define ${\frak k}$-forking, i.e.
replace linear orders other index set.
\bn
We define
\definition{\stag{dw5.1} Definition}  1) For any structure $I$ we say
that $\langle \bar a_t:t \in I\rangle$ is indiscernible (in ${\frak
C}$ over $A$) when: $\ell g(\bar a_t)$ depends only on the quantifier
type of $t$ in $I$ and:
\sn
\block
if $n < \omega$ and $\bar s = \langle
s_0,s_1,\dotsc,s_{n-1}\rangle,\bar t = \langle t_0,\dotsc,t_{n-1}\rangle$
realize the same quantifier-free type in $I$ \ub{then} $\bar a_{\bar t} := 
\bar a_{t_0} \char 94 \ldots \char 94 \bar a_{t_{n-1}}$ and $\bar a_{\bar
s} = \bar a_{s_0} \char 94 \ldots \char 94 \bar a_{s_{n-1}}$ realize
the same type (over $A$) in ${\frak C}$.
\endblock
\sn
2) We say that $\langle \bar b_u:u \in [I]^{< \aleph_0}\rangle$ is
indiscernible (in ${\frak C}$) (over $A$) similarly:
\sn
\block
if $n < \omega,w_0,\dotsc,w_{m-1} \subseteq \{0,\dotsc,n-1\}$ and $\bar s =
\langle s_\ell:\ell < n\rangle,\bar t = \langle t_\ell:\ell <
n\rangle$ realize the same quantifier-free types in $I$ and $u_\ell =
\{s_k:k \in w_\ell\},v_\ell = \{t_k:k \in w_\ell\}$ \ub{then} $\bar a_{u_0}
\char 94 \ldots \char 94 \bar a_{u_{n-1}},\bar a_{v_0} \char 94 \ldots
\bar a_{v_{n-1}}$ realize the same type in ${\frak C}$ (over $A$).
\endblock
\nl
3) We may use incr$(< \omega,I)$ instead of $[I]^{< \aleph_0}$ where
incr$({}^\alpha I) = \text{ incr}_\alpha(I) = \text{ incr}(\alpha,I) = 
\{\rho:\rho$ is an increasing sequence of length $\alpha$ of 
members of $I\}$; we can use $< \alpha$ or $\le \alpha$; clearly the
difference between incr$(< \omega,I)$ and $[I]^{<\aleph_0}$ is
notational only (when we have order).
\enddefinition
\bigskip

\definition{\stag{dw5.34} Definition}  1) We say that the $m$-type 
$p(\bar x)$ does $(\Delta,n)$-ict divide over $A$ (or 
$(\Delta,n)$-ict$^1$ divide over $A$) \ub{when}: there 
are an indiscernible sequence $\langle \bar a_t:
t \in I\rangle,I$ an infinite linear order and $s_0 <_I t_0 \le_I s_1
<_I t_1 <_I \ldots \le_I s_{n-1} <_I t_{n-1}$ such that 
\mr
\item "{$\circledast_1$}"  $p(\bar x) \vdash 
``\text{\rm tp}_\Delta(\bar x \char 94 \bar a_{s_\ell},A) \ne \text{\rm
tp}_\Delta(\bar x \char 94 \bar a_{t_\ell},A)"$ for $\ell < n$.
\ermn
2) We say that the $m$-type $p(\bar x)$ does
$(\Delta,n)$-ict$^2$-divides over $A$ \ub{when} above we replace
 $\circledast_1$ by:
\mr
\item "{$\circledast_2$}"  $p(\bar x) \vdash 
``\text{\rm tp}_\Delta(\bar x \char 94 \bar a_{s_\ell},\cup 
\{\bar a_{s_k}:k < \ell\} \cup A) \ne \text{\rm tp}_\Delta(\bar x 
\char 94 \bar a_{t_\ell},\cup\{\bar a_{s_k}:k < \ell\} \cup A)"$ 
for $\ell < n$.
\ermn
3) We say that the $m$-type $p(\bar x)$ does
$(\Delta,n)$-ict$^3$-divide over $A$ \ub{when} above $(\langle \bar a_t:t
\in I \cup \text{\rm incr}(< n,I)\rangle$ is indiscernible over $A$
and we replace $\circledast_1$ by
\mr
\item "{$\circledast_3$}"  $p(\bar x) \vdash 
 ``\text{\rm tp}_\Delta(\bar x \char 94 \bar a_{s_\ell},\bar a_{\langle
s_0,\dotsc,s_{\ell -1}\rangle} \cup A) \ne \text{\rm tp}_\Delta
(\bar x \char 94 \bar a_{t_\ell},\bar a_{\langle
s_0,\dotsc,s_{\ell-1}\rangle} \cup A)"$ for $\ell < n$.
\ermn
4) We say that the $m$-type $p(\bar x)$ does
$(\Delta,n)$-ict$^4$-divide over $A$ \ub{when} there are $n^* <
\omega$ and sequence $\langle \bar a_\eta:\eta \in 
\text{\rm inc}(\le n^*,I)\rangle$ indiscernible over $A$ 
such that (where comp$(I)$ is the completion of the linear order $I$):

if $\bar c$ realizes $p(\bar x)$ then for no set
$J \subseteq \text{ comp}(I)$ with $\le n$ members, the sequence $\langle
\bar a_\eta:\eta \in \text{\rm inc}(\le n^*,I^+)\rangle$ is
$\Delta$-indiscernible
over $A$ where $I^+ = (I,P_t)_{t \in J}$ and $P_t := \{s \in I:s < t\}$.  
Note that if $T$ is stable, we can 
equivalently require $J \subseteq I$.
\nl
5) For $k \in\{1,2,3,4\}$ we 
say that the $m$-type $p(\bar x)$ does $(\Delta,n)$-ict$^k$-forks
over $A$ \ub{when} for some sequence $\langle \psi_\ell(\bar x,\bar
a_\ell):\ell < \ell(*) < \omega \rangle$ we have
\mr
\item "{$(a)$}"  $p(\bar x) \vdash \dsize \bigvee_{\ell < \ell(*)}
\psi_\ell(\bar x,\bar a_i)$
\sn
\item "{$(b)$}"  $\psi_\ell(\bar x,\bar a_\ell)$ does
 $(\Delta,n)$-ict$^k$-divide over $A$.
\ermn
If $k=1$ we may omit it, if $\Delta = \Bbb L(\tau_T)$ we may omit it.
\nl
6) We define ict$^k-\text{rk}^m(p)$, an ordinal or $\infty$, as
 follows (easily well defined): ict$^k-\text{rk}^m(p) \ge
\alpha$ \ub{iff} $p$ is an $m$-type and for every finite $q \subseteq
p$, finite $A \subseteq \text{ Dom}(p)$ and $n < \omega$ and $\beta <
\alpha$ there is an $m$-type $r$ extending $q$ which 
$(\Bbb L(\tau_T),n)-\text{ict}^k$-forks over $A$ with ict$^k$-rk$^m(r) \ge
\beta$.  If ict$^k$-rk$^m(r) \ngeq \beta +1$; and we say that $n$
witnesses this if the demand above for this $n$ fails.  
If $n+1$ is the minimal witness
let $n = \text{ ict}^k-\text{wg}^n(r)$.
\nl
7) $\kappa^m_{k,\text{ict}}(T)$ is the first $\kappa \ge \aleph_0$ such
that for every $p \in \bold S^m(B),B \subseteq {\frak C}$ 
there is a set $A \subseteq B$ of cardinality $< \kappa$ such that
$p$ does not {\rm ict}$^k$-fork over $A$.  Omitting $m$ means for some
$m < \omega$; note that we write $\kappa_{k,\text{ict}}(T)$ to
distinguish it from Definition \scite{df2.3.5} of $\kappa_{\text{ict},2}$.
\nl
8) $T$ is strongly$_k$ dependent [stable] if $\kappa_{k,\text{ict}}(T)
 =\aleph_0$ [and $T$ is stable].
\nl
9) We define $\kappa_{k,\text{ict},*}(T)$ parallely i.e., now $p(\bar x)$
is the type of an indiscernible sequence of $m$-tuples and $T$ is
strongly$_{k,*}$ dependent [stable] if it is dependent [stable] and
$\kappa_{k,\text{ict},*}(T) = \aleph_0$.
\enddefinition
\bigskip

\proclaim{\stag{dw5.35} Claim}  1) For dependent $T$, the following
conditions are equivalent:
\mr
\item "{$(a)$}"  $\kappa_{4,\text{ict},*}(T) > \aleph_0$, see Definition
\scite{dw5.34}(4),(7),(9)
\sn
\item "{$(b)$}"   there are $m,\langle(\Delta_\ell,n_\ell):\ell <
\omega\rangle,I,{\bold J}$ such that
{\roster
\itemitem{ $(\alpha)$ }  $\Delta_\ell \subseteq \Bbb L(\tau_T)$ finite
and $n_\ell < \omega$ and $n_\ell > \ell$ for $\ell < \omega$
\sn  
\itemitem{ $(\beta)$ }   $I$ is an infinite linear order with
increasing $\omega$-sequence of members
\sn
\itemitem{ $(\gamma)$ }  ${\bold J} = 
\langle \bar a_\rho:\rho \in \text{\rm inc}_{< \omega}(I)\rangle$ is an 
indiscernible sequence with $\bar a_\rho \in {}^\omega {\frak C}$
\sn
\itemitem{ $(\delta)$ }  for $\eta \in {}^\omega I$ an increasing
sequence, for some $\bar c_\ell \in {}^m{\frak C}(\ell < \omega)$ we
have:
\sn  
\itemitem{ ${{}}$ }  $(i) \quad \langle \bar c_\ell:\ell < \omega
\rangle$ is an indiscernible sequence over 
\nl

\hskip35pt $\cup\{\bar a_\rho:\rho \in \text{\rm incr}(I,< \omega)\}$
\sn  
\itemitem{ ${{}}$ }  $(ii) \quad$ if $J$ is the completion of the
linear order $I$ \ub{then} for no
\nl

\hskip35pt  finite $J_0 \subseteq J$ do we
have: if $n < \omega$ and $\rho^\ell_0,\dotsc,\rho^\ell_{n-1} \in$
\nl

\hskip35pt $\text{\rm incr}(I,< \omega)$ for $\ell=1,2$ are such that $\rho^1_0
\char 94 \ldots \char 94 \rho^1_{n-1}$ and $\rho^2_0 \char 94 \ldots
\char 94 \rho^2_{n-1}$
\nl

\hskip35pt  realizes the same quantifier free type over
$J_0$ in $J$ then $\bar a_{\rho^1_0} \char 94 \ldots \char 94 \bar
a_{\rho^1_{n-1}}$,
\nl

\hskip35pt $\bar a_{\rho^2_0} \char 94 \ldots \char 94 \bar
a_{\rho^2_{n-1}}$ realize the same $\Delta_\ell$-type 
over $\cup\{\bar c_\ell:\ell < \omega\}$ in ${\frak C}$ 
\endroster}
\item "{$(c)$}"   the natural rank is always $< \infty$.
\ermn
2) For dependent $T$ the following conditions are equivalent
\mr
\item "{$(a)$}"  $\kappa^m_{4,\text{ict}}(T) > \aleph_0$
\sn
\item "{$(b)$}"  like (b) is part (1) only $\langle \bar c_\ell:\ell <
\omega\rangle$ is replaced by one $m$-tuple $\bar c$
\sn
\item "{$(c)$}"  \text{\rm ict}$^4-\text{\rm rk}^m(\bar x = \bar x) =
\infty$
\sn
\item "{$(d)$}"  \text{\rm ict}$^4-\text{\rm rk}^m(\bar x = \bar x)
\ge |T|^+$.
\ermn
3) Similarly (just simpler) for $k=1,2,3$ instead $4$.
\endproclaim
\bigskip

\demo{Proof}  Straight, but see details Cohen-Shelah
\cite{CoSh:919}.   \hfill$\square_{\scite{dw5.35}}$
\enddemo
\bn
\margintag{dw.39}\ub{\stag{dw.39} Question}:  1) Can we characterize the $T$ such that the
ict$^k$-rk$^1$ rank of the formula $x=x$ is $1$?
\nl
2) Do we have ict$^\ell$-rk$^m(\bar x = \bar x) = \infty$ iff
   ict$^\ell$-rk$^1(x=x) = \infty$, i.e. can we in part (2) say
   that the properties do not depend on $m$?  The positive answer will
appear in Cohen-Shelah \cite{CoSh:919}.
\bn
Now
\demo{\stag{dw5.19} Observation}  1) For $k=1,2,3$ if $p(\bar x)$
does $(\Delta,n)$-ict$^k$ forks over $A$ \ub{then} $p(\bar x)$ does
$(\Delta,n)$-ict$^{k+1}$ forks over $A$.
\nl
2) If $T$ is strongly$_{k +1}$ dependent/stable \ub{then} $T$ is
strongly$_k$ dependent/stable.
\nl
3) For $k \in \{1,2,3,4\}$ if $T$ is strongly$_k$ dependent/stable
\ub{then} $T$ is strongly dependent/stable; if $T_1$ is interpretable
in $T_2$ and $T_2$ is strongly$_k$ dependent/stable \ub{then} so
is $T_1$.
\nl
4) Assume $T$ is stable.
If $p \in \bold S^m(B)$ does not fork over $A \subseteq B$
\ub{then} ict$^k$-rk$^m(p) = \text{ ict}^k-\text{rk}^m(p \restriction A)$.
\enddemo
\bigskip

\remark{Remark}  Also the natural inequalities concerning
itc$_k$-rk$^n(-)$ follows by \scite{dw5.19}(1).
\endremark
\bigskip

\demo{Proof}  Straight.  \hfill$\square_{\scite{dw5.19}}$
\enddemo
\bn
\margintag{dw5.23}\ub{\stag{dw5.23} Example}:  1) There is a stable NDOP, NOTOP, not
multi-dimensional countable
complete theory which is not strongly$^2$ dependent.
\nl
2) $T = \text{ Th}({}^{\omega_1}(\Bbb Z_2),E_n)_{n < \omega}$ is as
   above where $\Bbb Z_2 = \Bbb Z/2 \Bbb Z$ as an additive group, $E_n
   = \{(\eta,\nu):\eta,\nu \in {}^{\omega_1}(\Bbb Z_2)$ are such that
$\eta \rest (\omega n) = \nu \rest (\omega n)$.
\nl
3) As in part (1) but $T$ is not strongly dependent.
\bigskip

\remark{Remark}  This is \cite[0.2=0z.5]{Sh:897}.  It shows that the
   theorem there adds more cases. 
\endremark
\bigskip

\demo{Proof}  1) By part (2).
\nl
2) So let $M_0$ be the additive group $({}^{\omega_1}(\Bbb Z_2),+)$
   where $+$ is coordinatewise addition and for $\alpha \le \omega$
   let $M_\alpha = ({}^{\omega_1}(\Bbb Z_2),P_n)_{n < \alpha}$, where
$P_n = \{\eta \in {}^{\omega_1}(\Bbb Z_2):\eta \rest (\omega n)\}$
   is constantly zero and $E_n = \{(\eta,\nu):\eta,\nu \in 
{}^{\omega_1}(\Bbb Z_2)$ are such that $\eta \rest (\omega n) = \nu
   \rest (\omega n)\}$ and $M'_\alpha = ({}^{\omega_1}(\Bbb Z_2),E_n)_{n
   < \alpha}$.  So $M'_\alpha,M_\alpha$ are bi-interpretable, so we
   shall use $M_\alpha$.  Let $T = \text{ Th}(M_\omega)$ and let
   $T_\alpha = \text{ Th}(M_\alpha)$.  So for a model $N$ of
   $T_\alpha$ is just an abelian group in which every element has
   order 2, with distinguished subgraph $P^N_n$ for $n < \alpha$ so a
   vector space over the field $\Bbb Z_2$.
\enddemo
\bn
\ub{$T$ is stable}:

For $n < \omega$, a model of $T_n$ is determined by finitely many
dimensions: $(P^N_k:P^N_{k+1})$ for $< n$ ($E^N_n$ is interpreted as
the equality), so $T_n$ is superstable not multi-dimensional.

Hence $T$ necessarily is stable.
\bn
\ub{$T$ is strongly dependent not strongly$^2$ dependent}:

As in \scite{dq1.8}, in fact it is strongly dependent.
\bn
\ub{$T$ is not multi-dimensional}:

If $N$ is an $\aleph_1$-saturated model of $T$ then it is determined
by the following dimension as vector spaces over $\Bbb Z_2$, for $n <
\omega$
\mr
\item "{$(*)_1$}"  $P^N_n/P^N_{n+1}$
\sn
\item "{$(*)_2$}"  $\dbca_{n < \omega} P^N_n$.
\ermn
Each corresponds to a regular type (in ${\frak C}^{\text{eq}}_T$).
\bn
\ub{$T$ has NDOP}:

Follows from uni-dimensionality.
\bn
\ub{$T$ has NOTOP}:

Assume $N_\ell \prec {\frak C}_T$ is $\aleph_1$-saturated, $N_0 \prec
N_\ell$ for $\ell-0,1,2$ such that tp$(N_1,N_2)$ does not fork over
$N_0$.  Let $A$ be the subgroup of ${\frak C}$ generated by $N_1 \cup
N_2$ and let $N_3 = {\frak C}_T \rest A$.  Easily $N_3 \prec {\frak
C}_T$, moreover $N_3$ is $\aleph_1$-saturated.

By \cite[XII]{Sh:c} this suffices.
\nl
3) Expand $M_\alpha$ by $Q_m = \{\eta \in {}^{\omega_1}(\Bbb Z_2):\eta
\rest [\omega m,\omega m + \omega)$ is constantly zero$\}$ for 
$m < n$.  \hfill$\square_{\scite{dw5.23}}$
\bn
(G) \quad \ub{strong$_3$ stable and primely minimal types}
\bigskip

\demo{\stag{6n.7} Hypothesis}  $T$ is stable (during \S5(G)).
\enddemo
\bigskip

\definition{\stag{6n.14} Definition}  [$T$ stable]  
We say $p \in \bold S^\alpha(A)$ is
primely regular (usually $\alpha < \omega$) 
\ub{when}: if $\kappa > |T| + |\alpha|$ is a regular cardinal, the model
$M$ is $\kappa$-saturated, the type tp$(\bar a,M)$ is parallel to $p$
(or just a stationarization of it)
and $N$ is $\kappa$-prime over $M+ \bar a$ and $\bar b \subseteq
{}^{\kappa >} N \backslash {}^{\kappa >} M$ \ub{then} tp$(\bar a,M +
\bar b)$ is $\kappa$-isolated, equivalently \footnote{because $N$ is
$\kappa$-prime over $M + \bar a + \bar c$ whenever $\bar c \in
{}^{\kappa >} N$} $N$ is $\kappa$-prime over $M + \bar b$.
\enddefinition
\bigskip

\proclaim{\stag{6n.21} Claim}  1) Definition \scite{6n.14} to
equivalent to: there are $\kappa,M,\bar a,N$ as there.
\nl
2) We can in part (1) replace ``$\kappa > |T| + |\alpha|$ regular,
$\kappa$-prime" by ``{\rm cf}$(\kappa) \ge \kappa(T),\bold
F^a_\kappa$-prime" respectively. 
\endproclaim
\bigskip

\demo{Proof}  Straight. \hfill$\square_{\scite{6n.21}}$
\enddemo
\bn
Now (recalling Definition \scite{ds5.1} and Observation \scite{dp5.11}).
\proclaim{\stag{6n.28} Claim}  [$T$ is strongly$_3$ stable]

If {\rm cf}$(\kappa) \ge \kappa_r(T)$ and $M \prec N$ are $\bold
F^a_\kappa$-saturated \ub{then} for some $a \in N \backslash M$ the
type {\rm tp}$(a,M)$ is primely regular.
\endproclaim
\bigskip

\demo{Proof}  The reader can note that by easy manipulations \wilog \,
$\kappa = \text{ cf}(\kappa)> |T|$; in fact, by this 
we can use {\rm tp} instead of {\rm stp}, etc.

Let $\alpha_* = \text{ min}\{\text{ict}^3-\text{rk}
(\text{tp}(a,M)):a \in N \backslash M\}$ and let
$a \in N \backslash M$ and $\varphi_*(x,\bar d_*) \in \text{
tp}(a,M)$ be such that $\alpha_* = 
\text{ ict}^3-\text{rk}(\{\varphi_*(x,\bar d_*)\})$.

Let $a \in N \backslash M$.  We try to choose $N_\ell,a_\ell,B_\ell$ by 
induction on $\ell < \omega$ such that
\mr
\item "{$\boxplus_\ell$}"  $(a) \quad M \prec N_\ell \prec N$ and $a_\ell
\in N_\ell \backslash M$
\sn
\item "{${{}}$}"  $(b) \quad N_\ell$ is $\bold F^a_\kappa$-primary over
$M + a_\ell$ and $a_0 = a$
\sn
\item "{${{}}$}"  $(c) \quad$ if $\ell = m+1$ \ub{then}
{\roster
\itemitem{ ${{}}$ }  $(\alpha) \quad N_\ell \prec N_m$ and 
{\rm tp}$(a_m,M+a_\ell)$ is not $\bold F^a_\kappa$-isolated
\sn
\itemitem{ ${{}}$ }  $(\beta) \quad N_m$ is $\bold F^a_\kappa$-primary 
over $N_\ell + a_m$
\sn
\itemitem{ ${{}}$ }  $(\gamma) \quad N_\ell$ is $\bold
F^a_\kappa$-constructible over $N_{\ell +1} + a_0$.
\endroster}
\item "{${{}}$}"  $(d)(\alpha) \quad B_\ell \subseteq N_\ell$
\sn
\item "{${{}}$}"  $\quad (\beta) \quad a_\ell \in B_\ell$
\sn
\item "{${{}}$}"  $\quad (\gamma) \quad |B_\ell| < \kappa$
\sn
\item "{${{}}$}"  $\quad (\delta) \quad$ every $\bold F^a_\kappa$-isolated 
type $q \in \bold S^{< \omega}(M \cup B_\ell)$ has no extension in
\nl

\hskip25pt  $\bold S^{< \omega}(M \cup \bigcup\{B_m:m \le \ell\})$
 which forks over $M \cup B_\ell$
\sn
\item "{${{}}$}"  $\quad (\varepsilon) \quad B_\ell$ is 
$\bold F^a_\kappa$-atomic over $M + b_\ell$.
\ermn
Let $(N_\ell,a_\ell)$ be defined iff $\ell < 1 + \ell(*) \le \omega$,
clearly $\ell(*) \ge 0$.  
\mr
\item "{$\boxtimes_1$}"  if $\ell(*) < \omega$ then
tp$(a_{\ell(*)},M)$ is primely regular.
\ermn
[Why?  If not, then for some $b \in N_{\ell(*)} \backslash M$ we have
tp$(a_{\ell(*)},M+b)$ is not $\bold F^a_\kappa$-isolated.  

We try to choose $\bar b'_\varepsilon$ by induction on $\varepsilon <
\kappa$ such that
\mr
\item "{$(\boxtimes_{1.1})$}"  $(\alpha) \quad \bar b_0 = \langle b \rangle$
\sn
\item "{${{}}$}"  $(\beta) \quad 
\bar b'_\varepsilon \in {}^{\omega >}(N_{\ell(*)})$
\sn
\item "{${{}}$}"  $(\gamma) \quad$ tp$(\bar b'_\varepsilon,M \cup \bigcup\{\bar
b'_\zeta:\zeta < \varepsilon\} \cup \{b\}\}$ is $\bold F^a_\kappa$-isolated
\sn
\item "{${{}}$}"  $(\delta) \quad$ 
tp$(\bar b'_\varepsilon, M \cup \bigcup\{\bar
b'_\zeta:\zeta < \varepsilon\} \cup \{b,a_k,\dotsc,a_{\ell(*)}\}$ is
$\bold F^a_\kappa$-isolated for
\nl

\hskip25pt  $k = \ell(*),\dotsc,0$
\sn
\item "{${{}}$}"  $(\varepsilon) \quad$ tp$(\bar a,M \cup 
\bigcup\{\bar b'_\zeta:\zeta \le \varepsilon\})$ 
forks over $M \cup \bigcup\{\bar b_\zeta:\zeta < \varepsilon\}$ 
for some
\nl

\hskip25pt  $\bar a \in {}^{\omega >}(B_{\ell(*)})$ when $\varepsilon > 0$.
\ermn
We are stuck for some $\varepsilon(*) < \kappa$ because
$|B_{\ell(*)}| < \kappa$ and let $B' = 
\cup\{\bar b'_\varepsilon:\varepsilon < \varepsilon(*)\}$.  Now we can find an
$\bold F^a_\kappa$-saturated $N'$ which is 
$\bold F^a_\kappa$-constructible over $M + B'$ and 
$\bold F^a_\kappa$-saturated $N''$ which is $\bold
F^a_\kappa$-constructible over $N' \cup B_{\ell(*)}$.  By the choice of
$B'$, the model $N'$ is $\bold F^a_\kappa$-constructible also over $M \cup
B_{\ell(*)} \cup B'$ (by the same construction)  hence $N''$ is $\bold
F^a_\kappa$-constructible over $M + B_{\ell(*)} + B'$.

Clearly $N''$ is $\bold F^a_\kappa$-prime over $M + B_{\ell(*)} + B'$
and $N_{\ell(*)}$ is $\bold F^a_\kappa$-prime over $M + B_{\ell(*)} +
B'$ (as $B' \subseteq N_{\ell(*)}$, see clause $(\beta)$ above) and
$B'$ has cardinality $< \kappa$.  So there is an isomorphism $f$ from
$N''$ onto $N_{\ell(*)}$ over $M \cup B_{\ell(*)} \cup B$.  Renaming
\wilog \, $f = \text{ id}_{N''}$ so $N'' = N_{\ell(*)}$.

Lastly, we shall show that $(N',b,B')$ is a legal choice for
$(N_{\ell(*)+1},a_{\ell(*)+1},B_{\ell(*)+1})$.
\mn
Why?  The non-obvious clauses are $(c)(\beta),(\gamma)$ and (d) of
$\boxplus$.

First, for clause (d) obviously $B' \subseteq |N'|,b \in f(N)$ and
 $|B'| < \kappa$, so $(d)(\alpha),(\beta),(\gamma)$ hold and clause
 $(d)(\varepsilon)$ holds by the clause $(\gamma)$.  As for
 $(d)(\delta)$ assume $q \in \bold S^{< \omega}(M \cup B')$ is
$\bold F^a_\kappa$-isolated let $\bar c \in {}^{\omega >}(N')$
realize $q$, and let $B_q \subseteq M \cup B'$ be of cardinality $<
 \kappa$ such that stp$(\bar c,B_q) \vdash \text{ stp}(\bar c,M \cup
 B')$.  Now we have
stp$(\bar c,M \cup B') \vdash \text{ stp}(\bar c,M
\cup B_{\ell(*)} \cup B')$
 as otherwise we can find $\bar c'_\ell$ in ${\frak C}$ realizing
 stp$(\bar c,B_q)$ hence stp$(\bar c,M \cup B')$ for $\ell=1,2$ such
 that stp$(\bar c_1,M \cup B_{\ell(*)} \cup B') \ne \text{ stp}(\bar
 c_2,M \cup B_{\ell(*)} \cup B')$; so for some finite $\bar a
 \subseteq B_{\ell(*)},\bar d \subseteq M$ we have stp$(\bar c,\bar d
 \cup \bar a \cup B') \ne \text{ stp}(\bar c_2,\bar d \cup \bar a \cup
 B')$ so \wilog \, $\bar c_1,\bar c_2$ are from $N_{\ell(*)}$
 contradicting the choice of $\varepsilon(*)$.  
Let $\bar{\bold b}$ list $B'$ without
repetitions, so by the induction hypothesis stp$(\bar{\bold b} \char 94
\bar c,M \cup B_{\ell(*)}) \vdash \text{ stp}(\bar{\bold b} \char 94
\bar c,M \cup B_0 \cup \ldots \cup B_{\ell(*)})$ hence stp$(\bar c,M
\cup B_{\ell(*)} \cup \bar{\bold b}) \vdash \text{ stp}(\bar c,M \cup
B_0 \cup \ldots \cup B_{\ell(*)} \cup \bar{\bold b})$ so by the choice of
$\bar b$ and the previous sentence really clause $(d)(\delta)$ holds for the
choice of $(N_{\ell(*)+1},a_{\ell(*)+1},B_{\ell(*)+1})$ above.  

Second, concerning clause $(c)(\beta)$ of $\boxplus$, by the sentence
after the choices of $B',N'$ above, we know that $N'$ is $\bold
F^a_\kappa$-constructively over $M \cup B_{\ell(*)} \cup B'$
 so clearly stp$(N',M \cup B') \vdash \text{ stp}(N',M \cup B' \cup
B_{\ell(*)})$ hence stp$(B_{\ell(*)},M \cup B') \vdash 
\text{ stp}(B_{\ell(*)},N')$, so easily stp$(B_{\ell(*)},M \cup B'
\cup\{a_{\ell(*)}\}) \vdash \text{ stp}(B_{\ell(*)},N')$.

Now $B_{\ell(*)} \cup B'$ is $\bold F^a_\kappa$-atomic over $M \cup
\{a_{\ell(*)}\}$ being $\subseteq N_{\ell(*)}$ recalling
$\boxplus_{\ell(*)}(b)$ holds;  hence $B_{\ell(*)}$ is
$\bold F^a_\kappa$-atomic over $M \cup B' \cup \{a_{\ell(*)}\}$ hence
by the previous sentence $B_{\ell(*)}$ is $\bold F^a_\kappa$-atomic
over $N' + a_{\ell(*)}$ but $|B_{\ell(*)}| < \kappa$ hence it is
$\bold F^a_\kappa$-constructible over $N' + a_{\ell(*)}$.
As $N''$ is $\bold F^a_\kappa$-constructible over $B_{\ell(*)} \cup
N'$ by its choice, (and $a_{\ell(*)} \in B_{\ell(*)}$ by
$\boxplus_{\ell(*)}(d)(\beta))$, clearly $N''$ 
is also $\bold F^a_\kappa$-constructible over $N' \cup
\{a_{\ell(*)}\}$ as required in $(c)(\beta)$.  

Clause $\boxplus_\ell(c)(\gamma)$ means that $N_{\ell(*)} = N''$ is $\bold
F^a_\kappa$-constructible over $N' + a_0$.  
Now $N_{\ell(*)} = N''$ is
$\bold F^a_\kappa$-constructible over $B_{\ell(*)} \cup N'$ and $\bar
a \in {}^{\omega >}(N_{\ell(*)})$ implies stp$(\bar a,B_{\ell(*)} \cup
N') \vdash \text{ stp}(\bar a,B_0 \cup \ldots \cup B_{\ell(*)} \cup
N')$ hence by monotonicity stp$(\bar a,B_{\ell(*)} \cup N') \vdash
\text{ stp}(\bar a,a_0 + B_{\ell(*)} + N')$, so by the same
construction, $N_{\ell(*)} = N''$ is $\bold F^a_\kappa$-constructible
over $a_0 + B_{\ell(*)} + N'$.  As $B_{\ell(*)} \subseteq
N_{\ell(*)},|B_{\ell(*)}| < \kappa$ it is enough to show that $B_{\ell(*)}$ is
$\bold F^a_\kappa$-atomic over $a_0 + N'$ and this is proved as in the
proof of clause $(d)(\delta)$ above.
So indeed $(N',b,B')$ is a legal choice for
$(N_{\ell(*)+1},a_{\ell(*)+1},B_{\ell(*)+1})$. 
\nl
But this contradicts the choice of $\ell(*)$, so we have finished
proving $\boxtimes_1$.]
\mr
\item "{$\boxtimes_2$}"  if $\ell=m+1 < 1 + \ell(*)$ 
\ub{then} {\rm tp}$(a_m,N_\ell)$ is
not orthogonal to $M$.
\ermn
[Why?  Toward contradiction assume {\rm tp}$(a_m,N_\ell) \bot M$.  So 
we can find $A_\ell \subseteq N_\ell$ of cardinality $< \kappa$ such that
{\rm tp}$(\langle a_0,\dotsc,a_m\rangle,A_\ell)$ is stationary,
tp$(\langle a_0,\dotsc,a_m\rangle,N_\ell)$ does not fork over $A_\ell$
and {\rm tp}$(A_\ell,M)$ does not fork over $C_\ell := A_\ell
\cap M$ and tp$(A_\ell,C_\ell)$ is stationary 
 and $a_\ell \in A_\ell$ and (recalling $N_\ell$ is $\bold
F^a_\kappa$-primary over $M+a_\ell$) we have stp$(A_\ell,C_\ell +
a_\ell) \vdash \text{ stp}(A_\ell,M+a_\ell)$; 
it follows that {\rm tp}$(M,A_\ell)$
does not fork over $C_\ell$.  As {\rm tp}$(a_m,M+A_\ell)$ is
parallel to {\rm tp}$(a_m,N_\ell)$ and to 
{\rm tp}$(a_m,A_\ell)$ and tp$(a_m,N_\ell) \perp M$ is assumed we get
that all three types are orthogonal to $M$.  It follows that {\rm
stp}$(a_m,A_\ell) \vdash \text{\rm stp}(a_m,M+A_\ell)$ but recall
$a_\ell \in A_\ell$ so {\rm stp}$(a_m,A_\ell) \vdash 
\text{\rm stp}(a_m,M+a_\ell)$.  As $|A_\ell| < \kappa$ this implies that
{\rm tp}$(a_m,M+A_\ell)$ is $\bold
F^a_\kappa$-isolated.  But recall {\rm stp}$(A_\ell,C_\ell + a_\ell) = 
\text{\rm stp}(A_\ell,(A_\ell \cap M) + a_\ell)
\vdash \text{\rm stp}(A_\ell,M+a_\ell)$.  Together {\rm stp}$(a_m +
A_\ell,C_\ell + a_\ell) \vdash \text{\rm stp}(a_m + A_\ell,
M +a_\ell)$ hence {\rm tp}$(a_m,M+ a_\ell)$ is
$\bold F^a_\kappa$-isolated, contradicting $\boxtimes_\ell(c)(\alpha)$.]

To complete the proof by $\boxtimes_1$
it suffices to show $\ell(*) < \omega$, so
toward contradiction assume:
\mr
\item "{$\boxtimes_3$}"  $\ell(*) = \omega$.
\ermn
As we are assuming $\boxtimes_3$, we can find 
$\langle N^+_\ell:\ell < \ell(*) = \omega\rangle$ such that
\mr
\item "{$\odot_1$}"  $(a) \quad N_\ell \prec N^+_\ell$
\sn
\item "{${{}}$}"  $(b) \quad N_\ell$ is saturated, e.g. of
cardinality $\|N\|^{|T|}$
\sn
\item "{${{}}$}"  $(c) \quad N^+_{\ell +1} \prec N^+_\ell$
\sn
\item "{${{}}$}"  $(d) \quad \text{\rm tp}(N^+_\ell,N)$ does not
fork over $N_\ell$
\sn
\item "{${{}}$}"  $(e) \quad (N^+_\ell,c)_{c \in N_\ell \cup
N^+_{\ell +1}}$ is saturated.
\ermn
[Why?  We can choose $N^+_\ell$ by induction on $\ell$.  For $\ell=0$
it is obvious and for $\ell=m+1$ we choose $N'_\ell$ and satisfying
the relevant demands in $\odot_1$ on $N_\ell$ and then choose $N'_m$
satisfying the relevant demands on $(N_\ell,N_m)$.  Lastly, by the
uniqueness of saturated model there is an isomorphism $f_\ell$ from
$N'_m$ onto $N_m$ over $N_m$ and let $N_\ell = f_\ell(N'_\ell)$.]

Next for $\ell < \ell(*)$ we can find $\bold I_\ell$ such that
\mr
\item "{$\odot_2$}"   $(a) \quad \bold I_\ell \subseteq N^+_\ell
\backslash N^+_{\ell +1}$ 
\sn
\item "{${{}}$}"  $(b) \quad \bold I_\ell$ is independent over
$(N^+_{\ell +1},M)$ 
\nl

\hskip25pt (i.e. $c \in \bold I_\ell \Rightarrow \text{
tp}(c,N^+_{\ell+1})$ does not fork over $M$ and $\bold I$ is
independent over $N^+_{\ell +1}$)
\sn
\item "{${{}}$}"  $(c) \quad$ {\rm tp}$(N^+_\ell,N^+_{\ell +1} 
\cup \bold I_\ell)$ is almost orthogonal to $M$
\sn
\item "{${{}}$}"  $(d) \quad$ if $c \in \bold I_\ell$ then either $c \in
\varphi_*({\frak C},\bar d_*)$ or tp$(c,M)$ is orthogonal to
\nl

\hskip25pt $\varphi_*(x,\bar d_*)$, i.e. to every $q \in \bold S(M)$ to
which $\varphi_*(x,\bar d_*)$ belongs
\sn
\item "{${{}}$}"  $(e) \quad$ if $q \in \bold S(N^+_{\ell +1})$ does
not fork over $M$ and $\varphi_*(x,\bar d_*) \in q$ or $q$ is
\nl

\hskip25pt orthogonal to $\varphi_*(x,\bar d_*)$ \ub{then} 
the set $\{c \in \bold I_\ell:c$ realizes $q\}$ 
\nl

\hskip25pt has cardinality $\|N_\ell\|$
\sn
\item "{${{}}$}"  $(f) \quad$ we let $\bold I'_\ell = \bold I_\ell
\cap \varphi_*({\frak C},\bar d_*)$.
\ermn
[Why possible?  As $(N^+_\ell,c)_{c \in N^+_{\ell +1}})$ is 
saturated.]  

Now for $\ell < \ell(*)$
\mr
\item "{$\odot_3$}"    $\bold I_\ell$ is not independent over
$(N^+_{\ell +1} +a,N^+_{\ell +1})$.
\ermn
[Why?  Recall $a = a_0$.  Assume toward contradiction that
\mr
\item "{$(*)_{3.1}$}"    $\bold I_\ell$ is
independent over $(N^+_{\ell +1} +a,N^+_{\ell +1})$.
\ermn
As by clause (b) of $\odot_2$ we have 
{\rm tp}$(\bold I_\ell,N^+_{\ell +1})$ does 
not fork over $M$, it follows that $\bold I_\ell$ is independent over 
$(N^+_{\ell +1} +a,M)$.  Also by $(*)_{3.1}$ we know that tp$(a,N^+_{\ell
+1} \cup \bold I_\ell)$ does not fork over $N^+_{\ell +1}$.  Also
tp$(a,N^+_{\ell +1})$ does not fork over $N_{\ell +1}$ (because
$a \in N$ and tp$(N^+_{\ell +1},N)$ does not fork over $N_{\ell +1}$
by $\odot_1(d)$), together it follows that 
\mr
\item "{$(*)_{3.2}$}"  {\rm tp}$(a,N^+_{\ell +1} + \bold I_\ell)$ 
does not fork over $N_{\ell +1}$.
\ermn
Recall that tp$(N_\ell,N^+_{\ell +1})$ does not fork 
over $N_{\ell+1}$ (by $\odot_1(d)$ because $N_\ell \prec N$ using
symmetry) and tp$(a,N_\ell \cup N^+_{\ell +1})$ does not fork over 
$N_\ell$ similarly hence tp$(N_\ell + a,N^+_{\ell +1})$
does not fork over $N_{\ell +1}$, hence
\mr
\item "{$\odot_{3.3}$}"   tp$(N_\ell,N^+_{\ell +1} + a)$
does not fork over $N_{\ell +1} + a$.
\ermn
Recall $N_\ell$ is $\bold F^a_\kappa$-constructible over $N_{\ell
+1}+a$ (by $\boxplus_{\ell +1}(c)(\gamma)),N_\ell$ is $\bold
F^a_\kappa$-saturated and {\rm tp}$(N^+_{\ell +1},N_\ell + a)$ does
not fork over $N_{\ell +1}$ clearly 
\mr
\item "{$(*)_{3.4}$}"  $N_\ell$ is also $\bold
F^a_\kappa$-constructible over $N^+_{\ell +1} +a$ (even by the same
construction).
\ermn
As {\rm tp}$(a,N^+_{\ell +1} + \bold I_\ell)$ does not fork over 
$N_{\ell +1}$ and $N^+_{\ell+1}$ is $\bold F^a_\kappa$-saturated, 
it follows that
\mr
\item "{$(*)_{3.5}$}"  {\rm tp}$(N_\ell,N^+_{\ell
+1} + \bold I_\ell)$ does not fork over $N^+_{\ell +1}$ hence over
$N_{\ell +1}$.
\ermn
But by $\odot_2$ clause (c), for every 
$\bar d \in {}^{\omega>}(N^+_\ell)$ the type
{\rm tp}$(\bar d,N^+_{\ell+1} + \bold I_\ell)$ is almost orthogonal to 
$M$ hence recalling $N_\ell \subseteq N^+_\ell$,
\mr
\item "{$(*)_{3.6}$}"  {\rm tp}$(N_\ell,N^+_{\ell +1} + \bold I_\ell)$ 
is almost orthogonal to $M$ (this does not depend on $\odot_{3.1} -
\odot_{3.5}$ so can be used later).
\ermn
Hence by $(*)_{3.5} + (*)_{3.6}$ we have
\mr
\item "{$(*)_{3.7}$}"  {\rm tp}$(N_\ell,N_{\ell+1})$ 
is almost orthogonal to $M$.
\ermn
But $N_{\ell +1}$ is $\bold F^a_\kappa$-saturated so this implies
\mr
\item "{$(*)_{3.8}$}"   {\rm tp}$(N_\ell,N_{\ell +1})$ is orthogonal to $M$.
\ermn
But by $\boxplus_\ell(b)$
\mr
\item "{$(*)_{3.9}$}"  $a_\ell \in N_\ell$.
\ermn
By $\boxtimes_2$ we have
\mr
\item "{$(*)_{3.10}$}"  tp$(a_\ell,N_{\ell +1})$ is not orthogonal to $M$.
\ermn
Together $(*)_{3.8} + (*)_{3.9} + (*)_{3.10}$ give a 
contradiction, so $(*)_{3.1}$ fails hence $\odot_3$ holds.]

Now (recalling clause (f) of $\odot_2$)
\mr
\item "{$\odot_4$}"  $\bold I'_\ell$ is not independent over $(N^+_{\ell
+1} +a,N^+_{\ell +1})$.
\ermn
[Why?  By $\odot_3$ + clauses (b)+(d) of $\odot_2$ recalling that $a
\in \varphi_*({\frak C},\bar d_*)$ by the choice of $a$ in the
beginning of the proof of \scite{6n.28}.]
\mr
\item "{$\odot_5$}"  for each $n$, tp$(a,N^+_n)$ does 
$(\Bbb L(\tau_T),n)$-ict$^3$-fork over $M$.
\ermn
[Why?  By \scite{6n.31} below with $\bold I_\ell,N^+_{n-\ell}$ 
here standing for $\bold I_{n-\ell-1},N_\ell$ there, clause (d) 
there holding by $\odot_3$ here.  $M,A$ there standing for $M,M$ here,
clause (a),(b),(c) there holds by $(*)_{3.6}$ here (recalling that
$(*)_{3.6}$ does not depend on $\odot_{3.1} - \odot_{3.5}$.]
\mr
\item "{$\odot_6$}"   $\alpha_* > 
\text{ ict}^3-\text{rk(tp}(a,N^+_n))$ for every $n < \omega$.
\ermn
[Why?  By the choice of $\varphi_*(x,\bar d_*),a,\alpha_*$ in the
beginning of the proof we have $\alpha^* = \text{
ict}^3-\text{rk}(\text{tp}(a,M))$ and by $\odot_5$ and the definition
of ict$^3-\text{rk}(-)$ this follows.]
\mr
\item "{$\odot_7$}"  for each $n$, tp$(a,N^+_n)$ is not orthogonal to
$M$.
\ermn
[Why?  By $\odot_2(b) + \odot_4$.]

Hence we can find $q \in \bold S(M)$ such that:
\mr
\item "{$\odot_8$}"  $(a) \quad$ some automorphism of ${\frak C}$ over
$\bar d_*$ maps tp$(a,N_n)$ to a type 
\nl

\hskip25pt parallel to $q$
\sn
\item "{${{}}$}"  $(b) \quad \text{ ict}^3-\text{rk}(q) < \alpha_*$
\sn  
\item "{${{}}$}"  $(c) \quad q$ and $\text{tp}(a,M)$ are not orthogonal
\sn
\item "{${{}}$}"  $(d) \quad$ if $q' \subseteq q,|q'| < \kappa$ then
$q'(N) \nsubseteq M$
\nl

\hskip25pt [actually clause (d) follows by (c)].
\ermn
This contradicts the choice of $\alpha_*$; so $\ell(*) < \omega$ and
so we are done.  \hfill$\square_{\scite{6n.28}}$
\enddemo
\bigskip

\proclaim{\stag{6n.31} Claim}  Assume $T$ is stable.  A sufficient
condition for ``{\rm tp}$(a,N_n)$ does $(\Delta,n)-\text{\rm ict}^3$-divide
over $A$"
is:
\mr
\item "{$\circledast$}"  $(a) \quad \langle N_\ell:\ell \le  n\rangle$ is
$\prec$-increasing
\sn
\item "{${{}}$}"  $(b) \quad A \subseteq M \prec N_0$
\sn  
\item "{${{}}$}"  $(c) \quad \bold I_\ell \subseteq N_{\ell +1}
\backslash N_\ell$ is independent over $(N_\ell,M)$ for $\ell < n$
\sn
\item "{${{}}$}"  $(d) \quad$ {\rm tp}$(a,N_\ell \cup \bold I_\ell)$
forks over $N_{\ell +1}$
\sn
\item "{${{}}$}"  $(e) \quad$ {\rm tp}$(N_{\ell +1},N_\ell + \bold
I_\ell)$ is almost orthogonal to $M$.
\endroster
\endproclaim
\bigskip

\demo{Proof}  Left to the reader noting that $\langle \bold
I_\ell:\ell < n\rangle$ are pairwise disjoint (by clauses (a) +(c))
and $\cup\{\bold I_\ell:\ell < n\}$ is independent).  
\hfill$\square_{\scite{6n.31}}$
\enddemo
\bigskip

\remark{\stag{6n.32} Remark}  1) We may give more details on the last
 proof and intend to continue the investigation of the theory of 
regular types (in order to get good theory of weight) 
in this context somewhere else.
\nl
2) We can use essentially \scite{6n.31} to define a variant of the
rank for stable theory.  So \scite{6n.31} can be written to use it
and so \scite{6n.35} connect the two ranks.
\endremark
\bigskip

\proclaim{\stag{6n.35} Claim}  Assume $k \in \{3,4\}$ 
and {\rm ict}$^k$-{\rm rk}$(T) < \infty$, see Definition \scite{dw5.34}(6).

If {\rm cf}$(\kappa) \ge |T|^+$ or less and $M \prec N$ are
$\kappa$-saturated \ub{then} for some $a,\varphi(x,\bar a),n^*$ we have:
\mr
\item "{$\circledast$}"  $(a) \quad a \in N \backslash M$
\sn
\item "{${{}}$}"  $(b) \quad$ if $T$ 
is stable, the type $p = \text{\rm tp}(a,M)$ is primely regular
\sn
\item "{${{}}$}"  $(c) \quad \bar a 
\in {}^{\omega >} M$ and $\varphi(x,\bar a) \in p$
\sn
\item "{${{}}$}"  $(d) \quad \omega \times 
(\text{\rm wict}^k$-{\rm rk}$(\varphi(x,\bar a))) + 
(\text{\rm ict}^k - \text{\rm wg}(\varphi(x,\bar a)))$ is minimal.
\endroster
\endproclaim
\bigskip

\demo{Proof}  We choose $a,\varphi_*(x,\bar d_*),\alpha,n_*$ such that
\mr
\item "{$\circledast$}"  $(a) \quad a \in N \backslash M$
\sn
\item "{${{}}$}"  $(b) \quad \bar d_* \subseteq M$
\sn
\item "{${{}}$}"  $(c) \quad {\frak C} \models \varphi[a,\bar d_*]$
\sn
\item "{${{}}$}"  $(d) \quad \alpha = 
\text{\rm ict}^k-\text{\rm rk}(\{\varphi_*(x,\bar d_*)\})$
\sn
\item "{${{}}$}"  $(e) \quad$ under clauses (a)-(d), the ordinal
$\alpha$ is minimal
\sn
\item "{${{}}$}"  $(f) \quad n_*$ witness $\alpha +1 \nleq 
\text{ ict}^k-\text{rk}(\{\varphi(x,\bar d_*)\})$
\sn
\item "{${{}}$}"  $(g) \quad$ under clauses (a)-(f) the number $n_*(<
\omega)$ is minimal.
\ermn
Clearly there are such $a,\varphi_*(x,\bar c),\alpha$ and $n_*$.  Then
we try to choose $(N_\ell,a_\ell)$ by induction on $\ell < \omega$
such that $\boxplus_\ell$ from the proof of \scite{6n.28} holds.  But
now we can prove similarly that $\ell(*) \le n_*$.  But still
tp$(a,N_{\ell(*)})$ is not orthogonal to $M$.
\nl
[Why?  We can choose $N^+_0,\dotsc,N^+_{\ell(*)},\bold
I_0,\dotsc,\bold I_{\ell(*)-1}$ as in $\odot_2 + \odot_3$ in the proof of
\scite{6n.21} and prove $\odot_3$ there which implies the statement
above.  As $\varphi_*(x,\bar d_*) \in \text{ tp}(a,N_{\ell(*)})$ it
follows that $\varphi(N_{\ell(*)},\bar c) \nsubseteq M$ and any $a'
\in \varphi(N_{\ell(*)},\bar c) \backslash M$ is as required.]
\nl 
This is enough.  \hfill$\square_{\scite{6n.35}}$ 
\enddemo
\bn
Similarly to Definition \scite{dw5.34}.
\definition{\stag{6n.63} Definition}  Let $T$ be stable.
\nl
1) For an $m$-type $p(\bar x)$ we define sict$^3$-rk$^m(p(\bar x))$ as
 an ordinal or $\infty $ by defining when ict$^3$-rk$^m(p(\bar x))
\ge \alpha$ for an ordinal $\alpha$ by induction on $\alpha$
\mr
\item "{$(*)^\alpha_{p(\bar x)}$}"   sict$^3$-rk$^m(p(\bar x)) \ge
 \alpha$ iff for every $\beta < \alpha$ and finite $q(\bar x)
 \subseteq p(x)$ and $n < \omega$ we have
\sn
\item "{$(**)^{\beta,n}_{q(\bar x)}$}"  we can find $\langle
 M_\ell:\ell \le n\rangle,\langle \bold I_\ell:\ell < n\rangle$ and
 $\bar a$
{\roster
\itemitem{ $(a)$ }  $M_\ell \prec {\frak C}$ is $\bold
 F^a_{\kappa_1(T)}$-saturated
\sn
\itemitem{ $(b)$ }  $M_\ell \prec M_{\ell +1}$
\sn
\itemitem{ $(c)$ }  $q(\bar x)$ is an $m$-type over $M_0$
\sn
\itemitem{ $(d)$ }  $\bar a$ realizes $q(\bar x)$ and $\beta \le
\text{ sict}^3-\text{rk}(\text{tp}(\bar a,M_n)) \ge \beta$
\sn
\itemitem{ $(e)$ }  $\bold I_\ell \subseteq {}^{\omega >}(M_{\ell
+1})$ is independent over $(M_\ell,M_0)$
\sn
\itemitem{ $(f)$ }  $\bold I_\ell$ is not independent over $(M_\ell +
\bar a,M_0)$
\nl
(clearly \wilog \, $\bold I_\ell$ is a singleton).
\endroster}
\ermn
2) If sict$^3$-rk$^m(p(\bar x)) = \alpha < \infty$ then we let
sict$^3$-wg$^m(p(\bar x))$ be the maximal $n$ such that for every
finite $q(\bar x) \subseteq p(\bar x)$ we have
$(**)^{\alpha,n}_{q(\bar x)}$.
\nl
3) Above instead sict$^3$-rk$(\text{tp}(\bar a,A))$ we may write
sict$^3$-rk$^m(\bar a,A)$; similarly for scit$^3$-wg$^m(\bar a,A)$; if
$m=1$ we may omit it.
\enddefinition
\bigskip

\proclaim{\stag{6n.70} Claim}  1) $T$ is strongly$_3$ stable \ub{iff}
$T$ is stable and {\rm sict}$^3$-{\rm rk}$^m(p(\bar x)) < \infty$ for every
$m$-type $p(\bar x)$.
\nl
2) For every type $p(\bar x)$ there is a finite $q(\bar x) \subseteq
p(\bar x)$ such that ({\rm sict}$^3$-{\rm rk}$(p(\bar x))$, 
{\rm sict}$^3$-{\rm wg}$(p(\bar x)) = \text{\rm sict}^3-\text{\rm rk}
(q(\bar x))$, {\rm sict}$^3$-{\rm wg}$(q(\bar x)))$.
\nl
3) If $p(\bar x) \vdash q(\bar x)$ then {\rm sict}$^3$-rk$(p(\bar x)) \le
\text{\rm sict}^3$-{\rm rk}$(q(\bar x))$ and if equality holds \ub{then}
{\rm sict}$^3$-{\rm wg}$^m(p(\bar x)) \le 
\text{\rm sict}^3$-{\rm wg}$^m(q(\bar x))$.
\nl
4) ($T$ stable) If $p(\bar x),q(\bar x)$ are stationary parallel
types, \ub{then} {\rm sict}$^3$-{\rm rk}$^m(p(\bar x)) = 
\text{\rm sict}^3$-{\rm rk}$^m(q(\bar x))$, etc.  If $\bar a_1,\bar
a_1$ realizes $p \in \bold S^m(A)$ then 
$\text{\rm sict}^3$-{\rm rk}$^m(\text{rm stp}(\bar a_1,A)) =
\text{\rm sict}^3$-{\rm rk}$^m(\text{\rm stp}(\bar a_2,A))$.
Similarly for $\text{\rm sict}^3$-{\rm wg}$^m$.  Also automorphisms of
${\frak C}$ preserve $\text{\rm sict}^3$-{\rm rk}$^m$ and
$\text{\rm sict}^3$-{\rm wg}.
\endproclaim
\bigskip

\proclaim{\stag{6n.49} Claim}  $p(\bar x)$ does $(\Delta,n)$-{\rm ict}$^3$
forks over $A$ for every $n$ \ub{when}:
\mr
\item "{$\odot$}"  $(a) \quad G$ is a definable group over $A$ (in
${\frak C}$)
\sn
\item "{${{}}$}"  $(b) \quad b \in G$ realizes a generic type of $G$
from $\bold S(A)$ as was proved to exist in
\nl

\hskip25pt  \cite[4.11]{Sh:783}, or $T$ stable
\sn
\item "{${{}}$}"  $(c) \quad p(\bar x) \in \bold S^{< \omega}(A + b)$
forks over $A$.
\endroster
\endproclaim
\bigskip

\remark{Remark}  We may have said it in \S5(F).
\endremark
\bigskip

\demo{Proof of \scite{6n.49}}  Straight.
\enddemo
\bigskip

\demo{\stag{6n.50} Conclusion}:  Assume $T$ is strongly$_3$ dependent.

If $G$ is a type-definable group in ${\frak C}_T$ \ub{then} there is no
decreasing sequence $\langle G_n:n < \omega\rangle$ of subgroups of
$G$ such that $(G_n:G_{n+1}) = \bar\kappa$ for every $n$.
\enddemo
\bigskip

\remark{\stag{6n.51} Remark}  1) In \scite{6n.49} we can replace
``ict$^3$": by ``ict$^4$" and also by suitable variants for stable
theories.
\nl
2) Similarly in \scite{6n.50}. 
\endremark
\bn
(H) \quad \ub{$T$ is $n$-dependent}

On related problems and background see \cite[2.9-2.20]{Sh:702},
(but, concerning indiscernibility, it speaks on finite tuples,
i.e. $\alpha < \omega$ in \scite{nd.49}, which affect the
definitions and the picture).  On a consequence of ``$T$ is $2$-dependent" for
definable subgroups in ${\frak C}$ (and more, e.g. concerning \scite{nd.7}), 
see \cite{Sh:886}.
\bigskip

\definition{\stag{nd.5} Definition}  1) A (complete first order) theory
$T$ is $n$-independent when clause $(a)^n$ in \scite{nd.7} below holds.
\nl
2) The negation is$n$-dependent.
\enddefinition
\bn
\margintag{nd.7}\ub{\stag{nd.7} Problem}  Sort out the relationships between the
following candidates for ``$T$ is $n$-independent" 
($T$ is order order complete, also we can fix $\varphi$; omitting
$m$ we mean 1)
\mr
\item "{$(a)^n$}"  some $\varphi(\bar x,\bar y_0,\bar y_1,\dotsc,\bar
y_{n-1})$ is $n$-independent, i.e. $(a)^n_m$ for some $m$
\sn
\item "{$(a)^n_m$}"  some $\varphi(\bar x,\bar y_0,\bar
y_1,\dotsc,\bar y_{n-1})$ is $n$-independent where $\ell g(\bar x) =
m$ \ub{where}
{\roster
\itemitem{ $\odot$ }   $\varphi(\bar x,\bar y_0,\bar
y_1,\dotsc,\bar y_{n-1})$ is $n$-independent when there are
$\bar a^\ell_\alpha \in {}^{\ell g(\bar y_\ell)}{\frak C}$
for $\alpha <\lambda,\ell < n$ and $\langle \varphi(\bar
x,\bar a^0_{\eta(0)},\dotsc,\bar a^{n-1}_{\eta(n-1)}):\eta \in {}^n
\lambda$ is increasing$\rangle$ is an independent (sequence of formulas)
\endroster}
\item "{$(b)^n_m$}"  there is an indiscernible sequence $\langle \bar
a_\alpha:\alpha <\lambda\rangle,\varphi = \varphi(\bar x,\bar
y_0,\dotsc,\bar y_{n-1}),m = \ell g(\bar x),
\ell g(\bar y_\ell) = \ell g(\bar a_\alpha)$
for $\ell < n,\alpha < \lambda$ and $\bar c \in 
{}^{\ell g(\bar x)}{\frak C}$ such that: 
\nl
if $k<n$ and $\langle R_\ell:\ell < \ell(*)\rangle$ is a finite
sequence of $k$-place relations on $\lambda$ \ub{then} for some 
sequence $\bar t,\bar s \in {}^n \lambda$ realizing the same
quantifier free type
\nl
in $(\lambda,<,R_0,R_1,\dotsc,R_{\ell(\alpha)})$ we have ${\frak C}
\models \varphi[\bar b,\bar a_{s_0},\dotsc,\bar a_{s_{n-1}}] \wedge
\neg \varphi[\bar b,\bar a_{t_0},\dotsc,\bar a_{t_{n-1}}]$
\sn
\item "{$(c)^n_m$}"  for some $\varphi = \varphi(\bar x,\bar
y_0,\dotsc,\bar y_{n-1}),\ell g(\bar x) = m$, for every $j \in
[1,\omega)$, for infinitely many $k$
there are $\bar a^\ell_i \in {}^{\ell g(\bar y)}{\frak C}$ for $i<k$
such that $|\{p \cap\{\varphi(\bar x,\bar a^0_{i_0},\dotsc,\bar
a^{n-1}_{i_{n-1}}):i_\ell < k$ for $\ell <n\}:p \in 
{\bold S}^m(\cup\{\bar a^\ell_i:\ell < n,i<k\}|\}| \ge 2^{k^{n-1} \times m}$.
\endroster
\bigskip

\remark{Remark}  We can phrase $(b)^n_m,(c)^n_m$ as alternative
definitions of ``$\varphi(\bar x,\bar y_0,\dotsc,\bar y_{n-1})$ is
$n$-independent".  So in $(b)^n_m$ better to have $n$ indiscernible sequences.
\endremark
\bigskip

\demo{\stag{nd.14} Observation}   If $\varphi(\bar x,\bar y_0,
\dotsc,\bar y_{n-1})$ satisfies clause $(a)^n$ \ub{then} it satisfies a strong
form of clause $(c)^n$ (for every $k$ and the number is $\ge 2^{k^n}$.
\enddemo
\bigskip

\remark{Remark}  Clearly Observation \scite{nd.14} can be read as a sufficient
condition for being $n$-dependent, e.g.
\endremark
\bigskip

\demo{\stag{nd.17} Conclusion}  $T$ is $n$-dependent \ub{when}: for
every $m,\ell$ and finite $\Delta \subseteq \Bbb L(\tau_T)$ for
infinitely many $k < \omega$ we have $|A| \le k \Rightarrow |\bold
S^m_\Delta(A)| < 2^{(k/\ell)^n}$.
\enddemo
\bn
\margintag{nd.21}\ub{\stag{nd.21} Question}:  1) Can we get clause (a) from clause (c)?
\nl
2) Can we use it to prove $(a)^n_1 \equiv (a)^n_m$?
\bigskip

\demo{\stag{nd.28} Observation}  In \scite{nd.7}, if clause (a) then clause (b).
\enddemo
\bn
\margintag{nd.35}\ub{\stag{nd.35} Question}:  Does (b) imply (a)?
\bigskip

\proclaim{\stag{nd.42} Claim}  If $T$ satisfies $(a)^n$ for every $n$
\ub{then}: if $\lambda \nrightarrow (\mu)^{< \omega}_2$ then $\lambda
\nrightarrow_T (\mu)_{\aleph_0}$ where  
\endproclaim
\bigskip

\definition{\stag{nd.49} Definition}  We say that $\lambda \rightarrow_T
(\mu)_\alpha$ when: if $\bar a_i \in {}^\alpha({\frak C}_T)$ for $i <
\lambda$ \ub{then} for some ${\Cal U} \in [\lambda]^\mu$ the sequence
$\langle \bar a_i:i \in {\Cal U}\rangle$ is an indiscernible sequence
in ${\frak C}_T$.
\enddefinition
\bigskip

\remark{Remark}  1) Note that for $\alpha < \omega$ this property behaves
differently.
\nl
2) Of course, if $\theta = 2^{|\alpha|+|T|}$ and $\lambda \rightarrow
   (\mu)^{< \omega}_\theta$ then $\lambda \rightarrow_T (\mu)_\alpha$.
\nl
3) See on the non-2-independent $T$ and definable groups in \cite{Sh:886}.
\endremark
\bigskip

\demo{\stag{nd.56} Conjecture}  Assume $\neg(a)^n$ (or another variant
of $n$-dependent).  \ub{Then} ZFC $\vdash \forall \alpha \forall \mu
\exists \lambda(\lambda \rightarrow_T (\mu)_\alpha)$.
\enddemo
\bn
\margintag{nd.63}\ub{\stag{nd.63} Question}:  Can we phrase and prove a generalization
of the type-decomposition theorems for dependent theories (\cite{Sh:900})
to $n$-dependent theories $T$,
e.g. when $(\lambda^{\lambda_\ell}_{\ell +1}) = \lambda_{\ell +1}$ for
$\ell < n,{\frak B}_\ell \prec ({\Cal H}(\bar \kappa^+),\in,<^*_{\bar
\kappa^+})$ has cardinality $\lambda_\ell,[{\frak B}_{\ell
+1}]^{\lambda_\ell} \subseteq {\frak B}_\ell,\{{\frak C}_T,{\frak
B}_{\ell +1},\dotsc,{\frak B}_n\} \in {\frak B}_\ell$.

\nocite{ignore-this-bibtex-warning} 
\newpage
    
REFERENCES.  
\bibliographystyle{lit-plain}
\bibliography{lista,listb,listx,listf,liste}

\enddocument